\documentclass[a4paper]{article} 
\usepackage{amssymb}
\usepackage{amscd}
\usepackage{amsmath} 
\usepackage{graphicx}
\usepackage{latexsym} 
\usepackage{enumerate} 
\usepackage[usenames]{color}

\usepackage{theorem}
\theoremstyle{plain}
\theorembodyfont{\itshape}
\newtheorem{theorem}{Theorem}[section]
\newtheorem{proposition}[theorem]{Proposition}
\newtheorem{lemma}[theorem]{Lemma}
\newtheorem{corollary}[theorem]{Corollary} 
\theorembodyfont{\rmfamily}
\newtheorem{definition}[theorem]{Definition}
\newtheorem{remark}[theorem]{Remark}

\newtheorem{algorithm}[theorem]{Algorithm}
\newtheorem{example}[theorem]{Example}
\newtheorem{recipe}[theorem]{Recipe} 
\newtheorem{observation}[theorem]{Observation}
\makeatletter
\@addtoreset{figure}{section}
\@addtoreset{table}{section}
\makeatother
\def\bC{\mathbb{C}}
\def\bD{\mathbb{D}}

\def\bP{\mathbb{P}}
\def\bQ{\mathbb{Q}}
\def\bR{\mathbb{R}} 
\def\bZ{\mathbb{Z}}
\def\cA{\mathcal{A}}
\def\cB{\mathcal{B}}
\def\cC{\mathcal{C}}
\def\cE{\mathcal{E}} 
\def\cK{\mathcal{K}}
\def\cO{\mathcal{O}}
\def\cX{\mathcal{X}} 
\def\rA{\mathrm{A}}
\def\ra{\mathrm{a}}
\def\rb{\mathrm{b}}

\def\rC{\mathrm{C}}
\def\ch{\mathrm{ch}}
\def\CT{\mathrm{CT}}
\def\rD{\mathrm{D}}

\def\disc{\mathrm{disc}}
\def\rE{\mathrm{E}}
\def\re{\mathrm{e}}

\def\Gal{\mathrm{Gal}}
\def\GL{\mathrm{GL}}
\def\rH{\mathrm{H}}
\def\ri{\mathrm{i}}
\def\idx{\mathrm{idx}} 
\def\Idx{\mathrm{Idx}}
\def\rI{\mathrm{I}} 
\def\Im{\mathrm{Im}} 
\def\rK{\mathrm{K}}

\def\rL{\mathrm{L}}
\def\LT{\mathrm{LT}} 
\def\rM{\mathrm{M}}
\def\MT{\mathrm{MT}} 
\def\NS{\mathrm{NS}}
\def\NT{\mathrm{NT}}
\def\rP{\mathrm{P}}
\def\Pic{\mathrm{Pic}} 
\def\Par{\mathrm{Par}}
\def\PL{\mathrm{PL}}

\def\rank{\mathrm{rank}} 

\def\Res{\mathrm{Res}}

\def\rt{\mathrm{t}}
\def\td{\mathrm{td}} 
\def\Tr{\mathrm{Tr}}
\def\rII{\mathrm{II}}

\def\ba{\mbox{\boldmath $a$}}
\def\sba{\mbox{\scriptsize \boldmath $a$}}
\def\bA{\mbox{\boldmath $A$}}
\def\sbA{\mbox{\scriptsize \boldmath $A$}}
\def\bb{\mbox{\boldmath $b$}}

\def\bB{\mbox{\boldmath $B$}}
\def\sbB{\mbox{\scriptsize \boldmath $B$}}
\def\bc{\mbox{\boldmath $c$}}
\def\bd{\mbox{\boldmath $d$}}
\def\be{\mbox{\boldmath $e$}}
\def\sbe{\mbox{\scriptsize \boldmath $e$}}
\def\bk{\mbox{\boldmath $k$}}
\def\sbk{\mbox{\scriptsize \boldmath $k$}}
\def\bl{\mbox{\boldmath $l$}}
\def\sbl{\mbox{\scriptsize \boldmath $l$}}

\def\br{\mbox{\boldmath $r$}}

\def\bs{\mbox{\boldmath $s$}}

\def\bu{\mbox{\boldmath $u$}}
\def\sbu{\mbox{\scriptsize \boldmath $u$}}
\def\bv{\mbox{\boldmath $v$}}

\def\bx{\mbox{\boldmath $x$}}
\def\by{\mbox{\boldmath $y$}}

\def\bdelta{\mbox{\boldmath $\delta$}}
\def\b0{\mbox{\boldmath $0$}}
\def\bX{\mbox{\boldmath $X$}}
\def\vD{\varDelta}
\def\ve{\varepsilon}

\def\ul{\underline}

\def\baon{\ba_{\mathrm{on}}} 
\def\bAon{\bA_{\mathrm{on}}} 
\def\bAoff{\bA_{\mathrm{off}}} 
\def\bAin{\bA_{\mathrm{in}}} 
\def\cAin{\cA_{\mathrm{in}}} 

\def\bbon{\bb_{\mathrm{on}}}
\def\bBon{\bB_{\mathrm{on}}} 
\def\cBon{\cB_{\mathrm{on}}} 
\def\bBoff{\bB_{\mathrm{off}}} 

\def\baoff{\ba_{\mathrm{off}}} 
\def\bboff{\bb_{\mathrm{off}}}

\def\defby{\stackrel{\scriptstyle \mathrm{def}}{\Longleftrightarrow}}
\def\mi{\phantom{-}}
\title{\bf Hypergeometric Groups and Dynamics on K3 Surfaces\footnote{
MSC(2020): 14J28, 14J50, 33C80. Keywords: hypergeometric groups; K3 surfaces;   
automorphisms; entropy; unimodular lattices; Salem numbers; Lehmer's number; Siegel disks.}} 
\author{Katsunori Iwasaki\thanks{Department of Mathematics, Faculty of Science, 
Hokkaido University, Kita 10, Nishi 8, Kita-ku, Sapporo 060-0810 Japan. 
{\tt iwasaki@math.sci.hokudai.ac.jp} (corresponding author).} \ and 
Yuta Takada\thanks{Department of Mathematics, Graduate School of Science, 
Hokkaido University, Kita 10, Nishi 8, Kita-ku, Sapporo 060-0810 Japan. 
{\tt takada@math.sci.hokudai.ac.jp}}}
\date{April 30, 2021}  
\begin{document}
%%%%%%%%%%%%%%%%%%%%%%%%%%%%%%%%%%%%%%%%%%%%%%%%%%%%%%%%%%%%%%%%%%%%%%
\maketitle
%%%%%%%%%%%%%%%%%%%%%%%%%%% abstract %%%%%%%%%%%%%%%%%%%%%%%%%%%%%%%%
\begin{abstract} 
A hypergeometric group is a matrix group modeled on the monodromy 
group of a generalized hypergeometric differential equation. 
This article presents a fruitful interaction between the theory of hypergeometric 
groups and dynamics on K3 surfaces by showing that a certain class of 
hypergeometric groups and related lattices lead to a lot of K3 surface automorphisms 
of positive entropy, especially such automorphisms with Siegel disks.  
\end{abstract} 
%%%%%%%%%%%%%%%%%%%%%%%%%%% sec:intro %%%%%%%%%%%%%%%%%%%%%%%%%%%%%%%%%%%
\section{Introduction} \label{sec:intro}
%%%%%%%%%%%%%%%%%%%%%%%%%%%%%%%%%%%%%%%%%%%%%%%%%%%%%%%%%%%%%%%%%%%%%%
This article originates from a simple question:  
What happens if we put the following two topics together ? 
One is the theory of hypergeometric groups due to Levelt \cite{Levelt},  
Beukers and Heckman \cite{BH}, and the other is dynamics on K3 surfaces 
due to McMullen \cite{McMullen1,McMullen3,McMullen4}; see also 
Gross and McMullen \cite{GM}.    
In this article we present a fruitful interaction between them by showing 
that a certain class of hypergeometric groups and related lattices produce 
a lot of K3 surface automorphisms of positive entropy, especially 
such automorphisms with Siegel disks.  
%%%%%%%
\par
%%%%%%%
A hypergeometric group is a group modeled on the monodromy group of 
a generalized hypergeometric differential equation. 
It is a matrix group $H = \langle A, B \rangle \subset \GL(n, \bC)$ generated 
by two invertible matrices $A$ and $B$ such that $\rank(A-B) = 1$, which   
is equivalent to the condition $\rank(I-C) = 1$ for the third matrix 
$C := A^{-1} B$. 
In the context of an $n$-th order hypergeometric equation the matrices 
$A$, $B$, $C$ are the local monodromy matrices around the regular singular 
points $z = \infty$, $0$, $1$, respectively.   
The rank condition for $C$ is then a consequence of the property that 
the differential equation has $n-1$ linearly independent holomorphic  
solutions around $z = 1$. 
Beukers and Heckman \cite{BH} established several fundamental properties 
of hypergeometric groups such as irreduciblity, invariant Hermitian form, 
signature, etc. and went on to classify finite hypergeometric groups. 
Along the way they also determined differential Galois groups of hypergeometric 
equations, that is, Zariski closures of hypergeomtric groups.         
%%%%%%
\par
%%%%%%
On the other hand, McMullen \cite{McMullen1} synthesized examples of 
K3 surface automorphisms with Siegel disks. 
His constructions were based upon  
(i) K3 lattices and K3 structures in Salem number fields, 
(ii) Lefschetz and Atiyah-Bott fixed point formulas, (iii) Siegel-Sternberg theory on 
linearizations of nonlinear maps and small divisor problems, and (iv) Gel'fond-Baker 
method in transcendence theory and Diophantine approximation. 
In his work the characteristic polynomials of the constructed automorphisms 
were Salem polynomials of degree $22$, so the topological entropies of them were 
logarithms of Salem numbers of degree $22$ by the Gromov-Yomdin theorem 
\cite{Gromov,Yomdin}.  
McMullen \cite{McMullen3,McMullen4} went on to construct K3 surface 
automorphisms with Salem numbers of lower degrees, especially 
ones with Lehmer's number $\lambda_{\rL}$ in \cite{Lehmer}, whose 
logarithm was the minimum of the positive entropy spectrum for all 
automorphisms on compact complex surfaces \cite{McMullen2}. 
He discussed the non-projective cases in \cite{McMullen3} and the projective 
ones in \cite{McMullen4}, though these papers did not touch on Siegel disks. 
Here we should recall from \cite[Theorem 7.2]{McMullen1} that an  
automorphism on a projective K3 surface never admits a Siegel disk.             
%%%%%
\par
%%%%%
Our chief idea in this article is to use hypergeometric groups and   
associated lattices, in place of Salem number fields mentioned in 
item (i) above. 
To outline our idea we need to review the minimal basics  
about K3 surfaces and their automorphisms (see Barth et al. 
\cite[Chap. VI\!I\!I]{BHPV}). 
The middle cohomology group $L = H^2(X, \bZ)$ of a K3 surface $X$ equipped 
with the intersection form is an even unimodular lattice of rank $22$ and  
signature $(3, 19)$.  
The geometry of $X$ then defines a triple, called the {\sl K3 structure} on 
$L$, consisting of Hodge structure 
$L \otimes \bC = H^{2,0} \oplus H^{1,1} \oplus H^{0,2}$, 
positive cone $\cC^+(X)$ and K\"ahler cone $\cK(X)$.  
It is accompanied by the related concepts of Picard lattice 
$\Pic(X)$ (or N\'eron-Severi lattice $\NS(X)$), root system $\varDelta(X)$ 
and Weyl group $W(X)$. 
The surface $X$ is projective if and only if $\Pic(X)$ contains an element 
of positive self-intersection.    
Any K3 surface automorphism $f : X \to X$ induces a lattice 
automorphism $f^* : L \to L$ preserving the K3 structure. 
Conversely, thanks to the Torelli theorem and surjectivity of the period map,  
any automorphism of an abstract K3 lattice preserving an abstract K3 structure is 
realized by a unique K3 surface automorphism up to isomorphisms.  
An automorphism $f : X \to X$ gives rise to two basic invariants: the {\sl entropy}    
$h(f) = \log \lambda(f)$ with $\lambda(f)$ being the spectral radius of $f^*|H^2(X)$, 
and the {\sl special eigenvalue} $\delta(f) \in S^1$ defined by $f^* \eta = \delta(f) \, \eta$ 
for a nowhere vanishing holomorphic $2$-form $\eta$ on $X$, 
or rather the {\sl special trace} $\tau(f) := \delta(f) + \delta(f)^{-1} = 
\Tr( f^*|H^{2,0}(X) \oplus H^{0,2}(X)) \in [-2, \, 2]$. 
Note that $\lambda(f)$ is either $1$ or a Salem number $\lambda > 1$, 
while $\delta(f)$ is either a root of unity or a conjugate to the Salem number $\lambda$. 
The spectral radius $\lambda(F)$, special eigenvalue $\delta(F)$ and special trace 
$\tau(F)$ are conceivable for a Hodge isometry $F : L \to L$ of an abstract 
K3 lattice $L$ with a Hodge structure.  
$F$ is said to be {\sl positive} if it preserves the positive cone.   
It falls into one of the elliptic, parabolic and hyperbolic types 
(see e.g. Cantat \cite[\S 1.4]{Cantat}). 
%%%%%
\par
%%%%% 
We turn our attention to hypergeometric groups. 
Let $\varphi(z)$ and $\psi(z)$ be a coprime pair of monic polynomials over $\bZ$ 
of an even degree $n = 2 N$ such that $\varphi(z)$ is anti-palindromic and 
$\psi(z)$ is palindromic, that is, 
$z^n \varphi(z^{-1}) = - \varphi(z)$ and $z^n \psi(z^{-1}) = \psi(z)$. 
Consider the hypergeometric group $H = \langle A, B \rangle \subset \GL(n, \bZ)$ 
generated by $A := Z(\varphi)$ and $B := Z(\psi)$, where $Z(f)$ is the companion matrix 
of a monic polynomial $f(z) = z^n + f_1 z^{n-1} + \cdots + f_n$,  
%%%%%%%%%%%%%%%%%%%%%%%%%% eqn:cm %%%%%%%%%%%%%%%%%%%%%%%%%%%%%%%%%%%%%
\begin{equation*} \label{eqn:cm}
Z(f) := 
\begin{pmatrix} 
0 &    &            &    &  -f_n      \\
1 & 0 &            &    &  -f_{n-1} \\
   & 1 & \ddots &    &  \vdots \\
   &    & \ddots & 0 &  -f_2      \\
   &    &           & 1 &  -f_1 
\end{pmatrix}.  
\end{equation*}
%%%%%%%%%%%%%%%%%%%%%%%%%%%%%%%%%%%%%%%%%%%%%%%%%%%%%%%%%%%%%%%%%%%%% 
%%%%%%%%%%%%%%%%%%%%%%%%%%%%% thm:main1 %%%%%%%%%%%%%%%%%%%%%%%%%%%%%%% 
\begin{theorem} \label{thm:main1} 
The matrix $C := A^{-1} B$ is a reflection, that is, it fixes a hyperplane in $\bQ^n$ 
pointwise and sends a nonzero vector $\br \in \bQ^n$ to its negative $-\br$.  
We have a free $\bZ$-module of rank $n$, 
%%%%%%%%%%%%%%%%%%%%%%%%%%%% eqn:AB %%%%%%%%%%%%%%%%%%%%%%%%%%%%%%%%%
\begin{equation} \label{eqn:AB}
L := \langle \br, A \br, \dots, A^{n-1} \br \rangle_{\bZ} = 
\langle \br, B \br, \dots, B^{n-1} \br \rangle_{\bZ}, 
\end{equation}
%%%%%%%%%%%%%%%%%%%%%%%%%%%%%%%%%%%%%%%%%%%%%%%%%%%%%%%%%%%%%%%%%%%
which is stable under the action of $H$. 
There exists a unique, non-degenerate, $H$-invariant, $\bZ$-valued symmetric 
bilinear form on $L$ such that $(\br, \br) = 2$, which makes $L$ into an even lattice.   
Its Gram matrix for the $A$-basis is given by $(A^{i-1} \br, A^{j-1} \br) = \xi_{|i-j|}$, 
where $\xi_0 := 2$ and $\{\xi_i\}_{i=1}^{\infty}$ is defined via the 
Taylor series expansion 
%%%%%%%%%%%%%%%%%%%%%%%%%% eqn:xi %%%%%%%%%%%%%%%%%%%%%%%%%%%%%%%%%%
\begin{equation} \label{eqn:xi}
\dfrac{\psi(z)}{\varphi(z)} = 1 + \sum_{i=1}^{\infty} \xi_i \, z^{-i} \qquad 
\mbox{around} \quad z = \infty.  
\end{equation}
%%%%%%%%%%%%%%%%%%%%%%%%%%%%%%%%%%%%%%%%%%%%%%%%%%%%%%%%%%%%%%%%%%
The Gram matrix for the $B$-basis is given in a similar manner by exchanging 
$\varphi(z)$ and $\psi(z)$ upside down in the Taylor expansion \eqref{eqn:xi}.   
The lattice $L$ is unimodular if and only if the resultant of $\varphi(z)$ and $\psi(z)$ is 
$\pm 1$.  
\end{theorem}
%%%%%%%%%%%%%%%%%%%%%%%%%%%%%%%%%%%%%%%%%%%%%%%%%%%%%%%%%%%%%%%%%%%%%%
\par
%%%%%%%%%
Since $\varphi(z)$ and $\psi(z)$ are anti-palindromic and palindromic of 
degree $n = 2N$, there exist monic polynomials $\Phi(w)$ and 
$\Psi(w)$ of degrees $N-1$ and $N$ over $\bZ$ such that      
$\varphi(z) = (z^2-1) z^{N-1} \Phi(z+z^{-1})$ and  
$\psi(z) = z^N \Psi(z+z^{-1})$. 
We refer to $\Phi(w)$ and $\Psi(w)$ as the {\sl trace polynomials} of 
$\varphi(z)$ and $\psi(z)$. 
Being coprime, $\varphi(z)$ and $\psi(z)$ have no roots in common 
and the same is true for $\Phi(w)$ and $\Psi(w)$.   
Let $\bA$ be the multi-set of all complex roots of $\Phi(w)$ counted with 
multiplicities. 
Let $\bAon$ and $\bAoff$ be those parts of $\bA$ which lie on and off 
the interval $[-2, \, 2]$ respectively. 
Define $\bB$, $\bBon$ and $\bBoff$ in a similar manner for $\Psi(w)$. 
Then $\bAon$ and $\bBon$ dissect each 
other into interlacing components, $\bA_1, \dots, \bA_{s+1}$ and 
$\bB_1, \dots, \bB_s$, called {\sl trace clusters}, such that 
%%%%%%%%%%%%%%%%%%%%%%%% eqn:trc %%%%%%%%%%%%%%%%%%%%%%%%%%%%%%%%%%%%%
\begin{equation} \label{eqn:trc}
-2 \le \bA_{s+1} < \bB_s < \bA_s < \dots < \bB_1 < \bA_1 \le 2, 
\end{equation}
%%%%%%%%%%%%%%%%%%%%%%%%%%%%%%%%%%%%%%%%%%%%%%%%%%%%%%%%%%%%%%%%%%%% 
where one or both of the end clusters $\bA_1$ and $\bA_{s+1}$ may be null, 
while any other cluster must be non-null.   
Let $\bA_{>2}$ and $\bB_{>2}$ be those parts of $\bA$ and $\bB$ which 
lie in $(2, \, \infty)$ respectively. 
We give a formula representing the index of $L$ in terms of the 
trace clusters \eqref{eqn:trc} as well as $\bA_{>2}$ and $\bB_{>2}$ 
(see Theorem \ref{thm:real}). 
It in particular says that the index up to sign $\pm$ depends only on the trace 
clusters \eqref{eqn:trc}, being independent of $\bA_{>2}$ and $\bB_{>2}$.      
%%%%%%
\par
%%%%%%
We focus on the specific rank $n = 22$, i.e. $N = 11$.  
It is natural to ask when a hypergeometric lattice $L$ or its negative $L(-1)$ 
becomes a K3 lattice with a Hodge structure such that the matrix $A$ 
(or $B$) is a positive Hodge isometry. 
Here we take $L$ or $L(-1)$ depending on whether the index of $L$ is 
negative or positive, because a K3 lattice has negative index $-16$. 
This procedure is called the {\sl renormalization} of $L$ and the renormalized 
$H$-invariant symmetric bilinear form is referred to as the intersection form on $L$.   
To state our theorems we introduce some notation and terminology:   
$|\bAon|$ stands for the cardinality of $\bAon$ counted with multiplicities; 
$[\bAon] = 0^{\nu_0} 1^{\nu_1} 2^{\nu_2} 3^{\nu_3}$ means that 
$\bAon$ consists of $\nu_0$ null clusters, $\nu_1$ simple clusters, 
$\nu_2$ double clusters, $\nu_3$ triple clusters, where $j^{\nu_j}$ is  
omitted if $\nu_j = 0$.   
The same rule applies to $\bBon$ and other related entities. 
By ``doubles adjacent'' we mean the situation in which there exist a unique 
double $\bAon$-cluster $\bA_i$ and a unique double $\bBon$-cluster $\bB_j$ 
and they are adjacent to each other. 
For such a pair if $\bA_i \cup \bB_j$ consists of four distinct elements 
$\lambda_1 < \lambda_2 < \lambda_3 < \lambda_4$, then $\lambda_2$ 
and $\lambda_3$ are referred to as the inner elements of the adjacent pair (AP).        
%%%%%%%%%%%%%%%%%%%%%%%%%% thm:main2  %%%%%%%%%%%%%%%%%%%%%%%%%%%%%%%%%
\begin{theorem} \label{thm:main2}
Let $L = L(\varphi, \psi)$ be a unimodular  
hypergeometric lattice of rank $22$. 
After renormalization, $L$ is a K3 lattice with a Hodge structure such that the 
matrix $A$ is a positive Hodge isometry, if and only if $\Phi(-2) \neq 0$ and 
the roots of $\Phi(w)$ and $\Psi(w)$ are all simple and have any one of the 
following configurations. 
%%%%%%%%%%%%%%%%%
\begin{enumerate}
\setlength{\itemsep}{-1pt}
\item[$(\rE)$] In elliptic case, any entry of Table $\ref{tab:ep-A}$ such that $\bA_1$ 
does not contain $2$, that is, $\Phi(2) \neq 0$.      
\item[$(\rP)$] In parabolic case, any entry of Table $\ref{tab:ep-A}$ 
such that $\bA_1$ does contain $2$, that is, $\Phi(2) = 0$.   
\item[$(\rH)$] In hyperbolic case, any entry of Table $\ref{tab:hyp-A}$ such that 
$\bA_1$ does not contain $2$, that is, $\Phi(2) \neq 0$.    
\end{enumerate}
%%%%%%%%%%%%%%%%  
The special trace $\tau(A)$ and the Hodge structure up to complex conjugation 
are uniquely determined by the pair $(\varphi, \psi)$.   
The location of $\tau(A)$ is exhibited in the last columns of Tables $\ref{tab:ep-A}$ 
and $\ref{tab:hyp-A}$, where we mean by ``middle of TC" that $\tau(A)$ is 
the middle element of the unique triple $\bAon$-cluster, and by ``inner of AP" 
that it is the inner element in $\bAon$ of the unique adjacent pair of double 
clusters in $\bAon \cup \bBon$.  
%%%%%%%%%%%%%%%%%%%%%%%%%%%% tab:ep-A %%%%%%%%%%%%%%%%%%%%%%%%%%%%
\begin{table}[hh]
\centerline{
\begin{tabular}{ccccccl|l}
\hline
\\[-4mm]
case & $s$ & $[\bAon]$ & $[\bBon]$ & $\bA_1$ & $|\bBoff|$ & constraints & ST $\tau(A)$ 
\\[1mm]
\hline
\\[-4mm]
$1$ & $8$ & $0^11^73^1$ & $1^8$ & non-null & $3$ & & middle of TC \\[1mm]
$2$ & $8$ & $0^11^73^1$ & $1^73^1$ & non-null & $1$ & & middle of TC \\[1mm]
$3$ & $9$ & $0^21^73^1$ & $1^82^1$ & null & $1$ & $|\bB_1| = 2$ & middle of TC \\[1mm]
$4$ & $8$ & $1^82^1$ & $1^8$ & non-null & $3$ & $|\bA_9|=2$ & $\min \bA_9$ \\[1mm]    
$5$ & $8$ & $1^82^1$ & $1^73^1$ & non-null & $1$ & $|\bA_9|=2$ & $\min \bA_9$ \\[1mm]
$6$ & $9$ & $0^11^82^1$ & $1^82^1$ & non-null & $1$ & doubles adjacent & inner of AP \\[1mm]
$7$ & $9$ & $0^11^82^1$ & $1^82^1$ & null & $1$ & $|\bA_{10}|=2$, $|\bB_1|=2$ & $\min \bA_{10}$ \\[1mm]
$8$ & $10$ & $0^21^82^1$ & $1^{10}$ & null & $1$ & $|\bA_2|=2$ & $\max \bA_2$ \\[1mm]
$9$ & $9$ & $1^{10}$ & $1^82^1$ & non-null & $1$ & $|\bB_9|=2$ & element of $\bA_{10}$ \\[1mm]
\hline   
\end{tabular}}
\caption{Conditions for $A$ to be a positive Hodge isometry of elliptic or parabolic type.} 
\label{tab:ep-A} 
\end{table}
%%%%%%%%%%%%%%%%%%%%%%%%%%%%%%%%%%%%%%%%%%%%%%%%%%%%%%%%%%%%%%%%%
%%%%%%%%%%%%%%%%%%%%%%%%%%%% tab:hy-A %%%%%%%%%%%%%%%%%%%%%%%%%%%%
\begin{table}[hh]
\centerline{
\begin{tabular}{ccccccl|l}
\hline
\\[-4mm]
case & $s$ & $[\bAon]$ & $[\bBon]$ & $|\bA_{>2}|$ & $|\bBoff|$ & constraints & ST $\tau(A)$ 
\\[1mm]
\hline
\\[-4mm]
$1$ & $8$ & $0^21^63^1$ & $1^8$ & $1$ & $3$ & & middle of TC \\[1mm]
$2$ & $8$ & $0^21^63^1$ & $1^73^1$ & $1$ & $1$ & & middle of TC \\[1mm]
$3$ & $8$ & $0^11^72^1$ & $1^8$ & $1$ & $3$ & $|\bA_1| = 2$ & $\max \bA_1$ \\[1mm]
$4$ & $8$ & $0^11^72^1$ & $1^8$ & $1$ & $3$ & $|\bA_9|=2$ & $\min \bA_9$ \\[1mm]    
$5$ & $8$ & $0^11^72^1$ & $1^73^1$ & $1$ & $1$ & $|\bA_1|=2$ & $\max \bA_1$ \\[1mm] 
$6$ & $8$ & $0^11^72^1$ & $1^73^1$ & $1$ & $1$ & $|\bA_9|=2$ & $\min \bA_9$ \\[1mm] 
$7$ & $9$ & $0^21^72^1$ & $1^82^1$ & $1$ & $1$ & doubles adjacent & inner of AP \\[1mm] 
$8$ & $9$ & $0^11^9$ & $1^82^1$ & $1$ & $1$ & $|\bA_1|=1$, $|\bB_1|=2$ & element of $\bA_1$ \\[1mm] 
$9$ & $9$ & $0^11^9$ & $1^82^1$ & $1$ & $1$ & $|\bA_{10}|=1$, $|\bB_9|=2$ & element of $\bA_{10}$ \\[1mm]
\hline   
\end{tabular}} 
\caption{Conditions for $A$ to be a positive Hodge isometry of hyperbolic type.} 
\label{tab:hyp-A} 
\end{table}
%%%%%%%%%%%%%%%%%%%%%%%%%%%%%%%%%%%%%%%%%%%%%%%%%%%%%%%%%%%%%%%%% 
\end{theorem}
%%%%%%%%%%%%%%%%%%%%%%%%%%%%%%%%%%%%%%%%%%%%%%%%%%%%%%%%%%%%%%%%%%%%%  
\par
%%%%%%%%%
Due to the asymmetry of $\varphi(z)$ and $\psi(z)$ the corresponding result 
for the matrix $B$ is somewhat different and a dominant role is played by 
$\bBon$ and $\bAin := \bA_2 \cup \cdots \cup \bA_s$ in place of $\bAon$.   
%%%%%%%%%%%%%%%%%%%%%%% thm:main3 %%%%%%%%%%%%%%%%%%%%%%%%%%%%%%%%%%
\begin{theorem}  \label{thm:main3} 
Let $L = L(\varphi, \psi)$ be a unimodular hypergeometric lattice of rank $22$. 
After renormalization, $L$ is a K3 lattice with a Hodge structure such that the 
matrix $B$ is a positive Hodge isometry, if and only if all roots of $\Psi(w)$ are 
simple, $\Psi(\pm2) \neq 0$, and the roots of $\Phi(w)$ and $\Psi(w)$ have 
any one of the configurations in Table $\ref{tab:hyp-B}$. 
The special trace $\tau(B)$ and the Hodge structure up to complex conjugation 
are uniquely determined by the pair $(\varphi, \psi)$.  
The location of $\tau(B)$ is exhibited in the last column of Table $\ref{tab:hyp-B}$,  
where we mean by ``middle of TC" that $\tau(B)$ is the middle element of the 
unique triple $\bBon$-cluster, and by ``inner of AP" that it is the inner element 
in $\bBon$ of the unique adjacent pair of double clusters in $\bAin \cup \bBon$.  
The Hodge isometry $B$ is always of hyperbolic type and any root of 
$\Phi(w)$ is simple except for at most one integer root of 
multiplicity $2$ or $3$. 
%%%%%%%%%%%%%%%%%%%%%%%%%% tab:hyp-B %%%%%%%%%%%%%%%%%%%%%%%%%%%%
\begin{table}[hh]
\centerline{
\begin{tabular}{ccccccl|l}
\hline
\\[-4mm]
case & $s$ & $[\bAin]$ & $[\bBon]$ & $|\bAin|$ & $|\bB_{>2}|$ & constraints & ST $\tau(B)$ 
\\[1mm]
\hline
\\[-4mm]
$1$ & $8$ & $1^7$ & $1^73^1$ & $7$ & $1$ &  & middle of TC \\[1mm]
$2$ & $8$ & $1^63^1$ & $1^73^1$ & $9$ & $1$ &  & middle of TC \\[1mm]
$3$ & $9$ & $1^72^1$ & $1^82^1$  & $9$ & $1$ & doubles adjacent & inner of AP \\[1mm]
$4$ & $9$ & $1^8$ & $1^82^1$ & $8$ & $1$ & $|\bB_1| = 2$ & $\max \bB_1$ \\[1mm] 
$5$ & $9$ & $1^8$ & $1^82^1$ & $8$ & $1$ & $|\bB_9| = 2$ & $\min \bB_9$ \\[1mm] 
$6$ & $9$ & $1^73^1$ & $1^82^1$ & $10$ & $1$ & $|\bB_1| = 2$ & $\max \bB_1$ \\[1mm] 
$7$ & $9$ & $1^73^1$ & $1^82^1$ & $10$ & $1$ & $|\bB_9| = 2$ & $\min \bB_9$ \\[1mm] 
$8$ & $10$ & $1^82^1$ & $1^{10}$ & $10$ & $1$ & $|\bA_2| = 2$ & element of $\bB_1$ \\[1mm] 
$9$ & $10$ & $1^82^1$ & $1^{10}$ & $10$ & $1$ & $|\bA_{10}| = 2$ & element of $\bB_{10}$ \\[1mm] 
\hline  
\end{tabular}} 
\caption{Conditions for $B$ to be a positive Hodge isometry, always of hyperbolic type.} 
\label{tab:hyp-B} 
\end{table}
%%%%%%%%%%%%%%%%%%%%%%%%%%%%%%%%%%%%%%%%%%%%%%%%%%%%%%%%%%%%%%  
\end{theorem}
%%%%%%%%%%%%%%%%%%%%%%%%%%%%%%%%%%%%%%%%%%%%%%%%%%%%%%%%%%%%%%%
%%%%%%%%%%%%%%%%%%%%%%% rem:main1 %%%%%%%%%%%%%%%%%%%%%%%%%%%%%% 
\begin{remark} \label{rem:main1} 
Put $\chi(z) = \varphi(z)$ and $F = A$ in Theorem \ref{thm:main2}, 
while $\chi(z) = \psi(z)$ and $F = B$ in Theorem \ref{thm:main3}. 
\begin{enumerate}
\setlength{\itemsep}{-1pt}
\item Any root of $\chi(z)$ is simple except for the triple root 
$z = 1$ in the parabolic case of Theorem \ref{thm:main2}. 
\item In the elliptic and parabolic cases of Theorem \ref{thm:main2},  
any irreducible component of $\chi(z)$ is a cyclotomic polynomial and 
the spectral radius $\lambda(F)$ is equal to $1$.   
\item In the hyperbolic case of Theorem \ref{thm:main2} and in the 
whole case of Theorem \ref{thm:main3}, $\chi(z)$ has a unique Salem 
polynomial factor with the associated Salem number $\lambda > 1$   
giving the spectral radius $\lambda(F)$, and any other irreducible 
component of $\chi(z)$ is a cyclotomic polynomial. 
\item In every case of Theorems \ref{thm:main2} and \ref{thm:main3} 
the special trace $\tau(F)$ lies in $(-2, \, 2)$, so that the special 
eigenvalue $\delta(F) \in S^1$ is either a root of unity different from 
$\pm 1$ or a conjugate to the Salem number $\lambda$.      
\end{enumerate}
\end{remark}
%%%%%%%%%%%%%%%%%%%%%%%%%%%%%%%%%%%%%%%%%%%%%%%%%%%%%%%%%%%%%%
\par
%%%%%%%%%%
Let $\Pic := L \cap H^{1,1}$ be the ``Picard lattice" arising from the 
constructions in Theorems $\ref{thm:main2}$ and $\ref{thm:main3}$. 
There are two cases: one is the {\sl projective} case where $\Pic$ contains 
a vector of positive self-intersection and the other is the {\sl non-projective} 
case where no such vector exists.      
By item (4) of Remark \ref{rem:main1} the minimal polynomial $\chi_0(z)$ 
of $\delta(F)$ is either a cyclotomic polynomial of degree $\ge 2$ or a 
Salem polynomial of degree $\ge 4$. 
Note that the degree of $\chi_0(z)$ is necessarily even in either case.   
%%%%%%%%%%%%%%%%%%%%%%%% thm:main4 %%%%%%%%%%%%%%%%%%%%%%%%%%%%
\begin{theorem} \label{thm:main4} 
The intersection form is non-degenerate on $\Pic$ and the Picard number 
$\rho := \rank \, \Pic$ is given by $\rho = 22 - \deg \chi_0(z)$, 
which is even and $\le 20$.    
The projective case occurs precisely when $\chi_0(z)$ is 
a cylcotomic polynomial of degree $\ge 2$, while  the 
non-projective case is exactly the case where $\chi_0(z)$ is 
a Salem polynomial of degree $\ge 4$ and the associated 
Salem number yields the spectral radius $\lambda(F)$. 
In particular both 
the elliptic and parabolic cases in 
Theorem $\ref{thm:main2}$ are projective.  
If $\br$ is the vector in Theorem $\ref{thm:main1}$ and 
$\bs := \chi_0(F) \br$, then the vectors 
$\bs, F \bs, \dots, F^{\rho-1} \bs$ form a free $\bZ$-basis 
of $\Pic$, which we call the standard basis.     
\end{theorem} 
%%%%%%%%%%%%%%%%%%%%%%%%%%%%%%%%%%%%%%%%%%%%%%%%%%%%%%%%%%%%%%
\par
%%%%%%%%%%
Theorem \ref{thm:main1} enables us to calculate the Gram matrix of $\Pic$ 
with respect to the standard basis and hence its discriminant.  
The root system is defined by $\vD := \{ \bu \in \Pic : (\bu, \bu) = -2 \}$ 
and the Weyl group $W$ is the group generated by all reflections in root 
vectors. 
Usually, $F$ does not preserve any Weyl chamber. 
Fix a Weyl chamber $C$.  
It may be sent to a different Weyl chamber $F(C)$, but there exists a unique 
element $w_F \in W$ that brings $F(C)$ back to the original chamber 
$C$, since $W$ acts on the set of Weyl chambers simply transitively.  
In this article we discuss the non-projective case only, in which the 
intersection form is negative definite on $\Pic$ so that  the root 
system $\vD$ and the Weyl group $W$ are finite.  
Lexicographical order with respect to the standard basis in Theorem 
\ref{thm:main4} leads to a set of positive roots $\vD^+$ and the 
corresponding Weyl chamber $\cK$.          
%%%%%%%%%%%%%%%%%%%%%%%%% thm:main5 %%%%%%%%%%%%%%%%%%%%%%%%%%%
\begin{theorem} \label{thm:main5} 
In the non-projective case there exists an algorithm to output the unique 
element $w_F \in W$ such that the modified matrix $\tilde{F} := w_F \circ F$ 
preserves the Weyl chamber $\cK$ (Algorithm $\ref{algorithm}$).   
This twist does not change the Hodge structure, spectral radius 
$\lambda(F) = \lambda(\tilde{F})$ and special trace $\tau(F) = \tau(\tilde{F})$.  
The modified Hodge isometry $\tilde{F} : L \to L$ lifts to a non-projective 
K3 surface automorphism $f : X \to X$ of positive entropy 
$h(f) = \log \lambda(F)$ with special trace $\tau(f) = \tau(F)$, Picard lattice 
$\Pic(X) \cong \Pic$ and Picard number $\rho(X) = \rho$ as in 
Theorem $\ref{thm:main4}$, while $\cK$ lifting to the 
K\"{a}hler cone of $X$. 
The map $f$ is uniquely determined by the pair $(\varphi, \psi)$ plus the choice 
of $F = A$ or $B$ up to biholomorphic conjugacy and complex conjugation.   
\end{theorem}
%%%%%%%%%%%%%%%%%%%%%%%%%%%%%%%%%%%%%%%%%%%%%%%%%%%%%%%%%%%%%%
\par
%%%%%%% 
The Weyl chamber $\cK$ is chosen just for the sake 
of algorithmic purpose.  
Taking any other Weyl chamber results in a modified matrix 
conjugate to the one in Theorem \ref{thm:main5} by a unique 
Weyl group element. 
It is an interesting future problem to consider what is going on in the 
projective case, in which the spectral radius $\lambda(\tilde{F})$ may 
change from the original one $\lambda(F)$. 
This case should produce automorphisms of projective K3 surfaces. 
%%%%%%%
\par
%%%%%%%
Back in the non-projective case we illustrate our method by working 
out some examples. 
Let $\lambda_1, \dots, \lambda_{10}$ be the ten Salem numbers of 
degree $22$ from McMullen \cite[Table 4]{McMullen1} reproduced in 
this article as Table \ref{tab:mcmullen}, and $S_1(z), \dots, S_{10}(z)$ 
be the associated Salem polynomials. 
Let $\lambda_{\rL}$ be Lehmer's number, the smallest Salem number 
ever known, and $\rL(z)$ be Lehmer's polynomial (see \eqref{eqn:lehmer1}).    
In this article we examine only three patterns: 
(i) $\varphi(z) = C(z)$, $\psi(z) = S_i(z)$ and matrix $B$;         
(ii) $\varphi(z) = \rL(z) \cdot C(z)$, $\psi(z) = S_i(z)$ and matrix $A$;    
(iii) $\varphi(z) = C(z)$, $\psi(z) = \rL(z) \cdot \tilde{C}(z)$ and matrix $B$,   
where $C(z)$ and $\tilde{C}(z)$ stand for products of cyclotomic polynomials. 
%%%%%%%%%%%%%%%%%%%%%%%%% thm:main6 %%%%%%%%%%%%%%%%%%%%%%%%%
\begin{theorem} \label{thm:main6} 
In each of the three patterns there exists a complete enumeration 
of all cases that produce non-projective K3 surface automorphisms 
$f : X \to X$, where we have $h(f) = \log \lambda_i$ and $\rho(X) = 0$ 
in pattern (i), while $h(f) = \log \lambda_{\rL}$ and $\rho(X) = 12$ in 
patterns (ii) and (iii); see Theorems $\ref{thm:SH}$, $\ref{thm:A-pattern}$  
and $\ref{thm:B-pattern}$, respectively.    
\end{theorem}
%%%%%%%%%%%%%%%%%%%%%%%%%%%%%%%%%%%%%%%%%%%%%%%%%%%%%%%%%%%% 
\par
%%%%%%%% 
Even from these restricted setups we can obtain a lot of K3 surface 
automorphisms of positive entropy, especially those with Siegel disks.  
As for the results on Siegel disks we refer to 
Theorems \ref{thm:SH2}--\ref{thm:A2}.  
With various other Salem polynomials and in various other settings 
our method allows us to construct a much larger number of 
K3 surface automorphisms of positive entropy.  
If a K3 surface admits an automorphism with special eigenvalue 
conjugate to a Salem number then its Picard number  
is an even integer between $0$ and $18$, and all possible 
Picard numbers can be covered by our constructions.    
This result will be reported upon separately.    
%%%%%%%
\par
%%%%%%%
It is usually difficult to decide whether two different choices 
of the pair $(\varphi, \psi)$ yield the same automorphism, or whether the 
construction in this article gives examples identical to those in previous works 
such as Oguiso \cite{Oguiso}, McMullen \cite{McMullen3} and others,  
unless some of the resulting invariants, e.g. entropies, special eigenvalues, 
Picard numbers, Dynkin types of root systems, etc. are distinct, 
in which case the answer is clearly negative.    
This is because, given a unimodular quadratic form over $\bZ$ with 
known Gram matrix with respect to some basis, there is no known theory 
(at least to the authors) to convert it to a canonical form. 
However, there exists a systematic way to compare the method of 
hypergeometric groups with that of Salem number fields due to 
McMullen \cite{McMullen1}, when $\psi(z)$ is an unramified Salem polynomial 
of degree $22$; see Theorem \ref{thm:LNF} and Corollary \ref{cor:LNF}. 
This comparison should be extended to Salem numbers of lower degrees.   
%%%%%%%%%%%%
\par
%%%%%%%%%%%%
We wonder whether our construction yields 
interesting examples of finite order automorphisms. 
By Theorem \ref{thm:main4} only automorphisms of 
positive entropy, necessarily of infinite order, can occur in the 
non-projective case. 
So we have to examine the projective case to answer this question.  
This is another problem yet to be discussed. 
%%%%%%%%%%%%%
\par
%%%%%%%%%%%%%
The organization of this article is as follows. 
The first part (\S \ref{sec:hgg}--\S \ref{sec:hgl}) is devoted to enriching 
infrastructure in the theory of hypergeometric groups in anticipation of its 
application in the second part as well as other applications elsewhere. 
In \S \ref{sec:hgg}, for a complex hypergeometric group 
$H = H(\varphi, \psi)$ we study the Gram matrix of an $H$-invariant 
Hermitian form and prove a large part of Theorem \ref{thm:main1}.  
Its index is expressed in terms of the clusters 
of roots of $\varphi(z)$ and $\psi(z)$ (Theorem \ref{thm:p-q}).  
When $H$ is a real hypergeometric group, the index is represented in 
terms of  the trace clusters of roots of $\Phi(w)$ and $\Psi(w)$ 
(Theorem \ref{thm:real}).  
A formula for local indices (Proposition \ref{prop:lsgn-r}) is also included for 
later use in Hodge structures. 
When $H$ is defined over $\bZ$, an $H$-invariant even lattice, 
called the hypergeometric lattice, is introduced in \S \ref{sec:hgl}.  
The second part (\S\ref{sec:K3}--\S\ref{sec:lnf}) is an 
application of hypergeometric groups to dynamics on K3 surfaces, whose 
main results are already stated above.    
A classification of all real hypergeometric groups of rank $22$ and 
index $\pm16$ is given in Theorem \ref{thm:pm16}. 
Theorems \ref{thm:main2} and \ref{thm:main3} are then proved in 
\S\ref{ss:st}. 
Theorems \ref{thm:main4}, \ref{thm:main5} and \ref{thm:main6} are 
established in \S \ref{ss:PL}, \S \ref{ss:algorithm} and \S \ref{sec:k3auto}
respectively.                   
%%%%%%%%%
\par
%%%%%%%%% 
As the final remark, Theorems \ref{thm:main2} and \ref{thm:main3} say 
that at least two eigenvalues of $B$ must be algebraic units off  the 
unit circle in order for a hypergeometric group to yield a K3 lattice.  
This offers a sharp contrast to such a situation as in the classification 
of finite hypergeometric groups by Beukers and Heckman \cite{BH} or  
in their treatment of Lorentzian hypergeometric groups by 
Fuchs, Meiri and Sarnak \cite{FMS}, in which all eigenvalues of 
$A$ and $B$ are roots of unity. 
Even if they are off the unit circle, hypergeometric groups can enjoy rich structures.    
%%%%%%%%%%%%%%%%%%%%%%%%%%% sec:hgg %%%%%%%%%%%%%%%%%%%%%%%%%%%%%%%%%%%%
\section{Hypergeometric Groups} \label{sec:hgg}
%%%%%%%%%%%%%%%%%%%%%%%%%%%%%%%%%%%%%%%%%%%%%%%%%%%%%%%%%%%%%%%%%%%%%%
The theory of hypergeometric groups is developed by Beukers and Heckman 
\cite{BH}. 
A lucid explanation of this concept can also be found in Heckman's lecture 
notes \cite{Heckman}.      
A hypergeometric group is a group $H = \langle A, \, B\rangle$ generated 
by two invertible matrices $A, B \in \GL(n, \bC)$ such that 
$\rank(A - B) = 1$. 
Let $\ba = \{a_1, \dots, a_n\}$ and $\bb = \{b_1, \dots, b_n \}$ be the 
eigenvalues of $A$ and $B$ respectively (they are multi-sets allowing 
repeated elements).  
Then $H$ acts on $\bC^n$ irreducibly if and only if 
%%%%%%%%%%%%%%%%%%%%%%%%%%% eqn:irr %%%%%%%%%%%%%%%%%%%%%%%%%%%%%%%%%%%%%
\begin{equation} \label{eqn:irr}
\ba \cap \bb = \emptyset 
\end{equation}
%%%%%%%%%%%%%%%%%%%%%%%%%%%%%%%%%%%%%%%%%%%%%%%%%%%%%%%%%%%%%%%%%%%%%% 
(see \cite[Theorem 3.8]{Heckman}).   
Hereafter we always assume condition \eqref{eqn:irr}.  
Let $a^{\dagger} := \bar{a}^{-1}$ for $a \in \bC^{\times}$. 
There then exists a non-degenerate $H$-invariant Hermitian form on $\bC^n$ 
if and only if 
%%%%%%%%%%%%%%%%%%%%%%%%% eqn:ihf %%%%%%%%%%%%%%%%%%%%%%%%%%%%%%%%%%%%%
\begin{equation} \label{eqn:ihf}
\ba^{\dagger} = \ba, \qquad \bb^{\dagger} = \bb,  
\end{equation} 
%%%%%%%%%%%%%%%%%%%%%%%%%%%%%%%%%%%%%%%%%%%%%%%%%%%%%%%%%%%%%%%%%%%%
where $\ba^{\dagger} := \{ a_1^{\dagger}, \dots, a_n^{\dagger} \}$  
(see \cite[Theorem 4.3]{BH} and \cite[Theorem 3.13]{Heckman}. 
We remark that the ``only if'' part is not mentioned there but it is an easy exercise.     
Hereafter we also assume condition \eqref{eqn:ihf}. 
%%%%%%%%
\par
%%%%%%%% 
Let $\varphi(z)$ and $\psi(z)$ be the characteristic polynomials of 
$A$ and $B$ respectively.  
It is interesting to describe conditions \eqref{eqn:irr} and \eqref{eqn:ihf} 
in terms of $\varphi(z)$ and $\psi(z)$.  
Note that $\varphi(0) = (-1)^n \det A$ and $\psi(0) = (-1)^n \det B$ are 
non-zero since $A$, $B \in \GL(n, \bC)$.    
Condition \eqref{eqn:irr} is equivalent to 
%%%%%%%%%%%%%%%%%%%%%%%% eqn:irr2 %%%%%%%%%%%%%%%%%%%%%%%%%%%%%%%%%%%%%
\begin{equation} \label{eqn:irr2} 
\Res(\varphi, \psi) \neq 0,   \tag{$\ref{eqn:irr}'$}
\end{equation}
%%%%%%%%%%%%%%%%%%%%%%%%%%%%%%%%%%%%%%%%%%%%%%%%%%%%%%%%%%%%%%%%%%%%
where $\Res(\varphi, \psi)$ denotes the resultant of $\varphi(z)$ and $\psi(z)$, 
while condition \eqref{eqn:ihf} is equivalent to 
%%%%%%%%%%%%%%%%%%%%%%%% eqn:ihf2 %%%%%%%%%%%%%%%%%%%%%%%%%%%%%%%%%%%%%
\begin{equation} \label{eqn:ihf2} 
z^n \bar{\varphi}(z^{-1}) = \bar{\varphi}(0) \cdot \varphi(z), 
\qquad 
z^n \bar{\psi}(z^{-1}) = \bar{\psi}(0) \cdot \psi(z), 
\tag{$\ref{eqn:ihf}'$}
\end{equation}
%%%%%%%%%%%%%%%%%%%%%%%%%%%%%%%%%%%%%%%%%%%%%%%%%%%%%%%%%%%%%%%%%%%%
where $\bar{f}(z) := \overline{f(\bar{z})}$ for $f(z) \in \bC[z]$. 
Comparing the constant terms in \eqref{eqn:ihf2} we have  
%%%%%%%%%%%%%%%%%%%%%%%% eqn:m1 %%%%%%%%%%%%%%%%%%%%%%%%%%%%%%%%%%%%
\begin{equation} \label{eqn:m1}
|\varphi(0)| = 1, \qquad |\psi(0)| = 1,  
\end{equation}  
%%%%%%%%%%%%%%%%%%%%%%%%%%%%%%%%%%%%%%%%%%%%%%%%%%%%%%%%%%%%%%%%%%%
because $\varphi(z)$ and $\psi(z)$ are monic polynomials. 
In the other way round,  start with any monic polynomials $\varphi(z)$ and $\psi(z)$ 
with nonzero constant terms that satisfy condition \eqref{eqn:irr2}.    
Then their companion matrices $Z(\varphi)$ and $Z(\psi)$ generate an irreducible 
hypergeometric group 
%%%%%%%%%%%%%%%%%%%%%%%%%% eqn:std %%%%%%%%%%%%%%%%%%%%%%%%%%%%%%%%%%%%
\begin{equation} \label{eqn:std} 
H(\varphi, \psi) := \langle Z(\varphi), Z(\psi) \rangle. 
\end{equation}
%%%%%%%%%%%%%%%%%%%%%%%%%%%%%%%%%%%%%%%%%%%%%%%%%%%%%%%%%%%%%%%%%%%%%  
Levelt's theorem \cite{Levelt} states that any irreducible hypergeometric 
group $H = \langle A, B \rangle$ is conjugate in $\GL(n, \bC)$ to this one 
(see \cite[Theorem 3.5]{BH} and \cite[Theorem 3.9]{Heckman}). 
In this sense $H$ is uniquely determined by the characteristic polynomials 
$\varphi(z)$ and $\psi(z)$ of $A$ and $B$, or equivalently by their eigenvalue sets 
$\ba$ and $\bb$. 
Thus we can write a general irreducible hypergeometric group 
$H = \langle A, B \rangle$ as $H(\varphi, \psi)$ or $H(\ba, \bb)$.       
%%%%%%
\par
%%%%%%
The structure of the invariant Hermitian form is discussed in \cite[\S 4]{BH} 
and \cite[\S 3.3]{Heckman} when $\ba$ and $\bb$ lie on the unit circle $S^1$. 
We can extend the discussion there to the general case without this restriction. 
Since $C := A^{-1} B$ is a complex reflection, that is, $\rank(I-C) = 1$, its 
determinant   
%%%%%%%%%%%%%%%%%%%%%%%% eqn:c %%%%%%%%%%%%%%%%%%%%%%%%%%%%%%%%%%%%%%%
\begin{equation} \label{eqn:c}
c := \det C = \frac{\psi(0)}{\varphi(0)} 
= \frac{b_1 \cdots b_n}{a_1 \cdots a_n} \in S^1      
\end{equation}
%%%%%%%%%%%%%%%%%%%%%%%%%%%%%%%%%%%%%%%%%%%%%%%%%%%%%%%%%%%%%%%%%%%%  
is an eigenvalue of $C$, called the distinguished eigenvalue.  
In this article we assume that 
%%%%%%%%%%%%%%%%%%%%%%% eqn:cneq1 %%%%%%%%%%%%%%%%%%%%%%%%%%%%%%%%%%%
\begin{equation} \label{eqn:cneq1}
c \neq 1. 
\end{equation}
%%%%%%%%%%%%%%%%%%%%%%%%%%%%%%%%%%%%%%%%%%%%%%%%%%%%%%%%%%%%%%%%%%%%
Let $\br$ be an eigenvector of $C$ corresponding to the eigenvalue $c$. 
As in the proof of \cite[Theorem 3.14]{Heckman} we have  
$(\br, \br) \in \bR^{\times}$ and 
%%%%%%%%%%%%%%%%%%%%%%%%%% eqn:CBA %%%%%%%%%%%%%%%%%%%%%%%%%%%%%%%%
\begin{subequations} \label{eqn:CBA}
\begin{align} 
C \bv &= \bv - \zeta (\bv, \br) \br, \qquad \bv \in \bC^n, 
\label{eqn:C} \\[1mm]
\frac{\psi(z)}{\varphi(z)}  
&= 1 + \zeta \, ((z I - A)^{-1} A \br, \, \br),   
\quad \mbox{where} \quad \zeta := \frac{1-c}{(\br, \br)}.
\label{eqn:cBA}
\end{align}
\end{subequations}
%%%%%%%%%%%%%%%%%%%%%%%%%%%%%%%%%%%%%%%%%%%%%%%%%%%%%%%%%%%%%%%%%%%%
\par
%%%%%
The invariant Hermitian form is uniquely determined up to scalar multiplications 
in $\bR^{\times}$. 
To eliminate this ambiguity we take the normalization 
%%%%%%%%%%%%%%%%%%%%%%%%% eqn:nmlz %%%%%%%%%%%%%%%%%%%%%%%%%%%%%%%%%%
\begin{equation} \label{eqn:nmlz}
(\br, \br) = |1-c| > 0 \qquad \mbox{so that} \qquad 
\zeta = \frac{1-c}{|1-c|} \in S^1.  
\end{equation} 
%%%%%%%%%%%%%%%%%%%%%%%%%%%%%%%%%%%%%%%%%%%%%%%%%%%%%%%%%%%%%%%%%%%
Let $\{\xi_i\}_{i=1}^{\infty}$ be the sequence defined by the Taylor series 
expansion  
%%%%%%%%%%%%%%%%%%%%%%%% eqn:taylor %%%%%%%%%%%%%%%%%%%%%%%%%%%%%%%%% 
\begin{equation} \label{eqn:taylor}
\dfrac{\psi(z)}{\varphi(z)} = 1 + \zeta \sum_{i=1}^{\infty} 
\xi_i \, z^{-i} \qquad \mbox{around} \quad z = \infty. 
\end{equation}
%%%%%%%%%%%%%%%%%%%%%%%%%%%%%%%%%%%%%%%%%%%%%%%%%%%%%%%%%%%%%%%%%%% 
%%%%%%%%%%%%%%%%%%%%%%%%% thm:ihf %%%%%%%%%%%%%%%%%%%%%%%%%%%%%%%%%%%%
\begin{theorem} \label{thm:ihf} 
The invariant Hermitian pairing $g_{i j} := (A^{i-1} \br, \, A^{j-1} \br)$ is given by 
%%%%%%%%%%%%%%%%%%%%%%%%% eqn:gram %%%%%%%%%%%%%%%%%%%%%%%%%%%%%%%%%
\begin{equation} \label{eqn:gram}
g_{i j} = 
\begin{cases}
\,\, \xi_{i-j} \quad & (i \ge j \ge 1), \\[1mm]  
\,\, \bar{\xi}_{j-i} \quad & (1 \le i < j),    
\end{cases} 
\end{equation}
%%%%%%%%%%%%%%%%%%%%%%%%%%%%%%%%%%%%%%%%%%%%%%%%%%%%%%%%%%%%%%%%%%
with convention $\xi_0 := |1-c|$.   
The determinant of $G := (g_{i j})_{i,j=1}^n$ has absolute value   
%%%%%%%%%%%%%%%%%%%%%%%% eqn:disc %%%%%%%%%%%%%%%%%%%%%%%%%%%
\begin{equation} \label{eqn:disc} 
|\det G| = | \Res(\varphi, \psi) |,     
\end{equation}
%%%%%%%%%%%%%%%%%%%%%%%%%%%%%%%%%%%%%%%%%%%%%%%%%%%%%%%%%%% 
which is non-zero by assumption \eqref{eqn:irr2}. 
In particular $\br, A \br, \dots, A^{n-1} \br$ form a basis of $\bC^n$. 
\end{theorem} 
%%%%%%%%%%%%%%%%%%%%%%%% begin proof %%%%%%%%%%%%%%%%%%%%%%%%%
{\it Proof}. 
First we show formula \eqref{eqn:gram}. 
Taylor expansion of \eqref{eqn:cBA} around $z = \infty$ reads   
%%%%%%%%%%%%%%%
$$
\frac{\psi(z)}{\varphi(z)}  
= 1 + \zeta \, ((I - z^{-1}A)^{-1} z^{-1} A \br, \, \br) 
= 1 + \zeta \sum_{i=1}^{\infty} (A^i \br, \, \br)  z^{-i}.  
$$
%%%%%%%%%%%%%%% 
Comparing this with expansion \eqref{eqn:taylor} together with 
convention $\xi_0 = |1-c|$ yields $(A^i \br, \, \br) = \xi_i$ for 
every $i \in \bZ_{\ge 0}$. 
Since the Hermitian form is $A$-invariant, we have    
%%%%%
$$
g_{i j} = (A^{i-1} \br, \, A^{j-1} \br) =
\begin{cases}
(A^{i-j} \br, \, \br) = \xi_{i-j}\quad & (i \ge j \ge 1), \\[1mm]
\overline{(A^{j-i} \br, \, \br)} = \bar{\xi}_{j-i} \quad & (1 \le i < j).  
\end{cases}
$$
%%%%%
This together with normalization \eqref{eqn:nmlz} leads to formula 
\eqref{eqn:gram}. 
%%%%%
\par
%%%%%
Next we show formula \eqref{eqn:disc}. 
Suppose for the time being that $a_1, \dots, a_n$ are mutually distinct.    
Let $\br = \br_1 + \cdots + \br_n$ be the decomposition of $\br$ into 
eigenvectors of $A$, where $\br_i$ correspond to the eigenvalue $a_i$. 
By condition $\ba^{\dagger} = \ba$ in \eqref{eqn:ihf} there exists  
a permutation $\sigma \in S_n$ such that $\sigma^2 = 1$ and 
$a_{\sigma(i)} = a_i^{\dagger}$ for $i = 1, \dots, n$. 
Note that $\sigma(i) = i$ if and only if $a_i \in S^1$. 
It follows from non-degeneracy and $A$-invariance of the Hermitian 
form that $(\br_i, \br_{\sigma(i)})$, $i = 1, \dots, n$, are non-zero 
while all the other Hermitian parings $(\br_i, \br_j)$ vanish. 
Thus equation \eqref{eqn:cBA} leads to 
%%%%%%
$$
\frac{\psi(z)}{\varphi(z)}           
= 1 + \zeta \sum_{i=1}^n \frac{a_i \, (\br_i, \, \br_{\sigma(i)})}{z-a_i}.        
$$
%%%%% 
Taking residue at $z = a_i$ we have for $i = 1, \dots, n$,  
%%%%%%%%%%%%%%%%%%%%% eqn:lam %%%%%%%%%%%%%%%%%%%%%%%%%%%%%%
\begin{equation} \label{eqn:lam}
\lambda_i := (\br_i, \, \br_{\sigma(i)}) = 
\frac{\psi(a_i)}{\zeta \, a_i \, \varphi_i(a_i)} \qquad \mbox{with} \quad 
\varphi_i(z) := \prod_{j \neq i} (z-a_j). 
\end{equation}
%%%%%%%%%%%%%%%%%%%%%%%%%%%%%%%%%%%%%%%%%%%%%%%%%%%%%%%%%%%
\par
%%%%%%
After rearranging $a_1, \dots, a_n$ if necessary, we may assume that 
$\sigma$ fixes $1, \dots, l$ and exchanges $l+2 i-1$ and $l+2 i$ 
for $i = 1, \dots, m$, where $n = l + 2 m$. 
Then $\br_1, \dots, \br_n$ form a basis of $\bC^n$ with respect to 
which the Gram matrix of the invariant Hermitian form is given by 
%%%%%%%%%%%%%%%%%%%%%% eqn:Lam %%%%%%%%%%%%%%%%%%%%%%%%%%%%%
\begin{equation} \label{eqn:Lam}
\varLambda = (\lambda_1) \oplus \cdots \oplus (\lambda_l) \oplus 
\varLambda_1 \oplus \cdots \oplus \varLambda_m 
\quad \mbox{with} \quad
\varLambda_i := 
\begin{pmatrix}
0 & \lambda_{l+ 2 i-1} \\[1mm]
\lambda_{l + 2 i} & 0 
\end{pmatrix}.  
\end{equation}
%%%%%%%%%%%%%%%%%%%%%%%%%%%%%%%%%%%%%%%%%%%%%%%%%%%%%%%%%%%
Note that $\lambda_1, \dots, \lambda_l \in \bR$ and 
$\lambda_{l+2 i-1} = \bar{\lambda}_{l + 2 i}$ for $i = 1, \dots, m$. 
From \eqref{eqn:lam} and \eqref{eqn:Lam} we have  
%%%%%%%%
$$
|\det \varLambda | = | (-1)^m \lambda_1 \cdots \lambda_n | 
= \left| \frac{ (-1)^m \psi(a_1) \cdots \psi(a_n) }{\zeta^n \,  
a_1 \cdots a_n \, \varphi_1(a_1) \cdots \varphi_n(a_n) } \right| 
= \frac{|\Res(\varphi, \psi)|}{\prod_{i < j} |a_i - a_j|^2}, 
$$
%%%%%%%
where $|\zeta| = 1 = |a_1 \cdots a_n|$ is used.  
Moreover we have 
$(\br, A \br, \dots, A^{n-1} \br) = (\br_1, \dots, \br_n) V$ 
with $V := (a_i^{j-1})_{i, j=1}^n$ being a Vandermonde matrix. 
This implies $G = {}^{\rt} V \varLambda \overline{V}$ and hence 
%%%%%%
$$
|\det G| = |\det \varLambda| |\det V|^2 =  |\det \varLambda| 
\prod_{i < j} |a_i-a_j|^2 = |\Res(\varphi, \psi)|, 
$$
%%%%%%
which proves formula \eqref{eqn:disc} when $a_1, \dots, a_n$ are 
distinct. 
The formula in the general case follows by a continuity argument, 
which works as far as conditions \eqref{eqn:irr2} and \eqref{eqn:ihf2} 
are fulfilled. \hfill $\Box$ \par\medskip
%%%%%%%%%%%%%%%%%%%%%%% end proof %%%%%%%%%%%%%%%%%%%%%%%%%%%
%%%%%%%%%%%%%%%%%%%% proof of thm:main1 %%%%%%%%%%%%%%%%%%%%%%
{\it Proof of Theorem $\ref{thm:main1}$}. 
A large part of the theorem is a special case of Theorem \ref{thm:ihf} 
where $\varphi(z)$ and $\psi(z)$ are polynomials over $\bZ$, 
$c = -1$ in \eqref{eqn:cneq1} so that  $\zeta = 1$ in 
\eqref{eqn:nmlz} and \eqref{eqn:taylor}. 
For the rest, see the beginning of \S \ref{sec:hgl}. 
\hfill $\Box$ 
%%%%%%%%%%%%%%%%%%% end proof %%%%%%%%%%%%%%%%%%%%%%%%%%%%%%%  
%%%%%%%%%%%%%%%%%%%%%%% rem:antip %%%%%%%%%%%%%%%%%%%%%%%%%%%
\begin{remark} \label{rem:antip}
It is obvious that if $H = \langle A, B \rangle \subset \GL(n, \bC)$ is 
a hypergeometric group then so is 
$H^{\ra} := \langle -A, -B \rangle$. 
We refer to $H^{\ra}$ as the {\sl antipode} of $H$. 
Note that $\varphi^{\ra}(z) = (-1)^n \varphi(-z)$ and 
$\psi^{\ra}(z) = (-1)^n \psi(-z)$, hence $\ba^{\ra} = - \ba$ and 
$\bb^{\ra} = - \bb$. 
It follows from $C^{\ra} = C$ and formula \eqref{eqn:gram} that 
$H$ and $H^{\ra}$ have the same invariant Hermitian form. 
This simple remark will be useful in discussing Hodge structures 
(see Remark \ref{rem:simple}).  
\end{remark} 
%%%%%%%%%%%%%%%%%%%%%%%%%%%%%%%%%%%%%%%%%%%%%%%%%%%%%%%%%%%
%%%%%%%%%%%%%%%%%%%%%%%%% sec:sgn %%%%%%%%%%%%%%%%%%%%%%%%%%%
\section{Index of the Invariant Hermitian Form} \label{sec:sgn}
%%%%%%%%%%%%%%%%%%%%%%%%%%%%%%%%%%%%%%%%%%%%%%%%%%%%%%%%%%%%
Let $\ba_{\mathrm{on}}$ and $\ba_{\mathrm{off}}$ be those components   
of $\ba$ whose elements lie on and off $S^1$ respectively.  
We define $\bb_{\mathrm{on}}$ and $\bb_{\mathrm{off}}$ in a similar 
manner for $\bb$.  
If both of $\baon$ and $\bbon$ are nonempty then they dissect 
each other into an equal number of components 
$\ba_1, \dots, \ba_t$ and $\bb_1, \dots, \bb_t$ so that 
$\ba_1, \bb_1$, $\ba_2, \bb_2$, $\dots$, $\ba_t, \bb_t$ are located 
consecutively on $S^1$ in the positive direction (anti-clockwise) as 
shown in Figure \ref{fig:clus}. 
%%%%%%%%%%%%%%%%%%%%%%% fig:clus %%%%%%%%%%%%%%%%%%%%%%%%%%%%%%
\begin{figure}[tt]
\centerline{
\input{cluster.tex}}
\caption{Clusters of $\ba_{\mathrm{on}}$ and $\bb_{\mathrm{on}}$ 
when $t = 5$.}  
\label{fig:clus}
\end{figure}
%%%%%%%%%%%%%%%%%%%%%%%%%%%%%%%%%%%%%%%%%%%%%%%%%%%%%%%%%%%%
Each $\ba_i$ is called a {\sl cluster} of $\ba_{\mathrm{on}}$ and 
is said to be {\sl simple, double, triple}, etc. if $|\ba_i| = 1$, $2$, $3$, 
and so on, where $|\bx|$ denotes the cardinality counted with 
multiplicities of a multi-set $\bx$.     
We write  
%%%%%%%
$$
[\baon] = 1^{\nu_1} 2^{\nu_2} 3^{\nu_3} \cdots   
$$ 
%%%%%%
if $\ba_{\mathrm{on}}$ consists of $\nu_1$ simple clusters, 
$\nu_2$ double clusters, $\nu_3$ triple clusters, etc.   
Note that $|\baon| = \nu_1 + 2 \nu_2 + 3 \nu_3 +\cdots$,  
$|\baoff|$ is even and $|\baon|+|\baoff| = n$;   
the same is true for $\bb$. 
Taking a branch-cut $\ell$ separating $\bb_t$ and $\ba_1$ as in 
Figure \ref{fig:clus} we define the argument of $z \in \bC^{\times}$ 
so that 
%%%%%%%%%%%%%%%%%%%%%%%% eqn:arg %%%%%%%%%%%%%%%%%%%%%%%%%%%%
\begin{equation} \label{eqn:arg}
\Theta \le \arg z < \Theta + 2 \pi, 
\end{equation}
%%%%%%%%%%%%%%%%%%%%%%%%%%%%%%%%%%%%%%%%%%%%%%%%%%%%%%%%%%%
where $\Theta \in [-\pi, \, \pi)$ is the angle of the ray $\ell$ 
to the positive real axis.   
%%%%%%
\par
%%%%%% 
Let $\arg a_i = 2 \pi \alpha_i$ and $\arg b_i = 2 \pi \beta_i$ for 
$i = 1, \dots, n$.  
Formula \eqref{eqn:c} allows us to write     
%%%%%%%%%%%%%%%%%%%%%%%%% eqn:cneq2 %%%%%%%%%%%%%%%%%%%%%%%%
\begin{equation} \label{eqn:cneq2}
c = \re^{2 \pi \ri \, \gamma} \in S^1 \qquad \mbox{with} \quad 
\gamma := \sum_{i=1}^n \beta_i - \sum_{i=1}^n \alpha_i \in \bR,  
\end{equation}
%%%%%%%%%%%%%%%%%%%%%%%%%%%%%%%%%%%%%%%%%%%%%%%%%%%%%%%%%%% 
where $\ri := \sqrt{-1}$, hence condition \eqref{eqn:cneq1} is equivalent 
to $\gamma \in \bR \setminus \bZ$, that is, $\sin \pi \gamma \in \bR^{\times}$.  
%%%%%%%%%%%%%%%%%%%%%%%% rem:branch %%%%%%%%%%%%%%%%%%%%%%%%
\begin{remark} \label{rem:branch} 
Taking another branch of $\arg$ has no effect on the contribution 
of $\baoff$ and $\bboff$ to the value of $\sin \pi \gamma$, because 
any pair $\lambda$, $\lambda^{\dagger} \in \baoff$ has a common 
argument so the sum $\arg \lambda + \arg \lambda^{\dagger}$ alters 
only by an even multiple of $2 \pi$; the same is true for  
$\lambda$, $\lambda^{\dagger} \in \bboff$.      
\end{remark}
%%%%%%%%%%%%%%%%%%%%%%%%%%%%%%%%%%%%%%%%%%%%%%%%%%%%%%%%%%% 
%%%%%%%%%%%%%%%%%%%%%%%% ss:sgn2 %%%%%%%%%%%%%%%%%%%%%%%%%%%%
\subsection{Clusters and Index} \label{ss:sgn2}
%%%%%%%%%%%%%%%%%%%%%%%%%%%%%%%%%%%%%%%%%%%%%%%%%%%%%%%%%%%% 
Let $(p, q)$ be the signature of the $H(\ba, \bb)$-invariant Hermitian 
form on $\bC^n$. 
Under the condition 
%%%%%%%%%%%%%%%%%%%%%%%%% eqn:nooff %%%%%%%%%%%%%%%%%%%%%%%%%%%%%
\begin{equation} \label{eqn:nooff}
\ba_{\mathrm{off}} = \bb_{\mathrm{off}} = \emptyset,   
\end{equation}
%%%%%%%%%%%%%%%%%%%%%%%%%%%%%%%%%%%%%%%%%%%%%%%%%%%%%%%%%%%%%%%  
Beukers and Heckman \cite[Theorem 4.5]{BH} gave  
a formula for the index $p-q$ (up to sign).   
We can state a refined version of it in terms of the clusters of $\baon$ 
and $\bbon$ without assuming \eqref{eqn:nooff}.        
%%%%%%%%%%%%%%%%%%%%%%%%%% thm:p-q %%%%%%%%%%%%%%%%%%%%%%%%%%%%%%%
\begin{theorem} \label{thm:p-q} 
If $\baon$ and $\bbon$ are nonempty then the invariant 
Hermitian form has index 
%%%%%%%%%%%%%%%%%%%%%%%%% eqn:p-q %%%%%%%%%%%%%%%%%%%%%%%%%%%%%%%%
\begin{equation} \label{eqn:p-q} 
p-q = \varepsilon \sum_{k \in \rK} (-1)^{\tau_k} \qquad 
\mbox{with} \quad 
\rK := \{\, k =1, \dots, t \, : \, |\ba_k| \equiv 1 \bmod 2 \,\},  
\end{equation} 
%%%%%%%%%%%%%%%%%%%%%%%%%%%%%%%%%%%%%%%%%%%%%%%%%%%%%%%%%%%%%%%%%
where $\varepsilon = \pm 1$ is the sign of $\sin \pi \gamma \in \bR^{\times}$ 
with $\gamma$ given in \eqref{eqn:cneq2} and $\tau_k$ is defined by 
%%%%%%
$$ 
\tau_1 := 0; \qquad 
\tau_k :=  |\ba_1| +|\bb_1| + \cdots + |\ba_{k-1}| +|\bb_{k-1}|, \quad 
k = 2, \dots, t. 
$$
%%%%%  
If at least one of $\baon$ and $\bbon$ is empty then the index $p-q$ is zero. 
\end{theorem}  
%%%%%%%%%%%%%%%%%%%%%%%%% begin proof %%%%%%%%%%%%%%%%%%%%%%%%%%%%%
{\it Proof}. 
We may assume that $a_1, \dots, a_n$ are mutually distinct, 
since the general case can be treated by a perturbation argument (see e.g. 
Kato \cite[Chapter I\!I, \S1.4]{Kato}).  
After rearranging the indices of $a_i$ and $b_j$ if necessary, 
we may further assume that  
%%%%%
\begin{alignat*}{4}
\ba_{\mathrm{on}} &= \{ a_1, \dots, a_l \}, \quad & 
\ba_{\mathrm{off}} &= \{ a_{l+1}, \dots, a_n \}, \quad & 
a_{l+2 i-1} &= a_{l+2 i}^{\dagger}, \quad & i &= 1, \dots, m, 
\\[1mm]  
\bb_{\mathrm{on}} &= \{ b_1, \dots, b_d \}, \quad &  
\bb_{\mathrm{off}} &= \{ b_{d+1}, \dots, b_n \}, \quad & 
b_{d+2 i-1} &= b_{d+2 i}^{\dagger}, \quad & i &= 1, \dots, e,   
\end{alignat*}
%%%%%
where $n = l + 2 m = d + 2 e$. 
Suppose that both of $\baon$ and $\bbon$ are nonempty, that is, 
$l \ge 1$ and $d \ge 1$.    
Since the Hermitian matrix $\varLambda_i$ in \eqref{eqn:Lam} has 
null index, we have    
%%%%%%%%%%%%%%%%%%%%%%%%% eqn:spm %%%%%%%%%%%%%%%%%%%%%%%%%%%%%%%%%%%
\begin{equation} \label{eqn:spm}
p-q = l_+ - l_-, \qquad
l_{\pm} := \# \{\, i = 1, \dots, l \, :\, \pm \lambda_i > 0 \, \}. 
\end{equation} 
%%%%%%%%%%%%%%%%%%%%%%%%%%%%%%%%%%%%%%%%%%%%%%%%%%%%%%%%%%%%%%        
\par
%%%%%%  
For $x$, $y \in \bR^{\times}$ we write $x \sim y$ if $x$ and $y$ have 
the same sign.   
We claim that 
%%%%%%%%%%%%%%%%%%%%%%%% eqn:sim %%%%%%%%%%%%%%%%%%%%%%%%%%%%%
\begin{equation} \label{eqn:sim}
\lambda_i \sim \sigma_i := \varepsilon \cdot 
\dfrac{\prod_{j=1}^d \sin \pi(\beta_j -\alpha_i)}{
\prod_{j=1}^{* \, l} \sin \pi(\alpha_j - \alpha_i)} \in 
\bR^{\times}, \qquad i = 1, \dots, l,    
\end{equation}
%%%%%%%%%%%%%%%%%%%%%%%%%%%%%%%%%%%%%%%%%%%%%%%%%%%%%%%%%%%%%
where $\prod_j^{*}$ is the product avoiding $j = i$. 
Indeed, equation \eqref{eqn:lam} together with \eqref{eqn:nmlz} yields 
%%%%%%%%%%%%
$$ 
\lambda_i = 
\dfrac{|1-c| \prod_{j=1}^d(a_i-b_j) \prod_{j=1}^e 
(a_i-b_{d+2 j}^{\dagger})(a_i-b_{d+2 j}) }{(1-c) \, a_i \prod_{j=1}^{* \, l} 
(a_i-a_j) \prod_{j=1}^m (a_i-a_{l+ 2 j}^{\dagger})(a_i - a_{l+2 j})}   
$$
%%%%%%%%%%%
for $i = 1, \dots, l$. 
To evaluate the right-hand side we use the following identities:  
%%%%%%%%%%%
\begin{align*}
n &= l + 2 m = d + 2 e, \\[1mm]  
1-c &= - c^{\frac{1}{2}} (c^{\frac{1}{2}} - c^{- \frac{1}{2}}) 
= - 2 \ri \cdot c^{\frac{1}{2}} \cdot \sin \pi \gamma, \\[1mm]
u - v &= - u^{\frac{1}{2}} v^{\frac{1}{2}} 
( u^{-\frac{1}{2}} v^{\frac{1}{2}} -  u^{\frac{1}{2}} v^{-\frac{1}{2}}) 
= - 2 \ri \cdot u^{\frac{1}{2}} v^{\frac{1}{2}} \cdot 
\sin \pi(\phi-\theta), \\[1mm]
(u -w^{\dagger})(u - w) &= - u w \cdot |u - w^{\dagger}|^2 
= -u (w^{\dagger} w)^{\frac{1}{2}} \cdot |w| |u - w^{\dagger}|^2,    
\end{align*}
%%%%%%%%%%%
for $u = \re^{2 \pi \ri \, \theta} \in S^1$, 
$v = \re^{2 \pi \ri \, \phi} \in S^1$ and $w \in \bC^{\times}$.  
Some calculations yield 
%%%%%%
$$
\lambda_i = \mu_i \cdot \sigma_i, 
\qquad 
\mu_i := 2^{d-l+1} \cdot 
\dfrac{\prod_{j=1}^e |b_{d+2 j}||a_i - b_{d+2 j}^{\dagger}|^2 }{
\prod_{j=1}^m |a_{l+2 j}||a_i-a_{l+ 2 j}^{\dagger}|^2} > 0  
$$
%%%%% 
for $i = 1, \dots, l$ and hence claim \eqref{eqn:sim} is proved. 
%%%%%
\par
%%%%%
Relation \eqref{eqn:sim} readily shows that  the sign  
$\varepsilon_i = \pm 1$ of $\lambda_i$ is determined by  
%%%%%%%%%%%%%%%%%%%%%%%%%%% eqn:sim2 %%%%%%%%%%%%%%%%%%%%%%%%%%%%%
\begin{equation} \label{eqn:sim2} 
\varepsilon_i = \varepsilon \cdot (-1)^{\delta_i}, \quad 
\delta_i := \# \{ j = 1, \dots, n \,:\, \alpha_j < \alpha_i \} + 
\# \{ j = 1, \dots, n \,:\, \beta_j < \alpha_i \},  
\end{equation}
%%%%%%%%%%%%%%%%%%%%%%%%%%%%%%%%%%%%%%%%%%%%%%%%%%%%%%%%%%%%%%%% 
where $\arg a_j = 2 \pi \alpha_j$ and $\arg b_j = 2 \pi \beta_j$. 
If $a_i$ is the $d_i$-th smallest element of $\ba_k$ in the argument  
then $\varepsilon = \varepsilon \cdot (-1)^{\tau_k + d_i -1}$. 
Since $d_i$ ranges over $1, \dots, |\ba_k|$ as $a_i$ runs through 
$\ba_k$, one has   
%%%%%%%
$$
\sum_{a_i \in \sba_k} (-1)^{d_i-1} = 
\begin{cases} 
1 \quad & (|\ba_k| \equiv 1 \bmod 2), \\[1mm]
0 \quad & (|\ba_k| \equiv 0 \bmod 2). 
\end{cases} 
$$
%%%%%%%     
Thus it follows from formula \eqref{eqn:spm} that  
%%%%%%%
$$
p-q = l_+ - l_- = \sum_{i=1}^l \varepsilon_i = \varepsilon 
\sum_{k=1}^t (-1)^{\tau_k} \sum_{a_i \in \sba_k} (-1)^{d_i-1}
= \varepsilon \sum_{k \in \rK} (-1)^{\tau_k},   
$$
%%%%%%
which establishes formula \eqref{eqn:p-q}.  
When $\baon$ is empty, the index is zero as we have $l = 0$ in \eqref{eqn:Lam}. 
When $\bbon$ is empty, replace $\ba$ with $\bb$ and proceed in 
a similar manner.  \hfill $\Box$ 
%%%%%%%%%%%%%%%%%%%%%%%%%% end proof %%%%%%%%%%%%%%%%%%%%%%%%%%%%%
%%%%%%%%%%%%%%%%%%%%%%%%%% rem:p-q %%%%%%%%%%%%%%%%%%%%%%%%%%%%%%
\begin{remark} \label{rem:p-q} 
The following remarks are helpful in applying Theorem \ref{thm:p-q}. 
\begin{enumerate}
\setlength{\itemsep}{-1pt}
\item Formula \eqref{eqn:p-q} is invariant under any cyclic permutation   
of the indices $k$ for $\ba_k$ and $\bb_k$. 
\item Note that $|p-q| \le |\rK| \le t \le n$ and $|p-q| \equiv |\rK| \equiv n \bmod 2$; 
moreover $|\rK| = t$ if and only if all $\ba_{\mathrm{on}}$-clusters 
$\ba_1, \dots, \ba_t$ have odd cardinalities.  
\item We may exchange the roles of $\ba$ and $\bb$ in Theorem \ref{thm:p-q}.  
\end{enumerate}   
\end{remark}
%%%%%%%%%%%%%%%%%%%%%%%%%%%%%%%%%%%%%%%%%%%%%%%%%%%%%%%%%%%%%%%%
%%%%%%%%%%%%%%%%%%%%%%%%% sec:L %%%%%%%%%%%%%%%%%%%%%%%%%%%%%%%%
\subsection{Lorentzian Hypergeometric Groups} \label{ss:L}
%%%%%%%%%%%%%%%%%%%%%%%%%%%%%%%%%%%%%%%%%%%%%%%%%%%%%%%%%%%%%%%
It is interesting to ask when the invariant Hermitian form is definite or 
Lorentzian. 
Under condition \eqref{eqn:nooff} Beukers and Heckman \cite[Corollary 4.7]{BH} 
obtained the interlacing criterion for definiteness, while Fuchs, Meiri and 
Sarnak \cite[\S 2.2]{FMS} derived the so-called almost interlacing criterion 
for the Lorentzian case.     
Even without assuming \eqref{eqn:nooff} a priori, Theorem \ref{thm:p-q} 
readily implies that the Hermitian form is definite if and only if 
$[\baon] = [\bbon] = 1^n$ and $\baoff = \bboff = \emptyset$.  
For the Lorentzian case we have the following classification, where    
types $1$ and $2$ appear in \cite{FMS}.        
%%%%%%%%%%%%%%%%%%%%%%%%%% thm:L %%%%%%%%%%%%%%%%%%%%%%%%%%%%%%%%
\begin{theorem} \label{thm:L}
The Lorentzian case is classified into five types in Table $\ref{tab:L}$.  
In type $1$ we mean by ``doubles adjacent'' that the unique double cluster 
in $\baon$ and the unique double cluster in $\bbon$ must be adjacent to each other. 
For the other types there are no constraints on the location of multiple clusters.    
%%%%%%%%%%%%%%%%%%%%%%%%%% tab:L %%%%%%%%%%%%%%%%%%%%%%%%%%%%
\begin{table}[hh]
\centerline{
\begin{tabular}{cccccc}
\hline
\\[-4mm]
type & $[\baon]$ & $[\bbon]$ & constraint & $|\baoff|$ & $|\bboff|$ \\[1mm]
\hline
\\[-4mm]
$1$ &  $1^{n-2} 2^1$ & $1^{n-2} 2^1$ & doubles adjacent & $0$ & $0$ \\[1mm]     
$2$ & $1^{n-3} 3^1$  & $1^{n-3} 3^1$ &                          & $0$ & $0$ \\[1mm]   
$3$ & $1^{n-3} 3^1$  & $1^{n-2}$      &                           & $0$ & $2$ \\[1mm]     
$4$ & $1^{n-2}$ & $1^{n-3} 3^1$       &                          & $2$ & $0$ \\[1mm]  
$5$ & $1^{n-2}$ & $1^{n-2}$             &                           & $2$ & $2$ \\[1mm]
\hline 
\end{tabular}} 
\caption{Lorentzian case.}
\label{tab:L}
\end{table} 
\end{theorem}
%%%%%%%%%%%%%%%%%%%%%%%%%% begin proof %%%%%%%%%%%%%%%%%%%%%%%%%%%%
{\it Proof}. 
It follows from item (2) of Remark \ref{rem:p-q} that 
$|p-q| = n-2 \le |\rK| \le t \le n$ and $|\rK| \equiv n \bmod 2$. 
So we have either $|\rK| = t = n$ or $|\rK| = n-2$, 
but $t = n$ is ruled out as it would lead to the definite case. 
Thus $|\rK| = n-2$ and $t = n-2$, $n-1$.   
%%%%%%
\par
%%%%%%
First we consider the case $t = n-1$. 
If follows from $n \ge |\baon| \ge t = n-1$;     
$n \equiv |\baon| \bmod 2$ and $|\rK| = n-2$ that 
$[\baon] = 1^{n-2}2^1$ and $|\baoff| = 0$. 
One has also $[\bbon] = 1^{n-2}2^1$ and $|\bboff| = 0$ 
by (3) of Remark \ref{rem:p-q}. 
Let $k$, $l \in \{1, \dots, n-1\}$ be the indices such that $|\ba_k| = 2$ 
and $|\bb_l| = 2$.    
By (1) of Remark \ref{rem:p-q} we may assume $k = n-1$ and 
hence $\rK = \{1, \dots, n-2\}$.  
Index formula \eqref{eqn:p-q} now reads $|\sum_{k=1}^{n-2} (-1)^{\tau_k}| = n-2$, 
which implies $\tau_k \equiv \tau_1 = 0 \bmod 2$ for $k = 2, \dots, n-2$ 
and hence $|\ba_k| + |\bb_k| = 1 + |\bb_k| \equiv 0 \bmod 2$ for 
$k = 1, \dots, n-3$, which in turn forces $l = n-2$ or $l = n-1$, that is, 
$\ba_k$ and $\bb_l$ must be adjacent.    
This case falls into type $1$ of Table \ref{tab:L}.  
%%%%%%
\par
%%%%%%
Next we proceed to the case $t = n-2$. 
Since $|\rK| = n-2 = t$, all of $|\ba_1|, \dots, |\ba_{n-2}|$ must be 
odd by item (2) of Remark \ref{rem:p-q}. 
It then follows from $n \ge |\baon| \ge n-2$ and 
$n \equiv |\ba_{\mathrm{on}}| \bmod 2$ that $\ba$ must satisfy 
either (A1) $[\baon] = 1^{n-3}3^1$, $|\baoff| = 0$; or    
(A2) $[\baon] = 1^{n-2}$, $|\baoff| = 2$. 
By item (3) of Remark \ref{rem:p-q}, $\bb$ must also satisfy either 
(B1) $[\bbon] = 1^{n-3}3^1$, $|\bboff| = 0$; or      
(B2) $[\bbon] = 1^{n-2}$, $|\bboff| = 2$.  
Then the combinations (A1)-(B1), (A1)-(B3), (A3)-(B1),   
(A3)-(B3) lead to types $2$, $3$, $4$, $5$ in Table \ref{tab:L}, 
respectively. 
The converse implication is easy to verify. \hfill $\Box$ 
%%%%%%%%%%%%%%%%%%%%%%%%%% end proof %%%%%%%%%%%%%%%%%%%%%%%%%%%%%
%%%%%%%%%%%%%%%%%%%%%%%%%% ss:lcs %%%%%%%%%%%%%%%%%%%%%%%%%%%%%%%%
\subsection{Local Index} \label{ss:lcs}
%%%%%%%%%%%%%%%%%%%%%%%%%%%%%%%%%%%%%%%%%%%%%%%%%%%%%%%%%%%%%%%%
Let $E(\nu)$ be the generalized eigenspace of $A$ corresponding 
to an eigenvalue $\nu \in \ba$.   
Note that $m(\nu) := \dim E(\nu)$ is the multiplicity of $\nu$ in the 
multi-set $\ba$.  
Put $E(\mu, \mu^{\dagger}) := E(\mu) \oplus E(\mu^{\dagger})$ for 
$\mu \in \baoff$. 
Some linear algebra shows that the non-degeneracy and 
$A$-invariance of the Hermitian form lead to an orthogonal direct 
sum decomposition 
%%%%%%%%%%%%%%%%%%%%%%%%%% eqn:esd %%%%%%%%%%%%%%%%%%%%%%%%%%%%%%%%%
\begin{equation} \label{eqn:esd}
\bC^n = \bigoplus_{|\lambda|=1} E(\lambda) \, \oplus \, 
\bigoplus_{|\mu| > 1} E(\mu, \mu^{\dagger}),    
\end{equation}
%%%%%%%%%%%%%%%%%%%%%%%%%%%%%%%%%%%%%%%%%%%%%%%%%%%%%%%%%%%%%%%%%%%
where $\lambda$ ranges over all distinct elements in $\baon$ and $\mu$ 
ranges over all distinct elements in $\baoff$ such that $|\mu| > 1$. 
The Hermitian form is non-degenerate on $E(\lambda)$ and  
$E(\mu, \mu^{\dagger})$, though it is null on $E(\mu)$ and  
$E(\mu^{\dagger})$ individually.  
Note that $m(\mu) = m(\mu^{\dagger})$.  
It is interesting to find the index of the Hermitian form restricted 
to $E(\lambda)$ or $E(\mu, \mu^{\dagger})$.  
The same problem also makes sense with $B$ and $\bb$ 
in place of $A$ and $\ba$, where if $\nu \in \bb$ then $E(\nu)$ 
is understood to be the generalized $\nu$-eigenspace of $B$. 
Thanks to \eqref{eqn:irr} or \eqref{eqn:irr2} the common notation $m(\nu)$ 
is allowed for $\nu \in \ba \cup \bb$, because $m(\nu)$ is the same as 
the multiplicity of $\nu$ in $\varphi(z) \cdot \psi(z)$.    
%%%%%% 
\par
%%%%%%
For two distinct elements $\lambda$, $\lambda' \in S^1$ we say that 
$\lambda'$ is {\sl smaller} than $\lambda$ if $\arg \lambda' < \arg \lambda$ with 
respect to the argument defined in \eqref{eqn:arg}. 
For any $\lambda \in S^1$ let 
%%%%%%%%%%%%%%%%%%%%%%%%%%% eqn:r(lam) %%%%%%%%%%%%%%%%%%%%%%%%%%%%%
\begin{equation} \label{eqn:r(lam)}
\mbox{$r(\lambda) := \#$ of all elements in $\baon \cup \bbon$ that are 
smaller than $\lambda$},  
\end{equation}
%%%%%%%%%%%%%%%%%%%%%%%%%%%%%%%%%%%%%%%%%%%%%%%%%%%%%%%%%%%%%%%%%
where $\#$ denotes the cardinality counted with multiplicities.       
%%%%%%%%%%%%%%%%%%%%%%%%%% prop:lsgn %%%%%%%%%%%%%%%%%%%%%%%%%%%%%%%
\begin{proposition} \label{prop:lsgn} 
For each $\lambda \in \baon \cup \bbon$ the Hermitian form restricted to 
$E(\lambda)$ has index 
%%%%%%%%%%%%%%%%%%%%%%%%%% eqn:lsgn1 %%%%%%%%%%%%%%%%%%%%%%%%%%%%%% 
\begin{equation} \label{eqn:lsgn1}
\idx(\lambda) := 
\begin{cases} 
\varepsilon \cdot (-1)^{r(\lambda)} \quad & 
\mbox{if $\lambda \in \baon$ and $m(\lambda)$ is odd}, \\[1mm]
\varepsilon \cdot (-1)^{r(\lambda) + 1} \quad & 
\mbox{if $\lambda \in \bbon$ and $m(\lambda)$ is odd}, \\[1mm]
\qquad 0 \quad & \mbox{if $m(\lambda)$ is even}, 
\end{cases}
\end{equation}
%%%%%%%%%%%%%%%%%%%%%%%%%%%%%%%%%%%%%%%%%%%%%%%%%%%%%%%%%%%%%%%%%
where $\varepsilon = \pm 1$ is the sign of $\sin \pi \gamma$ mentioned 
in Theorem $\ref{thm:p-q}$.  
For each $\mu \in \baoff \cup \bboff$ with $|\mu| > 1$ the Hermitian form 
restricted to $E(\mu, \mu^{\dagger})$ has null index, that is, 
%%%%%%%%%%%%%%%%%%%%%%%%%% eqn:lsgn2 %%%%%%%%%%%%%%%%%%%%%%%%%%%%%%
\begin{equation} \label{eqn:lsgn2}
\idx(\mu) = 0. 
\end{equation}
%%%%%%%%%%%%%%%%%%%%%%%%%%%%%%%%%%%%%%%%%%%%%%%%%%%%%%%%%%%%%%%%%
\end{proposition} 
%%%%%%%%%%%%%%%%%%%%%%%%%% begin proof %%%%%%%%%%%%%%%%%%%%%%%%%%%%
{\it Proof}. 
First we show claims \eqref{eqn:lsgn1} and \eqref{eqn:lsgn2} for 
$\lambda \in \baon$ and $\mu \in \baoff$.  
In the special case where $a_1, \dots, a_n$ are distinct and hence 
$\lambda$ and $\mu$ are simple, they are direct consequences of 
\eqref{eqn:sim2} and \eqref{eqn:Lam} respectively. 
In the general case they are then obtained by using perturbation theory 
of eigenprojections in Kato \cite[Chapter I\!I, \S 1.4]{Kato}, since 
$E(\lambda)$ and $E(\mu, \mu^{\dagger})$ are the total eigenspaces 
for the $\lambda$-group and  $\mu$-group (in Kato's terminology) 
respectively.  
%%%%%%
\par
%%%%%%
Next we can show the results for $\lambda \in \bbon$ and 
$\mu \in \bboff$ in a similar manner by exchanging the roles of 
$\ba$ and $\bb$. 
Notice that $r(\lambda) =: r_a(\lambda)$ and $\varepsilon =:  
\varepsilon_a$ in \eqref{eqn:lsgn1} are defined with respect to the 
branch cut $\ell = \ell_a$ lying between $\bb_t$ and $\ba_1$ as in 
Figure \ref{fig:clus}. 
Changing the roles of $\ba$ and $\bb$ we should replace $\ell_a$ by 
a new branch cut $\ell_b$ lying between $\ba_1$ and $\bb_1$ and 
consider the corresponding $r_b(\lambda)$ and $\varepsilon_b$. 
For $\lambda \in \bbon$ one has 
$r_b(\lambda) = r_a(\lambda) - |\ba_1|$ and 
$\varepsilon_b = - \varepsilon_a \cdot (-1)^{|\sba_1|}$, where 
Remark \ref{rem:branch} is used to obtain the latter relation.     
Thus $\varepsilon_b \cdot (-1)^{r_b(\lambda)} = \varepsilon_a 
\cdot (-1)^{r_a(\lambda)+1}$, which proves \eqref{eqn:lsgn1} for 
$\lambda \in \bbon$. 
The proof of \eqref{eqn:lsgn2} for $\mu \in \bboff$ is 
the same as that for $\mu \in \baoff$.  
\hfill $\Box$ 
%%%%%%%%%%%%%%%%%%%%%%%%%% end proof %%%%%%%%%%%%%%%%%%%%%%%%%%%%% 
%%%%%%%%%%%%%%%%%%%%%%%%%% rem:lsgn %%%%%%%%%%%%%%%%%%%%%%%%%%%%%%
\begin{remark} \label{rem:lsgn} 
Due to \eqref{eqn:lsgn2} the sum of $\idx(\lambda)$ over all distinct 
elements $\lambda$ of $\baon$ (or of $\bbon$) is equal to the global 
index $p-q$ of the Hermitian form on the whole space $\bC^n$. 
\end{remark}
%%%%%%%%%%%%%%%%%%%%%%%%%%%%%%%%%%%%%%%%%%%%%%%%%%%%%%%%%%%%%%%%  
%%%%%%%%%%%%%%%%%%%%%%%%%% sec:real %%%%%%%%%%%%%%%%%%%%%%%%%%%%%%%
\section{Real Hypergeometric Groups} \label{sec:real}
%%%%%%%%%%%%%%%%%%%%%%%%%%%%%%%%%%%%%%%%%%%%%%%%%%%%%%%%%%%%%%%%
A hypergeometric group $H = H(\varphi, \psi) = H(\ba, \bb) = 
\langle A, B \rangle$ is said to be {\sl real} if  
%%%%%
$$
\varphi(z) \in \bR[z], \quad \psi(z) \in \bR[z],   
\quad \mbox{or equivalently} \quad 
\bar{\ba} = \ba, \quad \bar{\bb} = \bb.   
$$ 
%%%%%
In the rest of this article we always assume that $H$ is a real hypergeometric group.  
Let $L_{\bR}$ be the $\bR$-linear span of $\br, A \br, \dots, A^{n-1} \br$. 
These vectors form an $\bR$-linear basis of $L_{\bR}$ as they are a 
$\bC$-linear basis of $\bC^n$ by Theorem \ref{thm:ihf}. 
Obviously $A$ preserves $L_{\bR}$. 
Since $\xi_i \in \bR$ for every $i \in \bZ_{\ge 0}$, formula 
\eqref{eqn:gram} shows that the Hermitian form is $\bR$-valued on $L_{\bR}$. 
The matrix $C = A^{-1} B$ acts on $L_{\bR}$ as a real reflection because   
formulas \eqref{eqn:C} and \eqref{eqn:nmlz} read
%%%%%%%%%%%%%%%%%%%%%%%%% eqn:C2 %%%%%%%%%%%%%%%%%%%%%%%%%%%%%%%%%
\begin{equation} \label{eqn:C2}
C \bv = \bv - (\bv, \br) \br, \quad \bv \in L_{\bR} \quad 
\mbox{with} \quad (\br, \br) = 2.    
\end{equation}
%%%%%%%%%%%%%%%%%%%%%%%%%%%%%%%%%%%%%%%%%%%%%%%%%%%%%%%%%%%%%%%%
So $L_{\bR}$ is also preserved by $B = A C$ and hence by the whole group 
$H$, thus $H \subset O(L_{\bR})$.  
%%%%%%%%%%%%%%%%%%%%%%%%%% ss:trp %%%%%%%%%%%%%%%%%%%%%%%%%%%%%%%%
\subsection{Trace Polynomials} \label{ss:trp}
%%%%%%%%%%%%%%%%%%%%%%%%%%%%%%%%%%%%%%%%%%%%%%%%%%%%%%%%%%%%%%%%
Conditions \eqref{eqn:m1} implies $\varphi(0) = \pm1$ and 
$\psi(0) = \pm 1$, while assumption \eqref{eqn:cneq1} forces 
%%%%%%%%%%%%%%%%%%%%%%%%%% eqn:rr2 %%%%%%%%%%%%%%%%%%%%%%%%%%%%%%%
\begin{equation} \label{eqn:rr2} 
c = -1, \qquad \zeta = 1, \qquad (\br, \, \br) = 2, 
\end{equation}    
%%%%%%%%%%%%%%%%%%%%%%%%%%%%%%%%%%%%%%%%%%%%%%%%%%%%%%%%%%%%%%%%
in \eqref{eqn:nmlz}, hence $\varphi(0)$ and $\psi(0)$ must have opposite signs. 
Thus after exchanging $\varphi(z)$ and $\psi(z)$ if necessary we may assume 
$\varphi(0) = -1$ and $\psi(0) = 1$ so that \eqref{eqn:ihf2} becomes 
$z^n \varphi(z^{-1}) = - \varphi(z)$ and $z^n \psi(z^{-1}) = \psi(z)$, that is, 
$\varphi(z)$ is anti-palindromic while $\psi(z)$ is palindromic. 
%%%%%%
\par
%%%%%%
In general a palindromic polynomial $f(z)$ of even degree $2 d$ can 
be expressed as $f(z) = z^d F(z+z^{-1})$ for a unique polynomial $F(w)$ of 
degree $d$, a palindromic polynomial $f(z)$ of odd degree factors as 
$f(z) = (z+1) g(z)$ with $g(z)$ being palindromic of even degree, and 
an anti-palindromic polynomial $f(z)$ factors as $f(z) = (z-1) g(z)$ with 
$g(z)$ being palindromic.     
Hence there exist unique monic real polynomials $\Phi(w)$ and $\Psi(w)$ 
such that 
%%%%%%%%%%%%%%%%%%%%%%% eqn:palin %%%%%%%%%%%%%%%%%%%%%%%%%%%%%%%%
\begin{subequations} \label{eqn:palin}
\begin{alignat}{3}
\varphi(z) &= (z^2-1) z^{N-1} \Phi(z+z^{-1}), \quad & 
\psi(z) &= z^N \Psi(z+z^{-1}), \quad & & \mbox{if $n = 2 N$}; 
\label{eqn:paline} \\[1mm]
\varphi(z) &= (z-1) z^N \Phi(z+z^{-1}), \quad & 
\psi(z) &= (z+1) z^N \Psi(z+z^{-1}), \quad & & \mbox{if $n = 2 N+1$}. 
\label{eqn:palino}
\end{alignat}
\end{subequations}
%%%%%%%%%%%%%%%%%%%%%%%%%%%%%%%%%%%%%%%%%%%%%%%%%%%%%%%%%%%%%%% 
We refer to $\Phi(w)$ and $\Psi(w)$ as the {\sl trace} polynomials of 
$\varphi(z)$ and $\psi(z)$. 
It is easily seen from \eqref{eqn:palin} that the resultant of $(\varphi, \psi)$ 
and that of $(\Phi, \Psi)$ are related by     
%%%%%%%%%%%%%%%%%%%%%%%% eqn:rest %%%%%%%%%%%%%%%%%%%%%%%%%%%%%%%
\begin{subequations} \label{eqn:rest}
\begin{alignat}{2} 
\Res(\varphi, \psi) &= \phantom{2} (-1)^N \cdot \Psi(2) \cdot \Psi(-2) 
\cdot \Res(\Phi, \Psi)^2, \qquad & & \mbox{if $n = 2 N$}; 
\label{eqn:reste} \\[1mm]
\Res(\varphi, \psi) &= 2 (-1)^N \cdot \Psi(2) \cdot \Phi(-2) \cdot 
\Res(\Phi, \Psi)^2, \qquad & & \mbox{if $n = 2 N+1$}. \label{eqn:resto}
\end{alignat}
\end{subequations}
%%%%%%%%%%%%%%%%%%%%%%%%%%%%%%%%%%%%%%%%%%%%%%%%%%%%%%%%%%%%%%%
The real hypergeometric group $H = H(\varphi, \psi)$ can also be 
expressed as $H = H(\Phi, \Psi)$.  
%%%%%%%
\par
%%%%%% 
It follows from \eqref{eqn:irr2} and \eqref{eqn:rest} that 
$\Res(\Phi, \Psi) \neq 0$ hence $\Phi(w)$ and $\Psi(w)$ have no root in 
common. 
By formulas \eqref{eqn:palin} the roots $\lambda \neq \pm 1$ 
of $\varphi(z)$ are in two-to-one correspondence with the roots 
$\tau \neq \pm 2$ of $\Phi(w)$ via the relation 
$\tau = \lambda + \lambda^{-1}$, since   
$w - \tau = z^{-1}(z-\lambda)(z-\lambda^{-1})$ with  
$w = z + z^{-1}$.   
The same statement is true for $\psi(z)$ and $\Psi(w)$. 
Moreover we have
%%%%%%%%%%%%%%%%%%%%%%%%% eqn:m-M %%%%%%%%%%%%%%%%%%%%%%%%%%%%%%
\begin{equation} \label{eqn:m-M}
m(\lambda) = 
\begin{cases}
M(\tau) & \mbox{if $\lambda \neq \pm 1$, i.e. $\tau \neq \pm 2$}, 
\\[1mm]
2 M(\tau) + 1 \quad & \mbox{if $\lambda = \pm 1$, i.e. $\tau = \pm 2$}, 
\end{cases}
\end{equation}
%%%%%%%%%%%%%%%%%%%%%%%%%%%%%%%%%%%%%%%%%%%%%%%%%%%%%%%%%%%%%%%%
under $\tau := \lambda + \lambda^{-1}$ where $M(\tau)$ is the multiplicity 
of $w = \tau$ in the equation $\Phi(w) \cdot \Psi(w) = 0$.  
%%%%%%%%%%%%%%%%%%%%%%%%%% ss:trace %%%%%%%%%%%%%%%%%%%%%%%%%%%%%%
\subsection{Trace Clusters and Index} \label{ss:trace}
%%%%%%%%%%%%%%%%%%%%%%%%%%%%%%%%%%%%%%%%%%%%%%%%%%%%%%%%%%%%%%%%%
In the real case the index formula \eqref{eqn:p-q} in Theorem \ref{thm:p-q} 
can be restated in terms of what we call trace clusters.        
In view of formulas \eqref{eqn:palin}, if $n$ is even then $\pm1 \in \baon$ 
 and hence $\baon$ is nonempty (but $\bbon$ may be empty), while if $n$ 
is odd then $1 \in \baon$, $-1 \in \bbon$ and hence both of $\baon$ and 
$\bbon$ are nonempty. 
In any case, as far as both of them are nonempty, the clusters 
$\ba_1, \bb_1, \dots, \ba_t, \bb_t$ can be indexed so that $1 \in \ba_1$.    
With this convention it is easy to see that 
%%%%%%%%%%%%%%%%%%%%%%%%% eqn:evod %%%%%%%%%%%%%%%%%%%%%%%%%%%%%%%
\begin{subequations}
\begin{align}
&\mbox{if $n$ is even then $t = 2 s$ is also even, $-1 \in \ba_{s+1}$ and  
$|\ba_1| \equiv |\ba_{s+1}| \equiv 1 \bmod 2$}, \label{eqn:ev}  \\[2mm]  
&\mbox{if $n$ is odd then $t = 2 s-1$ is also odd, $-1 \in \bb_s$ and  
$|\ba_1| \equiv |\bb_s| \equiv 1 \bmod 2$}.  \label{eqn:od}
\end{align} 
\end{subequations}
%%%%%%%%%%%%%%%%%%%%%%%%%%%%%%%%%%%%%%%%%%%%%%%%%%%%%%%%%%%%%%%%
In either case, with the convention $\ba_{t+1} = \ba_1$, $\bb_{t+1} = \bb_1$, 
we have    
%%%%%%%%%%%%%%%%%%%%%%%% eqn:bar %%%%%%%%%%%%%%%%%%%%%%%%%%%%%%%%%
\begin{equation} \label{eqn:bar}
\bar{\ba}_i = \ba_{t+2 -i}, \qquad 
\bar{\bb}_i = \bb_{t+1 -i}, \qquad i = 1, \dots, s. 
\end{equation}
%%%%%%%%%%%%%%%%%%%%%%%%%%%%%%%%%%%%%%%%%%%%%%%%%%%%%%%%%%%%%%%% 
Note that  $\ba_2, \dots, \ba_s$ are exactly those $\baon$-clusters which lie 
in the upper half-plane $\Im \, z > 0$.  
%%%%%%
\par
%%%%%%
Let $\bA$ be the multi-set of all complex roots of 
$\Phi(w)$ and $\bB$ be its $\Psi(w)$-counterpart. 
Let $\bAon$ resp. $\bAoff$ be the component  
of $\bA$ whose elements lie on resp. off $[-2, \, 2]$. 
Let $\bBon$ and $\bBoff$ be defined in a similar manner for $\bB$.     
Notice that both of $\baon$ and $\bbon$ are nonempty if and only if 
%%%%%%%%%%%%%%%%%%%%%%%% eqn:nemp %%%%%%%%%%%%%%%%%%%%%%%%%%%%%%%
\begin{equation} \label{eqn:nemp}
\mbox{either $n$ is even and $\bBon$ is nonempty, or $n$ is odd}, 
\end{equation}
%%%%%%%%%%%%%%%%%%%%%%%%%%%%%%%%%%%%%%%%%%%%%%%%%%%%%%%%%%%%%%%%
in which case $\bAon$ and $\bBon$ dissect each other into interlacing 
components called {\sl trace clusters}, 
%%%%%% 
$$
\bA_{s+1}, \bB_s, \bA_s, \dots, \bB_1, \bA_1 \quad \mbox{if $n$ is even}; 
\quad 
\bB_s, \bA_s, \dots, \bB_1, \bA_1 \quad \mbox{if $n$ is odd}, 
$$
%%%%%%
where one or both of the end clusters may be empty but all the other 
clusters must be nonempty. 
Put 
%%%%%%%%%%%%%%%%%%%%%%%% eqn:Ain %%%%%%%%%%%%%%%%%%%%%%%%%%%%%%%%%
\begin{equation} \label{eqn:Ain}
\bAin := \bA_2 \cup \cdots \cup \bA_s. 
\end{equation}
%%%%%%%%%%%%%%%%%%%%%%%%%%%%%%%%%%%%%%%%%%%%%%%%%%%%%%%%%%%%%%%%  
Finally, let $\bA_{>2}$ be the components of $\bAoff$ whose elements 
are real numbers greater than $2$ and $\bB_{>2}$ be defined in a 
similar manner for $\bB$.  
%%%%%%%%%%%%%%%%%%%%%%%%%% lem:gam %%%%%%%%%%%%%%%%%%%%%%%%%%%%%%
\begin{lemma} \label{lem:gam}  
Let $\gamma \in \bR \setminus \bZ$ be the number defined in 
\eqref{eqn:cneq2}. 
Under assumption \eqref{eqn:nemp} we have  
%%%%%%%%%%%%%%%%%%%%%%%%% eqn:gam %%%%%%%%%%%%%%%%%%%%%%%%%%%%%%%
\begin{equation} \label{eqn:gam}  
\ve = \sin \pi \gamma = (-1)^{|\sbA_1| + |\sbA_{>2}| + |\sbB_{>2}|}. 
\end{equation}
%%%%%%%%%%%%%%%%%%%%%%%%%%%%%%%%%%%%%%%%%%%%%%%%%%%%%%%%%%%%%%%%
\end{lemma}
%%%%%%%%%%%%%%%%%%%%%%%%%% begin proof %%%%%%%%%%%%%%%%%%%%%%%%%%%%
{\it Proof}. 
We consider how each component of $\ba$ contributes to the sum 
$2 \pi \alpha := 2 \pi \alpha_1 + \cdots + 2 \pi \alpha_n$, where 
$2 \pi \alpha_i := \arg a_i$. 
The $1$'s in $\baon$ has no contribution. 
The $-1$'s in $\baon$ has contribution $\pi m_a^-$ where  
$m_a^-$ is the multiplicity of $-1$ in $\ba$. 
For each non-real pair $\lambda, \bar{\lambda} \in \ba_{\mathrm{on}}$ 
with $\Im \, \lambda > 0$, the sum $\arg \lambda + \arg \bar{\lambda}$ 
is $0$ if $\lambda \in \ba_1$ and $2 \pi$ if $\lambda \not \in \ba_1$.    
So the total contribution of $\baon$ is given by 
%%%%%%
$$
\pi m_a^- + 2 \pi \cdot \frac{|\ba_{\mathrm{on}}|-|\ba_1| - m_a^{-}}{2} = 
\pi (|\ba_{\mathrm{on}}|-|\ba_1|). 
$$
%%%%%%
Let $\ba_{< -1}$ resp. $\ba_{>1}$ be the component of $\baoff$ 
whose elements are real numbers $< -1$ resp. $> 1$.    
For each real pair $\lambda$, $\lambda^{\dagger} \in \baoff$, 
the sum $\arg \lambda + \arg \lambda^{\dagger}$ is $0$ if 
$\lambda \in \ba_{> 1}$ and $2 \pi$ if $\lambda \in \ba_{< -1}$. 
For each non-real quartet $\lambda$, $\bar{\lambda}$, 
$\lambda^{\dagger}$, $\bar{\lambda}^{\dagger} \in 
\ba_{\mathrm{off}}$ with $|\lambda| > 1$ and $\Im \, \lambda > 0$ 
we have   
%%%%%%
$$
\arg \lambda + \arg \bar{\lambda} + \arg \lambda^{\dagger} + \arg 
\bar{\lambda}^{\dagger} = 
\begin{cases}
0 \quad & (0 < \arg \lambda \le |\Theta|), \\
4 \pi \quad & (|\Theta| < \arg \lambda < \pi),  
\end{cases}
$$
%%%%%%
where $\Theta$ is the number appearing in \eqref{eqn:arg}, which 
belongs to the interval $(-\pi, \, 0)$ due to assumption \eqref{eqn:nemp}. 
Therefore the contribution of $\ba_{\mathrm{off}}$ to the sum 
$2 \pi \alpha$ is $2 \pi |\ba_{<-1}| \bmod 4 \pi \bZ$.  
In total we have $2 \pi \alpha \equiv \pi (|\ba_{\mathrm{on}}|-|\ba_1|) + 
2 \pi |\ba_{<-1}| \bmod 4 \pi \bZ$, which yields a modulo $2$ congruence   
%%%%%%% 
$$ 
\alpha \equiv \frac{|\ba_{\mathrm{on}}|-|\ba_1|}{2} + 
|\ba_{<-1}| \mod 2. 
$$
%%%%%%
\par
%%%%%%
In a similar manner we consider how each component of $\bb$ contributes 
to the sum $2 \pi \beta := 2 \pi \beta_1 + \cdots + 2 \pi \beta_n$, where 
$2 \pi \beta_i := \arg b_i$.  
Taking $1 \not \in \bb$ into account we find that 
%%%%%%% 
$$ 
\beta \equiv \frac{|\bb_{\mathrm{on}}|}{2} + |\bb_{<-1}| \mod 2. 
$$
%%%%%%
Since $\gamma = \beta - \alpha$, we use relations   
$|\ba_{\mathrm{on}}| + |\ba_{\mathrm{off}}| = n$, 
$|\ba_{\mathrm{off}}| \equiv 2 |\ba_{<-1}| + 2 |\ba_{>1}| \bmod 4$ and 
their $\bb$-counterparts together with $|\ba_1| = 2 |\bA_1|+1$ to obtain  
a modulo $2$ congruence  
%%%%%%
\begin{align*}
\gamma 
&\equiv \frac{|\bb_{\mathrm{on}}| - |\ba_{\mathrm{on}}| + 
|\ba_1|}{2} + |\bb_{<-1}| - |\ba_{<-1}| 
\equiv \frac{|\ba_1|}{2} + |\ba_{>1}| - |\bb_{>1}| \\ 
&= \frac{1}{2} (2 |\bA_1|+1) + |\bA_{>2}| - |\bB_{>2}|
\equiv \frac{1}{2} + |\bA_1| + |\bA_{>2}| + |\bB_{>2}| \bmod 2.    
\end{align*}
%%%%%%
This establishes formula \eqref{eqn:gam}. \hfill $\Box$ \par\medskip
%%%%%%%%%%%%%%%%%%%%%%%%%% end proof %%%%%%%%%%%%%%%%%%%%%%%%%%%%% 
In the real case the index formula \eqref{eqn:p-q} in Theorem 
$\ref{thm:p-q}$ can be restated as follows.  
%%%%%%%%%%%%%%%%%%%%%%%%%% thm:real %%%%%%%%%%%%%%%%%%%%%%%%%%%%%%
\begin{theorem} \label{thm:real} 
Let $H = H(\Phi, \Psi) = \langle A, B \rangle$ be a real hypergeometric 
group of rank $n$. 
If condition \eqref{eqn:nemp} is satisfied then the index of the  
$H$-invariant Hermitian form is given by   
%%%%%%%%%%%%%%%%%%%%%%%%%% eqn:rsgn2 %%%%%%%%%%%%%%%%%%%%%%%%%%%%%
\begin{equation} \label{eqn:rsgn2}
p-q = \ve( 1+ \delta - 2 S),  
\end{equation}
%%%%%%%%%%%%%%%%%%%%%%%%%%%%%%%%%%%%%%%%%%%%%%%%%%%%%%%%%%%%%%%%
where $\ve = \pm 1$ is given in formula \eqref{eqn:gam} while 
$\delta$ and $S$ are defined by   
%%%%%%%%%%%%%%%%%%%%%%%%%% eqn:gam2a %%%%%%%%%%%%%%%%%%%%%%%%%%%%
\begin{subequations} \label{eqn:deltaS}
\begin{align} 
\delta &:= 
\begin{cases} 
(-1)^{ |\sbA_{\mathrm{in}}| + |\sbB_{\mathrm{on}} | +1}  & 
\mbox{if $n$ is even}, \\[1mm]
0 & \mbox{if $n$ is odd}, 
\end{cases} 
\label{eqn:delta} \\[2mm]
S &:= \sum_{i \in \rI} (-1)^{\sigma_i} \qquad \mbox{with} \quad  
\rI := \{\, i = 2, \dots, s : |\bA_i| \equiv 1 \bmod 2 \, \},   
\label{eqn:S} 
\end{align} 
\end{subequations} 
%%%%%%%%%%%%%%%%%%%%%%%%%%%%%%%%%%%%%%%%%%%%%%%%%%%%%%%%%%%%%%%% 
with $\bAin$ being as in \eqref{eqn:Ain} and $\sigma_i$ defined by 
$\sigma_2 := |\bB_1|$ and  
%%%%%%%
$$
\sigma_i := |\bB_1|+|\bA_2| + |\bB_2| + 
\cdots + |\bA_{i-1}|+|\bB_{i-1}|, \qquad i = 3, \dots, s. 
$$
%%%%%%%
If $n$ is even and $\bBon$ is empty then the index $p-q$ is zero. 
\end{theorem}
%%%%%%%%%%%%%%%%%%%%%%%%%% begin proof %%%%%%%%%%%%%%%%%%%%%%%%%%%%
{\it Proof}.  
If $n$ is even and $\bBon$ is nonempty, then \eqref{eqn:ev} and 
\eqref{eqn:bar} imply $|\ba_1| = 2 |\bA_1|+1$, $|\ba_{s+1}| = 2 |\bA_{s+1}|+1$, 
$|\ba_i| = |\ba_{2s+2-i}| = |\bA_i|$ for $i = 2, \dots, s$, and 
$|\bb_i| = |\bb_{2s+1-i}| = |\bB_i|$ for $i = 1, \dots, s$,  hence 
$\rI = \rK \cap \{2, \dots, s\}$ and  
$\rK = \{1, s+1\} \sqcup \rI \sqcup \{2s+2-i \,:\, i \in \rI \}$, 
where $\rK$ is defined in \eqref{eqn:p-q}.   
For each $i \in \rI$ we have 
$\tau_i = |\ba_1| + \sigma_i \equiv 1+ \sigma_i \bmod 2$ and   
%%%%%%%%%%%%
\begin{align*}
\tau_{2 s + 2 -i} 
&= |\ba_1|+ \cdots +|\ba_{2 s + 1-i}| + |\bb_1|+ \cdots + |\bb_{2 s + 1-i}| \\[1mm]
&\equiv |\ba_1|+ \cdots + |\ba_i| + |\ba_{s+1}| + |\bb_1|+ \cdots + |\bb_{i-1}| 
\bmod 2 \\[1mm]
&= \tau_i +  |\ba_i| + |\ba_{s+1}| \equiv \tau_i \bmod 2.  
\end{align*}
%%%%%%%%%%%%
Moreover, $\tau_{s+1} = |\ba_1| + |\bA_{\mathrm{in}}| + |\bBon| \equiv 
1 + |\bA_{\mathrm{in}}| + |\bBon| \bmod 2$. 
Formula \eqref{eqn:p-q} then yields  
%%%%%%%%%%%%
\begin{align*}
\ve(p-q) 
&= \sum_{i \in \rK} (-1)^{\tau_i} = (-1)^{\tau_1} + 
\sum_{i \in \rI} (-1)^{\tau_i} + (-1)^{\tau_{s+1}} + 
\sum_{i \in \rI} (-1)^{\tau_{2s+2-i}} \\[1mm]
&= 1+ (-1)^{\tau_{s+1}} + 2 \sum_{i \in \rI} (-1)^{\tau_i} = 
1 + \delta - 2 \sum_{i \in \rI} (-1)^{\sigma_i}.  
\end{align*}
%%%%%%%%%%%%
\par
%%%%%%%%%%%%
If $n$ is odd then \eqref{eqn:od} and \eqref{eqn:bar} imply  
$|\ba_1| = 2 |\bA_1|+1$, $|\ba_i| = |\ba_{2s+2-i}| = |\bA_i|$ for $i = 2, \dots, s$, 
$|\bb_s| = 2 |\bB_s|+1$, and $|\bb_i| = |\bb_{2s-i}| = |\bB_i|$ for $i = 1, \dots, s-1$, 
hence $\rI = \rK \cap \{2, \dots, s\}$ and  
$\rK = \{1 \} \sqcup \rI \sqcup \{2s+1-i \,:\, i \in \rI \}$.  
For each $i \in \rI$ we have 
$\tau_i = |\ba_1| + \sigma_i \equiv 1+ \sigma_i \bmod 2$ and     
%%%%%%%%%%%%
\begin{align*}
\tau_{2 s + 1 -i} 
&= |\ba_1|+ \cdots +|\ba_{2 s -i}| + |\bb_1|+ \cdots + |\bb_{2 s -i}| \\[1mm]
&\equiv |\ba_1|+ \cdots + |\ba_i|+ |\bb_1|+ \cdots + |\bb_{i-1}|  + |\bb_s|  
\bmod 2 \\[1mm]
&= \tau_i +  |\ba_i| + |\bb_s| \equiv \tau_i \bmod 2.  
\end{align*}
%%%%%%%%%%%%
Thus formula \eqref{eqn:p-q} in Theorem \ref{thm:p-q} then yields  
%%%%%%%%%%%%
\begin{align*}
\ve(p-q) 
&= \sum_{i \in \rK} (-1)^{\tau_i} = (-1)^{\tau_1} + 
\sum_{i \in \rI} (-1)^{\tau_i} + 
\sum_{i \in \rI} (-1)^{\tau_{2s+1-i}} \\[1mm]
&= 1 + 2 \sum_{i \in \rI} (-1)^{\tau_i} 
= 1 + \delta -2 \sum_{i \in \rI} (-1)^{\sigma_i}.  
\end{align*}
%%%%%%%%%%%%
\par
%%%%%%%%%%%%
In either case we have obtained formula \eqref{eqn:rsgn2}. 
If $n$ is even and $\bBon$ is empty, then $\bbon$ is empty 
and hence the index $p-q$ is zero by the last part of Theorem 
\ref{thm:p-q}.  \hfill $\Box$
%%%%%%%%%%%%%%%%%%%%%%%%%% end proof %%%%%%%%%%%%%%%%%%%%%%%%%%%%%
%%%%%%%%%%%%%%%%%%%%%%%%%% lem:estim %%%%%%%%%%%%%%%%%%%%%%%%%%%%%
\begin{lemma} \label{lem:estim}
There are the following numerical constraints 
%%%%%%%%%%%%%%%%%%%%%%%%%% eqn:estim %%%%%%%%%%%%%%%%%%%%%%%%%%%%%
\begin{equation} \label{eqn:estim} 
S \equiv |\rI| \equiv |\bA_{\mathrm{in}}| \bmod 2, \quad 
|S| \le |\rI| \le s-1 \le \frac{|\rI|+|\bA_{\mathrm{in}}|}{2} 
\le |\bA_{\mathrm{in}}|, \quad 
s \le |\bBon|.  
\end{equation}
%%%%%%%%%%%%%%%%%%%%%%%%%%%%%%%%%%%%%%%%%%%%%%%%%%%%%%%%%%%%%%%%%
\end{lemma}
%%%%%%%%%%%%%%%%%%%%%%%%%% begin proof %%%%%%%%%%%%%%%%%%%%%%%%%%%%
{\it Proof}. 
Congruence $S \equiv |\rI| \bmod 2$ follows from the definition of $S$ 
in \eqref{eqn:S} and $(-1)^{\sigma_i} \equiv 1 \bmod 2$.  
Since $|\bA_i| \ge 1$ for $i = 2, \dots, s$, the definition of $S$ and the 
inclusion $\rI \subset \{2, \dots, s\}$ imply $|S| \le |\rI| \le s-1 \le 
|\bA_2|+ \cdots + |\bA_s| = |\bA_{\mathrm{in}}|$. 
As $|\bA_i|$ is odd for $i \in \rI$ and even for $i \in \{2, \dots, s\} 
\setminus \rI$, we have $|\rI| \equiv |\bA_{\mathrm{in}}| \bmod 2$ 
and $|\rI| + 2(s-1-|\rI|) \le |\bA_{\mathrm{in}}|$, that is, 
$s-1 \le \frac{1}{2}(|\rI|+|\bA_{\mathrm{in}})$. 
Moreover we have $s \le |\bB_1|+\cdots+|\bB_s| = |\bBon|$ because 
$|\bB_i| \ge 1$ for $i = 1, \dots, s$.  
Putting all these together lead to the constraints 
\eqref{eqn:estim}. \hfill $\Box$ 
%%%%%%%%%%%%%%%%%%%%%%%%%% rem:refl %%%%%%%%%%%%%%%%%%%%%%%%%%%%%%%
\begin{remark} \label{rem:refl} 
If $n$ is even then the reflection of the trace clusters in $\bA_{\mathrm{in}} 
\cup \bBon$, 
%%%%%%%%%%%%%%%%%%%%%%%%%% eqn:refl %%%%%%%%%%%%%%%%%%%%%%%%%%%%%%%%
\begin{equation} \label{eqn:refl}
\begin{CD}
\bB_s, \bA_s, \bB_{s-1}, \dots, \bB_2, \bA_2, \bB_1 \quad 
@> \mathrm{reflection} >> \quad 
\bB_1, \bA_2, \bB_2, \dots, \bB_{s-1}, \bA_s, \bB_s
\end{CD}
\end{equation}
%%%%%%%%%%%%%%%%%%%%%%%%%%%%%%%%%%%%%%%%%%%%%%%%%%%%%%%%%%%%%%%%%
results in the change of signs $S \to \delta S$ and 
$p-q \to \delta (p-q)$, where $\delta$ is defined in \eqref{eqn:delta}.   
\end{remark}
%%%%%%%%%%%%%%%%%%%%%%%%%%%%%%%%%%%%%%%%%%%%%%%%%%%%%%%%%%%%%%%%
\subsection{Local Index in the Real Case} \label{ss:lsrc}
%%%%%%%%%%%%%%%%%%%%%%%%%%%%%%%%%%%%%%%%%%%%%%%%%%%%%%%%%%%%%%%%
Proposition \ref{prop:lsgn} gives formula \eqref{eqn:lsgn1} for the local 
index $\idx(\lambda)$ at $\lambda \in \baon \cup \bbon$. 
It can be restated in terms of $\tau := \lambda + \lambda^{-1} \in 
\bAon \cup \bBon$, where we are allowed to write $\idx(\lambda) = 
\idx(\lambda^{-1}) = \Idx(\tau)$ if $\lambda \neq \pm 1$, i.e. $\tau \neq \pm 2$.     
To state the result let $\rho : [-2, \, 2] \to \bZ_{\ge 0}$ be a function 
defined by  
%%%%%%%%%%%%%%%%%%%%%%%%%%% eqn:rho %%%%%%%%%%%%%%%%%%%%%%%%%%%%%%
\begin{equation} \label{eqn:rho} 
\rho(\tau) := \mbox{$\#$ of all real roots of $\Phi(w) \cdot \Psi(w)$ 
that are} 
\begin{cases}
> \tau \quad \mbox{if $\tau \in [-2, \, 2)$},  \\[1mm] 
\ge 2  \quad \mbox{if $\tau = 2$},       
\end{cases}
\end{equation} 
%%%%%%%%%%%%%%%%%%%%%%%%%%%%%%%%%%%%%%%%%%%%%%%%%%%%%%%%%%%%%%%%
where $\#$ denotes the cardinality counted with multiplicities. 
%%%%%%%%%%%%%%%%%%%%%%%%%% prop:lsgn-r %%%%%%%%%%%%%%%%%%%%%%%%%%%%
\begin{proposition} \label{prop:lsgn-r} 
For any $\tau \in \bAon \cup \bBon$ with $\tau \neq \pm 2$  
we have     
%%%%%%%%%%%%%%%%%%%%%%%%%% eqn:lsgn1r %%%%%%%%%%%%%%%%%%%%%%%%%%%%%% 
\begin{equation} \label{eqn:lsgn1r}
\Idx(\tau) := 
\begin{cases} 
(-1)^{\rho(\tau)+1} \quad & 
\mbox{if $\tau \in \bAon$ and $M(\tau)$ is odd}, \\[1mm]
(-1)^{\rho(\tau)} \quad & 
\mbox{if $\tau \in \bBon$ and $M(\tau)$ is odd}, \\[1mm]
\quad 0 \quad & \mbox{if $M(\tau)$ is even},     
\end{cases}
\end{equation}
%%%%%%%%%%%%%%%%%%%%%%%%%%%%%%%%%%%%%%%%%%%%%%%%%%%%%%%%%%%%%%%%%
where $M(\tau)$ is the multiplicity of $\tau$ in $\bAon \cup \bBon$.   
Moreover we have  
%%%%%%%%%%%%%%%%%%%%%%%%%% eqn:lsgn2r %%%%%%%%%%%%%%%%%%%%%%%%%%%%%
\begin{equation} \label{eqn:lsgn2r}
\idx(1) = (-1)^{\rho(2)}, \qquad \idx(-1) = (-1)^{\rho(-2)+n+1}. 
\end{equation}
%%%%%%%%%%%%%%%%%%%%%%%%%%%%%%%%%%%%%%%%%%%%%%%%%%%%%%%%%%%%%%%%%
\end{proposition} 
%%%%%%%%%%%%%%%%%%%%%%%%%% begin proof %%%%%%%%%%%%%%%%%%%%%%%%%%%% 
{\it Proof}. 
First we prove \eqref{eqn:lsgn1r}.  
Take an element $\lambda \in \baon \cup \bbon$ such that 
$\lambda + \lambda^{-1} = \tau$.  
We may assume $\Im \, \lambda > 0$ since 
$\idx(\lambda) = \idx(\lambda^{-1})$. 
If $\tau \in \bAon$, that is, $\lambda \in \baon$ then definition 
\eqref{eqn:r(lam)} and relation $|\ba_1| = 2 |\bA_1|+1$ give $r(\lambda) = 
|(\bAon \cup \bBon)_{> \tau}| + |\bA_1| + 1$, which together with  
\eqref{eqn:gam} and \eqref{eqn:rho} yields  
$\ve \cdot (-1)^{r(\lambda)} = 
(-1)^{|\sbA_1|+|\sbA_{>2}|+|\sbB_{>2}|} \cdot (-1)^{|(\sbA_{\mathrm{on}} \cup 
\sbB_{\mathrm{on}})_{> \tau}| + |\sbA_1| + 1} = (-1)^{\rho(\tau) +1}$.        
In a similar manner, if $\tau \in \bBon$ then $\ve \cdot (-1)^{r(\lambda) + 1} 
= (-1)^{\rho(\tau)}$. 
Thus \eqref{eqn:lsgn1r} follows from formula \eqref{eqn:lsgn1}. 
%%%%%
\par
%%%%%
Next we prove \eqref{eqn:lsgn2r}. 
By \eqref{eqn:r(lam)} and \eqref{eqn:rho} to gether with \eqref{eqn:m-M} 
we have 
%%%%%%%%%%%%%%%%
\begin{alignat*}{2}
r(1) &= |\bA_1|- M(2), \qquad & r(-1) &= |\bAon|+|\bBon|-M(-2) + |\bA_1|+1, 
\\[1mm]   
\rho(2) &= M(2) + |\bA_{>2}|+|\bB_{>2}|, \qquad &   
\rho(-2) &=  |\bAon|+|\bBon|-M(-2) + |\bA_{>2}|+|\bB_{>2}|. 
\end{alignat*}
%%%%%%%%%%%%%%% 
Putting these equations into \eqref{eqn:lsgn1} and using 
\eqref{eqn:m-M} and \eqref{eqn:gam} we establish 
\eqref{eqn:lsgn2r}, where we take into account that 
$1 \in \baon$ while $-1 \in \baon$ 
if $n$ is even and $-1 \in \bbon$ if $n$ is odd.   
\hfill $\Box$ \par\medskip
%%%%%%%%%%%%%%%%%%%%%%%%%% end proof %%%%%%%%%%%%%%%%%%%%%%%%%%%%%
Suppose that the rank $n$ is even; the odd case is omitted as it is not 
necessary in this article.  
Let $\bA_1^{\circ} := (\bA_1)_{< 2}$ and $\bA_{s+1}^{\circ} := (\bA_{s+1})_{>-2}$.  
For a multi-set $\bX$ given in boldface, its calligraphic style $\cX$ denotes 
the ordinary set of all {\sl distinct} elements of {\sl odd} multiplicity in $\bX$. 
We apply this rule to $\bX = \bA_1^{\circ}$, $\bA_{s+1}^{\circ}$, $\bAin$,  
$\bBon$ to define $\cX = \cA_1^{\circ}$, $\cA_{s+1}^{\circ}$, $\cAin$, 
$\cBon$ respectively. 
For these sets we put $\Idx(\cX) := \sum_{\tau \in \cX} \Idx(\tau)$.  
Note that $0 \le |\bX| - |\cX| \equiv 0 \bmod 2$.  
Moreover we denote by $\Par(m)$ the parity of $m \in \bZ$, that is, 
$\Par(m) = 0$ for $m$ even and $\Par(m) = 1$ for $m$ odd.    
%%%%%%%%%%%%%%%%%%%%%%%%%% thm:lc %%%%%%%%%%%%%%%%%%%%%%%%%%%%%%%%
\begin{theorem}  \label{thm:lc} 
In the situations of Theorem $\ref{thm:real}$ with the rank $n$ being even 
we have 
%%%%%%%%%%%%%%%%%%%%%%%%%% eqn:lc %%%%%%%%%%%%%%%%%%%%%%%%%%%%%%%%
\begin{subequations} \label{eqn:lc} 
\begin{alignat}{2}
\idx(1) &= \ve \cdot (-1)^{|\sbA_1^{\circ}|}, \qquad &  
\idx(-1) &= \ve \delta \cdot (-1)^{|\sbA_{s+1}^{\circ}|},  \label{eqn:lc1}
\\[2mm] 
\Idx(\cA_1^{\circ}) &= \ve \cdot \Par(|\bA_1^{\circ}|), \qquad & 
\Idx( \cA_{s+1}^{\circ}) 
&= \ve \delta \cdot \Par(|\bA_{s+1}^{\circ}|),  \label{eqn:lc2}
\\[2mm] 
\Idx(\cAin) &= - \ve S,  & \Idx(\cBon) &= (p-q)/2.   \label{eqn:lc3}
\end{alignat} 
\end{subequations}
%%%%%%%%%%%%%%%%%%%%%%%%%%%%%%%%%%%%%%%%%%%%%%%%%%%%%%%%%%%%%%%%% 
\end{theorem}
%%%%%%%%%%%%%%%%%%%%%%% begin proof %%%%%%%%%%%%%%%%%%%%%%%%%%%%%%%
{\it Proof}. 
Formulas in \eqref{eqn:lc1} are easy consequences of \eqref{eqn:lsgn2r}. 
The right-hand side of index formula \eqref{eqn:rsgn2} consists of three 
parts $\ve$, $\ve \delta$ and $-2 \ve S$. 
They are the contributions of $\ba_1$, $\ba_{s+1}$ and $\ba_{\mathrm{in}} := 
\ba_2 \cup \overline{\ba}_2 \cup \cdots \cup \ba_s \cup \overline{\ba}_s$, 
and hence equal to the sums of local indices of all distinct elements in the three 
respective sets.  
Since $\lambda \mapsto \tau = \lambda+\lambda^{-1}$ induces two-to-one 
maps $(\ba_1)_{\neq 1} \to \bA_1^{\circ}$,  
$(\ba_{s+1})_{\neq -1} \to \bA_{s+1}^{\circ}$, $\ba_{\mathrm{in}} \to \bAin$ 
and $\idx(\lambda) = \idx(\lambda^{-1}) = \Idx(\tau)$, we have  
%%%%%
$$
\ve = \idx(1) + 2 \Idx(\cA_1^{\circ}), 
\quad 
\ve \delta = \idx(-1) + 2 \Idx( \cA_{s+1}^{\circ}), 
\quad 
-2 \ve S = 2 \Idx(\cAin),  
$$
%%%%%
which in turn lead to the formulas in \eqref{eqn:lc2} and the first formula in 
\eqref{eqn:lc3}. 
Finally, the total index $p-q$ is the sum of local indices over all distinct 
elements in $\bbon$. 
In view of the two-to-one correspondence $\bbon \to \bBon$ we have 
the second formula in \eqref{eqn:lc3}. \hfill $\Box$
%%%%%%%%%%%%%%%%%%%%%%% end proof %%%%%%%%%%%%%%%%%%%%%%%%%%%%%%%%
%%%%%%%%%%%%%%%%%%%%%%%%%% sec:hgl %%%%%%%%%%%%%%%%%%%%%%%%%%%%%%%%
\section{Hypergeometric Lattices} \label{sec:hgl} 
%%%%%%%%%%%%%%%%%%%%%%%%%%%%%%%%%%%%%%%%%%%%%%%%%%%%%%%%%%%%%%%%
A hypergeometric group $H = H(\varphi, \psi) = H(\Phi, \Psi) = 
\langle A, B \rangle$ is said to be {\sl integral} if  
%%%%%%
$$
\varphi(z) \in \bZ[z], \quad \psi(z) \in \bZ[z], \qquad 
\mbox{or equivalently} \quad \Phi(w) \in \bZ[w], \quad \Psi(w) \in \bZ[w].    
$$ 
%%%%%%
In this case let  $L = L(\varphi, \psi) = L(\Phi, \Psi)$ be the $\bZ$-linear 
span of $\br, A \br, \dots, A^{n-1} \br$.  
Then $L$ is a free $\bZ$-module with free basis 
$\br, A \br, \dots, A^{n-1} \br$, the $H$-invariant Hermitian form is 
$\bZ$-valued on $L$ and $H \subset O(L)$. 
It follows from \eqref{eqn:C2} that $B \br = - A \br$ and 
$B^k \br = - A^k \br + \mbox{$\bZ$-linear combination of 
$A \br, \dots, A^{k-1} \br$}$ for every $k \ge 2$, hence 
$\br, B \br, \dots, B^{n-1} \br$ also form a $\bZ$-linear basis of $L$. 
Thus we have the equation in \eqref{eqn:AB}, that is, 
%%%%%%%%
$$
L = \langle \br, A \br, \dots, A^{n-1} \br \rangle_{\bZ} = 
\langle \br, B \br, \dots, B^{n-1} \br \rangle_{\bZ} 
$$
%%%%%%%
By the $A$-invariance of the Hermitian form and the normalization 
$(\br, \br) = 2$ we have $(\bv, \bv) \in 2 \bZ$ for all $\bv \in L$. 
Thus $L$ equipped with the $H$-invariant form becomes an even lattice, 
called a {\it hypergeometric lattice}. 
We say that the group $H$ is {\sl unimodular} if so is the lattice $L$. 
%%%%%%%%%%%%%%%%%%%%%%%%% ss:unim %%%%%%%%%%%%%%%%%%%%%%%%%%%%%%%%
\subsection{Unimodularity} \label{ss:unim} 
%%%%%%%%%%%%%%%%%%%%%%%%%%%%%%%%%%%%%%%%%%%%%%%%%%%%%%%%%%%%%%%%
Formula \eqref{eqn:disc} in Theorem \ref{thm:ihf} implies that the 
hypergeometric lattice $L$ is unimodular if and only if 
$\Res(\varphi, \psi) = \pm 1$. 
Then \eqref{eqn:reste} shows that when $n$ is even this condition is 
equivalent to 
%%%%%%%%%%%%%%%%%%%%%%%%%% eqn:unim %%%%%%%%%%%%%%%%%%%%%%%%%%%%%
\begin{equation} \label{eqn:unim} 
\Psi(2) = \pm 1, \quad \Psi(-2) = \pm 1; \quad \Res(\Phi, \Psi) = \pm 1,  
\end{equation}
%%%%%%%%%%%%%%%%%%%%%%%%%%%%%%%%%%%%%%%%%%%%%%%%%%%%%%%%%%%%%%%%
while \eqref{eqn:resto} shows that when $n$ is odd $L$ cannot be unimodular. 
We say that $\Psi$ is {\sl unramified} if it satisfies the first and second 
conditions in \eqref{eqn:unim}. 
Note that there is no unramified polynomial of degree one.     
With irreducible decompositions $\Phi(w) = \prod_i \Phi_i(w)$ and 
$\Psi(w) = \prod_j \Psi_j(w)$ in $\bZ[w]$ the conditions \eqref{eqn:unim} 
can be restated as 
%%%%%%%%%%%%%%%%%%%%%%%%%% eqn:unim2 %%%%%%%%%%%%%%%%%%%%%%%%%%%%
\begin{equation} \label{eqn:unim2} 
\Psi_j(2) = \pm 1, \quad \Psi_j(-2) = \pm 1; \quad 
\Res(\Phi_i, \Psi_j) = \pm 1   \qquad \mbox{for all} \quad i, j.  
\tag{$\ref{eqn:unim}'$}
\end{equation}
%%%%%%%%%%%%%%%%%%%%%%%%%%%%%%%%%%%%%%%%%%%%%%%%%%%%%%%%%%%%%%%%
This observation provides us with a recipe to construct a unimodular 
hypergeometric lattice.  
%%%%%%%%%%%%%%%%%%%%%%%%%% rec:unim %%%%%%%%%%%%%%%%%%%%%%%%%%%%%%
\begin{recipe} \label{rec:unim} 
Given an integer $N \in \bZ_{\ge 2}$, find a finite set of monic irreducible 
polynomials $\Phi_i(w)$, $\Psi_j(w) \in \bZ[w]$ that satisfies the 
unimodularity condition \eqref{eqn:unim2} and the degree 
condition    
%%%%%%%%%%
$$
\sum_{i} \deg \Phi_i = N-1, \qquad \sum_{j} \deg \Psi_j = N. 
$$
%%%%%%%%%%
Take the products $\Phi(w) := \prod_i \Phi_i(w)$ and 
$\Psi(w) := \prod_j \Psi_j(w)$ and consider the integral hypergeometric group 
$H = H(\Phi, \Psi)$.  
Then the associated lattice $L$ is an even unimodular lattice of 
rank $n = 2N$. 
\end{recipe}
%%%%%%%%%%%%%%%%%%%%%%%%%%%%%%%%%%%%%%%%%%%%%%%%%%%%%%%%%%%%%%%%
\par
%%%%%
We are especially interested in the settings where the irreducible factors 
$\Phi_i(w)$ and $\Psi_j(w)$ are cyclotomic trace polynomials or 
Salem trace polynomials.    
%%%%%%%%%%%%%%%%%%%%%%%% ss:cyclo %%%%%%%%%%%%%%%%%%%%%%%%%%%%%%%%%
\subsection{Cyclotomic Trace Polynomials} \label{ss:cyclotomic} 
%%%%%%%%%%%%%%%%%%%%%%%%%%%%%%%%%%%%%%%%%%%%%%%%%%%%%%%%%%%%%%%%%
For $k \in \bZ_{\ge 3}$ let $\rC_k(z) \in \bZ[z]$ denote the $k$-th cyclotomic 
polynomial. 
It is a monic irreducible polynomial of degree $\phi(k)$, 
where $\phi(k)$ is Euler's totient function.  
For $k = 1$, $2$, in view of our purpose, it is convenient to employ an 
unconventional definition 
%%%%%%%
$$
\rC_1(z) := (z-1)^2, \quad \rC_2(z) := (z+1)^2, \qquad 
\phi(1) = \phi(2) = 2.  
$$
%%%%%
Note that $\rC_k(x) \ge 0$ for every $x \in \bR$.  
We say that $\rC_k(z)$ is {\sl unramified} if 
$\rC_k(\pm1) = 1$. 
%%%%%
\par
%%%%%
For any $k \in \bZ_{\ge 1}$ Euler's totient $\phi(k)$ is an even integer and 
$\rC_k(z)$ is a monic palindromic polynomial of degree $\phi(k)$, so there 
exists a unique monic polynomial $\CT_k(w) \in \bZ[w]$ of degree 
$\phi(k)/2$, called the $k$-th {\sl cyclotomic trace polynomial}, such that 
%%%%%% 
$$
\rC_k(z) = z^{\phi(k)/2} \, \CT_k(z+z^{-1}).  
$$
%%%%%%
Note that $\CT_1(w) = w-2$, $\CT_2(w) = w+2$ and 
$\CT_k(w)$ is irreducible for every $k \ge 1$. 
We say that $\CT_k(w)$ is {\sl unramified} if so is $\rC_k(z)$, in which 
case $\CT_k(2) = 1$ and $\CT_k(-2) = (-1)^{\phi(k)/2}$. 
%%%%%%
\par
%%%%%%
The following lemma is helpful in checking the unimodularity condition 
\eqref{eqn:unim2} for cyclotomic trace factors of $\Phi(w)$ and $\Psi(w)$. 
%%%%%%%%%%%%%%%%%%%%%%%%% lem:apostol %%%%%%%%%%%%%%%%%%%%%%%%%%%%%
\begin{lemma} \label{lem:apostol} 
Let $k$ and $m$ be positive integers such that $k > m$.  
\begin{enumerate}
\setlength{\itemsep}{-1pt}
\item $\Res(\CT_k, \CT_m) = \pm 1$ if and only if the ratio $k/m$ is 
not a prime power,  
\item $\CT_m(w)$ is unramified if and only if neither $m$ nor $m/2$ 
is a prime power,    
\end{enumerate}
where a prime power is an integer of the form $p^l$ with a prime 
number $p$ and a positive integer $l$.  
\end{lemma}   
%%%%%%%%%%%%%%%%%%%%%%%%% begin proof %%%%%%%%%%%%%%%%%%%%%%%%%%%%%
{\it Proof}. 
Apostol \cite[Theorems 1, 3, 4]{Ap} evaluates the resultants of cyclotomic 
polynomials:   
%%%%%%%%%%%%%%%%%%%%%%%% eqn:apostol %%%%%%%%%%%%%%%%%%%%%%%%%%%%%%
\begin{equation*} 
\Res(\rC_k, \rC_m) = 
\begin{cases} 
p^{\phi(m)} & \mbox{if $k/m$ is a power of a prime $p$}, \\[1mm]
1              & \mbox{otherwise},    
\end{cases}
\end{equation*}
%%%%%%%%%%%%%%%%%%%%%%%%%%%%%%%%%%%%%%%%%%%%%%%%%%%%%%%%%%%%%%%%%
for $k > m$. 
In particular $\rC_k(z)$ is unramified if and only if neither $k$ nor $k/2$ is 
a prime power. 
Lemma \ref{lem:apostol} then readily follows from the relation 
$\Res(\rC_k, \rC_m) = \Res(\CT_k, \CT_m)^2$. \hfill $\Box$ 
%%%%%%%%%%%%%%%%%%%%%%%%% end proof %%%%%%%%%%%%%%%%%%%%%%%%%%%%%%
%%%%%%%%%%%%%%%%%%%%%%%%% lem:CTP %%%%%%%%%%%%%%%%%%%%%%%%%%%%%%%
\begin{lemma} \label{lem:CTP} 
There are exactly $41$ cyclotomic trace polynomials $\CT_k(w)$ of degree 
$\le 10$, among which $15$ are unramified. 
They are given in Table $\ref{tab:CTP}$ with unramified ones being  
underlined. 
%%%%%%%%%%%%%%%%%%%% tab:CTP %%%%%%%%%%%%%%%%%%%%%%%%%%%%%%%%%%%%
\begin{table}[hh]
\centerline{
\begin{tabular}{c|l||c|l}
\hline
\\[-4mm]
$\deg$ & $k$  & $\deg$ & $k$ \\[1mm]
\hline
\\[-4mm]
$1$ & $1$, $2$, $3$, $4$, $6$ & $6$ & $13$, $\ul{21}$, $26$, $\ul{28}$, $\ul{36}$, $\ul{42}$ \\[1mm]
$2$ & $5$, $8$, $10$, $\ul{12}$ & $8$ & $17$, $32$, $34$, $\ul{40}$, $\ul{48}$, $\ul{60}$ \\[1mm]
$3$ & $7$, $9$, $14$, $18$ & $9$ & $19$, $27$, $38$, $54$ \\[1mm]
$4$ & $\ul{15}$, $16$, $\ul{20}$, $\ul{24}$, $\ul{30}$ & $10$ & $25$, $\ul{33}$, $\ul{44}$, $50$, $\ul{66}$ 
\\[1mm]
$5$ & $11$, $22$ &  &  \\[1mm]  
\hline
\end{tabular}} 
\caption{Cyclotomic trace polynomials $\CT_k(z)$ of degree $\le 10$.} 
\label{tab:CTP} 
\end{table} 
%%%%%%%%%%%%%%%%%%%%%%%%%%%%%%%%%%%%%%%%%%%%%%%%%%%%%%%%%%%%%  
\end{lemma}
%%%%%%%%%%%%%%%%%%%%%%%% begin proof %%%%%%%%%%%%%%%%%%%%%%%%%%%%%
{\it Proof}. 
Indeed, the cases $k = 1$ and $k = 2$ are trivial. 
For $k \ge 3$ Euler's totient $\phi(k)$ admits a lower bound   
%%%%%%
$$
\phi(k) > \phi_0(k) := 
\dfrac{k}{\re^{\gamma} \log \log k + \frac{3}{\log \log k}},  
$$
%%%%%%
where $\gamma$ is Euler's gamma constant (see \cite[Theorem 8.8.7]{BS} 
and \cite[Theorem 15]{RS}). 
Thus the condition $20 \ge \phi(k)$ implies $20 \ge \phi_0(k)$ and a 
careful analysis of the function $\phi_0(x)$ in real variable $x \ge 3$ 
shows that the latter condition holds exactly for $k = 3, \dots, 93$. 
A case-by-case check in this range gives all solutions $k$ to the bound 
$\deg \CT_k \le 10$ as in Table \ref{tab:CTP}.   
\hfill $\Box$ \par\medskip
%%%%%%%%%%%%%%%%%%%%% end proof %%%%%%%%%%%%%%%%%%%%%%%%%%%%%%%%%%
%%%%%%%%%%%%%%%%%%%%%%%%% ss:salem %%%%%%%%%%%%%%%%%%%%%%%%%%%%%%%
\subsection{Salem Trace Polynomials} \label{ss:salem}
%%%%%%%%%%%%%%%%%%%%%%%%%%%%%%%%%%%%%%%%%%%%%%%%%%%%%%%%%%%%%%%%
A Salem number is an algebraic unit $\lambda > 1$ whose conjugates other 
than $\lambda^{\pm 1}$ lie on the unit circle $S^1$ (see Salem \cite{Salem}).  
The monic minimal polynomial of a Salem number is called a Salem 
polynomial. 
A Salem polynomial is a palindromic polynomial of even degree.   
A Salem trace is an algebraic integer $\tau > 2$ whose other conjugates 
lie in the interval $[-2, \, 2]$. 
The monic minimal polynomial of a Salem trace is called a Salem trace 
polynomial. 
Salem numbers $\lambda$ and Salem traces $\tau$ are in one-to-one 
correspondence via the relation $\tau = \lambda + \lambda^{-1}$. 
A Salem polynomial $S(z)$ and the associated Salem trace polynomial 
$R(w)$ are related by 
%%%%%%%
$$
S(z) = z^{d} R(z+z^{-1}) \qquad \mbox{with} \quad d := \deg R(w). 
$$
%%%%%%%% 
We say that $S(z)$ and $R(w)$ are {\sl unramified} if 
$|S(\pm 1)| = 1$ and $|R(\pm 2)| = 1$ respectively. 
%%%%%%
%%%%%%%%%%%%%%%%%%%% ex:mcmullen %%%%%%%%%%%%%%%%%%%%%%%%%%%%%%%%
\begin{example} \label{ex:mcmullen} 
McMullen \cite[Table 4]{McMullen1} gives a list of ten unramified 
Salem numbers $\lambda_i$ of degree $22$ and the associated 
Salem polynomials $S_i(z)$ and Salem trace polynomials $R_i(w)$. 
Approximate values of $\lambda_i$ and exact formulas for $R_i(w)$ are 
given in Table \ref{tab:mcmullen}.  
%%%%%%%%%%%%%%%%%%%%%%%%%%% tab:mcmullen %%%%%%%%%%%%%%%%%%%%%%%%
\begin{table}[hh]
\centerline{
\begin{tabular}{c|l|l}
\hline
\\[-4mm]
$i$ & $\lambda_i$ & Salem trace polynomial $R_i(w)$ 
\\[1mm]
\hline
\\[-4mm]
$1$ & $1.37289$ & $w (w - 1) (w + 1)^2 (w^2 - 4) (w^5 - 6 w^3 + 8 w - 2) - 1$ 
\\[1mm] 
$2$ & $1.45099$ & $(w + 1) (w^2 - 4) (w^8 - 8 w^6 - w^5 + 19 w^4 + 3 w^3 - 12 w^2 + 1) - 1$ 
\\[1mm] 
$3$ & $1.48115$ & $w (w - 1) (w + 1)^2 (w^2 - 4) (w^5 - 6 w^3 - w^2 + 8 w + 1) - 1$ 
\\[1mm] 
$4$ & $1.52612$ & $w^2 (w + 1) (w^2 - 4) (w^2 + w - 1) (w^4 - w^3 - 5 w^2 + 3 w + 5) - 1$ 
\\[1mm] 
$5$ & $1.55377$ & $(w - 1) (w + 1)^2 (w^2 - 4) (w^6 - 7 w^4 - w^3 + 12 w^2 + 2 w - 1) - 1$ 
\\[1mm] 
$6$ & $1.60709$ & $w (w + 1)^2 (w^2 - 4) (w^2 + w - 1) (w^4 - 2 w^3 - 3 w^2 + 6 w - 1) - 1$ 
\\[1mm] 
$7$ & $1.6298$  & $(w + 1) (w^2 - 4) (w^8 - 8 w^6 - w^5 + 19 w^4 + 2 w^3 - 14 w^2 + 2) - 1$ 
\\[1mm] 
$8$ & $1.6458$ & $w (w + 1)^2 (w^2 - 4) (w^6 - w^5 - 7 w^4 + 5 w^3 + 13 w^2 - 6 w - 2) - 1$ 
\\[1mm] 
$9$ & $1.66566$ & $w (w + 1) (w^2 - 4) (w^7 - 9 w^5 - 2 w^4 + 25 w^3 + 9 w^2 - 20 w - 8) - 1$ 
\\[1mm] 
$10$ & $1.69496$ & $w (w + 1)^2 (w^2 - 3) (w^2 - 4) (w^4 - w^3 - 4 w^2 + 2 w + 1) - 1$ 
\\[1mm]
\hline 
\end{tabular}} 
\caption{Salem trace polynomials of degree $11$ from McMullen \cite[Table 4]{McMullen1}.}
\label{tab:mcmullen} 
\end{table}
%%%%%%%%%%%%%%%%%%%%%%%%%%%%%%%%%%%%%%%%%%%%%%%%%%%%%%%%%%%%%%%%
\end{example}
%%%%%%%%%%%%%%%%%%%%%%%%%%%%%%%%%%%%%%%%%%%%%%%%%%%%%%%%%%%%%%%%
%%%%%%%%%%%%%%%%%%%%%% ex:lehmer %%%%%%%%%%%%%%%%%%%%%%%%%%%%%%%%
\begin{example} \label{ex:lehmer} 
Lehmer's number $\lambda_{\rL} \approx 1.17628$ discovered in \cite{Lehmer} 
is the smallest Salem number ever known (see e.g. Hironaka \cite{Hironaka}).  
The associated Salem polynomial and Salem trace polynomial, that is,  
Lehmer's polynomial and Lehmer's trace polynomial are given by 
%%%%%%%%%%%%%%%%%%%%%% eqn:lehmer %%%%%%%%%%%%%%%%%%%%%%%%%%%%%%%
\begin{subequations} \label{eqn:lehmer}
\begin{align}
\rL(z) &= z^{10} + z^9-z^7-z^6-z^5-z^4-z^3+z+1, \label{eqn:lehmer1}
\\[1mm]  
\LT(w) &= (w+1)(w^2-1)(w^2-4)-1, \label{eqn:lehmer2}   
\end{align}
\end{subequations}
%%%%%%%%%%%%%%%%%%%%%%%%%%%%%%%%%%%%%%%%%%%%%%%%%%%%%%%%%%%%%%%%% 
respectively. 
Note that they are unramified. 
Lehmer's trace $\tau_{\rL} := \lambda_{\rL} + \lambda_{\rL}^{-1}$ is 
approximately $2.02642$. 
\end{example}
%%%%%%%%%%%%%%%%%%%%%%%%%%%%%%%%%%%%%%%%%%%%%%%%%%%%%%%%%%%%%%%
%%%%%%%%%%%%%%%%%%%%%%%%%% sec:K3 %%%%%%%%%%%%%%%%%%%%%%%%%%%%%%%
\section{Hypergeometric K3 Lattices} \label{sec:K3}
%%%%%%%%%%%%%%%%%%%%%%%%%%%%%%%%%%%%%%%%%%%%%%%%%%%%%%%%%%%%%%%% 
An even unimodular lattice of rank $22$ and signature $(3, 19)$, that is, 
index $-16$ is called a {\sl K3 lattice}. 
It is well known that the second cohomology group $H^2(X, \bZ)$ of a K3 
surface $X$ equipped with the intersection form is a K3 lattice. 
We wonder whether a K3 lattice can be realized as a hypergeometric lattice. 
%%%%%%%%%%%%%%%%%%%%%%%%%% def:hgk3l %%%%%%%%%%%%%%%%%%%%%%%%%%%%%
\begin{definition} \label{def:hgk3l}
A {\sl hypergeometric K3 lattice} is a unimodular hypergeometric lattice of 
rank $22$ and index $\pm 16$, where the index is calculated with respect 
to the invariant Hermitian form normalized by $(\br, \br) = 2$ as in 
\eqref{eqn:rr2}, which we call {\sl hypergeometric normalization}. 
In the context of K3 lattice we should employ another normalization that 
makes the index always equal to $-16$, which we call {\sl K3 normalization}. 
When the index is $+16$ in the former normalization we can switch to 
the latter one by negating the invariant Hermitian form; this amounts to 
taking the reversed normalization $(\br, \, \br) = -2$. 
The invariant Hermitian form with K3 normalization 
is often referred to as the {\sl intersection form}.           
\end{definition}
%%%%%%%%%%%%%%%%%%%%%%%%%%%%%%%%%%%%%%%%%%%%%%%%%%%%%%%%%%%%%%%%
%%%%%%%%%%%%%%%%%%%%%%%%%% ss:rank22 %%%%%%%%%%%%%%%%%%%%%%%%%%%%%
\subsection{Real of Rank $\mbox{\boldmath $22$}$ and Index 
$\mbox{\boldmath $\pm 16$}$} \label{ss:rank22}
%%%%%%%%%%%%%%%%%%%%%%%%%%%%%%%%%%%%%%%%%%%%%%%%%%%%%%%%%%%%%%%%
\par
%%%%%%%
Putting aside the integral structure and unimodularity condition for the 
moment we shall classify all real hypergeometric group of rank $22$ 
and index $\pm 16$. 
We employ the hypergeometric normalization in Definition \ref{def:hgk3l} 
and use the notations in \S\ref{sec:real}.   
We now have $\deg \Phi(w) = |\bA| = 10$ and $\deg \Psi(w) = |\bB| = 11$.
%%%%%%%%%%%%%%%%%%%%%%%%%% thm:pm16 %%%%%%%%%%%%%%%%%%%%%%%%%%%%
\begin{theorem} \label{thm:pm16} 
A real hypergeometric group of rank $n = 22$ has index $p-q = \pm 16$,  
if and only if $\bAin$ and $\bBon$ have any one of the configurations in 
Table $\ref{tab:classify}$, where in case $5$ we mean by ``doubles adjacent" 
that the unique double cluster in $\bAin$ and the unique double cluster in 
$\bBon$ must be adjacent to each other.  
Each of cases $6$, $7$, $8$ divides into two subcases as indicated 
in the constraints column.  
%%%%%%%%%%%%%%%%%%%%%%%%%% tab:classify %%%%%%%%%%%%%%%%%%%%%%%%%%%%
\begin{table}[hh]
\centerline{
\begin{tabular}{cccccccc}
\hline
\\[-4mm]
case & $s$ & $|\bAin|$ & $|\bBon|$ & $[\bAin]$ & $[\bBon]$ & constraints & $\ve(p-q)$ \\[1mm]
\hline
\\[-4mm]
$1$ & $8$ & $7$ & $8$ & $1^7$ & $1^8$ &  & $16$ \\[1mm]
$2$ & $8$ & $9$ & $8$ & $1^63^1$ & $1^8$ &  & $16$ \\[1mm]
$3$ & $8$ & $7$ & $10$ & $1^7$ & $1^73^1$ &  & $16$ \\[1mm]
$4$ & $8$ & $9$ & $10$ & $1^63^1$ & $1^73^1$ &  & $16$ \\[1mm]
$5$ & $9$ & $9$ & $10$ & $1^72^1$ & $1^82^1$ & doubles adjacent & $16$ \\[1mm]
$6$ & $9$ & $8$ & $10$ & $1^8$ & $1^82^1$ & $|\bB_{5 \pm 4}| = 2$ & $\pm 16$ \\[1mm] 
$7$ & $9$ & $10$ & $10$ & $1^73^1$ & $1^82^1$ & $|\bB_{5 \pm 4}| = 2$ & $\pm 16$ \\[1mm] 
$8$ & $10$ & $10$ & $10$ & $1^82^1$ & $1^{10}$ & $|\bA_{6 \pm 4}| = 2$ & $\pm 16$ \\[1mm] 
\hline  
\end{tabular}} 
\caption{Real of rank $22$ and index $\pm 16$, where  
$\ve = \pm$ is defined in \eqref{eqn:gam}.} 
\label{tab:classify} 
\end{table}
%%%%%%%%%%%%%%%%%%%%%%%%%%%%%%%%%%%%%%%%%%%%%%%%%%%%%%%%%%%%%%
\end{theorem}
%%%%%%%%%%%%%%%%%%%%%%% begin proof %%%%%%%%%%%%%%%%%%%%%%%%%%%%%%
{\it Proof}. 
We use Theorem \ref{thm:real} and Lemma \ref{lem:estim} with $n = 22$. 
Equation \eqref{eqn:rsgn2} yields $|1+\delta -2S| = |p-q| = 16$, 
which is the case if and only if   
%%%%%%%%%%%%%%%%%%%%%%%%%% eqn:dS %%%%%%%%%%%%%%%%%%%%%%%%%%%%%%%
\begin{equation} \label{eqn:dS}
(\delta, S) = (1, -7), (1, 9), (-1, \mp 8).  
\end{equation}
%%%%%%%%%%%%%%%%%%%%%%%%%%%%%%%%%%%%%%%%%%%%%%%%%%%%%%%%%%%%%%%%
If $|\bBon|$ is odd then 
$\delta = (-1)^{|\sbA_{\mathrm{in}}|}$ by \eqref{eqn:delta} and hence 
$|\bAin| \not\equiv S \bmod 2$ by \eqref{eqn:dS}, which contradicts 
the congruence $S \equiv |\bAin| \bmod 2$ in \eqref{eqn:estim}.  
Hence $|\bBon|$ must be even and $\delta = -(-1)^{|\sbA_{\mathrm{in}}|}$.  
It follows from \eqref{eqn:estim} that   
$8 \le |S| +1 \le s \le |\bBon| \le |\bB| = 11$, so that we have either 
%%%%%%%%%%%%%%%%%%%%%%%%% eqn:Bs %%%%%%%%%%%%%%%%%%%%%%%%%%%%%%%
\begin{equation} \label{eqn:Bs}
|\bBon| = 8, \quad s = 8; 
\qquad \mbox{or} \qquad 
|\bBon| = 10, \quad s = 8, 9, 10.   
\end{equation}  
%%%%%%%%%%%%%%%%%%%%%%%%%%%%%%%%%%%%%%%%%%%%%%%%%%%%%%%%%%%%%%%
A careful inspection shows that the only ten cases in Table \ref{tab:cases} 
can meet the constraints \eqref{eqn:estim}, \eqref{eqn:dS}, 
\eqref{eqn:Bs} and $|\bAin| \le |\bA| = 10$. 
Let us make a case-by-case treatment. 
In what follows ``cases" refer to those in Table \ref{tab:cases}.  
%%%%%%%%%%%%%%%%%%%%% tab:cases %%%%%%%%%%%%%%%%%%%%%%%%%%%%%%%%
\begin{table}[hh]
\centerline{
\begin{tabular}{cccccrrcc}
\hline
\\[-4mm]
case & $s$    & $|\rI|$ & $|\bAin|$ & $|\bBon|$ & $\delta$ & $S$ & $[\bAin]$ & $[\bBon]$ \\[1mm]
\hline
\\[-4mm] 
$1$  & $8$    & $7$ & $7$   & $8$   & $1$ & $-7$        & $1^7$      & $1^8$ \\[1mm]
$2$  & $8$    & $7$ & $9$   & $8$   & $1$ & $-7$        & $1^63^1$ & $1^8$ \\[1mm]
$3$  & $8$    & $7$ & $7$   & $10$ & $1$ & $-7$        & $1^7$      & $1^73^1$, $1^62^2$ \\[1mm]
$4$  & $8$    & $7$ & $9$   & $10$ & $1$ & $-7$        & $1^63^1$ & $1^73^1$, $1^62^2$ \\[1mm]
$5$  & $9$    & $7$ & $9$   & $10$ & $1$ & $-7$        & $1^72^1$ & $1^82^1$ \\[1mm]
$6$  & $9$    & $8$ & $8$   & $10$ & $-1$ & $\mp 8$ & $1^8$      & $1^82^1$ \\[1mm] 
$7$  & $9$    & $8$ & $10$ & $10$ & $-1$ & $\mp 8$ & $1^73^1$ & $1^82^1$ \\[1mm] 
$8$  & $10$  & $8$ & $10$ & $10$ & $-1$ & $\mp 8$ & $1^82^1$ & $1^{10}$ \\[1mm] 
$9$  & $10$  & $9$ & $9$   & $10$ & $1$ & $-7$        & $1^9$     & $1^{10}$ \\[1mm] 
$10$ & $10$ & $9$ & $9$   & $10$ & $1$ & $9$          & $1^9$     & $1^{10}$ \\[1mm]
\hline  
\end{tabular}}
\caption{Ten cases.}
\label{tab:cases}
\end{table} 
%%%%%%%%%%%%%%%%%%%%%%%%%%%%%%%%%%%%%%%%%%%%%%%%%%%%%%%%%%%%%%%% 
\par
%%%%%%%
In cases $1$--$4$ we have $\rI = \{2, \dots, 8\}$.  
In cases $1$ and $2$, since $\sigma_i$ is odd for every $i = 2, \dots, 8$,  
$S = -7$ is actually realized. 
The same is true in cases $3$ and $4$ with $[\bBon] = 1^73^1$.   
In cases $3$ and $4$ with $[\bBon] = 1^62^2$, let $|\bB_k| = |\bB_l| = 2$ 
with $1 \le k < l \le 8$.  
Up to reflection \eqref{eqn:refl} we may assume 
(i) $k = 1$ and $l = 8$; or (ii) $2 \le k < l \le 8$. 
In case (i) $\sigma_i$ is even for every $i = 2, \dots, 8$, so that 
$S = 7$ contradicting $S = -7$. 
In case (ii) $\sigma_i$ is odd for $2 \le i \le k$ or $l+1 \le i \le 8$ 
and even for $k+1 \le i \le l$, so $S = -(k-1)-(8-l) + (l-k) = 2(l-k)-7 \ge -5$, 
which again contradicts $S = -7$. 
Thus cases $3$ and $4$ with $[\bBon] = 1^62^2$ cannot occur.  
%%%%%
\par
%%%%%
In case $5$ let $|\bA_k| = |\bB_l| = 2$ with   
$2 \le k \le 9$ and $1 \le l \le 9$, in which case    
$\rI = \{ 2, \dots, \hat{k}, \dots, 9\}$, 
Up to reflection \eqref{eqn:refl} we may assume $l < k$.  
Then $\sigma_i$ is odd for $2 \le i \le l$ or $k+1 \le i \le 9$ 
and even for $l+1 \le i \le k-1$, so that $S = -(l-1)-(9-k)+(k-l-1) 
= 2(k-l)-9 = -7$ implies $l = k-1$. 
Taking reflection \eqref{eqn:refl} we have also $l = k$. 
Thus $S = -7$ is realized if and only if  
$|\bA_k| = |\bB_{k-1}| = 2$ or $|\bA_k| = |\bB_k| = 2$ holds 
for some $2 \le k \le 9$, that is, the double cluster 
in $\bAin$ and the one in $\bBon$ must be adjacent to each other. 
%%%%%%
\par
%%%%%%
In cases $6$ and $7$ we have $\rI = \{2, \dots, 9\}$. 
Let $|\bB_k| = 2$ with $1 \le k \le 9$. 
Up to reflection \eqref{eqn:refl} we may assume $5 \le k \le 9$. 
Then $\sigma_i$ is odd for $2 \le i \le k$ and even for 
$k + 1 \le i \le 9$, so that $S = -(k-1)+(9-k) = 2(5-k) = \mp 8$ 
implies $(k, S) = (9, -8)$. 
Taking reflection \eqref{eqn:refl} we have also $(k, S) = (1, 8)$. 
In summary $S = \mp 8$ forces $|\bB_{5 \pm 4}| = 2$. 
%%%%%%
\par
%%%%%%
In case $8$ let $|\bA_k| = 2$ and $\rI = \{2, \dots, \hat{k}, \dots, 10\}$ 
with $2 \le k \le 10$. 
Then $\sigma_i$ is odd for $2 \le i \le k-1$ and even for 
$k+1 \le i \le 10$, so $S = -(k-2) + (10-k) = 2(6-k) = \mp 8$, 
which implies $k = 6 \pm 4$, that is, $|\bA_{6 \pm 4}| = 2$.  
%%%%%%
\par
%%%%%%
In cases $9$ and $10$ we have $[\bBon] = 1^{10}$ and 
$\rI = \{2, \dots, 10\}$, so $\sigma_i$ is odd for $i = 2, \dots, 10$, 
and hence $S = -9$, i.e. neither $S = -7$ nor $S = 9$ occurs.  
Thus these cases cannot happen. 
All these observations lead to the classification in Table \ref{tab:classify}.  
\hfill $\Box$ \par\medskip 
%%%%%%%%%%%%%%%%%%%%%%%%%% end proof %%%%%%%%%%%%%%%%%%%%%%%%%%%%%
Some important information about $\bA$ and $\bB$ can be extracted  
from Table \ref{tab:classify}. 
Recall that an end cluster in $\bAon$ may be empty. 
If this is the case then it is called a {\sl null} cluster.  
%%%%%%%%%%%%%%%%%%%%%%%%%% lem:Aon %%%%%%%%%%%%%%%%%%%%%%%%%%%%
\begin{lemma} \label{lem:Aon} 
The following are valid for $\bA$ and the same is true with $\bB$.   
\begin{enumerate}
\setlength{\itemsep}{-1pt}
\item Any non-null $\bAon$-cluster is simple except for at most one cluster of 
multiplicity $2$ or $3$.  
\item We have $|\bAoff| \le 3$ and if $|\bAoff| \ge 2$ then any non-null 
$\bAon$-cluster is simple. 
\item Any element of $\bA$ is simple except for at most one element 
of multiplicity $2$ or $3$.     
\end{enumerate}
%%%%%%%%%%%%%%%%%%%%%%%%%%%%%%%%%%%%%%%%%%%%%%%%%%%%%%%%%%%%%%%%
\end{lemma}
%%%%%%%%%%%%%%%%%%%%%%%% begin proof %%%%%%%%%%%%%%%%%%%%%%%%%%%%%%
{\it Proof}. 
Table \ref{tab:classify} tells us $|\bAin| \ge 7$ and so 
$|\bA_1|+|\bA_{s+1}|+|\bAoff| = 10 - |\bAin| \le 3$, in particular 
$|\bAoff| \le 3$, which gives the first part of assertion (2).  
We have also $|\bA_1|+|\bA_{s+1}| \le 10 - |\bAin|$.  
Using this inequality we consider how $\bAin$ can be extended to $\bAon$ 
by adding $\bA_1$ and $\bA_{s+1}$.  
A careful inspection of all cases in Table \ref{tab:classify} leads to assertion (1).  
Suppose $|\bAoff| \ge 2$ and hence $|\bAin| \le 8$. 
Then we must be in one of the cases 1, 3, 6 in Table \ref{tab:classify}, 
where $\bAin$ consists of simple clusters only. 
Then $\bAon$ can also contain simple or null clusters only, 
since $|\bA_1|+|\bA_{s+1}| \le 1$. 
This proves the second part of assertion (2).  
By assertion (1) any element of $\bAon$ is simple except for 
at most one element which is of multiplicity $2$ or $3$.  
Moreover it follows from $|\bAoff| \le 3$ that $\bAoff$ can contain 
at most one multiple element, which is of multiplicity $2$ or $3$.  
If there is such an element then the second part of assertion (2) implies 
that there is no multiple element in $\bAon$.  
This proves assertion (3). 
As for $\bB$ assertions (1) and (2) can be seen directly from 
Table \ref{tab:classify} and then assertion (3) is proved  
just in the same manner as in the case of $\bA$.  
\hfill $\Box$ \par\medskip
%%%%%%%%%%%%%%%%%%%%%%%% end proof %%%%%%%%%%%%%%%%%%%%%%%%%%%%%%%
\par
%%%%%%
We now take the integral structure and unimodularity condition into account. 
%%%%%%%%%%%%%%%%%%%%%%%%%% lem:me %%%%%%%%%%%%%%%%%%%%%%%%%%%%%%%
\begin{lemma} \label{lem:me}  
If $L(\Phi, \Psi)$ is a hypergeometric K3 lattice then any root of $\Phi(w)$ 
is simple except for at most one integer root of multiplicity 
$2$ or $3$, whereas any root of $\Psi(w)$ is simple.  
\end{lemma} 
%%%%%%%%%%%%%%%%%%%%%% begin proof %%%%%%%%%%%%%%%%%%%%%%%%%%%%%%%
{\it Proof}. 
It follows from (3) of Lemma \ref{lem:Aon} that $\Phi(w)$ admits 
at most one multiple root. 
If $\Phi(w)$ actually contains such a root $\tau$, then any conjugate $\tau'$ 
to $\tau$ is also a multiple root, so uniqueness forces $\tau' = \tau$, 
which means that $\tau$ must be an integer. 
For the same reason any root of $\Psi(w)$ is simple except for at most 
one integer root, but unramifiedness of $\Psi(w)$ rules out this exception 
because there is no unramified polynomial of degree one. \hfill $\Box$  
%%%%%%%%%%%%%%%%%%%%%% end proof %%%%%%%%%%%%%%%%%%%%%%%%%%%%%%%%%
%%%%%%%%%%%%%%%%%%%%%%%%% rem:me %%%%%%%%%%%%%%%%%%%%%%%%%%%%%%%%
\begin{remark} \label{rem:me} 
If $\Phi(w)$ admits a multiple root $\tau \in \bZ$ then the last condition 
in \eqref{eqn:unim} yields $\Psi(\tau) = \pm 1$. 
This equation help us find the multiple root of $\Phi(w)$ if it exists.   
\end{remark}
%%%%%%%%%%%%%%%%%%%%%% ss:hodge %%%%%%%%%%%%%%%%%%%%%%%%%%%%%%%%%%
\subsection{Hodge Isometry and Special Eigenvalue} \label{ss:hodge} 
%%%%%%%%%%%%%%%%%%%%%%%%%%%%%%%%%%%%%%%%%%%%%%%%%%%%%%%%%%%%%%%% 
Let $L$ be a K3 lattice and $L_{\bC} := L \otimes \bC$ be its 
complexification equipped with the induced Hermitian form.   
A {\sl Hodge structure} on $L$ is an orthogonal decomposition 
%%%%%%%%%%%%%%%%%%%%%% eqn:h-s %%%%%%%%%%%%%%%%%%%%%%%%%%%%%%%%%% 
\begin{equation} \label{eqn:h-s} 
L_{\bC} = H^{2,0} \oplus H^{1,1} \oplus H^{0,2}
\end{equation} 
%%%%%%%%%%%%%%%%%%%%%%%%%%%%%%%%%%%%%%%%%%%%%%%%%%%%%%%%%%%%%%%%
of signatures $(1, 0) \oplus (1, 19) \oplus (1, 0)$ such that 
$\overline{H^{i, j}} = H^{j, i}$.  
We remark that $H^{1,1}_{\bR} := H^{1,1} \cap L_{\bR}$ with 
$L_{\bR} := L \otimes \bR$ is a real Lorentzian space of signature 
$(1, 19)$ and the set of time-like vectors $\cC := 
\{ \bv \in H^{1,1}_{\bR} \,:\, (\bv, \bv) > 0\}$ consists of two disjoint 
connected cones, one of which is referred to as the {\sl positive cone} 
$\cC^+$ and the other as the negative cone $\cC^- = - \cC^+$.  
%%%%%%
\par
%%%%%%
A {\sl Hodge isometry} is a lattice automorphism $F : L \to L$ preserving 
the Hodge structure \eqref{eqn:h-s}.  
We then have either $F(\cC^{\pm}) = \cC^{\pm}$ or $F(\cC^{\pm}) = 
\cC^{\mp}$, according to which $F$ is said to be {\sl positive} or 
{\sl negative}, respectively.   
For any positive Hodge isometry $F$ there is a trichotomy 
(see e.g. Cantat \cite[\S 1.4]{Cantat}):    
%%%%%%
\begin{enumerate}
\setlength{\itemsep}{-1pt}
\item[(E)] There exists a line in $\cC \cup \{ 0 \}$ preserved by $F$.  
In this case the line is fixed pointwise by $F$ and all eigenvalues of 
$F$ lie on $S^1$. 
\item[(P)] There exists a unique line in $\overline{\cC}$ preserved by $F$ 
and this line is on the light-cone $\partial \cC$.  
In this case the line is fixed pointwise by $F$ and all eigenvalues of 
$F$ lie on $S^1$.   
\item[(H)] There exists a real number $\lambda > 1$ such that 
$\lambda^{\pm 1}$ are the only eigenvalues of $F$ outside $S^1$. 
In this case the eigenvalues $\lambda^{\pm 1}$ are simple 
and their eigen-lines are on the light-cone $\partial \cC$.     
\end{enumerate}  
%%%%%%
In the respective cases $F$ is said to be of {\sl elliptic}, {\sl parabolic} 
or {\sl hyperbolic type}.  
Since $F$ preserves the free $\bZ$-module $L$, any eigenvalue of $F$ 
must be a root of unity in the elliptic and parabolic cases, whereas 
it is either a root of unity or a conjugate to a unique Salem number 
$\lambda > 1$ in the hyperbolic case.    
In any case $F$ restricted to $H^{1,1}$ has a real eigenvalue 
$\lambda \ge 1$.  
%%%%%%%%%%%%%%%%%%%%%%%% rem:positive %%%%%%%%%%%%%%%%%%%%%%%%%
\begin{remark} \label{rem:positive} 
A Hodge isometry $F$ is positive and of hyperbolic type, if  
$F|H^{1,1}$ has a real eigenvalue $\lambda > 1$.  
\end{remark}
%%%%%%%%%%%%%%%%%%%%%%%%%%%%%%%%%%%%%%%%%%%%%%%%%%%%%%%%%%%%
\par
%%%%%
Any Hodge isometry $F : L \to L$ admits a number $\delta \in S^1$ 
such that $F|H^{2,0} = \delta I$ and $F|H^{0,2} = \delta^{-1} I$, 
where $I$ is an identity map.  
Note that if $\delta = \pm 1$ then $F|H^{2,0} \oplus H^{0,2} = \pm I$. 
We refer to $\delta^{\pm1}$ and $\tau := \delta + \delta^{-1} \in [-2, \, 2]$ 
as the {\sl special eigenvalues} and the {\sl special trace} of $F$ respectively. 
This observation leads us to the following.  
%%%%%%%%%%%%%%%%%%%%%% def:special %%%%%%%%%%%%%%%%%%%%%%%%%%%%%%%
\begin{definition} \label{def:special} 
Let $F : L \to L$ be an automorphism of a K3 lattice $L$. 
An eigenvalue $\delta \in S^1$ of $F$ is said to be {\sl special},  
if $\delta \neq \pm 1$ and there exists a line $\ell \subset L_{\bC}$ such that 
$F|\ell = \delta I$ and the induced Hermitian form is positive definite on $\ell$; 
or if $\delta = \pm 1$ and there exists a plane $P \subset L_{\bC}$ such that 
$F|P = \pm I$, $\overline{P} = P$ and the induced Hermitian form is 
positive definite on $P$. 
We refer to $\tau := \delta + \delta^{-1} \in [-2, \, 2]$ as the {\sl special trace} 
of $F$.     
\end{definition}
%%%%%%%%%%%%%%%%%%%%%%%%%%%%%%%%%%%%%%%%%%%%%%%%%%%%%%%%%%%%%%%
%%%%%%%%%%%%%%%%%%%%%% rem:special %%%%%%%%%%%%%%%%%%%%%%%%%%%%%%%
\begin{remark} \label{rem:special} 
Since the K3 lattice $L$ has signature $(3, 19)$, the special trace $\tau$ 
of $F$ is unique if it exists. 
Moreover, if $\tau \neq \pm 2$ then the pair $(\delta, \ell)$ in Definition 
\ref{def:special} is uniquely determined by $F$ 
up to exchange of $(\delta, \ell)$ for $(\xi^{-1}, \bar{\ell})$. 
\end{remark}
%%%%%%%%%%%%%%%%%%%%%%%%%%%%%%%%%%%%%%%%%%%%%%%%%%%%%%%%%%%%%%%
\par
%%%%%
If a lattice automorphism $F : L \to L$ admits a special eigenvalue 
$\delta \neq \pm 1$ with associated line $\ell \subset L_{\bC}$, then  
%%%%%%%%%%%%%%%%%%%%%% eqn:ell %%%%%%%%%%%%%%%%%%%%%%%%%%%%%%%%%%
\begin{equation} \label{eqn:ell}
L_{\bC} = H^{2,0} \oplus H^{1,1} \oplus H^{0,2} := 
\ell \oplus (\ell \oplus \overline{\ell})^{\perp} \oplus \overline{\ell}
\end{equation}
%%%%%%%%%%%%%%%%%%%%%%%%%%%%%%%%%%%%%%%%%%%%%%%%%%%%%%%%%%%%%%
gives a unique Hodge structure up to complex conjugation such that 
$F$ is a Hodge isometry. 
Here the description of the case $\delta = \pm 1$ is omitted as it is not   
used in this article.  
%%%%%
\par
%%%%%
Let $L = L(A, B)$ be a hypergeometric K3 lattice in K3 normalization 
(see Definition \ref{def:hgk3l}).  
It is natural to ask if $A$ (or $B$) admits a special eigenvalue 
or equivalently a special trace.   
First we focus on the matrix $A$. 
Recall from \S \ref{ss:lcs} that $E(\lambda)$ stands for the generalized 
eigenspace of $A$ corresponding to an eigenvalue $\lambda \in \baon$ and 
$m(\lambda) = \dim E(\lambda)$ is the multiplicity of $\lambda$.  
Denote by $V(\lambda)$ the $\lambda$-eigenspace of $A$ 
in the narrow sense.  
%%%%%%%%%%%%%%%%%%%%%% lem:simple %%%%%%%%%%%%%%%%%%%%%%%%%%%%%%%%
\begin{lemma} \label{lem:simple} 
For any eigenvalue $\lambda \in \baon$ one has $\dim V(\lambda) = 1$ and 
if $m(\lambda) \ge 2$ then $V(\lambda)$ must be isotropic with respect 
to the invariant Hermitian form.  
In particular if $A$ admits a special eigenvalue $\delta$ then it is simple 
i.e. $m(\delta) = 1$ and different from $\pm 1$, so that the special trace 
$\tau := \delta + \delta^{-1}$ is different from $\pm 2$.      
\end{lemma}  
%%%%%%%%%%%%%%%%%%%%%% begin proof %%%%%%%%%%%%%%%%%%%%%%%%%%%%%%%
{\it Proof}. 
The last part of Theorem \ref{thm:ihf} says that $\br$ is a cyclic 
vector for the matrix $A$. 
Let $\bv$ be its projection down to $E(\lambda)$ relative to 
the direct sum decomposition \eqref{eqn:esd}. 
Then $\bv$ is a cyclic vector of $A|E(\lambda)$, so if we put 
$\bv_j := (A-\lambda I)^{j-1} \bv$ for $j \in \bZ_{\ge1}$ then     
$\bv_1, \dots, \bv_m$ with $m := m(\lambda)$ form a basis of $E(\lambda)$.   
Since $\bv_{m+1} = \b0$, we have $A \bv_m = \lambda \bv_m$ 
and hence $V(\lambda) = \bC \bv_m$ is $1$-dimensional.  
If $m \ge 2$, using $\bv_m = A \bv_{m-1} - \lambda \bv_{m-1}$ 
we have 
%%%%%%%%%%%%%
\begin{align*}
\lambda (\bv_m, \bv_m) &= (A \bv_m, A \bv_{m-1} - \lambda \bv_{m-1}) 
= (A \bv_m, A \bv_{m-1}) - (\lambda \bv_m, \lambda \bv_{m-1}) \\
&= (\bv_m, \bv_{m-1}) - |\lambda|^2 (\bv_m, \bv_{m-1}) = 0, 
\end{align*}
%%%%%%%%%%%%% 
by the $A$-invariance of the Hermitian form and $\lambda \in S^1$, 
so the Hermitian form vanishes on $V(\lambda)$. 
Let $\delta$ be a special eigenvalue of $A$. 
If $\delta = \pm 1$ then the associated plane $P$ in Definition 
\ref{def:special} must be contained in the line $V(\delta)$, 
but this is impossible. 
So we have $\delta \neq \pm 1$ and $\ell = V(\delta)$. 
Since the Hermitian form is positive-definite on $V(\delta)$, 
we must have $m(\delta) = 1$. 
It is evident that the special trace $\tau$ is different from $\pm 2$.  
\hfill $\Box$
%%%%%%%%%%%%%%%%%%%%%% end proof %%%%%%%%%%%%%%%%%%%%%%%%%%%%%%%%
%%%%%%%%%%%%%%%%%%%%%% rem:simple %%%%%%%%%%%%%%%%%%%%%%%%%%%%%%%
\begin{remark} \label{rem:simple} 
If $A$ admits a special trace then the recipe \eqref{eqn:ell} yields a unique 
Hodge structure on $L$ up to complex conjugation such that $A$ is a 
Hodge isometry. 
Since $\pm A$ induce the same Hodge structure, it makes sense 
to speak of $A$ being a positive or negative Hodge isometry.  
The hypergeometric group $H = \langle A, B \rangle$ can then be normalized 
so that $A$ is positive, otherwise by replacing $H$ with its antipode 
$H^{\ra}$ mentioned in Remark \ref{rem:antip}.    
Lemma \ref{lem:simple} and these remarks are also valid for the matrix $B$. 
\end{remark}
%%%%%%%%%%%%%%%%%%%%%%%%%%%%%%%%%%%%%%%%%%%%%%%%%%%%%%%%%%%%%%%%%
%%%%%%%%%%%%%%%%%%%%%% ss:st %%%%%%%%%%%%%%%%%%%%%%%%%%%%%%%%%%%%
\subsection{Determination of Special Trace} \label{ss:st} 
%%%%%%%%%%%%%%%%%%%%%%%%%%%%%%%%%%%%%%%%%%%%%%%%%%%%%%%%%%%%%%%%%
To consider when a special trace exists and to determine it explicitly,   
we begin with the following.    
%%%%%%%%%%%%%%%%%%%%%%%% lem:endc %%%%%%%%%%%%%%%%%%%%%%%%%%%%%%%
\begin{lemma} \label{lem:endc} 
Suppose that $A$ admits a special trace and becomes a positive Hodge 
isometry with respect to the Hodge structure \eqref{eqn:ell}. 
Then the pair $(M(2), \, \idx(1))$ is $(0, 1)$ in the elliptic case; 
$(1, -1)$ in the parabolic case; and $(0, -1)$ in the hyperbolic case 
respectively, while one has $M(-2) = 0$ and $\idx(-1) = -1$ in every case, 
where the index is taken in K3 normalization in Definition $\ref{def:hgk3l}$. 
Moreover any root of $\Phi(w)$ is simple.   
\end{lemma}
%%%%%%%%%%%%%%%%%%%%%% begin proof %%%%%%%%%%%%%%%%%%%%%%%%%%%%%%
{\it Proof}. 
Let $\delta^{\pm1}$ be the special eigenvalues of $A$. 
Consider the associated Hodge structure \eqref{eqn:ell} and 
use Lemma \ref{lem:simple} repeatedly. 
Since $\delta \neq \pm 1$ we have $V(\pm1) \subset E(\pm1) \subset 
H^{1,1}$. 
Let $V_{\bR}(\pm 1) := V(\pm 1) \cap L_{\bR} \subset H^{1,1}_{\bR}$.  
If the line $V_{\bR}(-1)$ lies in $\overline{\cC}$ then $A$ sends 
$\cC^{\pm}$ to $\cC^{\mp}$ as $A = -I$ on $V_{\bR}(-1)$. 
This contradicts the positivity of $A$, so   
$V_{\bR}(-1)$ must be in the space-like region.    
Thus the Hermitian form $h$ is negative-definite on $V(-1)$ and we must 
have $m(-1) = 2 M(-2) + 1 = 1$, i.e. $M(-2) = 0$ and $\idx(-1) = -1$. 
In the elliptic and parabolic cases $V_{\bR}(1)$ is the unique line 
in $\overline{\cC}$ preserved by $A$. 
Accordingly, in the elliptic case $h$ is positive definite on $V(1)$, 
so we have $m(1) = 2 M(2) + 1 = 1$, i.e. $M(2) = 0$ and $\idx(1) = 1$. 
Similarly, in the parabolic case $h$ vanishes on $V(1)$ and has signature 
$(u, v) = (M(2), M(2)+1)$ or $(M(2)+1, M(2))$ on $E(1)$,  
which forces $M(2) = 1$ and $(u, v) = (1, 2)$, i.e. 
$\idx(1) = -1$, since $E(1) \subset H^{1,1}$ and 
$H^{1,1}$ has signature $(1, 19)$. 
In the hyperbolic case $h$ must be negative-definite on 
$E(1)$ because $h$ has signature $(1, 1)$ on 
$E(\lambda) \oplus E(\lambda^{-1}) \subset H^{1,1}$ where $\lambda$ 
is the real eigenvalue of $A$ strictly greater than $1$. 
This makes $M(2) = 0$ and $\idx(1) = -1$.    
Let $U$ be $E(1)$ in the elliptic or parabolic case and 
$E(\lambda) \oplus E(\lambda^{-1})$ in the hyperbolic case.    
Then $h$ is negative-definite on $U^{\perp} \cap H^{1,1}$. 
Hence by Propositions \ref{prop:lsgn} and \ref{prop:lsgn-r} 
any root of $\Phi(w)$ other than $\tau = \delta + \delta^{-1}$ and 
$2$ is simple. 
Simpleness of $\tau$ also follows from the same propositions, 
while if $2$ is a root then it is simple because $M(2) = 1$ as just 
shown.         
\hfill $\Box$ \par\medskip
%%%%%%%%%%%%%%%%%%%%%% end proof %%%%%%%%%%%%%%%%%%%%%%%%%%%%%%%%
We are now in a position to prove Theorems \ref{thm:main2} and 
\ref{thm:main3} stated in the Introduction.  
%%%%%%%%%%%%%%%%%%%%
\par\medskip
%%%%%%%%%%%%%%%%%%%%% begin proof %%%%%%%%%%%%%%%%%%%%%%%%%%%%%%%
{\it Proof of Theorem $\ref{thm:main2}$}.  
Suppose that $A$ has a special trace $\tau$ and is a positive Hodge isometry.     
The last part of Lemma \ref{lem:endc} implies  
$\cA_1^{\circ} = \bA_1^{\circ}$, $\cA_{s+1}^{\circ} = \bA_{s+1}^{\circ}$ 
and $\cAin = \bAin$ in the notation of Theorem \ref{thm:lc} . 
A careful inspection of Tables \ref{tab:classify} and \ref{tab:cases} allows us 
to rewrite \eqref{eqn:lc} in K3 normalization (see Definition \ref{def:hgk3l}).      
Formulas \eqref{eqn:lc1} now read    
%%%%%%%%%%%%%%%
\begin{alignat*}{3}
\idx(1) &= - (-1)^{|\sbA_1^{\circ}|}, \qquad & 
\idx(-1) &= - (-1)^{|\sbA_{s+1}^{\circ}|}  & 
\qquad & \mbox{in cases $1$--$5$}, \\
\idx(1) &= \mp (-1)^{|\sbA_1^{\circ}|}, \qquad & 
\idx(-1) &= \pm (-1)^{|\sbA_{s+1}^{\circ}|} & 
\qquad & \mbox{in cases $6$--$8$}     
\end{alignat*} 
%%%%%%%%%%%%%%% 
of Table \ref{tab:classify}. 
These equations are combined with Lemma \ref{lem:endc} to determine 
the parities of $|\bA_1^{\circ}|$ and $|\bA_{s+1}^{\circ}|$. 
Then $\Idx(\bA_1^{\circ})$ and 
$\Idx( \bA_{s+1}^{\circ})$ are evaluated by the formulas 
\eqref{eqn:lc2}, which now look like    
%%%%%%%%%%%%%%%
\begin{alignat*}{3}
\Idx(\bA_1^{\circ}) &= - \Par(|\bA_1^{\circ}|), \qquad & 
\Idx( \bA_{s+1}^{\circ}) &= - \Par(|\bA_{s+1}^{\circ}|) & \qquad &
\mbox{in cases $1$--$5$}, \\ 
\Idx(\bA_1^{\circ}) &= \mp \Par(|\bA_1^{\circ}|),  \qquad & 
\Idx( \bA_{s+1}^{\circ}) &= \pm \Par(|\bA_{s+1}^{\circ}|)  & \qquad &
\mbox{in cases $6$--$8$}     
\end{alignat*}
%%%%%%%%%%%%%% 
of Table \ref{tab:classify}. 
The information so obtained and Lemma \ref{lem:endc} confine 
the possibilities of $\bA_1$ and $\bA_{s+1}$. 
In particular, subcase $|\bB_9|=2$ of case $7$ and subcase 
$|\bA_{10}|=10$ of case $8$ are excluded in the elliptic case, whereas 
cases 7 and 8 are altogether ruled out in the parabolic and hyperbolic cases. 
%%%%%
\par
%%%%%    
In cases $1$--$5$ of Table \ref{tab:classify} the first formula in 
\eqref{eqn:lc3} reads $\Idx(\bAin) = - 7$.  
In cases $1$ and $3$ where $|\bAin| = 7$, this implies that all elements 
of $\bAin$ have local index $-1$, so the special root $\tau$ must lie 
either in $\bA_1$ or in $\bA_{s+1}$. 
From the data we have already had, it is easy to know which of them 
contains $\tau$.  
In cases $2$, $4$, $5$ where $|\bAin| = 9$, there is a unique element 
of $\bAin$ with local index $1$ (which is $\tau$)  and all the 
other elements have local index $-1$.  
Cases $6$--$8$ can be treated in a similar manner, where the first 
formula in \eqref{eqn:lc3} reads $\Idx(\bAin) = - 8$ with $|\bAin|$ 
being either $8$ or $10$. 
In any case, once the configuration of the clusters is fixed, formula 
\eqref{eqn:lsgn1r} tells us exactly where $\tau$ is located. 
In particular $\tau$ lies in the unique multiple cluster, if it exists, 
and $\tau$ is its middle element if it is triple.   
Exhausting all possibilities we have the assertions of 
Theorem \ref{thm:main2}.  \hfill $\Box$ 
%%%%%%%%%%%%%%%%%%%%% end proof %%%%%%%%%%%%%%%%%%%%%%%%%%%%%%%%
\par\medskip
%%%%%%%%%%%%%%%%%%%%%% begin proof %%%%%%%%%%%%%%%%%%%%%%%%%%%%%%%
{\it Proof of Theorem $\ref{thm:main3}$}. 
From Table \ref{tab:classify} we have $|\bBoff| = 3$ in cases $1$--$2$, 
while $|\bBoff| = 1$ in cases $3$--$8$. 
On the other hand, if $B$ admits a special trace $\tau$ then we have 
$|\bBoff| = 0$ in the elliptic or parabolic case, while $|\bBoff| = 1$ in 
the hyperbolic case. 
Thus cases $1$--$2$ are ruled out and we are exclusively in the 
hyperbolic case. 
Since $B$ is a positive Hodge isometry, we have $|\bBoff| = |\bB_{> 2}| = 1$.  
Assertions on simpleness of roots of $\Phi(w)$ and $\Psi(w)$  
follow from Lemma \ref{lem:me}. 
With $\cBon = \bBon$ the second formula in \eqref{eqn:lc3} adapted in 
K3 normalization reads $\Idx(\bBon) = -8$. 
Since $|\bBon|=10$ in cases $3$--$8$, this implies that 
there exists a unique element of $\bBon$ with local index $1$ 
(which is $\tau$) and all the other elements have index $-1$.  
The location of $\tau$ can be determined by formula \eqref{eqn:lsgn1r}. 
Exhausting all possibilities we have the assertions of 
Theorem \ref{thm:main3}.    
\hfill $\Box$   
%%%%%%%%%%%%%%%%%%%%%% end proof %%%%%%%%%%%%%%%%%%%%%%%%%%%%%%%%%
%%%%%%%%%%%%%%%%%%%%%% sec:k3str %%%%%%%%%%%%%%%%%%%%%%%%%%%%%%%%%
\section{K3 Structure} \label{sec:k3str} 
%%%%%%%%%%%%%%%%%%%%%%%%%%%%%%%%%%%%%%%%%%%%%%%%%%%%%%%%%%%%%%%%
Given a K3 lattice $L$ with a Hodge structure \eqref{eqn:h-s}, 
the Picard lattice and the root system are defined by 
$\Pic := L \cap H^{1,1}$ and  
$\vD := \{ \bu \in \Pic : (\bu, \bu) = -2 \}$ respectively.   
Given a positive cone $\cC^+ \subset \cC$, a {\sl set of positive roots} is 
by definition a subset $\vD^+ \subset \vD$ such that 
$\vD = \vD^+ \amalg (- \vD^+)$ and 
%%%%%%%%%%%%%%
\begin{equation*} 
\cK := \{ \bv \in \cC^+ : (\bv, \bu) > 0, \, 
{}^{\forall} \bu \in \vD^+ \}   
\end{equation*} 
%%%%%%%%%%%%%%
is nonempty, in which case $\cK$ is called the 
{\sl K\"ahler cone} associated with $\vD^+$. 
Note that $\vD^+$ determines a basis $\vD_{\rb}$ 
of the root system $\vD$ and vice versa. 
A {\sl Picard-Lefschetz reflection} is a lattice automorphism 
$\sigma_{\sbu} : L \to L$ defined by $\sigma_{\sbu}(\bv) := 
\bv + (\bv, \bu) \bu$ for each $\bu \in \vD^+$. 
The group $W := \langle \sigma_{\sbu} \,|\, \bu \in \vD^+ \rangle$ 
generated by those reflections is referred to as the {\sl Weyl group}. 
It acts on $\cC^+$ properly discontinuously and the closure in 
$\cC^+$ of the K\"ahler cone $\cK$ is a fundamental domain of 
this action. 
Each fundamental domain is called a {\sl Weyl chamber}.      
%%%%%%%%%%%%%%%%%%%%%% def:k3str %%%%%%%%%%%%%%%%%%%%%%%%%%%%%%%%%%%%
\begin{definition} \label{def:k3str} 
A {\sl K3 structure} on a K3 lattice $L$ is a specification of (i) a Hodge structure, 
(ii) a positive cone $\cC^+$ and (iii) a set of positive roots $\vD^+$.  
These data determines the associated K\"{a}hler cone $\cK$.  
Alternatively we can specify (i), (ii) and $(\mathrm{iii}')$ a Weyl chamber 
$\cK \subset \cC^+$ which we call the K\"ahler cone, in place of (iii).  
In this case the set of positive roots $\vD^+$ is associated 
afterwards.  
See McMullen \cite[\S 6]{McMullen3}.  
\end{definition}
%%%%%%%%%%%%%%%%%%%%%%%%%%%%%%%%%%%%%%%%%%%%%%%%%%%%%%%%%%%%%%%%%%
%%%%%%%%%%%%%%%%%%%%%% ss:?½@%%%%%%%%%%%%%%%%%%%%%%%%%%%% 
\subsection{Synthesis of Automorphisms}
%%%%%%%%%%%%%%%%%%%%%%%%%%%%%%%%%%%%%%%%%%%%%%%%%%%%%%%%%%%%%%%%%%
Any K3 surface $X$ induces a K3 structure on $H^2(X, \bZ)$; 
its Hodge structure is given by the Hodge-Kodaira decomposition of 
$H^2(X, \bC)$; its K\"ahler cone $\cK(X)$ is the set of all 
K\"ahler classes on $X$; its positive cone $\cC^+(X)$ is the 
connected component containing $\cK(X)$; the set of positive 
roots $\vD^+(X)$ is that of all effective $(-2)$-classes.      
Any automorphism $f : X \to X$ induces a lattice 
automorphism $f^* : H^2(X, \bZ) \to H^2(X, \bZ)$ that   
preserves the K3 structure. 
Conversely, the Torelli theorem and surjectivity of the period 
mapping tells us the following. 
%%%%%%%%%%%%%%%%%%%%%%%%% thm:great %%%%%%%%%%%%%%%%%%%%%%%%%%%%%%%%%
\begin{theorem} \label{thm:great} 
Let $L$ be a K3 lattice. 
Any K3 structure on $L$ is realized by a unique marked K3 surface 
$(X, \iota)$ up to isomorphism and any lattice automorphism 
$F : L \to L$  preserving the K3 structure is realized by a 
unique K3 surface automorphism $f : X \to X$ up to biholomorphic conjugacy, 
that is, there exists a commutative diagram:   
%%%%%%%%%%%%%%%%%%%%%%%%% cd:iF %%%%%%%%%%%%%%%%%%%%%%%%%%%%%%%%%
\begin{equation} \label{cd:iF}
\begin{CD}
H^2(X, \bZ) @> \iota >> L \\
@V f^* VV    @VV F V \\
H^2(X, \bZ) @> \iota >> L.  
\end{CD}
\end{equation} 
%%%%%%%%%%%%%%%%%%%%%%%%%%%%%%%%%%%%%%%%%%%%%%%%%%%%%%%%%%%%%%
\end{theorem}
%%%%%%%%%%%%%%%%%%%%%%%%%%%%%%%%%%%%%%%%%%%%%%%%%%%%%%%%%%%%%%
\par
%%%%%%% 
The K3 surface automorphism $f : X \to X$ so obtained is called the 
{\sl lift} of $F$.   
Thus constructing a K3 surface automorphism amounts to constructing 
an automorphism of a K3 lattice preserving a K3 structure. 
There may be curves in the surface $X$ and it will later be necessary 
to see how irreducible curves are permuted by $f$.     
%%%%%%%%%%%%%%%%%%%%%%%% lem:curve %%%%%%%%%%%%%%%%%%%%%%%%%%%%%%%
\begin{lemma} \label{lem:curve}
In Theorem $\ref{thm:great}$ suppose that the intersection form is 
negative-definite on $\Pic$.   
Then any irreducible curve $C \subset X$ is a $(-2)$-curve and hence the 
unique effective divisor representing the nodal class $[C] \in \Pic(X)$. 
Through the marking isomorphism $\iota$ in \eqref{cd:iF} the $(-2)$-curves 
in $X$ are in one-to-one correspondence with the simple roots in $\vD^+$, 
that is, the elements of $\vD_{\rb} \subset L$ and how these curves are 
permuted by $f$ is faithfully represented by the action of $F$ on $\vD_{\rb}$  
or equivalently on its Dynkin diagram.       
\end{lemma} 
%%%%%%%%%%%%%%%%%%%%%%%%% begin proof %%%%%%%%%%%%%%%%%%%%%%%%%%%%%%%
{\it Proof}. 
By assumption $\Pic \subset L$ is negative-definite and even,  
so is $\Pic(X) \subset H^2(X, \bZ)$.   
Since any irreducible curve $C$ represents a non-trivial effective 
class $[C] \in \Pic(X)$, its self-intersection number $(C, C)$ is an even 
negative integer. 
But the irreducibility of $C$ forces $(C, C) \ge -2$ and so  
$(C, C) = -2$, which implies that $C$ is a $(-2)$-curve. 
The remaining assertions follow from 
Barth et al. \cite[Chap. VI\!I\!I, (3.7) Proposition]{BHPV}.  
\hfill $\Box$ \par\medskip 
%%%%%%%%%%%%%%%%%%%%%%%%% end proof %%%%%%%%%%%%%%%%%%%%%%%%%%%%%%%%
The union of all $(-2)$-curves in $X$ is called the {\sl exceptional set} and is 
denoted by $\cE = \cE(X)$.  
The dual graph of its intersection relations is represented by 
the Dynkin diagram mentioned in Lemma \ref{lem:curve}.   
%%%%%%%
\par
%%%%%%%
If $F : L \to L$ is an automorphism of a K3 lattice $L$ admitting a special 
trace, then $F$ is a Hodge isometry with respect to the Hodge structure 
\eqref{eqn:ell}. 
We may assume that $F$ is a positive Hodge isometry, otherwise 
by replacing $F$ with $-F$.  
For any K3 structure having \eqref{eqn:ell} as its Hodge structure,  
$F$ obviously preserves the Hodge structure and the positive cone $\cC^+$, 
but it is not always true that $F$ preserves the K\"{a}hler cone $\cK$; it 
may be sent to a different Weyl chamber $F(\cK)$.  
However, there exists a unique element $w_F \in W$ that brings $F(\cK)$ 
back to $\cK$, since the Weyl group $W$ acts simply transitively on the set 
of Weyl chambers.    
The modified map 
%%%%%%%%%%%%%%%%%%%%%%%% eqn:modified %%%%%%%%%%%%%%%%%%%%%%%%%%%%%%%%
\begin{equation} \label{eqn:modified}
\tilde{F} := w_F \circ F 
\end{equation} 
%%%%%%%%%%%%%%%%%%%%%%%%%%%%%%%%%%%%%%%%%%%%%%%%%%%%%%%%%%%%%%%%%%%
preserves the K3 structure and hence lifts to a K3 surface automorphism via 
Theorem \ref{thm:great}. 
In \S \ref{ss:algorithm} we shall apply this scenario to a hypergeometric K3 
lattice with the map $F$ being the matrix $A$ or $B$.   
%%%%%%%%%%%%%%%%%%%%%%%%% ss:PL %%%%%%%%%%%%%%%%%%%%%%%%%%%%%%%%%%%%
\subsection{Picard Lattice}  \label{ss:PL}
%%%%%%%%%%%%%%%%%%%%%%%%%%%%%%%%%%%%%%%%%%%%%%%%%%%%%%%%%%%%%%%%%%%
Let $L$ be a K3 lattice with Hodge structure \eqref{eqn:h-s} and $F : L \to L$ 
be a Hodge isometry with special eigenvalues $\delta(F)^{\pm1}$.   
If $\chi_0(z)$ is their minimal polynomial then the characteristic polynomial 
$\chi(z)$ of $F$ factors as 
%%%%%%%%%%%%%%%%%%%%%%%%% eqn:chi01 %%%%%%%%%%%%%%%%%%%%%%%%%%%%%%%%
\begin{equation} \label{eqn:chi01}
\chi(z) = \chi_0(z) \cdot \chi_1(z)
\end{equation}
%%%%%%%%%%%%%%%%%%%%%%%%%%%%%%%%%%%%%%%%%%%%%%%%%%%%%%%%%%%%%%%% 
for some monic polynomial $\chi_1(z) \in \bZ[z]$.   
Consider the Picard lattice $\Pic := L \cap H^{1,1}$. 
It is projective or not depending on whether it contains a vector of 
positive self-intersection or not.       
%%%%%%%%%%%%%%%%%%%%%%%%% thm:pic %%%%%%%%%%%%%%%%%%%%%%%%%%%%%%%%%%
\begin{theorem} \label{thm:pic} 
Suppose that $\delta(F)$ is different from $\pm 1$ and that 
$\chi_0(z)$ and $\chi_1(z)$ are coprime.  
Then 
%%%%%%%%%%%%%%%%%%%%%%%% eqn:chi1 %%%%%%%%%%%%%%%%%%%%%%%%%%%%%%%%%
\begin{equation} \label{eqn:chi1}
\Pic = \{ \bv \in L : \chi_1(F) \bv = \b0 \}, 
\end{equation} 
%%%%%%%%%%%%%%%%%%%%%%%%%%%%%%%%%%%%%%%%%%%%%%%%%%%%%%%%%%%%%%%%
the intersection form is non-degenerate on $\Pic$ and the Picard number 
$\rho := \rank \, \Pic$ is given by $\rho = 22 - \deg \chi_0(z)$, 
which is necessarily even and $\le 20$. 
In the projective case $\chi_0(z)$ is a cyclotomic polynomial of 
degree $\ge 2$ and $\Pic$ has signature $(1, \rho-1)$. 
In the non-projective case $\chi_0(z)$ is a Salem polynomial of degree 
$\ge 4$, the associated Salem number gives the spectral radius 
$\lambda(F)$ and $\Pic$ has signature $(0, \rho)$. 
\end{theorem}
%%%%%%%%%%%%%%%%%%%% begin proof %%%%%%%%%%%%%%%%%%%%%%%%%%%%%%%%%%
{\it Proof}. 
There exist orthogonal direct sum decompositions 
%%%%%%%%%%%%%%%%%%%% eqn:LQ %%%%%%%%%%%%%%%%%%%%%%%%%%%%%%%%%%%%%
\begin{equation} \label{eqn:LQ}
L_{K} = L_0(K) \oplus L_1(K) \qquad \mbox{for} \quad 
K = \bQ, \, \overline{\bQ}, \, \bC,   
\end{equation}
%%%%%%%%%%%%%%%%%%%%%%%%%%%%%%%%%%%%%%%%%%%%%%%%%%%%%%%%%%%%%%%% 
where $L_K := L \otimes K$ and 
$L_i(K) := \{ \bv \in L_K : \chi_i(F) \bv = \b0 \}$ for $i = 0, 1$.    
Indeed, since $\chi_0(z)$ and $\chi_1(z)$ are coprime, we have  
$\phi_0(z) \cdot \chi_0(z) + \phi_1(z) \cdot \chi_1(z) = d$ for some 
$d \in \bZ_{\ge 1}$ and $\phi_0(z)$, $\phi_1(z) \in \bZ[z]$, so 
the orthogonal projection $P_i : L_{\bQ} \to L_i(\bQ)$ is given by 
$P_i := d^{-1} \phi_{1-i}(F) \cdot \chi_{1-i}(F)$ for $i = 0, 1$.  
This yields \eqref{eqn:LQ} for $K = \bQ$, which in turn gives  
\eqref{eqn:LQ} for $K = \overline{\bQ}, \, \bC$.   
Let $\bu \in L_0(\overline{\bQ})$ be a $\delta(F)$-eigenvector of $F$. 
Since $\bu \in H^{2, 0}$ and $\Pic \subset H^{1,1}$, we have 
$\bu \perp \Pic$ and so $\sigma(\bu) \perp \Pic$ for any  
$\sigma \in G := \Gal(\overline{\bQ}/\bQ)$. 
This implies $L_0(\overline{\bQ}) \perp \Pic$ and hence 
$\Pic \subset L_1(\overline{\bQ}) \subset H^{1,1}$ by 
\eqref{eqn:LQ} for $K = \overline{\bQ}$, because the vectors 
$\{ \sigma(\bu) \}_{\sigma \in G}$ span $L_0(\overline{\bQ})$. 
Therefore $\Pic = L \cap L_1(\overline{\bQ}) = 
\{ \bv \in L : \chi_1(F) \bv = \b0 \}$, which 
yields \eqref{eqn:chi1}. 
Moreover we have 
%%%%%%%%%%%%%%%%%%% eqn:L0L1 %%%%%%%%%%%%%%%%%%%%%%%%%%%%%%%% 
\begin{equation} \label{eqn:L0L1}
H^{2,0} \oplus H^{0,2} \subset L_0(\bC), \qquad \Pic_{\bQ} = L_1(\bQ), 
\end{equation}
%%%%%%%%%%%%%%%%%%%%%%%%%%%%%%%%%%%%%%%%%%%%%%%%%%%%%%%%%%%
where $\Pic_{\bQ} := \Pic \otimes \bQ$. 
Non-degeneracy of $\Pic$ and $\rho = 22 - \deg \chi_0(z)$ follow 
from \eqref{eqn:LQ}, the latter part of \eqref{eqn:L0L1}, 
non-degeneracy of $L$ and $\rank \, L = 22$.  
Recall that $L_{\bC}$ and $H^{2,0} \oplus H^{0,2}$ have signatures 
$(3, 19)$ and $(2, 0)$ respectively. 
Thus the former part of \eqref{eqn:L0L1} implies that the total signature 
$(3, 19)$ decomposes into  either (i) $(2, 20-\rho) \oplus (1, \rho-1)$; 
or (ii) $(3, 19-\rho) \oplus (0, \rho)$, along the orthogonal decomposition 
\eqref{eqn:LQ} for $K = \bC$. 
%%%%%%%
\par
%%%%%%%
On the other hand, $\chi_0(z)$ is either (C) a cyclotomic polynomial; 
or (S) a Salem polynomial. 
In case (C), since $\delta(F) \neq \pm 1$, the polynomial $\chi_0(z)$ 
is of degree even and has distinct roots  
$\lambda_1, \dots, \lambda_m, \bar{\lambda}_1, \dots, \bar{\lambda}_m 
\in S^1$ with $\deg \chi_0(z) = 2m \ge 2$. 
Then $L_0(\bC)$ admits an orthogonal decomposition 
$L_0(\bC) = \bigoplus_{j=1}^m E(\lambda_j) \oplus E(\bar{\lambda}_j)$ by 
eigenspaces, where $E(\lambda_j)$ and $E(\bar{\lambda}_j)$ 
have the same signature $(1, 0)$ or $(0, 1)$. 
Thus $L_0(\bC)$ must have a signature of the form 
$( \mbox{even}, \mbox{even})$, so we are in case (i) and 
$L_1(\bC)$ must have signature $(1, \rho-1)$.  
In case (S) the polynomial $\chi_0(z)$ has distinct roots 
$\lambda^{\pm 1}$ and 
$\lambda_1, \dots, \lambda_m, \bar{\lambda}_1, \dots, 
\bar{\lambda}_m \in S^1$ with $\lambda > 1$ being  
a Salem number and $\deg \chi_0(z) = 2(m+1) \ge 4$. 
Then $L_0(\bC)$ admits an orthogonal decomposition 
$L_0(\bC) = E(\lambda, \lambda^{-1}) \oplus 
\bigoplus_{j=1}^m E(\lambda_j) \oplus E(\bar{\lambda}_j)$ 
as in \eqref{eqn:esd}, where $E(\lambda, \lambda^{-1})$ has 
signature $(1, 1)$ while $E(\lambda_j)$ and $E(\bar{\lambda}_j)$ 
have the same signature $(1, 0)$ or $(0, 1)$. 
Thus $L_0(\bC)$ must have a signature of the form 
$( \mbox{odd}, \mbox{odd})$, so we are in case (ii) and 
$L_1(\bC)$ must have signature $(0, \rho)$.  
These observations and the latter part of \eqref{eqn:L0L1} 
then lead to the dichotomy in the theorem.  
Clearly cases (C) and (S) correspond to the projective and non-projective 
cases respectively.  \hfill $\Box$ \par\medskip
%%%%%%%%%%%%%%%%%%%% end proof %%%%%%%%%%%%%%%%%%%%%%%%%%%%%%%%%%%  
Let $L = L(\varphi, \psi)$ be a hypergeometric K3 lattice and let 
$F = A$ or $B$.   
Suppose that $F$ admits a special eigenvalue $\delta(F)$ such that 
$F$ is a positive Hodge isometry with respect to the Hodge 
structure \eqref{eqn:ell}. 
Let $\chi_0(z)$ be the minimal polynomial of $\delta(F)^{\pm 1}$. 
There then exists a factorization of polynomials \eqref{eqn:chi01}.   
%%%%%%%%%%%%%
\par\medskip
%%%%%%%%%%%%% 
{\it Proof of Theorem $\ref{thm:main4}$}. 
Recall from item (4) of Remark \ref{rem:main1} that 
$\delta(F)$ is different from $\pm 1$. 
By item (1) of the same remark any root of $\chi(z)$ 
is simple except for the triple root $z = 1$ in the parabolic case 
of Theorem \ref{thm:main2}.  
In any case $\chi_0(z)$ and $\chi_1(z)$ are coprime, since 
$z = 1$ cannot be a root of $\chi_0(z)$. 
Thus Theorem \ref{thm:pic} establishes all assertions in 
Theorem \ref{thm:main4} other than the one on standard basis.  
Formula \eqref{eqn:AB} shows that $\br$ is a cyclic vector 
for both $F = A$ and $B$. 
So any vector $\bv \in L$ can be expressed as 
$\bv = \chi_2(F) \br$ for a unique polynomial $\chi_2(z) \in \bZ[z]$ 
such that $\deg \chi_2(z) \le 21$.   
In view of \eqref{eqn:chi1} we have $\bv \in \Pic$ if and only if 
$\chi_1(F) \cdot \chi_2(F) \br = \b0$, which is the case exactly 
when $\chi(z)$ divides $\chi_1(z) \cdot \chi_2(z)$, that is, 
$\chi_0(z)$ divides $\chi_2(z)$. 
Upon writing $\chi_2(z) = \chi_3(z) \cdot \chi_0(z)$, any vector  
$\bv \in \Pic$ can be represented as $\bv = \chi_3(F) \bs$ 
with $\bs := \chi_0(F) \br$ for a unique polynomial $\chi_3(z) \in \bZ[z]$ 
such that $\deg \chi_3(z) \le \rho-1$.   
Hence $\bs_1, \dots, \bs_{\rho}$ form a free $\bZ$-basis of 
$\Pic$. \hfill $\Box$
%%%%%%%%%%%%%%%%%%%%%%% end proof %%%%%%%%%%%%%%%%%%%%%%%%%%%%%%%%%%
%%%%%%%%%%%%%%%%%%%%%%%%% ss:algorithm  %%%%%%%%%%%%%%%%%%%%%%%%%%%%%%
\subsection{Bringing-back Algorithm} \label{ss:algorithm} 
%%%%%%%%%%%%%%%%%%%%%%%%%%%%%%%%%%%%%%%%%%%%%%%%%%%%%%%%%%%%%%%%%%
In the non-projective case there exists an algorithm to output the Weyl group 
element $w_{F} \in W$ in \eqref{eqn:modified} for $F = A$ or $B$.  
Before stating it we make a remark about K\"{a}hler cone. 
There exists an orthogonal decomposition 
$H^{1,1}_{\bR} = V_{\bR} \oplus \Pic_{\bR}$ of signatures  
$(1, 19-\rho) \oplus (0, \rho)$, where $\Pic_{\bR} := \Pic \otimes \bR$ and the 
Weyl group $W$ acts on $V_{\bR}$ trivially.    
Thus the intersection form $(\bu, \bv)$ restricted to $\Pic_{\bR}$ 
is negative-definite. 
For the sake of convenience we turn it positive-definite by putting 
$\langle \bu , \bv \rangle := - (\bu , \bv)$. 
The essential part of the K\"ahler cone $\cK$ is given by 
%%%%%%%%%%%%%%%%%%%%%% eqn:K1 %%%%%%%%%%%%%%%%%%%%%%%%%%%%%%%%%%%%%
\begin{equation} \label{eqn:K1}
\cK_1 := \cK \cap \Pic_{\bR} = 
\{ \bv_1 \in \Pic_{\bR} :  \langle \bv_1, \bu \rangle > 0, \,  
{}^{\forall} \bu \in \vD^+ \}. 
\end{equation}
%%%%%%%%%%%%%%%%%%%%%%%%%%%%%%%%%%%%%%%%%%%%%%%%%%%%%%%%%%%%%%%%%
Indeed, if we put $\cC^+_0 := \cC^+ \cap V_{\bR}$ then the positive cone is 
given by 
%%%%%%%%
$$
\cC^+ = \{ \bv = \bv_0-\bv_1 : \bv_0 \in \cC^+_0,  \, \bv_1 \in \Pic_{\bR}, \,  
(\bv_0, \bv_0) > \langle \bv_1, \bv_1 \rangle \}, 
$$
%%%%%%%
and the vector $\bv = \bv_0-\bv_1 \in \cC^+$ belongs to $\cK$ if and only if 
$\bv_1 \in \cK_1$. 
%%%%%%%%%%%%%%%%%%%%%% algorithm %%%%%%%%%%%%%%%%%%%%%%%%%%%%%%%
\begin{algorithm} \label{algorithm} 
Take $A = Z(\varphi)$ and $B = Z(\psi)$ as in \eqref{eqn:std} and 
specify $F = A$ or $B$.  
%%%%%%%%%%%%%%%%%%%%
\par\medskip
%%%%%%%%%%%%%%%%%%%%
{\bf 1. Gram matrix}. The Gram matrix $\langle \bs_i, \bs_j \rangle$ 
for the standard basis $\bs_1, \dots, \bs_{\rho}$ of $\Pic$ 
(see Theorem \ref{thm:main4}) can be evaluated explicitly 
by using formula \eqref{eqn:xi} in Theorem \ref{thm:main1} 
(or its $B$-basis version).   
Note that $\langle \bs_i, \bs_j \rangle$ depends only on $|i-j|$. 
We give an algorithm to find all roots of the root system $\vD$. 
%%%%%%%%%%%%%%%%%%%%
\par\medskip
%%%%%%%%%%%%%%%%%%%%
{\bf 2. Finding all roots}. We have an even, positive-definite quadratic form  
%%%%%%%
$$
Q(t_1, \dots, t_{\rho}) := \langle \bu, \bu \rangle 
= \sum_{i,j=1}^{\rho} \langle \bs_i, \bs_j \rangle \, t_i t_j 
\qquad \mbox{in \quad $(t_1, \dots, t_{\rho}) \in \bZ^{\rho}$},  
$$
%%%%%%
by expressing each element $\bu \in \Pic$ as a $\bZ$-linear combination 
$\bu = t_1 \bs_1 + \cdots + t_{\rho} \bs_{\rho}$.   
Finding all roots in $\vD$ amounts to finding all integer solutions 
to the inequality $Q(t_1, \dots, t_{\rho}) \le 2$. 
Indeed, all but the trivial solution $\b0 = (0, \dots, 0)$ lead to the solutions 
of the equation $Q(t_1, \dots, t_{\rho}) = 2$, because $Q(t_1, \dots, t_{\rho})$ 
cannot take value $1$. 
As a positive-definite quadratic form, $Q(t_1, \dots, t_{\rho})$ can 
be expressed as 
%%%%%%%%%%%%%%%%%%%%%%% eqn:Q %%%%%%%%%%%%%%%%%%%%%%%%%%%%%%%%%%%%%
\begin{equation} \label{eqn:Q}
Q(t_1, \dots, t_{\rho}) = \sum_{j=1}^{\rho} c_j \{ t_j - p_j(t_1, \dots, t_{j-1}) \}^2, 
\qquad c_j \in \bQ_{> 0}, \quad j = 1, \dots, {\rho}, 
\end{equation}
%%%%%%%%%%%%%%%%%%%%%%%%%%%%%%%%%%%%%%%%%%%%%%%%%%%%%%%%%%%%%%%%%%
where  $p_1 = 0$ and for $j = 2, \dots, {\rho}$ the expression 
$p_j(t_1, \dots, t_{j-1})$ is a linear form over $\bQ$ in $t_1, \dots, t_{j-1}$. 
Note that $c_j$ and $p_j$ can be calculated explicitly in terms of the 
coefficients of $Q(t_1, \dots, t_{\rho})$.   
%%%%%
\par
%%%%%  
Let $Q_k(t_1, \dots t_{k-1}; t_k)$ be the partial sum of \eqref{eqn:Q} summed 
over $j = 1, \dots, k$; it is a quadratic function of $t_k$, once $t_1, \dots, t_{k-1}$ 
are given. 
We define a {\sl rooted forest} in the following manner. 
First, let all integer solutions $t_1$ of the inequality $Q_1(t_1) \le 2$ 
be the roots (in graph theory) of the forest.  
Next, for each root $t_1$, let all integer solutions $t_2$ of  
$Q_2(t_1; t_2) \le 2$ be the children of $t_1$. 
Inductively, given a parent $t_k$ with its ancestors $t_1, \dots, t_{k-1}$, 
let all integer solutions $t_{k+1}$ of $Q_{k+1}(t_1, \dots, t_k; t_{k+1}) \le 2$ 
be the children of $t_k$. 
Consider all paths from roots to leaves, say, $(t_1, \dots, t_k)$.  
Some of them may continue to the $\rho$-th generation, that is, 
$k = \rho$, while others may not. 
All paths $(t_1, \dots, t_{\rho})$ that continue to the $\rho$-th generation  
yield all integer solutions to the inequality $Q(t_1, \dots, t_{\rho}) \le 2$ and 
hence all elements $\bu \in \vD$ along with the origin $\b0$.     
%%%%%%%%%%%%
\par\medskip
%%%%%%%%%%%%
{\bf 3. Positive roots, simple roots and the K\"ahler cone}.  
Provide $\Pic$ with a lexicographic order in the following manner:  
for $\bu = t_1 \bs_1 + \cdots + t_{\rho} \bs_{\rho}$, 
$\bu' = t_1' \bs_1 + \cdots + t_{\rho}' \bs_{\rho} \in \Pic$,  
%%%%%%%%%%%%%%%%%%%%%%%% eqn:order %%%%%%%%%%%%%%%%%%%%%%%%%%%%%%%%  
\begin{equation} \label{eqn:order}
\mbox{$\bu \succ \bu'$ $\defby$ ${}^{\exists} i \in \{1, \dots, \rho \}$ 
such that $t_i > t_i'$ and $t_j = t_j'$ for all $j < i$}.   
\end{equation}
%%%%%%%%%%%%%%%%%%%%%%%%%%%%%%%%%%%%%%%%%%%%%%%%%%%%%%%%%%%%%%%%%   
Then in the previous step all solutions $\bu = t_1 \bs_1 + \cdots + 
t_{\rho} \bs_{\rho} \succ \b0$ give the set of positive roots $\vD^+$.  
An element $\bu \in \vD^+$ is a simple root if and only if 
$\bu - \bu' \not\in \vD^+$ for any $\bu' \in \vD^+$. 
This test can easily be carried out by computer and we obtain   
the set of simple roots, say $\vD_{\rb}$, that is, the basis of 
$\vD$ relative to $\vD^+$. 
Looking at $\vD_{\rb}$ we can draw the Dynkin diagram of 
$\vD$, which indicates the irreducible decomposition of $\vD$ 
as well as the Dynkin type of each irreducible component. 
The essential part $\cK_1$ of the K\"ahler cone $\cK$ is given by  
formula \eqref{eqn:K1}. 
It is the Weyl chamber $C(\bdelta)$ containing the regular vector 
$\bdelta := \frac{1}{2} \sum_{\sbu \in \vD^+} \bu \in \Pic_{\bR}$.  
%%%%%%%%%%%%
\par\medskip
%%%%%%%%%%%% 
{\bf 4. Bringing back}. The matrix $F$ sends $\cK_1 = C(\bdelta)$ to the 
Weyl chamber $C(\bd)$ containing the regular vector $\bd :=  F \bdelta 
\in \Pic_{\bR}$.    
Let $w_F \in W$ be an element maximizing     
$\langle w(\bd), \bdelta \rangle$ for $w \in W$.    
We claim that $w_F $ brings $C(\bd)$ back to $C(\bdelta)$. 
Indeed, for any $\bu \in \vD_{\rb}$ one has $\sigma_{\sbu}(\bdelta) = 
\bdelta - \bu$ by Humphreys \cite[\S10.2, Lemma B]{Humphreys} and 
hence the defining property of $w_F$ and $\sigma_{\sbu}$-invariance 
of the inner product yield   
%%%%%%%
$$
\langle w_F(\bd), \bdelta \rangle 
\ge \langle (\rho_{\sbu} w_F) (\bd), \bdelta \rangle = 
\langle w_F (\bd), \rho_{\sbu}(\bdelta) \rangle = 
\langle w_F (\bd), \bdelta-\bu \rangle = \langle w_F (\bd), \bdelta \rangle 
-\langle w_F(\bd), \bu \rangle,   
$$
%%%%%%% 
that is, $\langle w_F(\bd), \bu \rangle \ge 0$. 
This implies more strictly that $\langle w_F(\bd), \bu \rangle > 0$ for any 
$\bu \in \vD_{\rb}$, since $w_F(\bd)$ is a regular vector.   
Thus we have $w_F(\bd) \in C(\bdelta)$ and so $w_F (C(\bd)) = C(\bdelta)$. 
Note that the element $w_F \in W$ is unique because $W$ acts on the 
set of Weyl chambers simply transitively.  
%%%%%%%%%%%%
\par\medskip
%%%%%%%%%%%% 
{\bf 5. Product of Picard-Lefschetz reflections}. 
To determine $w_F$ explicitly, let $\PL$ be the set of all 
Picard-Lefschetz reflections together with the identity transformation $1$. 
Start with $\bd_0 := \bd$ and find an element $\sigma_1 \in \PL$ maximizing    
$\langle \sigma(\bd_0), \bdelta \rangle$ for $\sigma \in \PL$.  
If $\sigma_1 = 1$, stop here; otherwise, put $\bd_1 := \sigma_1(\bd_0)$ and 
find $\sigma_2 \in \PL$ maximizing $\langle \sigma(\bd_1), \bdelta \rangle$ 
for $\sigma \in \PL$.  
Inductively, given $\sigma_k \in \PL$ and $\bd_{k-1} \in \Pic_{\bR}$, 
if $\sigma_k =1$, stop at this stage; otherwise, put $\bd_k := \sigma_k(\bd_{k-1})$ 
and find yet another $\sigma_{k+1} \in \PL$ maximizing  
$\langle \sigma(\bd_k), \bdelta \rangle$ for $\sigma \in \PL$. 
We claim that 
%%%%%%%%%%%%%%%%%%%%%%%% eqn:claim %%%%%%%%%%%%%%%%%%%%%%%%%%%%%%
\begin{equation} \label{eqn:claim}
\mbox{if $\sigma_{k+1} = 1$ then $\bd_k \in C(\bdelta)$; otherwise,  
$\langle \sigma_{k+1} (\bd_k), \bdelta \rangle > \langle \bd_k, \bdelta \rangle$}. 
\end{equation}
%%%%%%%%%%%%%%%%%%%%%%%%%%%%%%%%%%%%%%%%%%%%%%%%%%%%%%%%%%%%%%%
\par
%%%%%%
Indeed, if $\sigma_{k+1} = 1$ then for any $\bu \in \vD_{\rb}$ we have from  
Humphreys \cite[\S10.2, Lemma B]{Humphreys},  
%%%%%%
$$
\langle \bd_k, \bdelta \rangle \ge \langle \sigma_{\sbu}(\bd_k), \bdelta \rangle 
= \langle \bd_k, \sigma_{\sbu}(\bdelta) \rangle = 
\langle \bd_k, \bdelta -\bu \rangle = 
\langle \bd_k, \bdelta \rangle -\langle \bd_k, \bu \rangle,  
$$
%%%%%
that is, $\langle \bd_k, \bu \rangle \ge 0$, which in turn implies 
$\langle \bd_k, \bu \rangle > 0$ for any $\bu \in \vD_{\rb}$, since 
$\bd_k$ is a regular vector.  
Thus $\bd_k \in C(\bdelta)$ and the first part of \eqref{eqn:claim} is proved. 
Next, suppose that $\sigma_{k+1} \neq 1$ is the Picard-Lefschetz reflection 
associated with a positive root $\bu_{k+1} \in \vD^+$. 
Then  
%%%%%%
$$
\langle \sigma_{k+1}(\bd_k), \bdelta \rangle = 
\langle \bd_k - \langle \bd_k, \bu_{k+1} \rangle \bu_{k+1}, \bdelta \rangle 
= \langle \bd_k, \bdelta \rangle - 
\langle \bd_k, \bu_{k+1} \rangle \langle \bu_{k+1}, \bdelta \rangle, 
$$
%%%%%%  
where $\langle \bd_k, \bu_{k+1} \rangle \langle \bu_{k+1}, \bdelta \rangle$ 
is nonzero, since $\bd_k$ and $\bdelta$ are regular vectors. 
On the other hand, by the defining property of $\sigma_{k+1}$ we have 
$\langle \sigma_{k+1} (\bd_k), \bdelta \rangle \ge 
\langle \bd_k, \bdelta \rangle$ and thus  
$\langle \sigma_{k+1} (\bd_k), \bdelta \rangle > \langle \bd_k, \bdelta \rangle$.  
%%%%%%
\par
%%%%%%
Since $\{ \langle w(\bd), \bdelta \rangle : w \in W \}$ is a finite set, it 
follows from \eqref{eqn:claim} that the step-by-step procedure mentioned 
above eventually terminates with $\sigma_{k+1} = 1$ and leads to the desired 
representation $w_F = \sigma_k \circ \sigma_{k-1} \circ \cdots \circ \sigma_1$ 
as a product of Picard-Lefschetz reflections. 
%%%%%%%%%%%%
\par\medskip
%%%%%%%%%%%% 
{\bf 6. Modified matrix}.  Let $w_F$ act on the whole lattice $L$ by 
extending it as identity to the orthogonal complement of $\Pic$.   
The modified matrix $\tilde{F} := w_F \circ F$ then preserves  
the K\"ahler cone $\cK$ and hence preserves 
the K3 structure constructed from $F$. 
Thanks to factorization \eqref{eqn:chi01} and the non-projective assumption 
the characteristic polynomial $\tilde{\chi}(z)$ of $\tilde{F}$ factors as 
%%%%%%%%%%%%%%%%%%%%%%%%% eqn:chi01t %%%%%%%%%%%%%%%%%%%%%%%%%%%%%%
\begin{equation} \label{eqn:chi01t}
\tilde{\chi}(z) = \chi_0(z) \cdot \tilde{\chi}_1(z), \qquad 
\deg \chi_0(z) = 22- \rho, \quad \deg \tilde{\chi}_1(z) = \rho,  
\end{equation}
%%%%%%%%%%%%%%%%%%%%%%%%%%%%%%%%%%%%%%%%%%%%%%%%%%%%%%%%%%%%%%%%%% 
where the Salem polynomial component $\chi_0(z)$ is the same as that 
in \eqref{eqn:chi01}, while $\tilde{\chi}_1(z)$ is the characteristic 
polynomial of $\tilde{F}|\Pic$.   
Thus $\tilde{F}$ has the same spectral radius as $F$.  
%%%%%%%%%%%%
\par\medskip
%%%%%%%%%%%% 
{\bf 7. Action on the Dynkin diagram}.   
Since $\tilde{F}$ preserves $\vD^+$, it also preserves $\vD_{\rb}$.  
Thus $\tilde{F}$ induces an automorphism of the corresponding Dynkin diagram. 
Describe this action explicitly. 
\end{algorithm} 
%%%%%%%%%%%%%%%%%%%%%%%%%%%%%%%%%%%%%%%%%%%%%%%%%%%%%%%%%%%%%%%%%%
\par
%%%%%%%%
Now we are able to obtain a K3 surface automorphism $f : X \to X$ as the lift of 
the modified Hodge isometry $\tilde{F} : L \to L$.  
This establishes Theorem \ref{thm:main5}.  
Step $7$ is used to apply Lemma \ref{lem:curve} to the constructed map $f$.  
%%%%%%%%%%%%%%%%%%%%%%%%% rem:pic %%%%%%%%%%%%%%%%%%%%%%%%%%%%%%%%%%
\begin{remark} \label{rem:pic}  
The Picard number $\rho$ is always positive for $F = A$, while it can be zero 
for $F = B$.    
If $\rho = 0$ then $\Pic$ is trivial, the K\"{a}hler cone $\cK$ coincides with 
the positive cone $\cC^+$, so there is no need to modify $B$ by a 
Weyl group element.   
In this case $\psi(z)$ is an unramified Salem polynomial of degree $22$.         
\end{remark}
%%%%%%%%%%%%%%%%%%%%%%%%%%%%%%%%%%%%%%%%%%%%%%%%%%%%%%%%%%%%%%%%%%  
%%%%%%%%%%%%%%%%%%%%%% sec:k3auto %%%%%%%%%%%%%%%%%%%%%%%%%%%%%%%%%%
\section{K3 Surface Automorphisms}  \label{sec:k3auto} 
%%%%%%%%%%%%%%%%%%%%%%%%%%%%%%%%%%%%%%%%%%%%%%%%%%%%%%%%%%%%%%%%%%
We illustrate the method of hypergeometric groups by constructing many 
examples of non-projective K3 surface automorphisms of positive entropy. 
In this article we present only a couple of cases involving ten Salem numbers 
of degree $22$ and Lehmer's number. 
Much wider applications of our method will be reported upon elsewhere.  
%%%%%%%%%%%%%%%%%%%%%%% ss:S22 %%%%%%%%%%%%%%%%%%%%%%%%%%%%%%%%%%%%%%%
\subsection{Salem Numbers of Degree $\mbox{\boldmath $22$}$} \label{ss:S22}
%%%%%%%%%%%%%%%%%%%%%%%%%%%%%%%%%%%%%%%%%%%%%%%%%%%%%%%%%%%%%%%%%%%%%
Things are simpler with Salem numbers of degree $22$ by Remark \ref{rem:pic}.     
McMullen \cite[Table 4]{McMullen1} gives a list of ten unramified Salem 
polynomials $S_i(z)$ of degree $22$, $i = 1, \dots, 10$, for each of which  
he constructs a K3 surface automorphism $f : X \to X$ with 
a Siegel disk such that the induced map 
$f^* : H^2(X, \bZ) \to H^2(X, \bZ)$ has $S_i(z)$ as its characteristic 
polynomial (see \cite[Theorem 10.1]{McMullen1}). 
The Salem trace polynomials $R_i(w)$ associated with $S_i(z)$ 
are given in Table \ref{tab:mcmullen}. 
Applying our method to them we are able to construct a much greater 
number of K3 surface automorphisms with a Siegel disk. 
We refer to \S \ref{sec:SD} for discussions about Siegel disks.  
%%%%%%
\par
%%%%%%
{\bf Setup and Tests}. 
In order to deal with Salem numbers of degree $22$ it is necessary to use the matrix $B$.  
Let $R(w)$ be an unramified Salem trace polynomial of degree $11$. 
We look for all hypergeomtric K3 lattices $L$ such that $\Psi(w) = R(w)$ 
and $\Phi(w)$ is a product of cyclotomic trace polynomials of the form 
%%%%%%%%%%%%%%%%%%%%%% eqn:CTk %%%%%%%%%%%%%%%%%%%%%%%%%%%%%%%%%%
\begin{equation} \label{eqn:CTk}
\Phi(w) = \CT_{\sbk}(w) := \prod_{k \in \sbk} \CT_k(w) \quad 
\mbox{with} \quad \sum_{k \in \sbk} \deg \CT_k(w) = 10, 
\end{equation}
%%%%%%%%%%%%%%%%%%%%%%%%%%%%%%%%%%%%%%%%%%%%%%%%%%%%%%%%%%%%%%%
where $\bk$ is a multi-set of positive integers whose elements $k$ come    
from Table \ref{tab:CTP}.   
By the last part of Theorem \ref{thm:main3} any element of $\bk$ is simple 
except for at most one element of multiplicity $2$ or $3$, where the multiple 
element must be one of $1$, $2$, $3$, $4$, $6$. 
An inspection of Table \ref{tab:CTP} shows that $\bk$ has at most 
seven distinct elements and this number reduces to six if $\bk$ has a triple 
element. 
With this setup our tests proceed as follows. 
%%%%%%%%%%%%%%%%%%
\begin{enumerate}
\setlength{\itemsep}{-1pt}
\item Find all $\bk$'s satisfying the unimodularity condition \eqref{eqn:unim2}, 
which reads $\Res(\CT_k, R) = \pm 1$ for every $k \in \bk$.  
\item Judge which of the $\bk$'s above are K3 lattices according to 
the criterion in Theorem \ref{thm:main3}.   
\item Identify the special trace $\tau \in (-2, \, 2)$ of the matrix $B$ 
by using Theorem \ref{thm:main3} again.  
\end{enumerate}
%%%%%%%%%%%%%%%%%
\par
%%%%%%%%% 
For each of the data in \eqref{eqn:CTk} passing two tests in steps (1) and (2), 
the information from step (3) allow us to provide $L$ with the Hodge 
structure \eqref{eqn:ell}, with respect to which $B$ is a positive Hodge isometry, 
since $\Psi(w) = R(w)$ is a Salem trace polynomial (see Remark \ref{rem:positive}).    
Specify one component of $\cC$ as the positive cone $\cC^+$ as well as  
the K\"{a}hler cone $\cK$. 
Thanks to Theorem \ref{thm:great} the Hodge isometry 
$B : L \to L$ then lifts to a K3 surface automorphism $f : X \to X$.   
%%%%%%%%%%%
\par
%%%%%%%%%%%
We carry out the above-mentioned tests for the ten Salem trace polynomials 
$R_i(w)$ in Table \ref{tab:mcmullen}. 
For each $i = 1, \dots, 10$, let $y_{10} < \cdots < y_2 < y_1$ be the roots of 
$R_i(w)$ in the interval $(-2, \, 2)$. 
Numerical values of them are given in Table \ref{tab:roots}, where the roots 
greater than $\tau_0 := 1-2\sqrt{2}$ are separated from those smaller than 
$\tau_0$ by a line; this information will be used to discuss Siegel disks in 
\S \ref{sec:SD}.  
Note that $y_{10}$ is always smaller than $\tau_0$. 
%%%%%%%%%%%%%%%%%%%%%%%% tab:roots %%%%%%%%%%%%%%%%%%%%%%%%%%%%%%%%%%
\begin{table}[hh]
\centerline{
\begin{tabular}{|c|c|c|c|c|c|}
\hline 
            & $R_1$            & $R_2$             & $R_3$             &  $R_4$           & $R_5$ \\
\hline
$y_1$    & $1.993294$   & $1.988033$     & $1.995401$     & $1.995839$    & $1.994419$   \\
$y_2$    & $1.205977$    & $1.76766$      & $1.339130$     & $1.495199$    & $1.594614$   \\
$y_3$    & $0.9297159$  & $0.7304257$   & $0.9728311$   & $ 0.5691842$  & $0.9860031$   \\
$y_4$    & $0.3600253$  & $0.3763138$   & $0.1183310$   & $0.2303372$   &  $0.2512725$  \\
$y_5$    & $-0.1005899$ & $-0.3628799$ & $-0.3009327$ & $ -0.2778793$ & $-0.5016395$  \\
$y_6$    & $-0.8098076$ & $-0.8420136$ & $ -0.736344$  & $-0.6300590$ & $-0.7481516$   \\
$y_7$    & $-1.19931$    & $-1.280397$   & $-1.297751$    & $-1.18174$    & $-1.209765$   \\
$y_8$    & $-1.667161$  & $-1.677966$  & $-1.492648$   & $-1.515876$   & $-1.746739$ \\
\cline{2-3}\cline{5-6}
$y_9$    & $-1.842436$   & $-1.891176$   & $-1.770639$  & $-1.897604$   & $-1.866712$   \\
\cline{4-4}
$y_{10}$ & $-1.970982$   & $-1.948177$   & $-1.983677$   &  $-1.96877$    & $-1.950669$   \\
\hline\hline
             & $R_6$             & $R_7$            & $R_8$ & $R_9$ & $R_{10}$ \\
\hline
$y_1$ & $1.997218$    & $1.996192$     & $1.992823$    & $1.988808$      & $1.995293$ \\
$y_2$     & $1.368950$    & $1.521855$    & $1.816402$     & $1.855143$     & $1.738103$ \\
$y_3$     & $0.526038$    & $0.8486504$   & $0.6683560$   & $1.144637$     & $0.773170$ \\
$y_4$     & $0.3786577$  & $0.4924207$   & $0.08744302$  & $0.02840390$ & $0.0659823$ \\
$y_5$     & $-0.1559579$ & $-0.4803898$ & $-0.3825477$ & $ -0.4829143$ & $-0.4868621$ \\
$y_6$     & $-0.7702764$ & $-0.8293021$ & $-0.7409521$ & $-0.8693841$  & $-0.6959829$ \\
$y_7$     & $-1.274529$   & $-1.297520$  & $-1.234218$   & $-1.419774$    & $ -1.265297$ \\
$y_8$     & $-1.464669$ & $-1.694755$ & $-1.698505$ & $-1.704027$ & $-1.548546$ \\
\cline{2-3}\cline{5-6}
$y_9$     & $ -1.85815$    & $-1.838895$  & $-1.790300$ & $-1.868712$    & $-1.904399$ \\
\cline{4-4}
$y_{10}$  & $-1.97660$     & $-1.961625$   & $ -1.971905$ & $-1.938202$    & $-1.956408$ \\
\hline
\end{tabular}} 
\caption{Roots smaller than $2$ of the Salem trace polynomials in Table \ref{tab:mcmullen}.} 
\label{tab:roots}
\end{table}
%%%%%%%%%%%%%%%%%%%%%%%%%%%%%%%%%%%%%%%%%%%%%%%%%%%%%%%%%%%%%%%%%%
%%%%%%%%%%%%%%%%%%%%%%%% thm:SH %%%%%%%%%%%%%%%%%%%%%%%%%%%%%%%%%%
\begin{theorem} \label{thm:SH} 
In the above setting we have a hypergeometric K3 lattice such 
that the matrix $B$ admits a special eigenvalue $\delta$ conjugate to the 
Salem number $\lambda_i$, if and only if $\Psi(w) = R_i(w)$ and $\bk$ are 
as in Tables $\ref{tab:SvH1}$--$\ref{tab:SvH3}$, where the special 
trace $\tau := \delta + \delta^{-1}$ is given in the ``ST" column, while 
the ``case" column refers to the case in Table $\ref{tab:hyp-B}$.  
For each entry of the tables the Hodge isometry $B : L \to L$ lifts to 
a K3 surface automorphism $f : X \to X$ of entropy $h(f) = \log \lambda_i$ 
with special trace $\tau(f) = \tau$ and Picard number $\rho(X) = 0$. 
The ``S/H" column is explained in Theorem $\ref{thm:SH2}$.      
%%%%%%%%%%%%%%%%%%%%%%%% tab:SvH1 %%%%%%%%%%%%%%%%%%%%%%%%%%%%%%%%%
\begin{table}[p]
\begin{minipage}{0.45\linewidth}
\centerline{
\begin{tabular}{|c|c|c|c|c|}\hline
$\Psi$ & case & $\bk$ & ST & S/H \\ \hline
$R_1$ & $1$ & $1,1,1,3,4,6,16$ & $y_8$ & S \\ \cline{2-5}
& $1$ & $1,1,2,3,4,6,16$ & $y_8$ & S \\ \cline{2-5}
& $1$ & $2,2,1,3,4,6,16$ & $y_8$ & S \\ \cline{2-5}
& $1$ & $2,2,2,3,4,6,16$ & $y_8$ & S \\ \cline{2-5}
& $2$ & $3,3,1,3,4,6,16$ & $y_8$ & S \\ \cline{2-5}
& $2$ & $3,3,2,3,4,6,16$ & $y_8$ & S \\ \cline{2-5}
& $2$ & $4,4,1,3,4,6,16$ & $y_8$ & S \\ \cline{2-5}
& $2$ & $4,4,2,3,4,6,16$ & $y_8$ & S \\ \cline{2-5}
& $2$ & $6,6,1,3,4,6,16$ & $y_8$ & S \\ \cline{2-5}
& $2$ & $6,6,2,3,4,6,16$ & $y_8$ & S \\ \cline{2-5}
& $3$ & $1,3,16,30$ & $y_2$ & S \\ \cline{2-5}
& $3$ & $2,3,16,30$ & $y_2$ & S \\ \cline{2-5}
& $3$ & $1,3,5,7,9$ & $y_2$ & S \\ \cline{2-5}
& $3$ & $2,3,5,7,9$ & $y_2$ & S \\ \cline{2-5}
& $3$ & $1,3,17$ & $y_2$ & S \\ \cline{2-5}
& $3$ & $2,3,17$ & $y_2$ & S \\ \cline{2-5}
& $5$ & $1,2,3,4,5,6,7$ & $y_{10}$ & H \\ \cline{2-5}
& $5$ & $1,1,3,4,5,6,7$ & $y_{10}$ & H \\ \cline{2-5}
& $5$ & $2,2,3,4,5,6,7$ & $y_{10}$ & H \\ \cline{2-5}
& $7$ & $3,3,3,4,5,6,7$ & $y_{10}$ & H \\ \cline{2-5}
& $7$ & $4,4,3,4,5,6,7$ & $y_{10}$ & H \\ \cline{2-5}
& $7$ & $6,6,3,4,5,6,7$ & $y_{10}$ & H \\ \cline{2-5}
& $8$ & $3,6,16,30$ & $y_1$ & S \\ \cline{2-5}
& $8$ & $4,5,7,20$ & $y_1$ & S \\ \cline{2-5}
& $8$ & $3,5,6,7,9$ & $y_1$ & S \\ \cline{2-5}
& $8$ & $3,6,17$ & $y_1$ & S \\ \hline 
$R_2$ & $1$ & $1,1,1,9,24$ & $y_7$ & S \\ \cline{2-5}
& $1$ & $1,1,2,9,24$ & $y_7$ & S \\ \cline{2-5}
& $1$ & $2,2,1,9,24$ & $y_7$ & S \\ \cline{2-5}
& $1$ & $2,2,2,9,24$ & $y_7$ & S \\ \cline{2-5}
& $3$ & $1,3,5,13$ & $y_8$ & S \\ \cline{2-5}
& $3$ & $1,3,5,42$ & $y_8$ & S \\ \cline{2-5}
& $3$ & $2,3,5,13$ & $y_8$ & S \\ \cline{2-5}
& $3$ & $2,3,5,42$ & $y_8$ & S \\ \cline{2-5}
& $3$ & $3,3,1,9,24$ & $y_7$ & S \\ \cline{2-5}
& $3$ & $3,3,2,9,24$ & $y_7$ & S \\ \cline{2-5}
& $3$ & $1,3,5,12,24$ & $y_4$ & S \\ \cline{2-5}
& $3$ & $2,3,5,12,24$ & $y_4$ & S \\ \cline{2-5} 
& $3$ & $1,3,12,13$ & $y_3$ & S \\ \cline{2-5} 
& $3$ & $1,3,12,42$ & $y_3$ & S \\ \cline{2-5} 
& $3$ & $2,3,12,13$ & $y_3$ & S \\ \cline{2-5} 
& $3$ & $2,3,12,42$ & $y_3$ & S \\ \cline{2-5} 
& $3$ & $1,3,24,30$ & $y_2$ & S \\ \cline{2-5} 
& $3$ & $2,3,24,30$ & $y_2$ & S \\ \cline{2-5} 
& $5$ & $66$ & $y_{10}$ & H \\ \cline{2-5} 
& $5$ & $1,2,3,5,22$ & $y_{10}$ & H \\ \cline{2-5} 
& $5$ & $1,1,3,5,22$ & $y_{10}$ & H \\ \hline 
\end{tabular}}
\end{minipage}
\begin{minipage}{0.45\linewidth}
\centerline{
\begin{tabular}{|c|c|c|c|c|}\hline
$\Psi$ & case & $\bk$ & ST & S/H \\\hline
$R_2$ & $5$ & $2,2,3,5,22$ & $y_{10}$ & H \\ \cline{2-5}
& $7$ & $14,16,18$ & $y_{10}$ & H \\ \cline{2-5} 
& $7$ & $3,3,3,5,22$ & $y_{10}$ & H \\ \cline{2-5} 
& $8$ & $3,12,18,24$ & $y_1$ & S \\ \hline
$R_3$ & $1$ & $1,1,1,4,36$ & $y_8$ & S \\ \cline{2-5} 
& $1$ & $1,1,2,4,36$ & $y_8$ & S \\ \cline{2-5} 
& $1$ & $2,2,1,4,36$ & $y_8$ & S \\ \cline{2-5} 
& $1$ & $2,2,2,4,36$ & $y_8$ & S \\ \cline{2-5} 
& $1$ & $1,1,1,4,13$ & $y_7$ & S \\ \cline{2-5} 
& $1$ & $1,1,2,4,13$ & $y_7$ & S \\ \cline{2-5} 
& $1$ & $2,2,1,4,13$ & $y_7$ & S \\ \cline{2-5} 
& $1$ & $2,2,2,4,13$ & $y_7$ & S \\ \cline{2-5} 
& $2$ & $3,3,1,4,36$ & $y_8$ & S \\ \cline{2-5} 
& $2$ & $3,3,2,4,36$ & $y_8$ & S \\ \cline{2-5} 
& $2$ & $4,4,1,4,36$ & $y_8$ & S \\ \cline{2-5} 
& $2$ & $4,4,2,4,36$ & $y_8$ & S \\ \cline{2-5} 
& $2$ & $6,6,1,4,36$ & $y_8$ & S \\ \cline{2-5} 
& $2$ & $6,6,2,4,36$ & $y_8$ & S \\ \cline{2-5} 
& $2$ & $4,4,1,4,13$ & $y_7$ & S \\ \cline{2-5} 
& $2$ & $4,4,2,4,13$ & $y_7$ & S \\ \cline{2-5} 
& $2$ & $6,6,1,4,13$ & $y_7$ & S \\ \cline{2-5} 
& $2$ & $6,6,2,4,13$ & $y_7$ & S \\ \cline{2-5} 
& $2$ & $1,3,8,42$ & $y_3$ & S \\ \cline{2-5} 
& $2$ & $2,3,8,42$ & $y_3$ & S \\ \cline{2-5} 
& $3$ & $1,3,7,11$ & $y_9$ & S \\ \cline{2-5} 
& $3$ & $2,3,7,11$ & $y_9$ & S \\ \cline{2-5} 
& $3$ & $1,3,4,7,30$ & $y_9$ & S \\ \cline{2-5} 
& $3$ & $2,3,4,7,30$ & $y_9$ & S \\ \cline{2-5} 
& $3$ & $3,3,1,4,13$ & $y_7$ & S \\ \cline{2-5} 
& $3$ & $3,3,2,4,13$ & $y_7$ & S \\ \cline{2-5}
 & $3$ & $1,3,11,18$ & $y_2$ & S \\ \cline{2-5} 
& $3$ & $2,3,11,18$ & $y_2$ & S \\ \cline{2-5} 
& $3$ & $1,3,4,18,30$ & $y_2$ & S \\ \cline{2-5} 
& $3$ & $2,3,4,18,30$ & $y_2$ & S \\ \cline{2-5} 
& $8$ & $3,4,8,13$ & $y_1$ & S \\ \cline{2-5} 
& $8$ & $3,6,11,18$ & $y_1$ & S \\ \cline{2-5} 
& $8$ & $3,4,6,18,30$ & $y_1$ & S \\ \hline
$R_4$ & $1$ & $1,1,1,4,5,24$ & $y_7$ & S \\ \cline{2-5} 
& $1$ & $1,1,2,4,5,24$ & $y_7$ & S \\ \cline{2-5} 
& $1$ & $2,2,1,4,5,24$ & $y_7$ & S \\ \cline{2-5} 
& $1$ & $2,2,2,4,5,24$ & $y_7$ & S \\ \cline{2-5} 
& $2$ & $4,4,1,4,5,24$ & $y_7$ & S \\ \cline{2-5} 
& $2$ & $4,4,2,4,5,24$ & $y_7$ & S \\ \cline{2-5} 
& $3$ & $3,3,1,4,5,24$ & $y_7$ & S \\ \cline{2-5} 
& $3$ & $3,3,2,4,5,24$ & $y_7$ & S \\ \cline{2-5}
& $3$ & $1,3,4,5,11$ & $y_3$ & S \\ \cline{2-5} 
& $3$ & $2,3,4,5,11$ & $y_3$ & S \\ \hline
\end{tabular}} 
\end{minipage}
\caption{K3 surface automorphisms from Salem numbers of degree $22$, Part 1.}
\label{tab:SvH1}
\end{table}
%%%%%%%%%%%%%%%%%%%%%% tab:SvH2 %%%%%%%%%%%%%%%%%%%%%%%%%%%%
\begin{table}[p]
\begin{minipage}{0.45\linewidth}
\centerline{
\begin{tabular}{|c|c|c|c|c|}\hline
$\Psi$ & case & $\bk$ & ST & S/H \\ \hline
$R_4$ & $3$ & $1,3,17$ & $y_3$ & S \\ \cline{2-5} 
& $3$ & $2,3,17$ & $y_3$ & S \\ \cline{2-5} 
& $3$ & $1,3,24,30$ & $y_2$ & S \\ \cline{2-5} 
& $3$ & $2,3,24,30$ & $y_2$ & S \\ \cline{2-5} 
& $8$ & $4,13,18$ & $y_1$ & S \\ \cline{2-5} 
& $8$ & $3,4,9,11$ & $y_1$ & S \\ \hline
$R_5$ & $1$ & $1,1,17$ & $y_9$ & H \\ \cline{2-5} 
& $1$ & $1,2,17$ & $y_9$ & H \\ \cline{2-5}
& $1$ & $2,2,17$ & $y_9$ & H \\ \cline{2-5}
& $1$ & $1,1,1,7,24$ & $y_7$ & S \\ \cline{2-5} 
& $1$ & $1,1,2,7,24$ & $y_7$ & S \\ \cline{2-5} 
& $1$ & $2,2,1,7,24$ & $y_7$ & S \\ \cline{2-5} 
& $1$ & $2,2,2,7,24$ & $y_7$ & S \\ \cline{2-5} 
& $1$ & $1,1,1,9,30$ & $y_6$ & S \\ \cline{2-5} 
& $1$ & $1,1,2,9,30$ & $y_6$ & S \\ \cline{2-5} 
& $1$ & $2,2,1,9,30$ & $y_6$ & S \\ \cline{2-5} 
& $1$ & $2,2,2,9,30$ & $y_6$ & S \\ \cline{2-5} 
& $2$ & $3,3,17$ & $y_9$ & H \\ \cline{2-5}
& $2$ & $6,6,17$ & $y_9$ & H \\ \cline{2-5}
& $2$ & $1,6,7,9,10$ & $y_7$ & S \\ \cline{2-5} 
& $2$ & $2,6,7,9,10$ & $y_7$ & S \\ \cline{2-5} 
& $2$ & $6,6,1,7,24$ & $y_7$ & S \\ \cline{2-5} 
& $2$ & $6,6,2,7,24$ & $y_7$ & S \\ \cline{2-5} 
& $2$ & $1,6,7,9,12$ & $y_6$ & S \\ \cline{2-5} 
& $2$ & $2,6,7,9,12$ & $y_6$ & S \\ \cline{2-5} 
& $2$ & $6,6,1,9,30$ & $y_6$ & S \\ \cline{2-5} 
& $2$ & $6,6,2,9,30$ & $y_6$ & S \\ \cline{2-5} 
& $3$ & $3,3,1,7,24$ & $y_7$ & S \\ \cline{2-5} 
& $3$ & $3,3,2,7,24$ & $y_7$ & S \\ \cline{2-5} 
& $3$ & $3,3,1,9,30$ & $y_6$ & S \\ \cline{2-5} 
& $3$ & $3,3,2,9,30$ & $y_6$ & S \\ \cline{2-5} 
& $3$ & $1,3,24,30$ & $y_2$ & S \\ \cline{2-5} 
& $3$ & $2,3,24,30$ & $y_2$ & S \\ \cline{2-5} 
& $8$ & $13,16$ & $y_1$ & S \\ \cline{2-5} 
& $8$ & $16,42$ & $y_1$ & S \\ \cline{2-5} 
& $8$ & $3,6,24,30$ & $y_1$ & S \\ \cline{2-5} 
& $8$ & $3,7,12,24$ & $y_1$ & S \\ \cline{2-5} 
& $8$ & $3,9,10,30$ & $y_1$ & S \\ \hline
$R_6$ & $1$ & $1,1,1,4,5,24$ & $y_7$ & S \\ \cline{2-5} 
& $1$ & $1,1,2,4,5,24$ & $y_7$ & S \\ \cline{2-5} 
& $1$ & $2,2,1,4,5,24$ & $y_7$ & S \\ \cline{2-5} 
& $1$ & $2,2,2,4,5,24$ & $y_7$ & S \\ \cline{2-5} 
& $2$ & $21,24$ & $y_7$ & S \\ \cline{2-5} 
& $2$ & $4,4,1,4,5,24$ & $y_7$ & S \\ \cline{2-5} 
& $2$ & $4,4,2,4,5,24$ & $y_7$ & S \\ \hline
\end{tabular}}
\end{minipage}
\begin{minipage}{0.45\linewidth}
\centerline{
\begin{tabular}{|c|c|c|c|c|}\hline
$\Psi$ & case & $\bk$ & ST & S/H \\ \hline
$R_6$ & $2$ & $1,3,5,36$ & $y_4$ & S \\ \cline{2-5} 
& $2$ & $2,3,5,36$ & $y_4$ & S \\ \cline{2-5} 
& $2$ & $1,3,8,42$ & $y_3$ & S \\ \cline{2-5}
& $2$ & $2,3,8,42$ & $y_3$ & S \\ \cline{2-5} 
& $3$ & $3,3,1,4,5,24$ & $y_7$ & S \\ \cline{2-5} 
& $3$ & $3,3,2,4,5,24$ & $y_7$ & S \\ \cline{2-5} 
& $3$ & $1,3,24,30$ & $y_2$ & S \\ \cline{2-5} 
& $3$ & $2,3,24,30$ & $y_2$ & S \\ \cline{2-5} 
& $7$ & $4,8,14,16$ & $y_{10}$ & H \\ \cline{2-5} 
& $8$ & $3,4,5,8,24$ & $y_1$ & S \\ \cline{2-5}
& $8$ & $4,19$ & $y_1$ & S \\ \hline
$R_7$ & $1$ & $1,1,1,7,16$ & $y_7$ & S \\ \cline{2-5} 
& $1$ & $1,1,1,7,24$ & $y_7$ & S \\ \cline{2-5} 
& $1$ & $1,1,2,7,16$ & $y_7$ & S \\ \cline{2-5} 
& $1$ & $1,1,2,7,24$ & $y_7$ & S \\ \cline{2-5} 
& $1$ & $2,2,1,7,16$ & $y_7$ & S \\ \cline{2-5} 
& $1$ & $2,2,1,7,24$ & $y_7$ & S \\ \cline{2-5} 
& $1$ & $2,2,2,7,16$ & $y_7$ & S \\ \cline{2-5} 
& $1$ & $2,2,2,7,24$ & $y_7$ & S \\ \cline{2-5} 
& $3$ & $3,3,1,7,16$ & $y_7$ & S \\ \cline{2-5} 
& $3$ & $3,3,1,7,24$ & $y_7$ & S \\ \cline{2-5} 
& $3$ & $3,3,2,7,16$ & $y_7$ & S \\ \cline{2-5} 
& $3$ & $3,3,2,7,24$ & $y_7$ & S \\ \cline{2-5} 
& $3$ & $1,3,7,11$ & $y_5$ & S \\ \cline{2-5} 
& $3$ & $2,3,7,11$ & $y_5$ & S \\ \cline{2-5} 
& $3$ & $1,3,5,7,9$ & $y_5$ & S \\ \cline{2-5} 
& $3$ & $2,3,5,7,9$ & $y_5$ & S \\ \cline{2-5} 
& $3$ & $1,3,7,11$ & $y_5$ & S \\ \cline{2-5} 
& $3$ & $2,3,7,11$ & $y_5$ & S \\ \cline{2-5} 
& $3$ & $1,3,5,7,9$ & $y_5$ & S \\ \cline{2-5} 
& $3$ & $2,3,5,7,9$ & $y_5$ & S \\ \cline{2-5} 
& $3$ & $1,3,16,30$ & $y_2$ & S \\ \cline{2-5} 
& $3$ & $1,3,24,30$ & $y_2$ & S \\ \cline{2-5} 
& $3$ & $2,3,16,30$ & $y_2$ & S \\ \cline{2-5} 
& $3$ & $2,3,24,30$ & $y_2$ & S \\ \cline{2-5} 
& $3$ & $1,3,16,30$ & $y_2$ & S \\ \cline{2-5} 
& $3$ & $1,3,24,30$ & $y_2$ & S \\ \cline{2-5} 
& $3$ & $2,3,16,30$ & $y_2$ & S \\ \cline{2-5} 
& $3$ & $2,3,24,30$ & $y_2$ & S \\ \hline
$R_8$ & $1$ & $1,1,1,4,12,30$ & $y_6$ & S \\ \cline{2-5} 
& $1$ & $1,1,2,4,12,30$ & $y_6$ & S \\ \cline{2-5} 
& $1$ & $2,2,1,4,12,30$ & $y_6$ & S \\ \cline{2-5} 
& $1$ & $2,2,2,4,12,30$ & $y_6$ & S \\ \cline{2-5} 
& $1$ & $1,1,1,3,12,30$ & $y_5$ & S \\ \cline{2-5} 
& $1$ & $1,1,2,3,12,30$ & $y_5$ & S \\ \hline
\end{tabular}}
\end{minipage}
\caption{K3 surface automorphisms from Salem numbers of degree $22$, Part 2.}
\label{tab:SvH2}
\end{table}
%%%%%%%%%%%%%%%%%%%%% tab:SvH3 %%%%%%%%%%%%%%%%%%%%%%%%%%%%
\begin{table}[p]
\begin{minipage}{0.45\linewidth}
\centerline{
\begin{tabular}{|c|c|c|c|c|}\hline
$\Psi$ & case & $\bk$ & ST & S/H \\\hline
$R_8$ & $1$ & $2,2,1,3,12,30$ & $y_5$ & S \\ \cline{2-5} 
& $1$ & $2,2,2,3,12,30$ & $y_5$ & S \\ \cline{2-5} 
& $2$ & $4,4,1,4,12,30$ & $y_6$ & S \\ \cline{2-5} 
& $2$ & $4,4,2,4,12,30$ & $y_6$ & S \\ \cline{2-5} 
& $2$ & $1,12,14,16$ & $y_5$ & S \\ \cline{2-5} 
& $2$ & $2,12,14,16$ & $y_5$ & S \\ \cline{2-5} 
& $2$ & $3,3,1,3,12,30$ & $y_5$ & S \\ \cline{2-5} 
& $2$ & $3,3,2,3,12,30$ & $y_5$ & S \\ \cline{2-5} 
& $2$ & $1,3,12,36$ & $y_4$ & S \\ \cline{2-5} 
& $2$ & $2,3,12,36$ & $y_4$ & S \\ \cline{2-5} 
& $3$ & $1,3,4,7,30$ & $y_9$ & S \\ \cline{2-5} 
& $3$ & $2,3,4,7,30$ & $y_9$ & S \\ \cline{2-5} 
& $3$ & $1,3,5,42$ & $y_8$ & S \\ \cline{2-5} 
& $3$ & $2,3,5,42$ & $y_8$ & S \\ \cline{2-5} 
& $3$ & $1,3,5,7,18$ & $y_8$ & S \\ \cline{2-5} 
& $3$ & $2,3,5,7,18$ & $y_8$ & S \\ \cline{2-5} 
& $3$ & $3,3,1,4,12,30$ & $y_6$ & S \\ \cline{2-5} 
& $3$ & $3,3,2,4,12,30$ & $y_6$ & S \\ \cline{2-5} 
& $3$ & $4,4,1,3,12,30$ & $y_5$ & S \\ \cline{2-5} 
& $3$ & $4,4,2,3,12,30$ & $y_5$ & S \\ \cline{2-5} 
& $3$ & $1,3,12,42$ & $y_3$ & S \\ \cline{2-5} 
& $3$ & $2,3,12,42$ & $y_3$ & S \\ \cline{2-5} 
& $3$ & $1,3,7,12,18$ & $y_3$ & S \\ \cline{2-5} 
& $3$ & $2,3,7,12,18$ & $y_3$ & S \\ \cline{2-5} 
& $3$ & $1,3,4,5,7,12$ & $y_2$ & S \\ \cline{2-5} 
& $3$ & $2,3,4,5,7,12$ & $y_2$ & S \\ \cline{2-5} 
& $6$ & $3,4,7,12,14$ & $y_1$ & S \\ \hline
$R_9$ & $1$ & $1,2,3,4,28$ & $y_9$ & H \\ \cline{2-5} 
& $1$ & $1,1,3,4,28$ & $y_9$ & H \\ \cline{2-5} 
& $1$ & $2,2,3,4,28$ & $y_9$ & H \\ \cline{2-5} 
& $1$ & $1,1,1,4,8,24$ & $y_8$ & S \\ \cline{2-5} 
& $1$ & $1,1,2,4,8,24$ & $y_8$ & S \\ \cline{2-5} 
& $1$ & $2,2,1,4,8,24$ & $y_8$ & S \\ \cline{2-5} 
& $1$ & $2,2,2,4,8,24$ & $y_8$ & S \\ \cline{2-5} 
& $1$ & $1,1,1,7,24$ & $y_7$ & S \\ \cline{2-5} 
& $1$ & $1,1,2,7,24$ & $y_7$ & S \\ \cline{2-5} 
& $1$ & $1,1,1,4,12,24$ & $y_7$ & S \\ \cline{2-5} 
& $1$ & $1,1,2,4,12,24$ & $y_7$ & S \\ \cline{2-5} 
& $1$ & $2,2,1,7,24$ & $y_7$ & S \\ \cline{2-5} 
& $1$ & $2,2,2,7,24$ & $y_7$ & S \\ \hline
\end{tabular}}
\end{minipage}
\begin{minipage}{0.45\linewidth}
\centerline{
\begin{tabular}{|c|c|c|c|c|}\hline
$\Psi$ & case & $\bk$ & ST & S/H \\ \hline
$R_9$ & $1$ & $2,2,1,4,12,24$ & $y_7$ & S \\ \cline{2-5} 
& $1$ & $2,2,2,4,12,24$ & $y_7$ & S \\ \cline{2-5} 
& $2$ & $3,3,3,4,28$ & $y_9$ & H \\ \cline{2-5}
& $2$ & $4,4,3,4,28$ & $y_9$ & H \\ \cline{2-5} 
& $2$ & $3,3,1,4,8,24$ & $y_8$ & S \\ \cline{2-5} 
& $2$ & $3,3,2,4,8,24$ & $y_8$ & S \\ \cline{2-5} 
& $2$ & $4,4,1,4,8,24$ & $y_8$ & S \\ \cline{2-5} 
& $2$ & $4,4,2,4,8,24$ & $y_8$ & S \\ \cline{2-5} 
& $2$ & $4,4,1,7,24$ & $y_7$ & S \\ \cline{2-5} 
& $2$ & $4,4,2,7,24$ & $y_7$ & S \\ \cline{2-5} 
& $2$ & $4,4,1,4,12,24$ & $y_7$ & S \\ \cline{2-5} 
& $2$ & $4,4,2,4,12,24$ & $y_7$ & S \\ \cline{2-5} 
& $3$ & $3,3,1,7,24$ & $y_7$ & S \\ \cline{2-5} 
& $3$ & $3,3,2,7,24$ & $y_7$ & S \\ \cline{2-5} 
& $3$ & $3,3,1,4,12,24$ & $y_7$ & S \\ \cline{2-5} 
& $3$ & $3,3,2,4,12,24$ & $y_7$ & S \\ \cline{2-5} 
& $3$ & $1,3,12,42$ & $y_3$ & S \\ \cline{2-5} 
& $3$ & $2,3,12,42$ & $y_3$ & S \\ \cline{2-5} 
& $5$ & $1,2,3,16,18$ & $y_{10}$ & H \\ \cline{2-5} 
& $5$ & $1,1,3,16,18$ & $y_{10}$ & H \\ \cline{2-5} 
& $5$ & $2,2,3,16,18$ & $y_{10}$ & H \\ \cline{2-5} 
& $7$ & $3,3,3,16,18$ & $y_{10}$ & H \\ \cline{2-5} 
& $7$ & $4,4,3,16,18$ & $y_{10}$ & H \\ \cline{2-5} 
& $8$ & $30,42$ & $y_1$ & S \\ \cline{2-5} 
& $8$ & $3,12,18,24$ & $y_1$ & S \\ \hline
$R_{10}$ & $1$ & $1,2,3,4,36$ & $y_9$ & H \\ \cline{2-5} 
& $1$ & $1,1,3,4,36$ & $y_9$ & H \\ \cline{2-5} 
& $1$ & $2,2,3,4,36$ & $y_9$ & H \\ \cline{2-5} 
& $1$ & $1,1,1,4,12,24$ & $y_7$ & S \\ \cline{2-5} 
& $1$ & $1,1,2,4,12,24$ & $y_7$ & S \\ \cline{2-5} 
& $1$ & $2,2,1,4,12,24$ & $y_7$ & S \\ \cline{2-5} 
& $1$ & $2,2,2,4,12,24$ & $y_7$ & S \\ \cline{2-5} 
& $2$ & $3,3,3,4,36$ & $y_9$ & H \\ \cline{2-5} 
& $2$ & $4,4,3,4,36$ & $y_9$ & H \\ \cline{2-5}
& $2$ & $4,4,1,4,12,24$ & $y_7$ & S \\ \cline{2-5} 
& $2$ & $4,4,2,4,12,24$ & $y_7$ & S \\ \cline{2-5} 
& $3$ & $3,3,1,4,12,24$ & $y_7$ & S \\ \cline{2-5} 
& $3$ & $3,3,2,4,12,24$ & $y_7$ & S \\ \cline{2-5} 
& $8$ & $3,12,18,24$ & $y_1$ & S \\ \cline{2-5} 
&       &                    &             &    \\ \hline 
\end{tabular}}
\end{minipage}
\caption{K3 surface automorphisms from Salem numbers of degree $22$, Part 3.}
\label{tab:SvH3}
\end{table} 
%%%%%%%%%%%%%%%%%%%%%%%%%%%%%%%%%%%%%%%%%%%%%%%%%%%%%%%%%%%%%%%%%
\end{theorem}  
%%%%%%%%%%%%%%%%%%%%%%% begin proof %%%%%%%%%%%%%%%%%%%%%%%%%%%%%%%
{\it Proof}. 
Carrying out the tests for each polynomial $R_i(w)$ in 
Table \ref{tab:mcmullen} by computer, we have all possible values of $\bk$ 
as in Tables $\ref{tab:SvH1}$--$\ref{tab:SvH3}$ along with the 
information about cases and special traces. 
Each of these data leads to a K3 surface automorphism $f : X \to X$ via 
Theorem \ref{thm:great}.      \hfill $\Box$ 
%%%%%%%%%%%%%%%%%%%%%%%%% end proof %%%%%%%%%%%%%%%%%%%%%%%%%%%%%%%    
%%%%%%%%%%%%%%%%%%%%%% ss:minimum %%%%%%%%%%%%%%%%%%%%%%%%%%%%%%%%%
\subsection{Minimum Entropy} \label{ss:minimum}
%%%%%%%%%%%%%%%%%%%%%%%%%%%%%%%%%%%%%%%%%%%%%%%%%%%%%%%%%%%%%%%%%%%
McMullen \cite[Theorem (A.1)]{McMullen2} shows that any automorphism 
$f$ of a compact complex surface with a positive entropy $h(f) > 0$ has a 
lower bound $h(f) \ge \log \lambda_{\rL}$, where $\lambda_{\rL}$ is 
Lehmer's number in Example \ref{ex:lehmer}. 
In \cite{McMullen3} he obtains some non-projective K3 examples that actually 
attain this lower bound and he goes on to construct projective ones in 
\cite{McMullen4}. 
Our method enables us to synthesize a goodly number of non-projective 
K3 surface automorphisms with minimum entropy $\log \lambda_{\rL}$, 
which necessarily have Picard number $12$. 
%%%%%%%%%%%%%%%%%%%%%%%%%%%%%%%%%% sss:first %%%%%%%%%%%%%%%%%%%%%%%%%%
\subsubsection{Use of Matrix $\mbox{\boldmath $A$}$} \label{sss:first}
%%%%%%%%%%%%%%%%%%%%%%%%%%%%%%%%%%%%%%%%%%%%%%%%%%%%%%%%%%%%%%%%%%%%%%%
Suppose that $\Phi(w) = \LT(w) \cdot \CT_{\sbk}(w)$, where $\LT(w)$ is 
Lehmer's trace polynomial in \eqref{eqn:lehmer2} and $\CT_{\sbk}(w)$ is 
a product of cyclotomic trace polynomials of the form 
%%%%%%%%%%%%%%%%%%%%%% eqn:CTk2 %%%%%%%%%%%%%%%%%%%%%%%%%%%%%%%%%%%%%%%
\begin{equation} \label{eqn:CTk2} 
\CT_{\sbk}(w) := \prod_{k \in \sbk} \CT_k(w) \qquad 
\mbox{with} \quad \sum_{k \in \sbk} \deg \CT_k(w) = 5, 
\end{equation}
%%%%%%%%%%%%%%%%%%%%%%%%%%%%%%%%%%%%%%%%%%%%%%%%%%%%%%%%%%%%%%%%%%%%%%%
where $\bk$ is a set of positive integers whose elements $k$ come    
from Table \ref{tab:CTP} with $\deg \le 5$, while $\Psi(w) = R_i(w)$ is one 
of the Salem trace polynomials in Table \ref{tab:mcmullen}. 
We remark that $\bk$ is not a multi-set but an ordinary set because  
an assertion of Theorem \ref{thm:main2} rules out the occurrence of 
multiple elements. 
To construct K3 surface automorphisms we utilize the matrix $A$ and 
its modification $\tilde{A}$. 
Let $x_4 < x_3 < x_2 < x_1$ denote the roots in $[-2, \, 2]$ of Lehmer's trace 
polynomial $\LT(w)$ in \eqref{eqn:lehmer2}, whose numerical values are given by   
%%%%%%%%%%%%%%%%%%%%%%%%% eqn:lehconj %%%%%%%%%%%%%%%%%%%%%%%%%%%%%%%%%%%
\begin{equation} \label{eqn:lehconj}
x_4\approx -1.88660, \quad x_3\approx -1.46887, \quad 
x_2\approx -0.584663, \quad x_1\approx 0.913731. 
\end{equation}
%%%%%%%%%%%%%%%%%%%%%%%%%%%%%%%%%%%%%%%%%%%%%%%%%%%%%%%%%%%%%%%%%%%%%%%
%%%%%%%%%%%%%%%%%%%%% tab:minA %%%%%%%%%%%%%%%%%%%%%%%%%%%%%%%%%%%%%%%%
\begin{table}[hh]
\centerline{
\begin{tabular}{|c|c|c|c|c|c|c|c|}
\hline
           &          &      &                  &                                &                          &  \\[-3mm]
$\Psi$ & $\bk$ & ST & Dynkin type & $\tilde{\varphi}_1(z)$ & $\Tr \, \tilde{A}$ & S/H \\ 
\hline
$R_1$ & $4,20$ & $x_4$ &  $\rE_6 \oplus \rE_6$ & $(z-1)^4(z+1)^4(z^2+1)^2$ & $-1$ & H \\ 
\hline
$R_1$ & $4,6,7$ & $x_2$ & $\rD_{10}$ & $(z-1)^9(z+1)(z^2+1)$ & $7$ & S \\ 
\hline
$R_3$ & $3,15$ & $x_4$ & $\rA_2$ & $(z-1)^2(z^2+z+1) \, \rC_{15}(z)$ & $1$ & S \\ 
\hline
$R_3$ & $3,4,6,8$ & $x_3$ & $\rE_6 \oplus \rE_6$ & $(z-1)^4(z+1)^4(z^2+1)^2$ & $-1$ & S \\ 
\hline
$R_3$ & $4,6,18$ & $x_2$ & $\rD_{10}$ & $(z-1)^9(z+1)(z^2+1)$ & $7$ & S \\ 
\hline
$R_4$ & $3,4,9$ & $x_4$ & $\rD_{10}$ & $(z-1)^9(z+1)(z^2+1)$ & $7$ & H \\ 
\hline
$R_4$ & $4,24$ & $x_2$ & $\rE_8 \oplus \rA_2 \oplus \rA_2$ & $(z-1)^9(z+1)(z^2+1)$ & $7$ & S \\ 
\hline
$R_4$ & $4,20$ & $x_1$ & $\rE_6 \oplus \rE_6$ & $(z-1)^4(z+1)^4(z^2+1)^2$ & $-1$ & S \\ 
\hline
$R_5$ & $7,12$ & $x_3$ & $\rE_8 \oplus \rA_2 \oplus \rA_2$ & $(z-1)^9(z+1)(z^2+1)$ & $7$ & H \\ 
\hline
$R_9$ & $4,20$ & $x_4$ & $\rE_6 \oplus \rE_6$ & $(z-1)^4(z+1)^4(z^2+1)^2$ & $-1$ & H \\
\hline
$R_9$ & $12,18$ & $x_3$ & $\rE_8 \oplus \rA_2 \oplus \rA_2$ & $(z-1)^9(z+1)(z^2+1)$ & $7$ & H \\ 
\hline
$R_9$ & $4,30$ & $x_1$ & $\rE_8 \oplus \rA_2 \oplus \rA_2$ & $(z-1)^9(z+1)(z^2+1)$ & $7$ & H \\ 
\hline
$R_{10}$ & $4,16$ & $x_4$ & $\rE_6 \oplus \rE_6$ & $(z-1)^4(z+1)^4(z^2+1)^2$ & $-1$ & H \\ 
\hline
$R_{10}$ & $4,24$ & $x_2$ & $\rE_8 \oplus \rA_2 \oplus \rA_2$ & $(z-1)^9(z+1)(z^2+1)$ & $7$ & S \\ 
\hline
$R_{10}$ & $3,15$ & $x_1$ & $\rA_2$ & $(z-1)^2(z^2+z+1) \, \rC_{15}(z)$ & $1$ & S \\ 
\hline
\end{tabular}}
\caption{K3 surface automorphisms of minimum entropy from matrix $\tilde{A}$.} 
\label{tab:minA}
\end{table}
%%%%%%%%%%%%%%%%%%%%%%%%%%%%%%%%%%%%%%%%%%%%%%%%%%%%%%%%%%%%%%%%%%%%%%%
%%%%%%%%%%%%%%%%%%%%%%%%% thm:A-pattern %%%%%%%%%%%%%%%%%%%%%%%%%%%%%%%%%
\begin{theorem} \label{thm:A-pattern} 
In the above setting we have a hypergeometric K3 lattice such that the matrix $A$ 
admits a special eigenvalue $\delta$ conjugate to Lehmer's number  
$\lambda_{\rL}$, if and only if $\Psi(w) = R_i(w)$ and $\bk$ are as in 
Table $\ref{tab:minA}$, where the special trace $\tau := \delta + \delta^{-1}$ is 
given in the ``ST" column. 
The modified Hodge isometry $\tilde{A} := w_A \circ A : L \to L$ lifts to a K3 surface 
automorphism $f : X \to X$ of entropy $h(f) = \log \lambda_{\rL}$ with special 
trace $\tau(f) = \tau$ and Picard number $\rho(X) = 12$, having exceptional set 
$\cE(X)$ of Dynkin type indicated in Table $\ref{tab:minA}$. 
The matrix $\tilde{A}$ or the map $f^*|H^2(X, \bZ)$ has characteristic polynomial  
$\tilde{\varphi}(z) = \rL(z) \cdot \tilde{\varphi}_1(z)$, where $\rL(z)$ is Lehmer's 
polynomial in \eqref{eqn:lehmer1} and $\tilde{\varphi}_1(z)$ is given in Table $\ref{tab:minA}$.    
The value of $\Tr \, \tilde{A} = \Tr \, f^*|H^2(X, \bZ)$ is also included there.  
The ``S/H" column is explained in Theorems $\ref{thm:min-Siegel1}$--$\ref{thm:A2}$.    
\end{theorem} 
%%%%%%%%%%%%%%%%%%%%%%%%%% begin proof %%%%%%%%%%%%%%%%%%%%%%%%%%%%%%%%%%
{\it Proof}. 
First, pick out all pairs $(i, \bk)$ such that 
$\Phi(w) = \LT(w) \cdot \CT_{\sbk}(w)$ and $\Psi(w) = R_i(w)$ satisfy 
unimodularity condition \eqref{eqn:unim}, where $\bk$ must be subject to the 
degree constraint in \eqref{eqn:CTk2}. 
Secondly, determine all the correct solutions as well as their special traces 
by using the hyperbolic case of Theorem \ref{thm:main2}.    
Since the special trace comes from Lehmer's polynomial factor, the matrix 
$F = A$ is a positive Hodge isometry by Remark \ref{rem:positive}.    
Thirdly, for each solution run Algorithm \ref{algorithm} to find the 
set of positive roots $\vD^+$, its basis $\vD_{\rb}$, its Dynkin type  
and the modified matrix $\tilde{A}$. 
The characteristic polynomial $\tilde{\varphi}(z) = \varphi_0(z) \cdot \tilde{\varphi}_1(z)$ 
of $\tilde{A}$ is determined by the formulas \eqref{eqn:chi01} and \eqref{eqn:chi01t}  
together with the non-projectivity assumption $\varphi_0(z) = \rL(z)$, where 
$\chi(z) = \varphi(z)$ and $\chi_i(z) = \varphi_i(z)$.   
\hfill $\Box$ \par\medskip  
%%%%%%%%%%%%%%%%%%%%%%%%%% end proof %%%%%%%%%%%%%%%%%%%%%%%%%%%%%%%%%%
We illustrate how Algorithm \ref{algorithm} works for an entry in 
Table \ref{tab:minA}.  
%%%%%%%%%%%%%%%%%%%%%%%% ex:algorithm %%%%%%%%%%%%%%%%%%%%%%%%%%%%%%%%%%%%
\begin{example} \label{ex:algorithm} 
Consider the case where $\Phi(w) = \LT(w) \cdot \CT_{\sbk}(w)$ 
with $\bk = \{3, 4, 6, 8\}$ and $\Psi(w) = R_3(w)$ in Table \ref{tab:minA}.  
Steps $1$--$3$ of Algorithm \ref{algorithm} return us $72$ elements 
$\b0 \prec \bu_1 \prec \bu_2 \prec \cdots \prec \bu_{72}$ in the   
lexicographic order \eqref{eqn:order} for the set of positive roots $\vD^+$  
and then $12$ elements for the basis $\vD_{\rb}$, implying that the  
root system is of type $\rE_6 \oplus \rE_6$.  
If the basis $\vD_{\rb} = \{ \be_1, \dots, \be_6, \be_1', \dots, \be_6' \}$ is 
labeled as in Figure \ref{fig:E6E6} then the simple roots are given by  
%%%%%%%%%%%%%%% 
\begin{alignat*}{6}
\be_1 &= \bu_{23}, \quad & \be_2 &= \bu_{8}, \quad & \be_3 &= \bu_{1}, \quad & 
\be_4 &= \bu_{3}, \quad & \be_5 &= \bu_{16}, \quad & \be_6 &= \bu_{25}, \\ 
\be_1' &= \bu_{7}, \quad & \be_2' &= \bu_{24}, \quad & \be_3' &= \bu_{2}, \quad & 
\be_4' &= \bu_{5}, \quad & \be_5' &= \bu_{9}, \quad & \be_6' &= \bu_{35}. 
\end{alignat*}
%%%%%%%%%%%%%%
Steps $4$--$6$ tell us that the Weyl group element $w_A \in W$ 
bringing $A(\cK)$ back to $\cK$ is
%%%%%%% 
$$
w_A = \sigma_{5} \circ \sigma_{23} \circ \sigma_{35} \circ \sigma_{41} \circ 
\sigma_{62} \circ \sigma_{57} \circ \sigma_{72},
$$
%%%%%%% 
where $\sigma_j$ denotes the Picard-Lefschetz reflection corresponding to the 
$j$-th positive root $\bu_j$. 
Step $7$, that is, how the modified matrix $\tilde{A}$ acts on $\vD_{\rb}$ will 
be mentioned in \S \ref{sss:action}.  
%%%%%%%%%%%%%%%%%%%%% fig:E6E6 %%%%%%%%%%%%%%%%%%%%%%%%%%%%%%%%%%%%%%%%
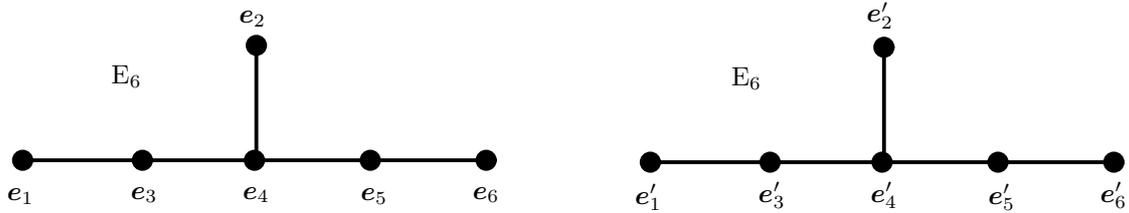
\begin{figure}[hh] 
\begin{center}
%WinTpicVersion4.32a
{\unitlength 0.1in%
\begin{picture}(57.7500,9.6000)(5.2000,-13.4000)%
% CIRCLE 0 0 0 0 Black White  
% 4 600 1210 620 1250 620 1250 620 1250
% 
\special{sh 1.000}%
\special{ia 600 1210 45 45 0.0000000 6.2831853}%
\special{pn 20}%
\special{ar 600 1210 45 45 0.0000000 6.2831853}%
% CIRCLE 0 0 0 0 Black White  
% 4 1220 1210 1240 1250 1240 1250 1240 1250
% 
\special{sh 1.000}%
\special{ia 1220 1210 45 45 0.0000000 6.2831853}%
\special{pn 20}%
\special{ar 1220 1210 45 45 0.0000000 6.2831853}%
% CIRCLE 0 0 0 0 Black White  
% 4 1800 1210 1820 1250 1820 1250 1820 1250
% 
\special{sh 1.000}%
\special{ia 1800 1210 45 45 0.0000000 6.2831853}%
\special{pn 20}%
\special{ar 1800 1210 45 45 0.0000000 6.2831853}%
% CIRCLE 0 0 0 0 Black White  
% 4 2400 1210 2420 1250 2420 1250 2420 1250
% 
\special{sh 1.000}%
\special{ia 2400 1210 45 45 0.0000000 6.2831853}%
\special{pn 20}%
\special{ar 2400 1210 45 45 0.0000000 6.2831853}%
% CIRCLE 0 0 0 0 Black White  
% 4 1810 610 1830 650 1830 650 1830 650
% 
\special{sh 1.000}%
\special{ia 1810 610 45 45 0.0000000 6.2831853}%
\special{pn 20}%
\special{ar 1810 610 45 45 0.0000000 6.2831853}%
% CIRCLE 0 0 0 0 Black White  
% 4 3000 1210 3020 1250 3020 1250 3020 1250
% 
\special{sh 1.000}%
\special{ia 3000 1210 45 45 0.0000000 6.2831853}%
\special{pn 20}%
\special{ar 3000 1210 45 45 0.0000000 6.2831853}%
% LINE 0 0 3 0 Black White  
% 2 620 1210 3000 1210
% 
\special{pn 20}%
\special{pa 620 1210}%
\special{pa 3000 1210}%
\special{fp}%
% LINE 0 0 3 0 Black White  
% 2 1810 610 1810 1220
% 
\special{pn 20}%
\special{pa 1810 610}%
\special{pa 1810 1220}%
\special{fp}%
% STR 2 0 3 0 Black White  
% 4 520 1340 520 1440 2 0 0 0
% $\be_1$
\put(5.2000,-14.4000){\makebox(0,0)[lb]{$\be_1$}}%
% STR 2 0 3 0 Black White  
% 4 1160 1330 1160 1430 2 0 0 0
% $\be_3$
\put(11.6000,-14.3000){\makebox(0,0)[lb]{$\be_3$}}%
% STR 2 0 3 0 Black White  
% 4 1740 1330 1740 1430 2 0 0 0
% $\be_4$
\put(17.4000,-14.3000){\makebox(0,0)[lb]{$\be_4$}}%
% STR 2 0 3 0 Black White  
% 4 2350 1340 2350 1440 2 0 0 0
% $\be_5$
\put(23.5000,-14.4000){\makebox(0,0)[lb]{$\be_5$}}%
% STR 2 0 3 0 Black White  
% 4 2930 1330 2930 1430 2 0 0 0
% $\be_6$
\put(29.3000,-14.3000){\makebox(0,0)[lb]{$\be_6$}}%
% STR 2 0 3 0 Black White  
% 4 1720 410 1720 510 2 0 0 0
% $\be_2$
\put(17.2000,-5.1000){\makebox(0,0)[lb]{$\be_2$}}%
% CIRCLE 0 0 0 0 Black White  
% 4 3850 1220 3870 1260 3870 1260 3870 1260
% 
\special{sh 1.000}%
\special{ia 3850 1220 45 45 0.0000000 6.2831853}%
\special{pn 20}%
\special{ar 3850 1220 45 45 0.0000000 6.2831853}%
% CIRCLE 0 0 0 0 Black White  
% 4 4470 1220 4490 1260 4490 1260 4490 1260
% 
\special{sh 1.000}%
\special{ia 4470 1220 45 45 0.0000000 6.2831853}%
\special{pn 20}%
\special{ar 4470 1220 45 45 0.0000000 6.2831853}%
% CIRCLE 0 0 0 0 Black White  
% 4 5050 1220 5070 1260 5070 1260 5070 1260
% 
\special{sh 1.000}%
\special{ia 5050 1220 45 45 0.0000000 6.2831853}%
\special{pn 20}%
\special{ar 5050 1220 45 45 0.0000000 6.2831853}%
% CIRCLE 0 0 0 0 Black White  
% 4 5650 1220 5670 1260 5670 1260 5670 1260
% 
\special{sh 1.000}%
\special{ia 5650 1220 45 45 0.0000000 6.2831853}%
\special{pn 20}%
\special{ar 5650 1220 45 45 0.0000000 6.2831853}%
% CIRCLE 0 0 0 0 Black White  
% 4 5060 620 5080 660 5080 660 5080 660
% 
\special{sh 1.000}%
\special{ia 5060 620 45 45 0.0000000 6.2831853}%
\special{pn 20}%
\special{ar 5060 620 45 45 0.0000000 6.2831853}%
% CIRCLE 0 0 0 0 Black White  
% 4 6250 1220 6270 1260 6270 1260 6270 1260
% 
\special{sh 1.000}%
\special{ia 6250 1220 45 45 0.0000000 6.2831853}%
\special{pn 20}%
\special{ar 6250 1220 45 45 0.0000000 6.2831853}%
% LINE 0 0 3 0 Black White  
% 2 3870 1220 6250 1220
% 
\special{pn 20}%
\special{pa 3870 1220}%
\special{pa 6250 1220}%
\special{fp}%
% LINE 0 0 3 0 Black White  
% 2 5060 620 5060 1230
% 
\special{pn 20}%
\special{pa 5060 620}%
\special{pa 5060 1230}%
\special{fp}%
% STR 2 0 3 0 Black White  
% 4 3770 1370 3770 1470 2 0 0 0
% $\be_1'$
\put(37.7000,-14.7000){\makebox(0,0)[lb]{$\be_1'$}}%
% STR 2 0 3 0 Black White  
% 4 4410 1360 4410 1460 2 0 0 0
% $\be_3'$
\put(44.1000,-14.6000){\makebox(0,0)[lb]{$\be_3'$}}%
% STR 2 0 3 0 Black White  
% 4 4990 1360 4990 1460 2 0 0 0
% $\be_4'$
\put(49.9000,-14.6000){\makebox(0,0)[lb]{$\be_4'$}}%
% STR 2 0 3 0 Black White  
% 4 5600 1370 5600 1470 2 0 0 0
% $\be_5'$
\put(56.0000,-14.7000){\makebox(0,0)[lb]{$\be_5'$}}%
% STR 2 0 3 0 Black White  
% 4 6180 1360 6180 1460 2 0 0 0
% $\be_6'$
\put(61.8000,-14.6000){\makebox(0,0)[lb]{$\be_6'$}}%
% STR 2 0 3 0 Black White  
% 4 4970 420 4970 520 2 0 0 0
% $\be_2'$
\put(49.7000,-5.2000){\makebox(0,0)[lb]{$\be_2'$}}%
% STR 2 0 3 0 Black White  
% 4 1060 730 1060 830 2 0 0 0
% $\rE_6$
\put(10.6000,-8.3000){\makebox(0,0)[lb]{$\rE_6$}}%
% STR 2 0 3 0 Black White  
% 4 4270 740 4270 840 2 0 0 0
% $\rE_6$
\put(42.7000,-8.4000){\makebox(0,0)[lb]{$\rE_6$}}%
\end{picture}}%
\end{center} 
\caption{Dynkin diagram of type $\rE_6 \oplus \rE_6$.}
\label{fig:E6E6}
\end{figure}
%%%%%%%%%%%%%%%%%%%%%%%%%%%%%%%%%%%%%%%%%%%%%%%%%%%%%%%%%%%%%%%%%%%%%   
\end{example}
%%%%%%%%%%%%%%%%%%%%%%%%%%%%%%%%%%%%%%%%%%%%%%%%%%%%%%%%%%%%%%%%%%%%%%
%%%%%%%%%%%%%%%%%%%%%%%%%%% rem:A-pattern %%%%%%%%%%%%%%%%%%%%%%%%%%%%%%%
\begin{remark} \label{rem:A-pattern} 
A careful inspection shows that all entries in Table \ref{tab:minA} fall into case 
$7$ of Table \ref{tab:hyp-A}. 
However, other cases can occur in other settings, although empirically they are 
rather rare.  
The following two examples are in cases $1$ and $9$ of Table \ref{tab:hyp-A}, 
respectively: 
%%%%%%%%%
\begin{alignat*}{2}
\Phi(w) &= \LT(w) \cdot \CT_{18}(w) \cdot \CT_5(w), \quad &  
\Psi(w) &= \{ (w+1)(w^2-4) -1\} \cdot \CT_{60}(w), \\[1mm]
\Phi(w) &= \LT(w) \cdot \CT_{15}(w) \cdot \CT_4(w), \quad &  
\Psi(w) &= \{ w(w^2-1)(w^2-3)(w^2-4) -1\} \cdot \CT_{24}(w).    
\end{alignat*}
%%%%%%%%%%
\end{remark}
%%%%%%%%%%%%%%%%%%%%%%%%%%% sss:second %%%%%%%%%%%%%%%%%%%%%%%%%%%%%%%%
\subsubsection{Use of Matrix $\mbox{\boldmath $B$}$} \label{sss:second}
%%%%%%%%%%%%%%%%%%%%%%%%%%%%%%%%%%%%%%%%%%%%%%%%%%%%%%%%%%%%%%%%%%%%%
Suppose that $\Phi(w)$ is a product of cyclotomic trace polynomials of the 
form \eqref{eqn:CTk}, while $\Psi(w)$ is a product of Lehmer's trace 
polynomial and cyclotomic trace polynomials 
%%%%%%%
$$
\Psi(w) = \LT(w) \cdot \CT_{\sbl}(w) \qquad \mbox{such that} \qquad 
\sum_{l \in \sbl} \deg \CT_l(w) = 6. 
$$
%%%%%%%
To construct K3 surface automorphisms we utilize the matrix 
$B$ and its modification $\tilde{B}$. 
By Theorem \ref{thm:main3} any element of $\bl$ must be simple and  
the unramifiedness condition in \eqref{eqn:unim} reduces the possibilities 
of $\bl$ considerably, confining $\Psi(w)$ into one of the following possibilities. 
%%%%%%%%%%%%%%
\begin{alignat*}{2}
L_1(w)&=\LT(w) \cdot \CT_{21}(w), \qquad & 
L_2(w)&=\LT(w) \cdot \CT_{28}(w), \\[1mm]
L_3(w)&=\LT(w) \cdot \CT_{36}(w), \qquad &
L_4(w)&=\LT(w) \cdot \CT_{42}(w), \\
L_5(w)&=\LT(w) \cdot \CT_{12}(w) \cdot \CT_{15}(w), \qquad & 
L_6(w)&=\LT(w) \cdot \CT_{12}(w) \cdot \CT_{20}(w), \\[1mm]       
L_7(w)&=\LT(w) \cdot \CT_{12}(w) \cdot \CT_{24}(w), \qquad &        
L_8(w)&=\LT(w) \cdot \CT_{12}(w) \cdot \CT_{30}(w).        
\end{alignat*}
%%%%%%%%%%%%%
%%%%%%%%%%%%%%%%%%%%%%%%% thm:B-pattern %%%%%%%%%%%%%%%%%%%%%%%%%%%%%%%%%
\begin{theorem} \label{thm:B-pattern} 
In the above setting we have a hypergeometric K3 lattice such that the matrix $B$ 
admits a special eigenvalue $\delta$ conjugate to Lehmer's number $\lambda_{\rL}$, 
if and only if $\Psi(w)$ and $\bk$ are as in Table $\ref{tab:minB}$, where the special 
trace $\tau := \delta + \delta^{-1}$ is given in the ``ST" column while the ``case" 
column refers to the case in Table $\ref{tab:hyp-B}$.  
The modified Hodge isometry $\tilde{B} := w_B \circ B : L \to L$ lifts to a K3 
surface automorphism $f : X \to X$ of entropy $h(f) = \log \lambda_{\rL}$ 
with special trace $\tau(f) = \tau$ and Picard number $\rho(X) = 12$, 
having exceptional set $\cE(X)$ of Dynkin type indacated in Table $\ref{tab:minB}$. 
The matrix $\tilde{B}$ or the map $f^*|H^2(X, \bZ)$ has characteristic polynomial   
$\tilde{\psi}(z) = \rL(z) \cdot \tilde{\psi}_1(z)$ with $\tilde{\psi}_1(z)$ given in 
Table $\ref{tab:minB}$.  
The value of $\Tr \, \tilde{B} = \Tr \, f^*|H^2(X, \bZ)$ is also included there.  
The ``S/H" column is explained in Theorems $\ref{thm:min-Siegel1}$ and 
$\ref{thm:min-Siegel2}$.    
%%%%%%%%%%%%%%%%%%%%%%% tab:minB %%%%%%%%%%%%%%%%%%%%%%%%%%%%%%%%%%%%%%%
\begin{table}[hh]
\centerline{
\begin{tabular}{|c|c|c|c|c|c|c|c|}
\hline
           &         &          &      &                  &                            &                          &       \\[-3mm]
$\Psi$ & case & $\bk$ & ST & Dynkin type & $\tilde{\psi}_1(z)$ & $\Tr \, \tilde{B}$ & S/H \\
\hline
$L_3$ & $2$ & $3,6,10,21$ & $x_4$ & $\rE_6 \oplus \rE_6$ & $(z-1)^4(z+1)^4(z^2+1)^2$ & $-1$ & H \\
\hline
$L_3$ & $2$ & $1,6,8,28$ & $x_2$ & $\rE_6 \oplus \rE_6$ & $(z-1)^4(z+1)^4(z^2+1)^2$  & $-1$ & S \\
\hline
$L_3$ & $2$ & $2,6,8,28$ & $x_2$ & $\rE_6 \oplus \rE_6$ & $(z-1)^4(z+1)^4(z^2+1)^2$ & $-1$ & S \\
\hline
$L_3$ & $2$ & $8,10,42$ & $x_1$ & $\rE_6 \oplus \rE_6$ & $(z-1)^4(z+1)^4(z^2+1)^2$ & $-1$ & S \\
\hline
$L_3$ & $3$ & $1,3,5,6,11$ & $x_2$ & $\rE_6 \oplus \rE_6$ & $(z-1)^4(z+1)^4(z^2+1)^2$ & $-1$ & S \\
\hline
$L_3$ & $3$ & $1,3,7,11$ & $x_2$  & $\rE_6 \oplus \rE_6$ & $(z-1)^4(z+1)^4(z^2+1)^2$ & $-1$ & S \\
\hline
$L_3$ & $3$ & $2,3,5,6,11$ & $x_2$  & $\rE_6 \oplus \rE_6$ & $(z-1)^4(z+1)^4(z^2+1)^2$ & $-1$ & S \\ 
\hline
$L_3$ & $3$ & $2,3,7,11$ & $x_2$  & $\rE_6 \oplus \rE_6$ & $(z-1)^4(z+1)^4(z^2+1)^2$ & $-1$ & S \\
\hline
$L_3$ & $3$ & $50$ & $x_2$ & $\rE_6 \oplus \rE_6$ & $(z-1)^4(z+1)^4(z^2+1)^2$ & $-1$ & S \\
\hline
$L_6$ & $1$ & $1,1,8,13$ & $x_4$  & $\rE_8 \oplus \rA_2 \oplus \rA_2$ & $(z-1)^9(z+1)(z^2+1)$ & $7$ & S \\
\hline
$L_6$ & $1$ & $1,2,8,13$ & $x_4$  & $\rE_8 \oplus \rA_2 \oplus \rA_2$ & $(z-1)^9(z+1)(z^2+1)$ & $7$ & S \\
\hline
$L_6$ & $1$ & $2,2,8,13$ & $x_4$  & $\rE_8 \oplus \rA_2 \oplus \rA_2$ & $(z-1)^9(z+1)(z^2+1)$ & 7 & S \\
\hline
$L_7$ & $1$ & $1,1,17$ & $x_4$  & $\rE_6 \oplus \rE_6$ & $(z-1)^4(z+1)^4(z^2+1)^2$ & $-1$ & H \\
\hline
$L_7$ & $1$ & $1,2,17$ & $x_4$  & $\rE_6 \oplus \rE_6$ & $(z-1)^4(z+1)^4(z^2+1)^2$ & $-1$ & H \\
\hline
$L_7$ & $1$ & $2,2,17$ & $x_4$  & $\rE_6 \oplus \rE_6$ & $(z-1)^4(z+1)^4(z^2+1)^2$ & $-1$ & H \\
\hline
$L_7$ & $1$ & $1,27$ & $x_4$  & $\rE_6 \oplus \rE_6$ & $(z-1)^4(z+1)^4(z^2+1)^2$ & $-1$ & H \\
\hline
$L_7$ & $1$ & $2,27$ & $x_4$ & $\rE_6 \oplus \rE_6$ & $(z-1)^4(z+1)^4(z^2+1)^2$  & $-1$ & H \\
\hline
$L_7$ & $3$ & $16,42$ & $x_2$  & $\rE_6 \oplus \rE_6$ & $(z-1)^4(z+1)^4(z^2+1)^2$ & $-1$ & S \\
\hline
$L_7$ & $3$ & $5,20,30$ & $x_2$  & $\rE_6 \oplus \rE_6$ & $(z-1)^4(z+1)^4(z^2+1)^2$ & $-1$ & S \\
\hline
$L_8$ & $1$ & $1,1,1,7,16$ & $x_3$ & $\rE_8 \oplus \rA_2 \oplus \rA_2$ & $(z-1)^9(z+1)(z^2+1)$  & $7$ & H \\
\hline
$L_8$ & $1$ & $1,1,2,7,16$ & $x_3$  & $\rE_8 \oplus \rA_2 \oplus \rA_2$ & $(z-1)^9(z+1)(z^2+1)$ & $7$ & H \\
\hline
$L_8$ & $1$ & $2,2,1,7,16$ & $x_3$ & $\rE_8 \oplus \rA_2 \oplus \rA_2$ & $(z-1)^9(z+1)(z^2+1)$  & $7$ & H \\
\hline
$L_8$ & $1$ & $2,2,2,7,16$ & $x_3$  & $\rE_8 \oplus \rA_2 \oplus \rA_2$ & $(z-1)^9(z+1)(z^2+1)$ & $7$ & H \\
\hline
$L_8$ & $2$ & $7,9,20$ & $x_3$ & $\rE_8 \oplus \rA_2 \oplus \rA_2$ & $(z-1)^9(z+1)(z^2+1)$  & $7$ & H \\
\hline
\end{tabular}}
\caption{K3 surface automorphisms of minimum entropy from matrix $\tilde{B}$.}
\label{tab:minB}
\end{table}
%%%%%%%%%%%%%%%%%%%%%%%%%%%%%%%%%%%%%%%%%%%%%%%%%%%%%%%%%%%%%%%%%%%%%%%%%%%
\end{theorem} 
%%%%%%%%%%%%%%%%%%%%%%%% begin proof %%%%%%%%%%%%%%%%%%%%%%%%%%%%%%%%%%
{\it Proof}. 
First, pick out all pairs $(\bk, i)$ such that 
$\Phi(w) = \CT_{\sbk}(w)$ and $\Psi(w) = L_i(w)$ satisfy 
unimodularity condition \eqref{eqn:unim}, where $\bk$ must be subject to the 
degree constraint in \eqref{eqn:CTk}. 
Secondly, determine all the correct solutions as well as their special traces 
by using Theorem \ref{thm:main3}.    
Since the special trace comes from Lehmer's polynomial factor, 
the matrix $F = B$ is a positive Hodge isometry by Remark \ref{rem:positive}.    
Thirdly, for each solution run Algorithm \ref{algorithm} to find the set of 
positive roots $\vD^+$, its basis $\vD_{\rb}$, its Dynkin type  
and the modified matrix $\tilde{B}$.  
The characteristic polynomial $\tilde{\psi}(z) = \psi_0(z) \cdot \tilde{\psi}_1(z)$ 
of $\tilde{B}$ is determined by the formulas \eqref{eqn:chi01} and \eqref{eqn:chi01t} 
together with the non-projectivity assumption $\psi_0(z) = \rL(z)$, where 
$\chi(z) = \psi(z)$ and $\chi_i(z) = \psi_i(z)$.   
\hfill $\Box$ \par\medskip 
%%%%%%%%%%%%%%%%%%%%%%%% end proof %%%%%%%%%%%%%%%%%%%%%%%%%%%%%%%%%%%%
In Table \ref{tab:minB} only four polynomials $L_3(w)$, $L_6(w)$, $L_7(w)$, $L_8(w)$ 
appear as $\Psi(w)$.     
%%%%%%%%%%%%%%%%%%%%%%%% ss:action %%%%%%%%%%%%%%%%%%%%%%%%%%%%%%%%%%%
\subsubsection{Action on Dynkin Diagram} \label{sss:action}
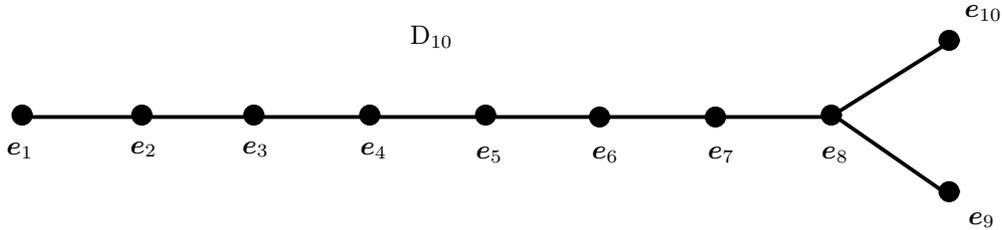
\begin{figure}[hh] 
\centerline{
%WinTpicVersion4.32a
{\unitlength 0.1in%
\begin{picture}(49.8000,10.8000)(5.2000,-16.7000)%
% CIRCLE 0 0 0 0 Black White  
% 4 600 1210 620 1250 620 1250 620 1250
% 
\special{sh 1.000}%
\special{ia 600 1210 45 45 0.0000000 6.2831853}%
\special{pn 20}%
\special{ar 600 1210 45 45 0.0000000 6.2831853}%
% CIRCLE 0 0 0 0 Black White  
% 4 1220 1210 1240 1250 1240 1250 1240 1250
% 
\special{sh 1.000}%
\special{ia 1220 1210 45 45 0.0000000 6.2831853}%
\special{pn 20}%
\special{ar 1220 1210 45 45 0.0000000 6.2831853}%
% CIRCLE 0 0 0 0 Black White  
% 4 1800 1210 1820 1250 1820 1250 1820 1250
% 
\special{sh 1.000}%
\special{ia 1800 1210 45 45 0.0000000 6.2831853}%
\special{pn 20}%
\special{ar 1800 1210 45 45 0.0000000 6.2831853}%
% CIRCLE 0 0 0 0 Black White  
% 4 2400 1210 2420 1250 2420 1250 2420 1250
% 
\special{sh 1.000}%
\special{ia 2400 1210 45 45 0.0000000 6.2831853}%
\special{pn 20}%
\special{ar 2400 1210 45 45 0.0000000 6.2831853}%
% CIRCLE 0 0 0 0 Black White  
% 4 4790 1210 4810 1250 4810 1250 4810 1250
% 
\special{sh 1.000}%
\special{ia 4790 1210 45 45 0.0000000 6.2831853}%
\special{pn 20}%
\special{ar 4790 1210 45 45 0.0000000 6.2831853}%
% CIRCLE 0 0 0 0 Black White  
% 4 3000 1210 3020 1250 3020 1250 3020 1250
% 
\special{sh 1.000}%
\special{ia 3000 1210 45 45 0.0000000 6.2831853}%
\special{pn 20}%
\special{ar 3000 1210 45 45 0.0000000 6.2831853}%
% STR 2 0 3 0 Black White  
% 4 520 1340 520 1440 2 0 0 0
% $\be_1$
\put(5.2000,-14.4000){\makebox(0,0)[lb]{$\be_1$}}%
% STR 2 0 3 0 Black White  
% 4 1160 1330 1160 1430 2 0 0 0
% $\be_2$
\put(11.6000,-14.3000){\makebox(0,0)[lb]{$\be_2$}}%
% STR 2 0 3 0 Black White  
% 4 1740 1330 1740 1430 2 0 0 0
% $\be_3$
\put(17.4000,-14.3000){\makebox(0,0)[lb]{$\be_3$}}%
% STR 2 0 3 0 Black White  
% 4 2350 1340 2350 1440 2 0 0 0
% $\be_4$
\put(23.5000,-14.4000){\makebox(0,0)[lb]{$\be_4$}}%
% STR 2 0 3 0 Black White  
% 4 2950 1360 2950 1460 2 0 0 0
% $\be_5$
\put(29.5000,-14.6000){\makebox(0,0)[lb]{$\be_5$}}%
% STR 2 0 3 0 Black White  
% 4 4740 1360 4740 1460 2 0 0 0
% $\be_8$
\put(47.4000,-14.6000){\makebox(0,0)[lb]{$\be_8$}}%
% CIRCLE 0 0 0 0 Black White  
% 4 5400 1610 5420 1650 5420 1650 5420 1650
% 
\special{sh 1.000}%
\special{ia 5400 1610 45 45 0.0000000 6.2831853}%
\special{pn 20}%
\special{ar 5400 1610 45 45 0.0000000 6.2831853}%
% CIRCLE 0 0 0 0 Black White  
% 4 5400 820 5420 860 5420 860 5420 860
% 
\special{sh 1.000}%
\special{ia 5400 820 45 45 0.0000000 6.2831853}%
\special{pn 20}%
\special{ar 5400 820 45 45 0.0000000 6.2831853}%
% CIRCLE 0 0 0 0 Black White  
% 4 4190 1220 4210 1260 4210 1260 4210 1260
% 
\special{sh 1.000}%
\special{ia 4190 1220 45 45 0.0000000 6.2831853}%
\special{pn 20}%
\special{ar 4190 1220 45 45 0.0000000 6.2831853}%
% CIRCLE 0 0 0 0 Black White  
% 4 3590 1220 3610 1260 3610 1260 3610 1260
% 
\special{sh 1.000}%
\special{ia 3590 1220 45 45 0.0000000 6.2831853}%
\special{pn 20}%
\special{ar 3590 1220 45 45 0.0000000 6.2831853}%
% STR 2 0 3 0 Black White  
% 4 5480 620 5480 720 2 0 0 0
% $\be_{10}$
\put(54.8000,-7.2000){\makebox(0,0)[lb]{$\be_{10}$}}%
% STR 2 0 3 0 Black White  
% 4 5500 1700 5500 1800 2 0 0 0
% $\be_9$
\put(55.0000,-18.0000){\makebox(0,0)[lb]{$\be_9$}}%
% STR 2 0 3 0 Black White  
% 4 3550 1360 3550 1460 2 0 0 0
% $\be_6$
\put(35.5000,-14.6000){\makebox(0,0)[lb]{$\be_6$}}%
% STR 2 0 3 0 Black White  
% 4 4150 1360 4150 1460 2 0 0 0
% $\be_7$
\put(41.5000,-14.6000){\makebox(0,0)[lb]{$\be_7$}}%
% STR 2 0 3 0 Black White  
% 4 2610 760 2610 860 2 0 0 0
% $\rD_{10}$
\put(26.1000,-8.6000){\makebox(0,0)[lb]{$\rD_{10}$}}%
% LINE 0 0 3 0 Black White  
% 2 580 1220 4790 1220
% 
\special{pn 20}%
\special{pa 580 1220}%
\special{pa 4790 1220}%
\special{fp}%
% LINE 0 0 3 0 Black White  
% 2 4800 1210 5400 830
% 
\special{pn 20}%
\special{pa 4800 1210}%
\special{pa 5400 830}%
\special{fp}%
% LINE 0 0 3 0 Black White  
% 2 4800 1210 5400 1630
% 
\special{pn 20}%
\special{pa 4800 1210}%
\special{pa 5400 1630}%
\special{fp}%
\end{picture}}%
}
\caption{Dynkin diagram of type $\rD_{10}$.} 
\label{fig:D10}
\end{figure}

%%%%%%%%%%%%%%%%%%%%%%%%%%%%%%%%%%%%%%%%%%%%%%%%%%%%%%%%%%%%%%%%%%%%%
Denote by $\tilde{F}$ the modified matrix $\tilde{A}$ or $\tilde{B}$.  
We are interested in how $\tilde{F}$ acts on the simple roots $\vD_{\rb}$. 
Let $\tilde{\chi}_1(z)$ denote the polynomial $\tilde{\varphi}_1(z)$ in 
Table \ref{tab:minA} or $\tilde{\psi}_1(z)$ in Table \ref{tab:minB}. 
Then $\tilde{\chi}_1(z)$ is divisible by the characteristic polynomial of 
$\tilde{F}|\mathrm{Span} \, \vD_{\rb}$.  
Observe that Tables \ref{tab:minA} and \ref{tab:minB} contain only root 
systems of types 
%%%%%%%%%%%%%%%%%%%%%%%% eqn:dynkin %%%%%%%%%%%%%%%%%%%%%%%%%%%%%%%%%%
\begin{equation} \label{eqn:dynkin}
\rE_6 \oplus \rE_6, \qquad \rE_8 \oplus \rA_2 \oplus \rA_2, \qquad 
\rD_{10}, \qquad \rA_2.   
\end{equation}
%%%%%%%%%%%%%%%%%%%%%%%%%%%%%%%%%%%%%%%%%%%%%%%%%%%%%%%%%%%%%%%%%%%% 
With the polynomials $\tilde{\chi}_1(z)$ in Tables \ref{tab:minA} and \ref{tab:minB}
we can identify the action of $\tilde{F}$ on $\vD_{\rb}$ in a unique way by 
enumerating all Dynkin automorphisms of types \eqref{eqn:dynkin} and their 
characteristic polynomials.  
This may not be true if the root system is different from \eqref{eqn:dynkin} 
or if $\tilde{\chi}_1(z)$ is a different polynomial.    
Moreover, it is a consequence of direct calculations that all entries of the 
same Dynkin type have the same polynomial $\tilde{\chi}_1(z)$. 
We do not know if this is just by accident or with any reason. 
In any case, within Tables \ref{tab:minA} and \ref{tab:minB}, 
how $\tilde{F}$ acts on $\vD_{\rb}$ depends only on the Dynkin type 
of $\vD_{\rb}$.  
%%%%%%%%%%%%%%%%%%%%%%%%% fig:E8A2A2 %%%%%%%%%%%%%%%%%%%%%%%%%%%%%%%%%%
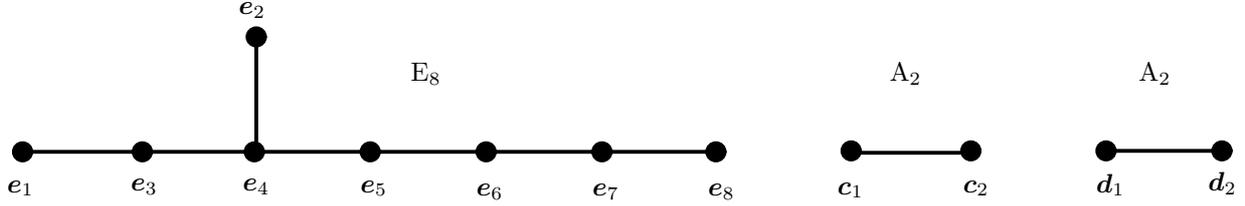
\begin{figure}[hh] 
\begin{center}
%WinTpicVersion4.32a
{\unitlength 0.1in%
\begin{picture}(63.3500,9.5000)(5.2000,-13.3000)%
% CIRCLE 0 0 0 0 Black White  
% 4 600 1210 620 1250 620 1250 620 1250
% 
\special{sh 1.000}%
\special{ia 600 1210 45 45 0.0000000 6.2831853}%
\special{pn 20}%
\special{ar 600 1210 45 45 0.0000000 6.2831853}%
% CIRCLE 0 0 0 0 Black White  
% 4 1220 1210 1240 1250 1240 1250 1240 1250
% 
\special{sh 1.000}%
\special{ia 1220 1210 45 45 0.0000000 6.2831853}%
\special{pn 20}%
\special{ar 1220 1210 45 45 0.0000000 6.2831853}%
% CIRCLE 0 0 0 0 Black White  
% 4 1800 1210 1820 1250 1820 1250 1820 1250
% 
\special{sh 1.000}%
\special{ia 1800 1210 45 45 0.0000000 6.2831853}%
\special{pn 20}%
\special{ar 1800 1210 45 45 0.0000000 6.2831853}%
% CIRCLE 0 0 0 0 Black White  
% 4 2400 1210 2420 1250 2420 1250 2420 1250
% 
\special{sh 1.000}%
\special{ia 2400 1210 45 45 0.0000000 6.2831853}%
\special{pn 20}%
\special{ar 2400 1210 45 45 0.0000000 6.2831853}%
% CIRCLE 0 0 0 0 Black White  
% 4 1810 610 1830 650 1830 650 1830 650
% 
\special{sh 1.000}%
\special{ia 1810 610 45 45 0.0000000 6.2831853}%
\special{pn 20}%
\special{ar 1810 610 45 45 0.0000000 6.2831853}%
% CIRCLE 0 0 0 0 Black White  
% 4 3000 1210 3020 1250 3020 1250 3020 1250
% 
\special{sh 1.000}%
\special{ia 3000 1210 45 45 0.0000000 6.2831853}%
\special{pn 20}%
\special{ar 3000 1210 45 45 0.0000000 6.2831853}%
% LINE 0 0 3 0 Black White  
% 2 1810 610 1810 1220
% 
\special{pn 20}%
\special{pa 1810 610}%
\special{pa 1810 1220}%
\special{fp}%
% STR 2 0 3 0 Black White  
% 4 520 1340 520 1440 2 0 0 0
% $\be_1$
\put(5.2000,-14.4000){\makebox(0,0)[lb]{$\be_1$}}%
% STR 2 0 3 0 Black White  
% 4 1160 1330 1160 1430 2 0 0 0
% $\be_3$
\put(11.6000,-14.3000){\makebox(0,0)[lb]{$\be_3$}}%
% STR 2 0 3 0 Black White  
% 4 1740 1330 1740 1430 2 0 0 0
% $\be_4$
\put(17.4000,-14.3000){\makebox(0,0)[lb]{$\be_4$}}%
% STR 2 0 3 0 Black White  
% 4 2350 1340 2350 1440 2 0 0 0
% $\be_5$
\put(23.5000,-14.4000){\makebox(0,0)[lb]{$\be_5$}}%
% STR 2 0 3 0 Black White  
% 4 2950 1360 2950 1460 2 0 0 0
% $\be_6$
\put(29.5000,-14.6000){\makebox(0,0)[lb]{$\be_6$}}%
% STR 2 0 3 0 Black White  
% 4 1720 410 1720 510 2 0 0 0
% $\be_2$
\put(17.2000,-5.1000){\makebox(0,0)[lb]{$\be_2$}}%
% CIRCLE 0 0 0 0 Black White  
% 4 4890 1205 4910 1245 4910 1245 4910 1245
% 
\special{sh 1.000}%
\special{ia 4890 1205 45 45 0.0000000 6.2831853}%
\special{pn 20}%
\special{ar 4890 1205 45 45 0.0000000 6.2831853}%
% CIRCLE 0 0 0 0 Black White  
% 4 5510 1205 5530 1245 5530 1245 5530 1245
% 
\special{sh 1.000}%
\special{ia 5510 1205 45 45 0.0000000 6.2831853}%
\special{pn 20}%
\special{ar 5510 1205 45 45 0.0000000 6.2831853}%
% CIRCLE 0 0 0 0 Black White  
% 4 4190 1210 4210 1250 4210 1250 4210 1250
% 
\special{sh 1.000}%
\special{ia 4190 1210 45 45 0.0000000 6.2831853}%
\special{pn 20}%
\special{ar 4190 1210 45 45 0.0000000 6.2831853}%
% CIRCLE 0 0 0 0 Black White  
% 4 6210 1205 6230 1245 6230 1245 6230 1245
% 
\special{sh 1.000}%
\special{ia 6210 1205 45 45 0.0000000 6.2831853}%
\special{pn 20}%
\special{ar 6210 1205 45 45 0.0000000 6.2831853}%
% CIRCLE 0 0 0 0 Black White  
% 4 3600 1210 3620 1250 3620 1250 3620 1250
% 
\special{sh 1.000}%
\special{ia 3600 1210 45 45 0.0000000 6.2831853}%
\special{pn 20}%
\special{ar 3600 1210 45 45 0.0000000 6.2831853}%
% CIRCLE 0 0 0 0 Black White  
% 4 6810 1205 6830 1245 6830 1245 6830 1245
% 
\special{sh 1.000}%
\special{ia 6810 1205 45 45 0.0000000 6.2831853}%
\special{pn 20}%
\special{ar 6810 1205 45 45 0.0000000 6.2831853}%
% STR 2 0 3 0 Black White  
% 4 4820 1350 4820 1450 2 0 0 0
% $\bc_1$
\put(48.2000,-14.5000){\makebox(0,0)[lb]{$\bc_1$}}%
% STR 2 0 3 0 Black White  
% 4 5470 1340 5470 1440 2 0 0 0
% $\bc_2$
\put(54.7000,-14.4000){\makebox(0,0)[lb]{$\bc_2$}}%
% STR 2 0 3 0 Black White  
% 4 6160 1340 6160 1440 2 0 0 0
% $\bd_1$
\put(61.6000,-14.4000){\makebox(0,0)[lb]{$\bd_1$}}%
% STR 2 0 3 0 Black White  
% 4 6740 1330 6740 1430 2 0 0 0
% $\bd_2$
\put(67.4000,-14.3000){\makebox(0,0)[lb]{$\bd_2$}}%
% LINE 0 0 3 0 Black White  
% 2 590 1210 4190 1210
% 
\special{pn 20}%
\special{pa 590 1210}%
\special{pa 4190 1210}%
\special{fp}%
% STR 2 0 3 0 Black White  
% 4 3550 1360 3550 1460 2 0 0 0
% $\be_7$
\put(35.5000,-14.6000){\makebox(0,0)[lb]{$\be_7$}}%
% STR 2 0 3 0 Black White  
% 4 4150 1360 4150 1460 2 0 0 0
% $\be_8$
\put(41.5000,-14.6000){\makebox(0,0)[lb]{$\be_8$}}%
% LINE 0 0 3 0 Black White  
% 2 4920 1215 5550 1215
% 
\special{pn 20}%
\special{pa 4920 1215}%
\special{pa 5550 1215}%
\special{fp}%
% LINE 0 0 3 0 Black White  
% 2 6210 1205 6820 1205
% 
\special{pn 20}%
\special{pa 6210 1205}%
\special{pa 6820 1205}%
\special{fp}%
% STR 2 0 3 0 Black White  
% 4 2610 760 2610 860 2 0 0 0
% $\rE_8$
\put(26.1000,-8.6000){\makebox(0,0)[lb]{$\rE_8$}}%
% STR 2 0 3 0 Black White  
% 4 5090 760 5090 860 2 0 0 0
% $\rA_2$
\put(50.9000,-8.6000){\makebox(0,0)[lb]{$\rA_2$}}%
% STR 2 0 3 0 Black White  
% 4 6380 760 6380 860 2 0 0 0
% $\rA_2$
\put(63.8000,-8.6000){\makebox(0,0)[lb]{$\rA_2$}}%
\end{picture}}%
\end{center}
\caption{Dynkin diagram of type $\rE_8 \oplus \rA_2 \oplus \rA_2$.} 
\label{fig:E8A2A2}
\end{figure}
%%%%%%%%%%%%%%%%%%%%%%%%%%%%%%%%%%%%%%%%%%%%%%%%%%%%%%%%%%%%%%%%%%%%     
%%%%%%%%%%%%%%%%%%%%%%%%%%% obs %%%%%%%%%%%%%%%%%%%%%%%%%%%%%%%%%%%%%
\begin{observation} \label{obs} 
According to Dynkin types we have the following observations, where 
$(c_1, c_2, \dots, c_k)$ stands for the cyclic permutation 
$c_1 \to c_2 \to \cdots \to c_k \to c_1$.   
\begin{enumerate} 
\setlength{\itemsep}{-1pt}
\item In case of type $\rE_6 \oplus \rE_6$ the matrix $\tilde{F}$ acts on the 
simple roots in Figure \ref{fig:E6E6} by 
%%%%%%%%%%%%%%%%%%%%%%%%%%%%%% eqn:cyclic1 %%%%%%%%%%%%%%%%%%%%%%%%%%%%%%%%%%% 
\begin{equation*}
(\be_1, \be_1', \be_6, \be_6')(\be_3, \be_3', \be_5, \be_5')(\be_2, \be_2')(\be_4, \be_4').  
\end{equation*}
%%%%%%%%%%%%%%%%%%%%%%%%%%%%%%%%%%%%%%%%%%%%%%%%%%%%%%%%%%%%%%%%%%%%%%%%%%%
In particular $\tilde{F}$ exchanges the two connected $\rE_6$-components.  
\item In case of type $\rE_8 \oplus \rA_2 \oplus \rA_2$ the matrix $\tilde{F}$ fixes    
$\be_1, \dots, \be_8$ in the $\rE_8$-component while it acts on the simple roots in the 
$\rA_2 \oplus \rA_2$-component by  
$(\bc_1, \bd_1, \bc_2, \bd_2)$ in Figure \ref{fig:E8A2A2}.  
\item In case of type $\rD_{10}$ the matrix $\tilde{F}$ fixes $\be_1, \dots, \be_8$ 
and exchanges $\be_9$ and $\be_{10}$ in Figure \ref{fig:D10}. 
\item In case of type $\rA_2$ the matrix $\tilde{F}$ fixes all simple roots.   
\end{enumerate} 
\end{observation}
%%%%%%%%%%%%%%%%%%%%%%%%%%%%%%%%%%%%%%%%%%%%%%%%%%%%%%%%%%%%%%%%%%%%%%
%%%%%%%%%%%%%%%%%%%%%%%%% fig:D10 %%%%%%%%%%%%%%%%%%%%%%%%%%%%%%%%%%%%%%
%%%%%%%%%%%%%%%%%%%%%%%%%%%%%%%%%%%%%%%%%%%%%%%%%%%%%%%%%%%%%%%%%%%%%%%
\par
%%%%%%%%%%%
Observation \ref{obs} shows how the $(-2)$-curves in $X$ are permuted 
by the K3 surface automorphism $f : X \to X$ arising from the Hodge isometry 
$\tilde{F} : L \to L$ (see Lemma \ref{lem:curve}). 
%%%%%%%%%%%%%%%%%%%%%%% sec:SD %%%%%%%%%%%%%%%%%%%%%%%%%%%%%%%%%%%%%%%
\section{Siegel Disks} \label{sec:SD}
%%%%%%%%%%%%%%%%%%%%%%%%%%%%%%%%%%%%%%%%%%%%%%%%%%%%%%%%%%%%%%%%%%
Let $\bD$ and $S^1$ be the unit open disk and the unit circle in $\bC$ 
respectively.   
Given $(\alpha_1, \alpha_2) \in T := S^1 \times S^1$, the map 
$g : \bD^2 \to \bD^2$,  $(z_1, z_2) \mapsto (\alpha_1 z_1, \alpha_2 z_2)$ is 
said to be an {\sl irrational rotation} if $g : T \to T$ has dense orbits; this 
condition is equivalent to saying that $\alpha_1$ and $\alpha_2$ are 
multiplicatively independent, meaning that  
$\alpha_1^{m_1} \alpha_2^{m_2} = 1$, $m_1, m_2 \in \bZ$, implies 
$m_1 = m_2 = 0$.  
Let $f : X \to X$ be an automorphism of a complex surface $X$, and  
$U$ be an open neighborhood of a point $p \in X$.  
We say that $f$ has a {\sl Siegel disk} $U$ centered at $p$ if $f(U, p) = (U, p)$  
and $f|(U, p)$ is biholomorphically conjugate to an irrational rotation 
$g|(\bD^2, 0)$.    
Namely, a Siegel disk is an invariant subset modeled on an irrational rotation 
around the origin. 
%%%%%%%%%%%%%%%%%%%%%%%%%%% ss:exist %%%%%%%%%%%%%%%%%%%%%%%%%%%%%%%%
\subsection{Existence of a Siegel Disk} \label{ss:existence}
%%%%%%%%%%%%%%%%%%%%%%%%%%%%%%%%%%%%%%%%%%%%%%%%%%%%%%%%%%%%%%%%%%
McMullen \cite{McMullen1} synthesizes examples of K3 surface automorphisms 
with a Siegel disk from Salem numbers of degree $22$. 
To this end he establishes a sufficient condition for the existence of a Siegel 
disk (see \cite[Theorems 7.1 and 9.2]{McMullen1}).  
We shall relax his criterion so that it may be applied to Salem numbers of lower 
degrees in the presence of exceptional sets.  
Let $F : L \to L$ be an automorphism of a K3 lattice $L$ admitting a special 
eigenvalue $\delta \in S^1$ and preserving a K3 structure. 
Let $f : X \to X$ be the K3 surface automorphism obtained as the lift of $F$. 
In \cite{McMullen1} the number $\delta$ is called the {\sl determinant} of $f$  
because $\delta$ is equal to the determinant of the holomorphic tangent map 
$(d f)_p : T_p X \to T_p X$ at any fixed point $p$ of $f$. 
Thus the eigenvalues of $(d f)_p$ can be expressed as 
$\alpha_1 = \delta^{1/2} \, \alpha$ and $\alpha_2 = \delta^{1/2} \, \alpha^{-1}$ 
for some $\alpha \in \bC^{\times}$.       
%%%%%%%%%%%%%%%%%%%%%%% prop:siegel %%%%%%%%%%%%%%%%%%%%%%%%%%%%%%%%
\begin{proposition} \label{prop:siegel} 
Suppose that the special eigenvalue $\delta$ is conjugate to a Salem number 
and there exists a rational function $q(w) \in \bQ(w)$ such that 
$(\alpha+ \alpha^{-1})^2 = q(\tau)$ for the special trace 
$\tau := \delta + \delta^{-1} \in (-2, \, 2)$.    
If $\tau$ satisfies $0 \le q(\tau) \le 4$ and admits a conjugate 
$\tau' \in (-2, \, 2)$ such that $q(\tau') > 4$, then the fixed point 
$p \in X$ is the center of a Siegel disk for the map $f$. 
If $\tau$ satisfies $q(\tau) > 4$ then $p$ is a hyperbolic fixed point of $f$. 
\end{proposition} 
%%%%%%%%%%%%%%%%%%%%%%% begin proof %%%%%%%%%%%%%%%%%%%%%%%%%%%%%%%
{\it Proof}. 
Since $\tau \in (-2, \, 2)$, we have $\delta \in S^1$ and   
$\delta^{1/2} \in S^1$. 
It follows from $0 \le q(\tau) \le 4$ that $\alpha \in S^1$ and hence 
$\alpha_1$, $\alpha_2 \in S^1$. 
Note that $\alpha$ is an algebraic number.  
To show $\alpha_1$, $\alpha_2 \in \overline{\bQ}$ are multiplicatively 
independent we mimic the proof of \cite[Lemma 7.5]{McMullen1}. 
Let $\delta'$ and $\alpha'$ be the conjugates of 
$\delta$ and $\alpha$ corresponding to $\tau'$. 
We have $\delta' \in S^1$ and $\alpha' \not\in S^1 $ by 
$\tau' \in (-2, \, 2)$ and $\{\alpha' +(\alpha')^{-1} \}^2 = q(\tau') > 4$. 
If $\alpha_1^m \alpha_2^n = \delta^{(m+n)/2} \alpha^{m-n} =1$ with 
$m$, $n \in \bZ$ then $(\delta')^{(m+n)/2} (\alpha')^{m-n} =1$ and hence 
$|\alpha'|^{m-n} = 1$, which forces $m = n$ and $\delta^m = 1$, 
but the latter equation implies $m = 0$ because $\delta$ is 
not a root of unity.   
The Gel'fond-Baker method then shows that they are jointly Diophantine 
and the Siegel-Sternberg theory tells us that there exists a Siegel 
disk centered at $p$ (see \cite[Theorems 5.2 and 5.3]{McMullen1}).  
On the other hand, if $q(\tau) > 4$ then $\delta^{1/2} \in S^1$ and 
$\alpha \not\in S^1$, so we have $|\alpha_i| < 1 < |\alpha_j|$ for 
$(i, j) = (1, 2)$ or $(2, 1)$. 
Thus $p$ is a hyperbolic fixed point of $f$. \hfill $\Box$ \par\medskip
%%%%%%%%%%%%%%%%%%%%%%% end proof %%%%%%%%%%%%%%%%%%%%%%%%%%%%%%%%
When $f$ has no fixed curves in $X$, McMullen \cite{McMullen1} uses  
Proposition \ref{prop:siegel} with the following setting     
%%%%%%%%%%%%%%%%%%%%%%% eqn:trq %%%%%%%%%%%%%%%%%%%%%%%%%%%%%%%%%%
\begin{equation} \label{eqn:trq}
(\mathrm{a}) \quad \Tr \, [F : L \to L] = -1; \qquad  
(\mathrm{b}) \quad q(w) = \dfrac{(w+1)^2}{w+2},  
\end{equation}
%%%%%%%%%%%%%%%%%%%%%%%%%%%%%%%%%%%%%%%%%%%%%%%%%%%%%%%%%%%%%%%%
by showing that under condition (a) the Lefschetz fixed point formula implies   
the existence of a unique transverse fixed point $p \in X$ of $f$ and  
the Atiyah-Bott formula determines $q(w)$ in the form (b).   
To guarantee the non-existence of fixed curves, he takes Salem numbers 
of degree $22$ so that the Picard lattice and root system in it are empty.  
His conclusion is that the unique fixed point $p$ is the center of a Siegel 
disk, provided  
%%%%%%%%%%%%%%%%%%%%%%% eqn:tau0 %%%%%%%%%%%%%%%%%%%%%%%%%%%%%%%%%
\begin{equation} \label{eqn:tau0}
\mbox{$\tau > \tau_0 := 1 - 2 \sqrt{2} \approx -1.8284271$ 
and $\tau$ has a conjugate $\tau' < \tau_0$}.   
\end{equation}
%%%%%%%%%%%%%%%%%%%%%%%%%%%%%%%%%%%%%%%%%%%%%%%%%%%%%%%%%%%%%%%% 
\par
%%%%%%
To any simple root $\bu \in \vD_{\rb}$ the corresponding $(-2)$-curve 
in $X$ is denoted by the same symbol $\bu$. 
If $F$ fixes a simple root $\bu$ then $f$ preserves the curve 
$\bu \cong \bP^1$, inducing a M\"{o}bius transformation on it, so that 
$\bu$ is a fixed curve of $f$ precisely when the induced map is identity.  
In the possible occurrence of fixed curves we have to use fixed point formulas 
stronger than the classical Lefschetz and Atiyah-Bott formulas 
(see \S \ref{ss:min-Siegel}).  
Any version of  Lefschetz-type formula involves the value of $\Tr \, f^*|H^2(X)$ 
and the following remark is helpful in this respect.  
%%%%%%%%%%%%%%%%%%%%%%%%% rem:trace %%%%%%%%%%%%%%%%%%%%%%%%%%%%%%%%
\begin{remark} \label{rem:trace} 
The {\sl trace} of a monic polynomial $P(z) = z^d + c_1 z^{d-1} + \cdots + c_d$ 
is defined by $\Tr \, P := -c_1$, the sum of its roots. 
A palindromic polynomial of even degree and its trace polynomial have the same 
trace. 
Let $\chi(z)$ be the characteristic polynomial of $F : L \to L$.   
It is a palindromic polynomial of degree $22$. 
We can calculate $\Tr \, f^*|H^2(X) = \Tr \, F$ as the trace 
of $\chi(z)$ or equivalently as the trace of its trace polynomial.  
\end{remark}
%%%%%%%%%%%%%%%%%%%%%%%%%%%% ss:deg22 %%%%%%%%%%%%%%%%%%%%%%%%%%%%%%
\subsection{Salem Numbers of Degree $\mbox{\boldmath $22$}$} 
\label{ss:deg22}
%%%%%%%%%%%%%%%%%%%%%%%%%%%%%%%%%%%%%%%%%%%%%%%%%%%%%%%%%%%%%%%%%%%
Our result on Siegel disks for Salem numbers of degree $22$ is stated 
as follows. 
%%%%%%%%%%%%%%%%%%%%%%%% thm:SH2 %%%%%%%%%%%%%%%%%%%%%%%%%%%%%%%%%%
\begin{theorem} \label{thm:SH2} 
In Theorem $\ref{thm:SH}$ each K3 surface automorphism $f : X \to X$ has 
a unique fixed point $p \in X$, which is either the center of a Siegel disk or 
a hyperbolic fixed point, where the former case is indicated by ``S" while the 
latter case by ``H" in the last columns of 
Tables $\ref{tab:SvH1}$--$\ref{tab:SvH3}$.       
\end{theorem} 
%%%%%%%%%%%%%%%%%%%%%%% begin proof %%%%%%%%%%%%%%%%%%%%%%%%%%%%%%%%
{\it Proof}. 
Observe from Table \ref{tab:mcmullen} that $\Tr \, R_i = -1$ for $i = 1, \dots, 10$. 
Condition (a) in \eqref{eqn:trq} is checked by Remark \ref{rem:trace}. 
Following the framework \eqref{eqn:trq}-\eqref{eqn:tau0} 
the theorem can be established by Proposition \ref{prop:siegel}.  
\hfill $\Box$ \medskip\par
%%%%%%%%%%%%%%%%%%%%%%%% end proof %%%%%%%%%%%%%%%%%%%%%%%%%%%%%%%%%
We remark that Tables \ref{tab:SvH1}--\ref{tab:SvH3} contain a total of $263$ 
entries, among which $230$ are ``S" and $33$ are ``H". 
A majority of the entries have $\bk$'s with a multiple element.        
%%%%%%%%%%%%%%%%%%%%%%%%% ss:min-Siegel %%%%%%%%%%%%%%%%%%%%%%%%%%%%
\subsection{Minimum Entropy} \label{ss:min-Siegel}
%%%%%%%%%%%%%%%%%%%%%%%%%%%%%%%%%%%%%%%%%%%%%%%%%%%%%%%%%%%%%%%%%
We discuss the existence of Siegel disks for K3 surface automorphisms 
constructed in Theorems \ref{thm:A-pattern} and \ref{thm:B-pattern}.    
Let $x_4 < x_3 < x_2 < x_1$ be the roots in $[-2, \, 2]$ of Lehmer's trace 
polynomial $\LT(w)$ as in \eqref{eqn:lehconj}.      
We denote by $\cE$ the exceptional set in $X$ (see Lemma \ref{lem:curve}). 
Even when $\cE$ is nonempty, the framework \eqref{eqn:trq}-\eqref{eqn:tau0} 
remains valid if $f$ has no irreducible fixed curve in $\cE$. 
In particular this is the case if $\tilde{F}$ acts on $\vD_{\rb}$ freely.   
%%%%%%%%%%%%%%%%%%%%%%%%%% thm:min-Siegel1 %%%%%%%%%%%%%%%%%%%%%%%%%%%%%%
\begin{theorem} \label{thm:min-Siegel1} 
For each entry of Dynkin type $\rE_6 \oplus \rE_6$ in Tables $\ref{tab:minA}$ and 
$\ref{tab:minB}$, the K3 surface automorphism $f : X \to X$ has no fixed 
point on $\cE$ and a unique fixed point $p \in X \setminus \cE$. 
If the special trace $\tau$ is any one of $x_1$, $x_2$, $x_3$ then $p$ is the center of 
a Siegel disk, while if $\tau = x_4$ then $p$ is a hyperbolic fixed point; this 
information is indicated in the ``S/H" columns of the tables.       
\end{theorem} 
%%%%%%%%%%%%%%%%%%%%%%%%%% begin proof %%%%%%%%%%%%%%%%%%%%%%%%%%%%%%%%%%% 
{\it Proof}. 
It follows from item (1) of Observation \ref{obs} that $f$ exchanges the two  
$\rE_6$-components of $\cE$, having no fixed point there.  
Thus $f$ has no fixed curve on $X$ and so \eqref{eqn:trq}-\eqref{eqn:tau0} 
can be applied to show the existence of a unique 
transverse fixed point $p \in X \setminus \cE$. 
In view of $x_4 < \tau_0 < x_3$ Proposition \ref{prop:siegel} leads to the theorem.  
\hfill $\Box$ \par\medskip
%%%%%%%%%%%%%%%%%%%%%%%%%% end proof %%%%%%%%%%%%%%%%%%%%%%%%%%%%%%%%%%%%
To discuss the cases where the automorphism $f : X \to X$ may have 
fixed curves, we use S. Saito's fixed point formula \cite[(0.2)]{Saito}. 
Originally it was stated for projective surfaces, but it works for compact 
K\"ahler surfaces as well (see Dinh et al. \cite[Theorem 4.3]{DNT}).  
In the notation of Iwasaki and Uehara \cite[Theorem 1.2]{IU} the formula reads      
%%%%%%%%%%%%%%%%%%%%%%%%% eqn:Saito %%%%%%%%%%%%%%%%%%%%%%%%%%%%%%%%%%%%%
\begin{equation} \label{eqn:Saito}
L(f) := \sum_{i=0}^4 (-1)^i \Tr f^*|H^i(X)   
= \sum_{p \in X_0(f)} \nu_p(f) + \sum_{C \in X_{\rI}(f)} \chi_{C} \cdot \nu_C(f) + 
\sum_{C \in X_{\rII}(f)} \tau_{C} \cdot \nu_C(f), 
\end{equation}
%%%%%%%%%%%%%%%%%%%%%%%%%%%%%%%%%%%%%%%%%%%%%%%%%%%%%%%%%%%%%%%%%%%%%%% 
where $X_0(f)$ is the set of fixed points of $f$ while $X_{\rI}(f)$ and $X_{\rII}(f)$ are 
the sets of irreducible fixed curves of types I and I\!I respectively, $\chi_C$ is the 
Euler number of the normalization of a curve $C$ and $\tau_C$ is the 
self-intersection number of $C$. 
We refer to \cite[\S 3]{IU} for the definitions of the indices $\nu_p(f)$ and 
$\nu_C(f)$ as well as for a detailed account of the terminology used here.    
%%%%%
\par
%%%%%
We also use the Toledo-Tong fixed point formula \cite[Theorem (4.10)]{TT} in the 
special case where $X$ is a compact complex surface, $f : X \to X$ is a 
holomorphic map, $E$ is a holomorphic line bundle on $X$ and $\phi : f^* E \to E$ 
is a holomorphic bundle map. 
Suppose that any isolated fixed point $p$ is transverse and any connected 
component of the $1$-dimensional fixed point set  is a smooth curve $C$  
such that the induced differential map $d^Nf$ on the normal line bundle 
$N = N_C$ to $C$ does not have eigenvalue $1$. 
The formula is then stated as    
%%%%%%%%%%%%%%%%%%%%%%%%%%% eqn:TT %%%%%%%%%%%%%%%%%%%%%%%%%%%%%%%%%%%%%
\begin{equation} \label{eqn:TT}
L(f, \phi) := \sum_{i=0}^2 (-1)^i \, \Tr \, (f, \phi)^*|H^i(X, \cO(E))  
= \sum_{p} \nu_p(f, \phi) + \sum_{C} \nu_C(f, \phi), 
\end{equation}
%%%%%%%%%%%%%%%%%%%%%%%%%%%%%%%%%%%%%%%%%%%%%%%%%%%%%%%%%%%%%%%%%%%%%%%
where the sums are taken over all isolated fixed points $p$ and all 
connected fixed curves $C$. 
If $d^Nf$ has eigenvalue $\lambda_C$ on $N_C$ while $\phi$ has eigenvalues 
$\mu_p$ and $\mu_C$ on $E_p$ and $E|_{C}$ respectively, then the indices 
$\nu_p(f, \phi)$ and $\nu_C(f, \phi)$ are given by   
%%%%%%%%%%%%%%%%%%%%%%%%%% eqn:TTnu %%%%%%%%%%%%%%%%%%%%%%%%%%%%%%%%%%%%
\begin{subequations} \label{eqn:TTnu} 
\begin{align}
\nu_p(f, \phi) &= \frac{\mu_p}{1 - \Tr (d f)_p + \det (d f)_p}, 
\label{eqn:TTnu1} \\[1mm] 
\nu_C(f, \phi) &= \int_C \td(C) \cdot \{ 1 - \lambda_C \, \ch(\check{N}) \}^{-1} 
\cdot \mu_C \, \ch(E), \label{eqn:TTnu2}
\end{align} 
\end{subequations}
%%%%%%%%%%%%%%%%%%%%%%%%%%%%%%%%%%%%%%%%%%%%%%%%%%%%%%%%%%%%%%%%%%%%%%%
where $\td(C)$ is the Todd class of $C$, 
$\ch(E)$ is the Chern character of $E$, $\check{N}$ is the dual line bundle to 
$N$ and the integral sign stands for evaluation on the fundamental cycle of $C$.  
%%%%%
\par
\begin{figure}[hh]
\begin{center}
\input{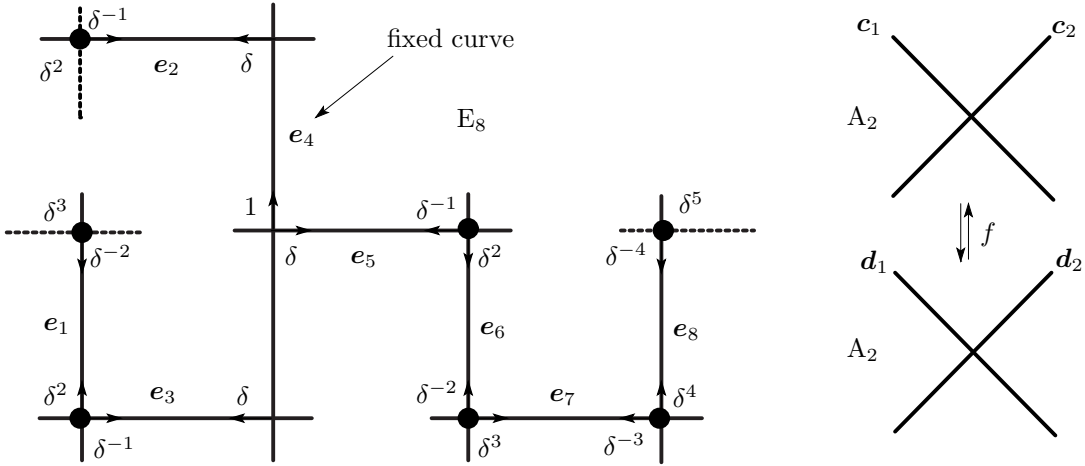}
\end{center}
\caption{Eigenvalues at the isolated fixed points on the exceptional set $\cE$.}
\label{fig:e-curve} 
\end{figure}

%%%%%
A {\it three-cycle of Siegel disks} for $f$ is a sequence of open subsets 
$U$, $f(U)$, $f^2(U)$ in $X$ such that $U$ is a Siegel disk for $f^3$ centered 
at a point $p \in U$ which is a periodic point of period $3$, that is,   
$p$, $f(p)$, $f^2(p)$ are mutually distinct and $f^3(p) = p$. 
We remark that $U$, $f(U)$, $f^2(U)$ are Siegel disks for $f^3$ centered at 
the points $p$, $f(p)$, $f^2(p)$, respectively.         
%%%%%%%%%%%%%%%%%%%%%%%%%% thm:min-Segel2 %%%%%%%%%%%%%%%%%%%%%%%%%%%%%%%
\begin{theorem} \label{thm:min-Siegel2} 
For each entry of type $\rE_8 \oplus \rA_2 \oplus \rA_2$  
in Tables $\ref{tab:minA}$ and $\ref{tab:minB}$, the K3 surface automorphism  
$f : X \to X$ has a unique periodic orbit $p$, $f(p)$, $f^2(p) \in X$ of period $3$, 
which lie in $X \setminus \cE$. 
If the special trace $\tau$ is either $x_2$ or $x_4$ then the three points are 
the centers of a three-cycle of Siegel disks, while if $\tau$ is either 
$x_1$ or $x_3$ then they are hyperbolic periodic points; this information is 
indicated in the ``S/H" columns of the tables.        
\end{theorem}
%%%%%%%%%%%%%%%%%%%%%%%%% begin proof %%%%%%%%%%%%%%%%%%%%%%%%%%%%%%%%%%%
{\it Proof}. 
It follows from (2) of Observation \ref{obs} that $f$ has no fixed point on the 
$\rA_2 \oplus \rA_2$-component of $\cE$. 
On the other hand $f$ preserves each $(-2)$-curve on the 
$\rE_8$-component, inducing a M\"obius transformation on it. 
A M\"obius transformation falls into one of the three categories according to 
the number $n$ of its fixed points; parabolic for $n = 1$, non-parabolic for $n = 2$, 
and the identity for $n \ge 3$. 
In the parabolic case the derivative at the unique fixed point is $1$.     
In the non-parabolic case, if the derivatives at one fixed points is $c$   
then the derivative at the other fixed point is $c^{-1}$, where $c \neq 1$.   
Notice that $\be_4$ is a fixed curve of $f$ since $\be_4$ is fixed at the three 
intersections $\be_4 \cap \be_2$, $\be_4 \cap \be_3$, $\be_4 \cap \be_5$ 
(see Figure \ref{fig:e-curve}). 
%%%%%%%%%%%%%%%%%%%%%%%%%% fig:e-curve %%%%%%%%%%%%%%%%%%%%%%%%%%%%%%%%%%%
%%%%%%%%%%%%%%%%%%%%%%%%%%%%%%%%%%%%%%%%%%%%%%%%%%%%%%%%%%%%%%%%%%%%%%%% 
%%%%
\par
%%%%
Consider the arm  $\ba := \be_5 \cup \be_6 \cup \be_7 \cup \be_8$ 
emanating from $\be_4$. 
Let $q_0, q_1, q_2, q_3$ be the points at $\be_4 \cap \be_5$, 
$\be_5 \cap \be_6$, $\be_6 \cap \be_7$, $\be_7 \cap \be_8$. 
Note that they are fixed points of $f$.  
We use the fact that $\det (d f)_q = \delta$ at any fixed point $q$ of $f$. 
Since $(d f)_{q_0} = 1$ along $\be_4$ we have $(d f)_{q_0} = \delta \neq 1$ 
along $\be_5$.  
Thus $f$ induces a non-parabolic M\"obius transformation on $\be_5$ 
having derivative $\delta^{-1}$ at $q_1$. 
Then $f$ induces a non-parabolic M\"obius transformation on $\be_6$ 
having derivatives $\delta^2$ at $q_1$ and $\delta^{-2}$ at $q_2$. 
Repeating this argument shows that $f$ has four fixed points on $\ba$, 
three of which are $q_1$, $q_2$, $q_3$ and the final one $q_4$ is on $\be_8$, 
and that $(d f)_{q_j}$ has eigenvalues $\delta^{-j}$ and 
$\delta^{j+1}$ for $j = 1, 2, 3, 4$.  
Similar statements can be made for the shorter arms 
$\bb := \be_1 \cup \be_3$ and $\be_2$. 
In total there are seven isolated fixed points indicated by $\bullet$  
and one fixed curve $\be_4$ on  $\cE$. 
%%%%%%
\par
%%%%%% 
We apply Saito's formula \eqref{eqn:Saito} to $f$.  
For each entry of type $\rE_8 \oplus \rA_2 \oplus \rA_2$ in 
Tables \ref{tab:minA} and \ref{tab:minB} we have $\Tr \, \tilde{F} = 7$  
where $\tilde{F} = \tilde{A}$ or $\tilde{B}$, and so 
$L(f) = 2 + \Tr \, \tilde{F} = 9$.  
The seven isolated fixed point on $\cE$ are transverse and hence 
of index $1$. 
The fixed curve $\be_4$ is of type I and has index $1$, since $d f$ 
has eigenvalue $\delta \neq 1$ in its normal direction. 
Any point on $\be_4$ has index $0$. 
The Euler number of $\be_4 \cong \bP^1$ is $2$. 
There is no fixed curve of type I\!I. 
Thus formula \eqref{eqn:Saito} reads 
%%%%%%%%%%%%%%%%%%% eqn:Saito2 %%%%%%%%%%%%%%%%%%%%%%%%%%%%
\begin{equation} \label{eqn:Saito2}
9 = 7 + \sum_{p \in X_0(f) \setminus \cE} \nu_p(f) 
+ 2, \qquad \mbox{i.e.} \qquad 
\sum_{p \in X_0(f) \setminus \cE} \nu_p(f) = 0,   
\end{equation}
%%%%%%%%%%%%%%%%%%%%%%%%%%%%%%%%%%%%%%%%%%%%%%%%%%%%%%%%%
where $7$ and $2$ on the RHS are the contributions of the seven 
isolated fixed point on $\cE$ and the fixed curve $\be_4$ respectively.  
This implies that $f$ has no fixed point on $X \setminus \cE$.   
%%%%%%
\par
%%%%%%
Next we apply the Toledo-Tong formula \eqref{eqn:TT} to $f$ upon setting 
$E = K_X := \wedge^2 T^* X$ and $\phi = f^* : f^* K_X \to K_X$, 
the pull-back mapping.   
Notice that $\mu_p = \lambda_C = \mu_C = \delta$ with $C = \be_4$ in 
\eqref{eqn:TTnu}. 
Let $\eta$ be a nowhere vanishing holomorphic $2$-form on $X$. 
At each $x \in C$ the map $T_x X \to T_x^* C$, 
$v \mapsto \eta_x(v, \,\cdot \,)$ induces an isomorphism 
$N_x := T_x X/T_x C \cong T_x^* C$, hence $\check{N} \cong T C$. 
A little calculation yields $\nu_C(f, f^*) = \delta (\delta+1)/(\delta-1)^2$ 
with $C = \be_4$ in \eqref{eqn:TTnu2}.  
It is easy to see that $L(f, f^*) = 1 + \delta$.  
Formula \eqref{eqn:TT} then asserts that $z = \delta$ is a solution 
to the equation 
%%%%%%%
\begin{align*}
1+z 
&= \sum_{j=1}^4 \dfrac{z}{1-(z^{-j}+z^{j+1})+z} + 
\sum_{j=1}^2 \dfrac{z}{1-(z^{-j}+z^{j+1})+z}  \\
&\phantom{=} + \dfrac{z}{1-(z^{-1}+z^2)+z} + \dfrac{z(z+1)}{(z-1)^2},   
\end{align*}
%%%%%%%
where the first three terms in the RHS are the contributions of the fixed points 
on the arms $\ba$, $\bb$, $\be_2$ and the last term is that of the fixed curve 
$\be_4$.  
A careful inspection shows that the difference $D(z)$ of the LHS from the RHS 
above admits a clean factorization  
%%%%%%%%%%%%%%%%%%%%%%%%%%% eqn:D(z) %%%%%%%%%%%%%%%%%%%%%%%%%%%%%%%%%%%
\begin{equation} \label{eqn:D(z)}
D(z) = \frac{\rL(z)}{(z+1) \cdot \rC_1(z) \cdot \rC_3(z) \cdot \rC_5(z)} 
= \dfrac{z \cdot \LT(w)}{(z+1) \cdot \CT_1(w) \cdot \CT_3(w) \cdot \CT_5(w)}, 
\end{equation}
%%%%%%%%%%%%%%%%%%%%%%%%%%%%%%%%%%%%%%%%%%%%%%%%%%%%%%%%%%%%%%%%%%%%%%
where $w := z+ z^{-1}$. 
So the Toledo-Tong formula for $f$ is nothing other than the tautological fact that 
$\delta$ is a root of Lehmer's polynomial, but the formula \eqref{eqn:D(z)} itself   
will be useful later. 
%%%%%%
\par
%%%%%%
We again apply Saito's formula \eqref{eqn:Saito} this time to $f^3$. 
If $\tilde{F}$ has characteristic polynomial $\tilde{\chi}(z) = 
z^{22} - e_1 z^{21} + e_2 z^{20} - e_3 z^{19} + \cdots$ then 
$\Tr (\tilde{F}^3) = p_3 = 3(e_3 - e_1 e_2) + e_1^3$ by the relation between 
power sums and elementary symmetric polynomials.  
We have $e_1 = 7$, $e_2 = 20$, $e_3 = 29$ and $p_3 = 10$ in the cases of type 
$\rE_8 \oplus \rA_2 \oplus \rA_2$.   
This implies $L(f) = 2 + \Tr(\tilde{F}^3) = 12$, hence the equation 
\eqref{eqn:Saito2} turns into 
%%%%%%%%%%%%%%%%%%% eqn:Saito3 %%%%%%%%%%%%%%%%%%%%%%%%%%%%
\begin{equation*}
12 = 7 + \sum_{p \in X_0(f^3) \setminus \cE} \nu_p(f^3) 
+ 2, \qquad \mbox{i.e.} \qquad 
\sum_{p \in X_0(f^3) \setminus \cE} \nu_p(f^3) = 3.   
\end{equation*}
%%%%%%%%%%%%%%%%%%%%%%%%%%%%%%%%%%%%%%%%%%%%%%%%%%%%%%%%%
Thus $f^3$ has three fixed points on $X \setminus \cE$ 
counted with multiplicity. 
Let $p$ be one of them. 
Since $f$ has no fixed point on 
$X \setminus \cE$, we have $p \neq f(p)$ and hence 
$f(p) \neq f^2(p)$ and $f^2(p) \neq f^3(p) = p$. 
Accordingly, $p$, $f(p)$, $f^2(p)$ must be distinct and 
$\nu_{f^j(p)}(f^3) = 1$ for $j = 0, 1, 2$.  
%%%%%%
\par
%%%%%%
The tangent maps $(d f^3)_{f^j(p)} : T_{f^j(p)} X \to T_{f^j(p)} X$, $j = 0, 1, 2$, 
have common determinant $\delta^3$ and common eigenvalues 
of the form $\delta^{3/2} \alpha^{\pm 1}$ with $\alpha \in \bC^{\times}$.  
The Toledo-Tong formula \eqref{eqn:TT} for $f^3$ is then 
expressed as  
%%%%%%
$$
D(\delta^3) = \dfrac{3 \delta^3}{1- \delta^{3/2} (\alpha + \alpha^{-1}) + \delta^3},  
$$
%%%%%%
in terms of the rational function $D(z)$ in \eqref{eqn:D(z)}. 
Solving this equation we have  
%%%%%%
\begin{align*} 
q(\tau) &:= (\alpha + \alpha^{-1})^2 = \delta^{-3} 
\left\{ 1 + \delta^3 - \frac{3 \delta^3}{D(\delta^3)} \right\}^2 
= \frac{(1+\delta^3)^2 \rM(\delta^3)^2}{\delta^3 \, \rL(\delta^3)^2} \\[2mm]
&= \frac{(\delta^3 + \delta^{-3}+2) \, \MT(\delta^3+\delta^{-3})^2}{
\LT(\delta^3 + \delta^{-3})^2} 
= \frac{(\tau+2) (\tau-1)^2 \, \MT(\tau^3-3 \tau)^2}{\LT(\tau^3 - 3 \tau)^2},  
\end{align*}
%%%%%%
where $\rM(z)$ and $\MT(z)$ are Salem polynomial and its trace polynomial 
defined by   
%%%%%%%
\begin{align*}
\rM(z) &:= z^{10} - 2 z^9 - z^7 + 2 z^6 - z^5 + 2z^4 -z^3 - 2 z+ 1, 
\\[1mm]
\MT(w) &:= (w+1)(w-2)(w^3-w^2-4 w+ 1) -1. 
\end{align*}
%%%%%%%
Here we observe that $0 < q(x_j) < 4$ for $j = 2$, $4$ and $q(x_j) > 4$ for $j = 1$, $3$.   
Thus Proposition \ref{prop:siegel} leads to the conclusion of the theorem.  
\hfill $\Box$ 
%%%%%%%%%%%%%%%%%%%%%%%%% end proof %%%%%%%%%%%%%%%%%%%%%%%%%%%%%%%%%%%%%
%%%%%%%%%%%%%%%%%%%%%%%%% thm:D10 %%%%%%%%%%%%%%%%%%%%%%%%%%%%%%%%%%%%
\begin{theorem} \label{thm:D10} 
For each entry of type $\rD_{10}$ in Table $\ref{tab:minA}$ the K3 
surface automorphism $f : X \to X$ has a unique periodic orbit $p$, $f(p)$, 
$f^2(p) \in X$ of period $3$ which lies in $X \setminus \cE$. 
If the special trace $\tau$ is $x_2$ then the three points are the centers of a 
three-cycle of Siegel disks, while if $\tau = x_4$ then they are hyperbolic 
periodic points; see the ``S/H" column of Table $\ref{tab:minA}$.        
\end{theorem}
%%%%%%%%%%%%%%%%%%%%%%%%% begin proof %%%%%%%%%%%%%%%%%%%%%%%%%%%%%%%%%%
{\it Proof}. 
The proof is similar to that of Theorem \ref{thm:min-Siegel2} or even simpler 
because Saito's formula is not necessary and the classical Lefschetz 
formula is sufficient for dealing with $f$ and $f^3$.  
%%%%%%
\par
%%%%%%
For each $j = 1, \dots, 7$, let $p_j$ be the intersection of the $(-2)$-curves 
$\be_{9-j}$ and $\be_{8-j}$ in Figure \ref{fig:D10}. 
Moreover let $p_{\pm}$ be the intersection of $\be_8$ and $\be_{(19 \pm 1)/2}$.   
It follows from (3) of Observation \ref{obs} that $f$ fixes $p_1, \dots, p_7$ 
and exchanges $p_{\pm}$. 
In particular the M\"{o}bius transformation $f|_{\sbe_8}$ has derivative $-1$ at $p_1$, 
and hence admits another fixed point $p_0 \in \be_8 \setminus \{p_1, p_{\pm}\}$ 
with derivative $-1$. 
Thus $(d f)_{p_0}$ and $(d f)_{p_1}$ have eigenvalues $-1$ and $-\delta$.    
Inductively, $(d f)_{p_j}$ has eigenvalues $- \delta^{1-j}$ and $-\delta^j$ for 
$j = 1, \dots, 7$. 
This is also true for $j = 8$ with another fixed point $p_8$ of $f|_{\sbe_1}$. 
So the map $f$ has nine transverse fixed points $p_0, p_1, \dots, p_8$ on 
$\cE$ and no fixed curve on $X$.  
By the Lefschetz formula the Lefschetz number $L(f) = 2 + \Tr \, \tilde{A} 
= 2 + 7 = 9$ is represented as 
%%%%%%%
$$
L(f) = \sum_{j=0}^8 \nu_{p_j}(f) + \sum_{p \in X_0(f) \setminus \cE} \nu_p(f) 
= 9 + \sum_{p \in X_0(f) \setminus \cE} \nu_p(f),  
$$
%%%%%%%
which implies $\sum_{p \in X_0(f) \setminus \cE} \nu_p(f) = 0$, that is,  
$f$ has no fixed point in $X \setminus \cE$. 
%%%%%%%
\par
%%%%%%%
Similarly the map $f^3$ has no fixed curve and nine transverse fixed points 
$p_0, \dots, p_8 \in \cE$, so the Lefschetz number 
$L(f^3) = 2 + \Tr((f^3)^*|H^2(X)) = 2 + 10 = 12$ is equal to  
%%%%%%%
$$
L(f^3) = \sum_{j=0}^8 \nu_{p_j}(f^3) + 
\sum_{p \in X_0(f^3) \setminus \cE} \nu_p(f^3) 
= 9 + \sum_{p \in X_0(f^3) \setminus \cE} \nu_p(f^3),  
$$
%%%%%%%
which implies $\sum_{p \in X_0(f^3) \setminus \cE} \nu_p(f^3) = 3$, that is,  
$f^3$ has three fixed points in $X \setminus \cE$ counted with multiplicities.   
As in the proof of Theorem \ref{thm:min-Siegel2}, they are distinct 
and transverse by the fact that $f$ has no fixed point in $X \setminus \cE$. 
Thus $f$ has a unique transverse periodic cycle $p$, $f(p)$, $f^2(p) \in X 
\setminus \cE$ of period $3$. 
If $\delta^{3/2} \alpha^{\pm 1}$ are the common eigenvalues of 
$(d f^3)_{f^j(p)}$, $j = 0, 1, 2$, then the Atiyah-Bott formula for the 
canonical line bundle $K_X$ tells us that $z = \delta^3$ is a root of  
%%%%%%%
$$
1 + z = \dfrac{z}{1 + (1+z) + z} + \sum_{j=1}^8 \dfrac{z}{1 +(z^{1-j} + z^j) + z} 
+ \dfrac{3 z}{1- z^{\frac{1}{2}}(\alpha + \alpha^{-1}) + z}, 
$$
%%%%%%%
where the LHS comes from the holomorphic Lefschetz number $L(f^3, (f^3)^*)$,  
while three terms on the RHS are the contributions of 
$p_0$; $p_1, \dots, p_8$; and $p$, $f(p)$, $f^2(p)$, respectively.  
Solving the above equation with respect to $\alpha + \alpha^{-1}$, 
we have $(\alpha + \alpha^{-1})^2 = q(\tau)$, where   
%%%%%%
$$
q(w) := \frac{(w + 2)(w-1)^2 \NT(w^3-3 w)^2}{\LT(w^3-3 w)^2}, 
$$
%%%%%%
with $\NT(w) := (w^2-2w-2)(w^3-3w+1) -1$ being a Salem trace 
polynomial of degree $5$.  
%%%%%%
\par
%%%%%%
Since $0 < q(x_j) < 4$ for $j = 1$, $2$, $3$ and $q(x_4) > 4$, Proposition \ref{prop:siegel} 
shows that if $\tau = x_1$, $x_2$, $x_3$ then $p$, $f(p)$, $f^2(p)$ are the centers 
of a three-cycle of Siegel disks for $f$, and if $\tau = x_4$ they are hyperbolic 
periodic points of period $3$. 
Actually we have $\tau = x_2$, $x_4$ for the entries of type $\rD_{10}$ in 
Table \ref{tab:minA}. 
We remark that there is no entry of type $\rD_{10}$ in Table \ref{tab:minB}.     
\hfill $\Box$
%%%%%%%%%%%%%%%%%%%%%%%% end proof %%%%%%%%%%%%%%%%%%%%%%%%%%%%%%%%%%%%%%%
%%%%%%%%%%%%%%%%%%%%%%%% thm:A2 %%%%%%%%%%%%%%%%%%%%%%%%%%%%%%%%%%%%%%%%
\begin{theorem} \label{thm:A2} 
For each entry of type $\rA_2$ in Table $\ref{tab:minA}$ the K3 surface automorphism 
$f : X \to X$ has exactly three transverse fixed points $p_0$ and $p_{\pm}$ on the exceptional 
set $\cE$ and no fixed point outside $\cE$, where $p_0$ is at the intersection 
of two $(-2)$-curves $\be_{\pm} \subset \cE$ while $p_{\pm}$ lie on 
$\be_{\pm} \setminus \{ p_0 \}$ as in Figure $\ref{fig:A2a}$.  
Moreover $f$ has three Siegel disks centered at $p_0$ and $p_{\pm}$ respectively; 
see the ``S/H" column of Table $\ref{tab:minA}$.      
\end{theorem}
%%%%%%%%%%%%%%%%%%%%%%%%%%%%%% fig:A2 %%%%%%%%%%%%%%%%%%%%%%%%%%%%%%%%%%%
\begin{figure}[h]
\begin{minipage}{0.45\linewidth}
\centerline{%WinTpicVersion4.32a
{\unitlength 0.1in%
\begin{picture}(28.7000,17.5000)(5.3000,-26.1000)%
% LINE 0 0 3 0 Black White  
% 2 1810 2590 3400 1010
% 
\special{pn 20}%
\special{pa 1810 2590}%
\special{pa 3400 1010}%
\special{fp}%
% LINE 0 0 3 0 Black White  
% 2 2200 2610 610 1020
% 
\special{pn 20}%
\special{pa 2200 2610}%
\special{pa 610 1020}%
\special{fp}%
% CIRCLE 0 0 0 0 Black White  
% 4 2000 2410 2040 2440 2040 2440 2040 2440
% 
\special{sh 1.000}%
\special{ia 2000 2410 50 50 0.0000000 6.2831853}%
\special{pn 20}%
\special{ar 2000 2410 50 50 0.0000000 6.2831853}%
% CIRCLE 0 0 0 0 Black White  
% 4 900 1330 940 1360 940 1360 940 1360
% 
\special{sh 1.000}%
\special{ia 900 1330 50 50 0.0000000 6.2831853}%
\special{pn 20}%
\special{ar 900 1330 50 50 0.0000000 6.2831853}%
% LINE 1 2 3 0 Black White  
% 2 740 1510 1370 860
% 
\special{pn 13}%
\special{pa 740 1510}%
\special{pa 1370 860}%
\special{dt 0.045}%
% STR 2 0 3 0 Black White  
% 4 530 890 530 990 2 0 0 0
% $\be_-$
\put(5.3000,-9.9000){\makebox(0,0)[lb]{$\be_-$}}%
% STR 2 0 3 0 Black White  
% 4 3340 890 3340 990 2 0 0 0
% $\be_+$
\put(33.4000,-9.9000){\makebox(0,0)[lb]{$\be_+$}}%
% VECTOR 0 0 3 0 Black White  
% 2 2150 2250 2250 2160
% 
\special{pn 20}%
\special{pa 2150 2250}%
\special{pa 2250 2160}%
\special{fp}%
\special{sh 1}%
\special{pa 2250 2160}%
\special{pa 2187 2190}%
\special{pa 2210 2196}%
\special{pa 2214 2219}%
\special{pa 2250 2160}%
\special{fp}%
% VECTOR 0 0 3 0 Black White  
% 4 1850 2260 1770 2170 1770 2170 1770 2170
% 
\special{pn 20}%
\special{pa 1850 2260}%
\special{pa 1770 2170}%
\special{fp}%
\special{sh 1}%
\special{pa 1770 2170}%
\special{pa 1799 2233}%
\special{pa 1805 2210}%
\special{pa 1829 2207}%
\special{pa 1770 2170}%
\special{fp}%
\special{pa 1770 2170}%
\special{pa 1770 2170}%
\special{fp}%
% VECTOR 0 0 3 0 Black White  
% 6 1000 1420 1120 1530 1130 1530 1130 1530 1130 1530 1130 1530
% 
\special{pn 20}%
\special{pa 1000 1420}%
\special{pa 1120 1530}%
\special{fp}%
\special{sh 1}%
\special{pa 1120 1530}%
\special{pa 1084 1470}%
\special{pa 1081 1494}%
\special{pa 1057 1500}%
\special{pa 1120 1530}%
\special{fp}%
\special{pa 1130 1530}%
\special{pa 1130 1530}%
\special{fp}%
\special{pa 1130 1530}%
\special{pa 1130 1530}%
\special{fp}%
% VECTOR 0 2 3 0 Black White  
% 2 950 1290 1100 1140
% 
\special{pn 20}%
\special{pa 950 1290}%
\special{pa 1100 1140}%
\special{dt 0.045}%
\special{sh 1}%
\special{pa 1100 1140}%
\special{pa 1039 1173}%
\special{pa 1062 1178}%
\special{pa 1067 1201}%
\special{pa 1100 1140}%
\special{fp}%
% STR 2 0 3 0 Black White  
% 4 1940 2560 1940 2660 2 0 0 0
% $p_0$
\put(19.4000,-26.6000){\makebox(0,0)[lb]{$p_0$}}%
% STR 2 0 3 0 Black White  
% 4 2260 2250 2260 2350 2 0 0 0
% $\delta^{\frac{1}{2}} \alpha$
\put(22.6000,-23.5000){\makebox(0,0)[lb]{$\delta^{\frac{1}{2}} \alpha$}}%
% STR 2 0 3 0 Black White  
% 4 1290 2250 1290 2350 2 0 0 0
% $\delta^{\frac{1}{2}} \alpha^{-1}$
\put(12.9000,-23.5000){\makebox(0,0)[lb]{$\delta^{\frac{1}{2}} \alpha^{-1}$}}%
% STR 2 0 3 0 Black White  
% 4 1120 1200 1120 1300 2 0 0 0
% $\delta^{\frac{3}{2}} \alpha^{-1}$
\put(11.2000,-13.0000){\makebox(0,0)[lb]{$\delta^{\frac{3}{2}} \alpha^{-1}$}}%
% STR 2 0 3 0 Black White  
% 4 740 1550 740 1650 2 0 0 0
% $\delta^{-\frac{1}{2}} \alpha$
\put(7.4000,-16.5000){\makebox(0,0)[lb]{$\delta^{-\frac{1}{2}} \alpha$}}%
% STR 2 0 3 0 Black White  
% 4 650 1290 650 1390 2 0 0 0
% $p_-$
\put(6.5000,-13.9000){\makebox(0,0)[lb]{$p_-$}}%
% CIRCLE 0 0 0 0 Black White  
% 4 3050 1340 3090 1370 3090 1370 3090 1370
% 
\special{sh 1.000}%
\special{ia 3050 1340 50 50 0.0000000 6.2831853}%
\special{pn 20}%
\special{ar 3050 1340 50 50 0.0000000 6.2831853}%
% LINE 1 2 3 0 Black White  
% 4 2650 880 3240 1580 3240 1580 3240 1580
% 
\special{pn 13}%
\special{pa 2650 880}%
\special{pa 3240 1580}%
\special{dt 0.045}%
\special{pa 3240 1580}%
\special{pa 3240 1580}%
\special{dt 0.045}%
% STR 2 0 3 0 Black White  
% 4 3180 1310 3180 1410 2 0 0 0
% $p_+$
\put(31.8000,-14.1000){\makebox(0,0)[lb]{$p_+$}}%
% VECTOR 0 0 3 0 Black White  
% 2 3010 1390 2860 1550
% 
\special{pn 20}%
\special{pa 3010 1390}%
\special{pa 2860 1550}%
\special{fp}%
\special{sh 1}%
\special{pa 2860 1550}%
\special{pa 2920 1515}%
\special{pa 2896 1511}%
\special{pa 2891 1488}%
\special{pa 2860 1550}%
\special{fp}%
% VECTOR 0 2 3 0 Black White  
% 2 3040 1360 2860 1140
% 
\special{pn 20}%
\special{pa 3040 1360}%
\special{pa 2860 1140}%
\special{dt 0.045}%
\special{sh 1}%
\special{pa 2860 1140}%
\special{pa 2887 1204}%
\special{pa 2894 1181}%
\special{pa 2918 1179}%
\special{pa 2860 1140}%
\special{fp}%
% STR 2 0 3 0 Black White  
% 4 2570 1190 2570 1290 2 0 0 0
% $\delta^{\frac{3}{2}} \alpha$
\put(25.7000,-12.9000){\makebox(0,0)[lb]{$\delta^{\frac{3}{2}} \alpha$}}%
% STR 2 0 3 0 Black White  
% 4 2840 1620 2840 1720 2 0 0 0
% $\delta^{-\frac{1}{2}} \alpha^{-1}$
\put(28.4000,-17.2000){\makebox(0,0)[lb]{$\delta^{-\frac{1}{2}} \alpha^{-1}$}}%
\end{picture}}%
}
\caption{Case of type $\rA_2$ (a). }
\label{fig:A2a}
\end{minipage}
\hspace{5mm} 
\begin{minipage}{0.45\linewidth}
\centerline{%WinTpicVersion4.32a
{\unitlength 0.1in%
\begin{picture}(28.7000,17.5000)(5.3000,-26.1000)%
% LINE 0 0 3 0 Black White  
% 2 1810 2590 3400 1010
% 
\special{pn 20}%
\special{pa 1810 2590}%
\special{pa 3400 1010}%
\special{fp}%
% LINE 0 0 3 0 Black White  
% 2 2200 2610 610 1020
% 
\special{pn 20}%
\special{pa 2200 2610}%
\special{pa 610 1020}%
\special{fp}%
% CIRCLE 0 0 0 0 Black White  
% 4 2000 2410 2040 2440 2040 2440 2040 2440
% 
\special{sh 1.000}%
\special{ia 2000 2410 50 50 0.0000000 6.2831853}%
\special{pn 20}%
\special{ar 2000 2410 50 50 0.0000000 6.2831853}%
% CIRCLE 0 0 0 0 Black White  
% 4 900 1330 940 1360 940 1360 940 1360
% 
\special{sh 1.000}%
\special{ia 900 1330 50 50 0.0000000 6.2831853}%
\special{pn 20}%
\special{ar 900 1330 50 50 0.0000000 6.2831853}%
% LINE 1 2 3 0 Black White  
% 2 740 1510 1370 860
% 
\special{pn 13}%
\special{pa 740 1510}%
\special{pa 1370 860}%
\special{dt 0.045}%
% STR 2 0 3 0 Black White  
% 4 530 890 530 990 2 0 0 0
% $\be_-$
\put(5.3000,-9.9000){\makebox(0,0)[lb]{$\be_-$}}%
% STR 2 0 3 0 Black White  
% 4 3340 890 3340 990 2 0 0 0
% $\be_+$
\put(33.4000,-9.9000){\makebox(0,0)[lb]{$\be_+$}}%
% VECTOR 0 0 3 0 Black White  
% 2 2150 2250 2250 2160
% 
\special{pn 20}%
\special{pa 2150 2250}%
\special{pa 2250 2160}%
\special{fp}%
\special{sh 1}%
\special{pa 2250 2160}%
\special{pa 2187 2190}%
\special{pa 2210 2196}%
\special{pa 2214 2219}%
\special{pa 2250 2160}%
\special{fp}%
% VECTOR 0 0 3 0 Black White  
% 4 1850 2260 1770 2170 1770 2170 1770 2170
% 
\special{pn 20}%
\special{pa 1850 2260}%
\special{pa 1770 2170}%
\special{fp}%
\special{sh 1}%
\special{pa 1770 2170}%
\special{pa 1799 2233}%
\special{pa 1805 2210}%
\special{pa 1829 2207}%
\special{pa 1770 2170}%
\special{fp}%
\special{pa 1770 2170}%
\special{pa 1770 2170}%
\special{fp}%
% VECTOR 0 0 3 0 Black White  
% 6 1000 1420 1120 1530 1130 1530 1130 1530 1130 1530 1130 1530
% 
\special{pn 20}%
\special{pa 1000 1420}%
\special{pa 1120 1530}%
\special{fp}%
\special{sh 1}%
\special{pa 1120 1530}%
\special{pa 1084 1470}%
\special{pa 1081 1494}%
\special{pa 1057 1500}%
\special{pa 1120 1530}%
\special{fp}%
\special{pa 1130 1530}%
\special{pa 1130 1530}%
\special{fp}%
\special{pa 1130 1530}%
\special{pa 1130 1530}%
\special{fp}%
% VECTOR 0 2 3 0 Black White  
% 2 950 1290 1100 1140
% 
\special{pn 20}%
\special{pa 950 1290}%
\special{pa 1100 1140}%
\special{dt 0.045}%
\special{sh 1}%
\special{pa 1100 1140}%
\special{pa 1039 1173}%
\special{pa 1062 1178}%
\special{pa 1067 1201}%
\special{pa 1100 1140}%
\special{fp}%
% STR 2 0 3 0 Black White  
% 4 1940 2560 1940 2660 2 0 0 0
% $p_0$
\put(19.4000,-26.6000){\makebox(0,0)[lb]{$p_0$}}%
% STR 2 0 3 0 Black White  
% 4 2260 2250 2260 2350 2 0 0 0
% $1$
\put(22.6000,-23.5000){\makebox(0,0)[lb]{$1$}}%
% STR 2 0 3 0 Black White  
% 4 1630 2260 1630 2360 2 0 0 0
% $\delta$
\put(16.3000,-23.6000){\makebox(0,0)[lb]{$\delta$}}%
% STR 2 0 3 0 Black White  
% 4 1120 1200 1120 1300 2 0 0 0
% $\delta^2$
\put(11.2000,-13.0000){\makebox(0,0)[lb]{$\delta^2$}}%
% STR 2 0 3 0 Black White  
% 4 830 1600 830 1700 2 0 0 0
% $\delta^{-1}$
\put(8.3000,-17.0000){\makebox(0,0)[lb]{$\delta^{-1}$}}%
% STR 2 0 3 0 Black White  
% 4 1860 1180 1860 1280 2 0 0 0
% fixed curve or 
\put(18.6000,-12.8000){\makebox(0,0)[lb]{fixed curve or }}%
% STR 2 0 3 0 Black White  
% 4 1880 1390 1880 1490 2 0 0 0
% parabolic
\put(18.8000,-14.9000){\makebox(0,0)[lb]{parabolic}}%
% STR 2 0 3 0 Black White  
% 4 1880 1580 1880 1680 2 0 0 0
% M\"{o}bius t. 
\put(18.8000,-16.8000){\makebox(0,0)[lb]{M\"{o}bius t. }}%
% VECTOR 2 0 3 0 Black White  
% 2 2650 1510 2740 1610
% 
\special{pn 8}%
\special{pa 2650 1510}%
\special{pa 2740 1610}%
\special{fp}%
\special{sh 1}%
\special{pa 2740 1610}%
\special{pa 2710 1547}%
\special{pa 2704 1570}%
\special{pa 2681 1574}%
\special{pa 2740 1610}%
\special{fp}%
% STR 2 0 3 0 Black White  
% 4 620 1290 620 1390 2 0 0 0
% $p_-$
\put(6.2000,-13.9000){\makebox(0,0)[lb]{$p_-$}}%
\end{picture}}%
}
\caption{Case of type $\rA_2$ (b). }
\label{fig:A2b}
\end{minipage}
\end{figure}
%%%%%%%%%%%%%%%%%%%%%%%%%%%%%%%%%%%%%%%%%%%%%%%%%%%%%%%%%%%%%%%%%%%%%%%    
%%%%%%%%%%%%%%%%%%%%%%% begin proof %%%%%%%%%%%%%%%%%%%%%%%%%%%%%%%%%%%%%
{\it Proof}. 
First we show that neither $\be_+$ nor $\be_-$ is fixed by $f$. 
Otherwise, one of them say $\be_+$ is fixed by $f$.  
It is a transverse fixed curve and hence of type I, so any point $p \in \be_+$ 
has index $\nu_p(f) = 0$ and $f$ has a unique fixed point $p_- \in \cE \setminus \be_+$, 
which is transverse, as in Figure \ref{fig:A2b}. 
Thus by Saito's formula \eqref{eqn:Saito} the Lefschetz number 
$L(f) = 2 + \Tr \, \tilde{A} = 2 + 1 = 3$ is expressed as   
%%%%%%%%
$$
L(f) = \nu_{p_-}(f) + \sum_{p \in X_0(f) \setminus \cE} \nu_p(f) + \chi_{\sbe_+} \cdot \nu_{\sbe_+}(f) 
= 1 + \sum_{p \in X_0(f) \setminus \cE} \nu_p(f) + 2 \cdot 1, 
$$
%%%%%%%
which implies $\sum_{p \in X_0(f) \setminus \cE} \nu_p(f) = 0$, that is, 
$f$ has no fixed point in $X \setminus \cE$.  
Since $\be_+$ is a transverse fixed curve, the Toledo-Tong formula \eqref{eqn:TT} 
shows that the holomorphic Lefschetz number $L(f, f^*) = 1 + \delta^{-1}$ for 
the structure sheaf $\cO_X$ is equal to 
%%%%%%%%%%%%%%%%%%%%%%%%%%%%%% eqn:LfO %%%%%%%%%%%%%%%%%%%%%%%%%%%%%%%%%
\begin{equation} \label{eqn:LfO}
L(f, f^*) = \dfrac{1}{1-(\delta^{-1} + \delta^2) + \delta} + \dfrac{1 + \delta}{(1-\delta)^2}, 
\end{equation}
%%%%%%%%%%%%%%%%%%%%%%%%%%%%%%%%%%%%%%%%%%%%%%%%%%%%%%%%%%%%%%%%%%%%%%
where the first and second terms on the RHS are the contributions of $p_-$ and 
$\be_+$ respectively. 
This shows that $\delta$ is a root of the quartic equation 
$z^4-z^3-3z^2-z+1 = 0$, which contradicts the fact that 
$\delta$ is a root of Lehmer's equation $\rL(z) = 0$. 
Hence $\be_+$ cannot be fixed by $f$.   
%%%%%%%%
\par
%%%%%%%%
Next we show that $p_0$ is a transverse fixed point of $f$.   
Otherwise, we may suppose that $f|_{\be_+}$ is a parabolic M\"{o}bius transformation 
while $f|_{\be_-}$ is a non-parabolic M\"{o}bius transformation. 
Since $f$ has a fixed point $p_0$ of multiplicity $\nu_{p_0}(f) \ge 2$, a fixed point 
$p_- \in \be_- \setminus \{p_0\}$ of multiplicity $\nu_{p_-}(f) = 1$ and no fixed curve  
as in Figure \ref{fig:A2b}, Lefschetz formula shows that  
%%%%%%%%
$$
3 = L(f) = \nu_{p_0}(f) + \nu_{p_-}(f) + \sum_{p \in X_0(f) \setminus \cE} \nu_p(f) 
\ge 2 + 1 + 0,  
$$
%%%%%%%%
which implies that $\nu_{p_0}(f) = 2$ and $f$ has no fixed point in $X \setminus \cE$. 
Let us consider the holomorphic Lefschetz number $L(f, f^*) = 1 + \delta^{-1}$ for the 
structure sheaf $\cO_X$ again.  
By a version of the Atiyah-Bott formula allowing isolated multiple fixed points we have 
%%%%%%%%%%%%%%%%%%%%%%%%%%%%%%%%%%%%%%%%%%%%%%%%%%%%%%%%%%%%%%%%%%%%%%%
\begin{equation*}
L(f, f^*) = \dfrac{1}{1-(\delta^{-1} + \delta^2) + \delta} + \Res_{p_0} \omega 
\quad \mbox{with} \quad 
\omega := \dfrac{d z_1 \wedge d z_2 }{(z_1 - f_1(z))(z_2 - f_2(z))},   
\end{equation*}
%%%%%%%%%%%%%%%%%%%%%%%%%%%%%%%%%%%%%%%%%%%%%%%%%%%%%%%%%%%%%%%%%%%%%%
where the second term on the RHS is a Grothendieck residue 
(see Toledo \cite[formula (6.3)]{Toledo}) with $z = (z_1, z_2)$ being local coordinates 
around $p_0 \leftrightarrow z = (0, 0)$ such that $e_+ = \{z_2 = 0\}$ and  
$e_- = \{z_1 = 0\}$ along with a local representation $(f_1(z), f_2(z))$ for $f$. 
We have 
%%%%%%%%%%%%%%%%%%%%%%%%%%%% eqn:GR %%%%%%%%%%%%%%%%%%%%%%%%%%%%%%%%%%%
\begin{equation} \label{eqn:GR}
\Res_{p_0} \omega = \dfrac{1 + \delta}{(1-\delta)^2}. 
\end{equation}
%%%%%%%%%%%%%%%%%%%%%%%%%%%%%%%%%%%%%%%%%%%%%%%%%%%%%%%%%%%%%%%%%%%%%%  
To see this notice that a suitable choice of coordinates makes  
$f_1(z_1, 0) = z_1/(1+z_1)$, $f_2(z_1, 0) = 0$, $f_1(0, z_2) = 0$ and $f_2(0, z_2) = z_2 
\{\delta  + O(z_2) \}$. 
Consider the perturbation $f^t(z_1, z_2) = f(t z_1, t^{-1} z_2)$ 
for a real $t \in (1/2, \, 2)$ with $t \neq 1$.  
Near $p_0 = (0, 0)$ the map $f^t$ has exactly two fixed points;   
$p_0$ itself and $q^t = (\ve, 0)$ with $\ve := (t-1)/t$, both transverse. 
It is easy to see that $(d f^t)_{p_0}$ has eigenvalues 
$t$ and $\delta t^{-1}$. 
It is also easy to see that one eigenvalue of $(d f^t)_{q^t}$ is $t^{-1}$ 
and the other eigenvalue is $t \delta^t$ with   
$\delta^t : = \det( d f^t)_{q^t}$. 
It now follows from $f^* \eta = \delta \eta$ for a nowhere vanishing holomorphic $2$-form 
$\eta$ on the K3 surface $X$ that $\delta^t = \delta(1 + a \ve^2)$ with 
$a = O(1)$ as $\ve \to 0$.  
So
%%%%%%%
\begin{align*}
\Res_{p_0} \omega^t + \Res_{q^t} \omega^t 
&= \dfrac{1}{(1-t)(1-\delta t^{-1})} + \dfrac{1}{(1-t^{-1})(1- t \delta^t)}  \\[2mm]
&= \dfrac{(1-\ve)(1+\delta + a \delta \ve)}{(1-\delta+ \delta \ve)(1-\delta-\ve-a \delta \ve^2)} 
\to \dfrac{1+\delta}{(1-\delta)^2} \quad 
\mbox{as} \quad \ve \to 0, 
\end{align*}
%%%%%%%
where $\omega^t$ is the corresponding perturbation of $\omega$. 
Continuity principle then yields \eqref{eqn:GR}, which in turn  
leads to \eqref{eqn:LfO} and hence a contradiction.  
Thus $p_0$ must be a transverse fixed point of $f$.    
%%%%%%%
\par
%%%%%%%
Thirdly, since $p_0$ is a transverse fixed point, $f|_{\be_{\pm}}$ are 
nontrivial, non-parabolic M\"{o}bius transformations, having exactly two isolated 
fixed points $p_{\pm} \in \be_{\pm} \setminus \{ p_0\}$ as in Figure \ref{fig:A2a}.  
Since $f$ has no fixed curve, the Lefschetz formula reads 
%%%%%%%
$$
3 = L(f) = \nu_{p_0}(f) + \nu_{p_+}(f) + \nu_{p_-}(f) + 
\sum_{p \in X_0(f) \setminus \cE} \nu_p(f) \ge 1 + 1 + 1 + 0, 
$$ 
%%%%%%%
which implies that $p_{\pm}$ must be transverse and there 
is no fixed point in $X \setminus \cE$. 
If the eigenvalues of $(d f)_{p_0}$ are $\delta^{\frac{1}{2}} \alpha^{\pm1}$ then 
those of $(d f)_{p_{\ve}}$ are $\delta^{\frac{1}{2}}(\delta^{\ve} \alpha)^{\pm 1}$ for 
$\ve = \pm 1$. 
Note that $\alpha \neq \delta^{\pm \frac{1}{2}}$, $\delta^{\pm \frac{3}{2}}$, 
because $p_0$ and $p_{\pm}$ are all transverse. 
By the Atiyah-Bott formula the holomorphic Lefschetz number 
$L(f, f^*) = 1 + \delta^{-1}$ for the structure sheaf $\cO_X$ is represented as 
%%%%%%%%
$$
L(f, f^*) = \dfrac{1}{(1-\delta^{\frac{1}{2}} \alpha)(1- \delta^{\frac{1}{2}} \alpha^{-1})} 
+ \dfrac{1}{(1-\delta^{\frac{3}{2}} \alpha)(1- \delta^{-\frac{1}{2}} \alpha^{-1}) } 
+  \dfrac{1}{(1-\delta^{-\frac{1}{2}} \alpha)(1- \delta^{\frac{3}{2}} \alpha^{-1}) }.  
$$
%%%%%%%%
A manipulation of this equation yields $(\alpha + \alpha^{-1})^2 = q(\tau)$ with 
$q(w) := (w^2-3)^2/(w+2)$. 
%%%%%%%
\par
%%%%%%%
Finally we have $0 < q(x_j) < 4$ for $j =1$, $3$, $4$ and $q(x_2) > 4$, hence 
Proposition \ref{prop:siegel} is applicable. 
The multiplicative independence of the eigenvalues of $(d f)_{p_{\ve}}$, $\ve = \pm1$, 
can be verified in the same manner as that of the eigenvalues of $(d f)_{p_0}$ 
(see the proof of Proposition \ref{prop:siegel}).   
We can conclude that $p_0$ and $p_{\pm}$ are the centers of Siegel disks, 
since we have $\tau = x_1$, $x_4$ for the entries of type $\rA_2$ in 
Table \ref{tab:minA}. 
Note that there is no entry of type $\rA_2$ in Table \ref{tab:minB}. 
\hfill $\Box$
%%%%%%%%%%%%%%%%%%%%%%% end proof %%%%%%%%%%%%%%%%%%%%%%%%%%%%%%%%%  
%%%%%%%%%%%%%%%%%%%%%%% sec:comparison %%%%%%%%%%%%%%%%%%%%%%%%%%%
\section{Lattices in Number Fields} \label{sec:lnf} 
%%%%%%%%%%%%%%%%%%%%%%%%%%%%%%%%%%%%%%%%%%%%%%%%%%%%%%%%%%%%%%%%
It is interesting to discuss the relationship between the method of hypergeometric 
groups and that of Salem number fields by McMullen 
\cite{McMullen1} and Gross and McMullen \cite{GM}. 
We briefly review the construction in \cite{McMullen1}. 
Let $S(z)$ be an unramified Salem polynomial of degree $22$ with the 
associated Salem trace polynomial $R(w)$ where $w = z + z^{-1}$. 
Consider the $\bZ$-algebra $L_S := \bZ[z]/(S(z))$ together with its field 
of fractions $K := \bQ[z]/(S(z))$.  
Let $U(w)$ be a unit in $\bZ[w]/(R(w))$ such that the polynomial $R(w)$ 
admits a unique root 
%%%%%%%%%%%%%%%%%%%%%%%% eqn:comp %%%%%%%%%%%%%%%%%%%%%%%%%%%%%%%
\begin{equation} \label{eqn:comp}
\tau \in (-2, \, 2) \quad \mbox{that satisfies} \quad U(\tau) \, R'(\tau) > 0,  
\qquad \mbox{see \cite[Theorem 8.3]{McMullen1}}. 
\end{equation}
%%%%%%%%%%%%%%%%%%%%%%%%%%%%%%%%%%%%%%%%%%%%%%%%%%%%%%%%%%%%%%%%% 
This root $\tau$ corresponds to the ``special trace'' in our  
method.    
By defining the inner product 
%%%%%%%%%%%%%%%%%%%%%%% eqn:ips %%%%%%%%%%%%%%%%%%%%%%%%%%%%%%%%%
\begin{equation} \label{eqn:ips}
(g_1, \, g_2 )_S 
:= \Tr^K_{\bQ} \left( \dfrac{U(w) \, g_1(z) \, g_2(z^{-1})}{R'(w)} \right), 
\qquad g_1(z), \, g_2(z) \in \bZ[z] \bmod S(z), 
\end{equation} 
%%%%%%%%%%%%%%%%%%%%%%%%%%%%%%%%%%%%%%%%%%%%%%%%%%%%%%%%%%%%%%%%
one can make $L_S$ into a K3 lattice equipped with a Hodge structure such 
that multiplication by $z$, that is, $M_z : L_S \to L_S$, $g(z) \mapsto z g(z)$ 
induces a Hodge isometry. 
The map $M_z$ lifts to a K3 surface automorphism $f : X \to X$ such that 
$f^* : H^2(X, \bZ) \to H^2(X, \bZ)$ has characteristic polynomial $S(z)$. 
%%%%%%%
\par
%%%%%%%
%%%%%%%%%%%%%%%%%%%%%%% tab:vs22 %%%%%%%%%%%%%%%%%%%%%%%%%%%%%%%%%
\begin{table}[hh]
\centerline{
\begin{tabular}{c|c|cc|l}
\hline
\\[-4mm] 
case & $R$ &  ST in \cite{McMullen1} & $(1, 1)_S$ &  
ST in this article \\[1mm]
\hline
\\[-4mm] 
1 & $R_1$ & $\by_8$ & $\mi 0$ & $y_1$, $y_2$, $\by_8$, $y_{10}$ \\[1mm]
2 & $R_2$ & $\by_4$ & $\mbox{\boldmath $-2$}$ & $y_1$, $y_2$, $y_3$, $\by_4$, 
$y_7$, $y_8$, $y_{10}$ \\[1mm]
3 & $R_3$ & $\by_7$ & $\mi 0$ & $y_1$, $y_2$, $y_3$, $\by_7$, $y_8$, $y_9$  
\\[1mm]
4 & $R_4$ & $y_5$ & $\mbox{\boldmath $-2$}$ & $y_1$, $y_2$, $y_3$, $y_7$  
\\[1mm]
5 & $R_5$ & $y_4$ & $\mbox{\boldmath $-2$}$ & $y_1$, $y_2$, $y_6$, $y_7$  
\\[1mm]
6 & $R_6$ & $y_8$ & $\mi 0$ & $y_1$, $y_2$, $y_3$, $y_4$, $y_7$, $y_{10}$  \\[1mm]
7 & $R_7$ & $\by_7$ & $\mi 0$ & $y_2$, $y_5$, $\by_7$  \\[1mm]
8 & $R_8$ & $\by_6$ & $\mi 0$ & $y_1$, $y_2$, $y_3$, $y_4$, $y_5$, $\by_6$, $y_8$, $y_9$ 
\\[1mm]
9 & $R_9$    &   $y_6$ & $\mbox{\boldmath $-2$}$ & $y_1$, $y_3$, $y_7$, 
$y_8$, $y_9$, 
$y_{10}$ \\[1mm]
10 & $R_{10}$ &  $y_5$ & $\mbox{\boldmath $-2$}$ & $y_1$, $y_7$, $y_9$  
\\[1mm]
\hline
\end{tabular}}
\caption{Table 4 in McMullen \cite{McMullen1} vs.  
Tables \ref{tab:SvH1}--\ref{tab:SvH3} in this article.} 
\label{tab:vs22} 
\end{table}
%%%%%%%%%%%%%%%%%%%%%%%%%%%%%%%%%%%%%%%%%%%%%%%%%%%%%%%%%%%%%%%%
McMullen \cite[Table 4]{McMullen1} gives a list of ten triples 
$(S_i(z), R_i(w), U_i(w))$, $i = 1, \dots, 10$, to which his method is applied.  
Table \ref{tab:vs22} gives a comparison between his Table 4 and our 
Tables \ref{tab:SvH1}--\ref{tab:SvH3} (in Theorem \ref{thm:SH}), 
showing to what extent the special traces are common or not, 
where the common ones are shown in boldface.  
%%%%%%%%
\par
%%%%%%%
An interesting thought occurs to us when we notice a similarity 
between the lattice in number field $L_S$ and the hypergeometric 
lattice $L_H = L(A, B) = L(\Phi, \Psi)$, that is,   
%%%%%%%
$$
L_S = \langle 1, z, z^2, \dots, z^{21} \rangle_{\bZ} \circlearrowleft M_z  
\quad \mbox{and} \quad   
L_H = \langle \br, B \br, \dots, B^{21} \br \rangle_{\bZ} \circlearrowleft B. 
$$
%%%%%%%
Recall from \eqref{eqn:rr2} and Definition \ref{def:hgk3l} that the vector 
$\br$ is normalized as $( \br, \, \br )_H = \pm 2$. 
So we wonder whether when $(1, \, 1)_S = \pm 2$ we can 
go back and forth between $L_S$ and $L_H$ via the correspondences 
$1 \leftrightarrow \br$ and $M_z \leftrightarrow B$.  
The value of $(1, \, 1)_S$ is twice the coefficient of 
$w^{10}$ in $U(w)$, so cases 2, 4, 5, 9, 10 in Table \ref{tab:vs22} are 
relevant to this question.  
In these cases we can recover the matrix $A = B C$ and hence the 
polynomail  $\Phi(w)$ from McMullen's data, since the reflection in 
the vector $1 \in L_S$ corresponds to the matrix $C$, that is, the 
reflection in $\br \in L_H$.     
The results are given in Table \ref{tab:merhm}, where $\Phi(w)$ contains a 
non-cyclotomic trace factor, which is Salem in cases 4 and 10, but not 
Salem in cases 2, 5 and 9.    
These are not covered by Tables \ref{tab:SvH1}--\ref{tab:SvH3}, 
because there $\Phi(w)$ is restricted to a product of cyclotomic 
trace polynomials (see \eqref{eqn:CTk}). 
%%%%%%%%%%%%%%%%%%%%%%% tab:merhm %%%%%%%%%%%%%%%%%%%%%%%%%%%%%%%% 
\begin{table}[hh] 
\centerline{
\begin{tabular}{cll}
\hline
\\[-4mm] 
case & $R$ & $\Phi(w)$ \\[1mm]
\hline
\\[-4mm] 
2 & $R_2$ & $\CT_3(w) \cdot (w^9-9w^7+25w^5-2w^4-21w^3+7w^2+2w-2)$ \\[1mm]
\hline
\\[-4mm] 
4 & $R_4$ & $\CT_4(w) \cdot \CT_{42}(w) \cdot (w^3-w^2-3w+1)$ \\[1mm]
\hline
\\[-4mm] 
5 & $R_5$ & $w^{10}-10 w^8-2w^7+33w^6+12w^5-37 w^4-16 w^3+6w^2-3w-3$ 
\\[1mm]
\hline
\\[-4mm] 
9 & $R_9$ & $w^{10}-11 w^8-3w^7+42w^6+22w^5-62w^4-49w^3+23w^2+33w+8$ 
\\[1mm]
\hline
\\[-4mm] 
10 & $R_{10}$ & $\CT_4(w) \cdot  
(w^9-w^8-10w^7+7w^6+35w^5-14w^4-48w^3+7w^2+18w-1)$ 
\\[1mm]
\hline
\end{tabular}}
\caption{Examples in McMullen \cite[Table 4]{McMullen1} recovered by hypergeometric method.}  
\label{tab:merhm}
\end{table}
%%%%%%%%%%%%%%%%%%%%%%%%%%%%%%%%%%%%%%%%%%%%%%%%%%%%%%%%%%%%%%%%   
%%%%%%
\par 
%%%%%%
In the other way round some hypergeometric lattices can be realized 
as lattices in number fields. 
Thinking a little bit more generally, let $S(z) \in \bZ[z]$ be a monic, 
irreducible, palindromic polynomial of  degree $2N$ with the associated 
trace polynomial $R(w) \in \bZ[w]$.  
Let $L_H$ be an even unimodular lattice of rank $2N$ equipped with an 
isometry $F : L_H \to L_H$ whose characteristic polynomial is $S(z)$. 
Suppose that $F$ admits a {\sl cyclic vector} $\br$ over $\bZ$, that is, 
%%%%%%%%%%%%%%%%%%%%%%% eqb:cyclic-v %%%%%%%%%%%%%%%%%%%%%%%%%%%%%%
\begin{equation} \label{eqn:cyclic-v}
L_H = \langle \br, F \br, \dots, F^{2 N-1} \br \rangle_{\bZ}. 
\end{equation}
%%%%%%%%%%%%%%%%%%%%%%%%%%%%%%%%%%%%%%%%%%%%%%%%%%%%%%%%%%%%%%%%
\par
%%%%%% 
Let $P_j(w) := z^j + z^{-j}$ for $j \in \bZ_{\ge 0}$ and think of them 
as polynomials in $w = z+ z^{-1}$. 
Then they satisfy three-term recurrence relation 
%%%%%%
$$
P_0(w) = 2, \qquad P_1(w) = w, \qquad P_{j+1}(w) - w \, P_j(w) + P_{j-1}(w) = 0, 
\quad j \ge 1,    
$$
%%%%%% 
which shows that if $j \ge 1$ then $P_j(w)$ is a monic polynomial of degree 
$j$ in $\bZ[w]$.  
Given a polynomial $g(w) \in \bZ[w]$, we denote by $[g(w)]_R \in \bZ$ 
the coefficient of $w^{N-1}$ in the remainder of $g(w)$ divided 
by $R(w)$. 
We define the integers $u_1, \dots, u_N \in \bZ$ inductively by  
%%%%%%%%%%%%%%%%%%%%%%% eqn:uj %%%%%%%%%%%%%%%%%%%%%%%%%%%%%%%%%%%
\begin{equation} \label{eqn:uj}
u_1 = \dfrac{1}{2} (\br, \, \br)_H, \qquad 
u_j = (F^{j-1} \br, \, \br)_H - \sum_{k=1}^{j-1} c_{j k} \, u_k, \quad 
j = 2, \dots, N,  
\end{equation}
%%%%%%%%%%%%%%%%%%%%%%%%%%%%%%%%%%%%%%%%%%%%%%%%%%%%%%%%%%%%%%%% 
where $c_{j k} := [P_{j-1}(w) \cdot w^{N-k}]_R \in \bZ$. 
Note that $u_1 \in \bZ$ since $L_H$ is an even lattice.  
%%%%%%%%%%%%%%%%%%%%%%% thm:LNF %%%%%%%%%%%%%%%%%%%%%%%%%%%%%%%%%
\begin{theorem} \label{thm:LNF} 
Under condition \eqref{eqn:cyclic-v} the pair $(L_H, \, F)$ is isomorphic to 
$(L_S, \, M_z)$ with $L_S := \bZ[z]/(S(z))$ and $M_z$ being multiplication 
by $z$, where $L_S$ carries the inner product \eqref{eqn:ips} with 
%%%%%%%%%%%%%%%%%%%%%%% eqn:U(w) %%%%%%%%%%%%%%%%%%%%%%%%%%%%%%%%
\begin{equation} \label{eqn:U(w)}
U(w) = u_1 \, w^{N-1} + u_2 \, w^{N-2} + \cdots + u_N \in \bZ[w], 
\end{equation}
%%%%%%%%%%%%%%%%%%%%%%%%%%%%%%%%%%%%%%%%%%%%%%%%%%%%%%%%%%%%%%% 
whose coefficients are determined by the recurrence \eqref{eqn:uj}.  
It is a unit in $\bZ[w]/(R(w))$.  
\end{theorem}
%%%%%%%%%%%%%%%%%%%%%%% begin proof %%%%%%%%%%%%%%%%%%%%%%%%%%%%%%
{\it Proof}. 
Identify $(L_H, \, F)$ and $(L_S, \, M_z)$ via $\br \leftrightarrow 1$ and 
$F \leftrightarrow M_z$. 
Export the inner product on $L_H$ to $L_S$ isometrically.  
The existence in $\bQ[w]/(R(w))$ of $U(w)$ that makes \eqref{eqn:ips} 
valid is mentioned in Gross and McMullen \cite[\S4, Remark]{GM}, where 
they suppose $\cO_K = \bZ[z]/(S(z))$ but the existence can be proved   
without this assumption.    
To determine the coefficients of $U(w)$, 
substitute $g_1(z) = z^{j-1}$, $j = 1, \dots, N$, and 
$g_2(z) = 1$ into \eqref{eqn:ips} and 
use \eqref{eqn:U(w)} to have 
%%%%%%%%
\begin{align*}
(F^{j-1} \br, \, \br)_H &= (z^{j-1}, \, 1)_S = 
\Tr^K_{\bQ} \left( \dfrac{U(w) \cdot z^{j-1}}{R'(w)} \right) 
= \Tr^{J}_{\bQ} \circ \Tr^K_J \left( \dfrac{U(w) \cdot z^{j-1}}{R'(w)} \right)   
\\[2mm]
&= \Tr^{J}_{\bQ} \left( \dfrac{U(w) \, P_{j-1}(w)}{R'(w)} \right) 
= [U(w) \, P_{j-1}(w)]_R = \sum_{k=1}^N c_{j k} \, u_k = \sum_{k=1}^j c_{j k} \, u_k,  
\end{align*}
%%%%%%%%  
where $J := \bQ[w]/(R(w))$, the fifth equality is by residue calculus 
as in \cite[page 222]{McMullen1} and the final equality follows from 
$c_{j k} = 0$ for $j < k \leq N$. 
For $j = 1$ we have $(\br, \, \br)_H = c_{11} \, u_1 = 2 u_1$, which yields 
the first equality in \eqref{eqn:uj}. 
For $j \ge 2$ we have the second equality in \eqref{eqn:uj}, since 
$c_{jj} = 1$. 
Thus $U(z) \in \bQ[w]$ belongs to $\bZ[w]$.   
The discriminant of $L_S$ is $\disc(L_S) = \pm \det^2 M_{U(w)}$, where 
$M_{U(w)} : \bZ[w]/(R(w)) \to \bZ[w]/(R(w))$ is multiplication by $U(w)$. 
Since $L_S \cong L_H$ is unimodular, i.e. $\disc(L_S) = \pm 1$, we have  
$\det M_{U(w)} = \pm 1$ and hence $U(w)$ is a unit in $\bZ[w]/(R(w))$.  
\hfill $\Box$ \par\medskip 
%%%%%%%%%%%%%%%%%%%%%%% end proof %%%%%%%%%%%%%%%%%%%%%%%%%%%%%%%
Theorem \ref{thm:LNF} can be applied to the unimodular hypergeometric 
lattice $L_H= L(A, B)$ that arises from an irreducible hypergeometric group 
$H = H(A, B)$, since it satisfies condition \eqref{eqn:cyclic-v} for 
$F = B$.   
In particular for hypergeometric K3 lattices we have the following.  
%%%%%%%%%%%%%%%%%%%%%%% cor:LNF %%%%%%%%%%%%%%%%%%%%%%%%%%%%%%%%
\begin{corollary} \label{cor:LNF} 
In the situation of Theorem $\ref{thm:main3}$ suppose that the polynomial 
$\psi(z)$ is an unramified Salem polynomial $S(z)$ of degree $22$.  
Then the pair $(L_H, \, B)$ admits McMullen's construction with unit 
$U(w) \in \bZ[w]/(R(w))$ determined by \eqref{eqn:uj} and \eqref{eqn:U(w)}, 
and the special trace $\tau$ satisfies the compatibility condition 
\eqref{eqn:comp}.    
\end{corollary}
%%%%%%%%%%%%%%%%%%%%%%%%%%%%%%%%%%%%%%%%%%%%%%%%%%%%%%%%%%%%%% 
\par
%%%%%%%%
As an illustration we give an example from the top entry of 
Table \ref{tab:SvH1}. 
The $B$-basis version of formula \eqref{eqn:xi} gives the values of 
$(B^{j-1} \br, \, \br)_H$ for $j = 1, \dots, 11$, then formulas 
\eqref{eqn:uj} and \eqref{eqn:U(w)} yield   
%%%%%%%
$$
U(w) = -w^{10} + 6 w^9 - 7 w^8 - 22 w^7 + 54 w^6 - 4 w^5-70 w^4 
+ 36 w^3 + 24 w^2 - 16 w. 
$$ 
%%%%%%
\par
%%%%%%
As for non-projective K3 surface automorphisms with minimum entropy, 
McMullen \cite{McMullen3} gives only one such example. 
It has characteristic polynomial $\rL(z) \cdot (z-1)^9 (z+1)(z^2+1)$, 
special trace $x_4$ in \eqref{eqn:lehconj}, and root system of 
Dynkin type $\rE_8 \oplus \rA_2 \oplus \rA_2$.  
Our Table \ref{tab:minB} contains three entries with the same data. 
It is interesting to compare them in a deeper level, but it needs  
to discuss gluing of lattices. 
%%%%%%%%%%%%%%%%%%%%%%% end main text %%%%%%%%%%%%%%%%%%%%%%%%%%%%%
\par\vspace{3mm} \noindent
%%%%%%%%%%%%%%%%%%%%%% acknowledgment %%%%%%%%%%%%%%%%%%%%%%%%%%%%
{\bf Acknowledgments}. 
This work is supported by Grant-in-Aid for Scientific Research, 
JSPS, 19K03575 (C).  
The authors thank Takato Uehara for fruitful discussions. 
Their appreciations are also due to anonymous reviewers whose questions 
and comments led to the discussion in \S \ref{sec:lnf} as well as to 
a better presentation of this article.                   
%%%%%%%%%%%%%%%%%%%%%%%%%% References %%%%%%%%%%%%%%%%%%%%%%%%%%%%%%%%%
 
%%%%%%%%%%%%%%%%%%%%%%%%%%%%%%%%%%%%%%%%%%%%%%%%%%%%%%%%%%%%%%%%%%%%%%
\end{document}